\documentclass[12pt]{scrartcl}    

\usepackage{a4} 
\usepackage{amsmath}    
\usepackage{paralist}
\usepackage{latexsym} 
\usepackage{amssymb}  
\usepackage{amsfonts}  
\usepackage{mathrsfs}  
\usepackage{dsfont}
\usepackage{xcolor}
\usepackage{bbm,exscale}
\definecolor{Myblue}{rgb}{0,0,0.6}  
\usepackage[colorlinks,citecolor=Myblue,linkcolor=Myblue,urlcolor=Myblue,pdfpagemode=UseNone]{hyperref}
\usepackage{amsthm}
\usepackage{accents}
\usepackage[square,numbers,sort&compress]{natbib} 
\usepackage[all,cmtip]{xy}
\usepackage{ifthen} 
\usepackage{bbding}
\usepackage{stmaryrd}  
\usepackage{verbatim}
\usepackage{bbding} 
\usepackage{wasysym}  
\usepackage{soul}  
\usepackage[yyyymmdd,hhmmss]{datetime}
\usepackage{booktabs}
\usepackage{enumitem}
\usepackage{color}
\usepackage{textcomp}
\usepackage{gensymb}
\usepackage{pdfpages}
\usepackage{tcolorbox}
\usepackage[bb=boondox]{mathalfa} 
\usepackage{tikz}
\usepackage{tikz-cd}
\usetikzlibrary{calc}
\usetikzlibrary{decorations.markings}
\usetikzlibrary{fadings,decorations.pathreplacing}
\usetikzlibrary{matrix,arrows}
\usetikzlibrary{patterns}
\usetikzlibrary{arrows,calc,decorations.pathreplacing,decorations.markings,shapes.geometric,shadows}

\DeclareSymbolFont{bbold}{U}{bbold}{m}{n}
\DeclareSymbolFontAlphabet{\mathbbold}{bbold}

\tikzset{
	string/.style={draw=#1, postaction={decorate}, decoration={markings,mark=at position .51 with {\arrow[draw=#1]{>}}}},
	costring/.style={draw=#1, postaction={decorate}, decoration={markings,mark=at position .51 with {\arrow[draw=#1]{<}}}},
	ostring/.style={draw=#1, postaction={decorate}, decoration={markings,mark=at position .47 with {\arrow[draw=#1]{>}}}},
	ustring/.style={draw=#1, postaction={decorate}, decoration={markings,mark=at position .56 with {\arrow[draw=#1]{>}}}},
	oostring/.style={draw=#1, postaction={decorate}, decoration={markings,mark=at position .43 with {\arrow[draw=#1]{>}}}},
	uustring/.style={draw=#1, postaction={decorate}, decoration={markings,mark=at position .59 with {\arrow[draw=#1]{>}}}},
	directed/.style={string=blue!50!black}, 
	odirected/.style={ostring=blue!50!black}, 
	udirected/.style={ustring=blue!50!black}, 
	oodirected/.style={oostring=blue!50!black}, 
	uudirected/.style={uustring=blue!50!black},     
	redirected/.style={costring= blue!50!black},
	redirectedgreen/.style={costring= green!50!black},
	directedgreen/.style={string= green!50!black},
}

\tikzset{-dot-/.style={decoration={
			markings,
			mark=at position 0.5 with {\fill circle (2pt);}},postaction={decorate}}}

\tikzset{
	Fdot/.style={circle, draw, fill, inner sep=0pt}, 
	Odot/.style={circle, draw, inner sep=0.1pt, minimum size=0.1cm}
}

\newcommand\tikzzbox[1]
{#1}

\usepackage{tikzit}

\tikzstyle{black dot}=[fill=black, draw=black, shape=circle, minimum size=3pt, inner sep=0pt]
\tikzstyle{blue cloud}=[fill={rgb,255: red,130; green,135; blue,200}, draw={rgb,255: red,130; green,135; blue,200}, shape=circle, minimum size=25pt, inner sep=0pt, tikzit fill={rgb,255: red,130; green,135; blue,200}, tikzit draw={rgb,255: red,130; green,135; blue,200}, fill opacity=0.3, draw opacity=0.3]
\tikzstyle{green cloud}=[fill={rgb,255: red,111; green,200; blue,151}, draw={rgb,255: red,111; green,200; blue,151}, shape=circle, minimum size=25pt, inner sep=0pt, fill opacity=0.3, draw opacity=0.3, tikzit fill={rgb,255: red,111; green,200; blue,151}, tikzit draw={rgb,255: red,111; green,200; blue,151}]

\tikzstyle{blue mid arrow}=[-, draw={rgb,255: red,0; green,0; blue,128}, postaction={on each segment={mid-arrow}}, line width=1.5pt, tikzit draw={rgb,255: red,0; green,0; blue,128}]
\tikzstyle{blue mid arrow bg}=[-, draw={rgb,255: red,0; green,0; blue,128}, tikzit draw={rgb,255: red,0; green,0; blue,128}, line width=1.5pt, postaction={on each segment={mid-arrow}}, opacity=0.4]
\tikzstyle{blue}=[-, draw={rgb,255: red,0; green,0; blue,128}, line width=1.5pt, tikzit draw={rgb,255: red,0; green,0; blue,128}]
\tikzstyle{blue bg}=[-, draw={rgb,255: red,0; green,0; blue,128}, tikzit draw={rgb,255: red,0; green,0; blue,128}, line width=1.5pt, opacity=0.4]
\tikzstyle{blue dashed}=[-, draw={rgb,255: red,0; green,0; blue,128}, line width=1.5pt, dashed, tikzit draw={rgb,255: red,0; green,0; blue,128}]
\tikzstyle{black}=[-, draw=black, line width=0.75pt, tikzit draw=black, fill=none]
\tikzstyle{black bg}=[-, draw=black, line width=0.75pt, tikzit draw=black, opacity=0.4]
\tikzstyle{black dashed}=[-, draw=black, line width=0.5pt, dashed, tikzit draw=black]
\tikzstyle{black dashed bg}=[-, draw=black, line width=0.5pt, dashed, tikzit draw=black, opacity=0.4]
\tikzstyle{red mid arrow}=[-, postaction={on each segment={mid-arrow}}, line width=1.5pt, draw={rgb,255: red,204; green,0; blue,0}, tikzit draw={rgb,255: red,204; green,0; blue,0}]
\tikzstyle{red mid arrow bg}=[-, postaction={on each segment={mid-arrow}}, line width=1.5pt, draw={rgb,255: red,204; green,0; blue,0}, tikzit draw={rgb,255: red,204; green,0; blue,0}, opacity=0.4]
\tikzstyle{red}=[-, draw={rgb,255: red,204; green,0; blue,0}, line width=1.5pt, tikzit draw={rgb,255: red,204; green,0; blue,0}]
\tikzstyle{red bg}=[-, draw={rgb,255: red,204; green,0; blue,0}, line width=1.5pt, tikzit draw={rgb,255: red,204; green,0; blue,0}, opacity=0.4]
\tikzstyle{grey}=[-, draw={rgb,255: red,128; green,128; blue,128}, draw opacity=0.3, line width=10pt, tikzit draw={rgb,255: red,128; green,128; blue,128}]
\tikzstyle{orange front}=[-, draw={rgb,255: red,255; green,225; blue,193}, opacity=0.8, line width=20pt, tikzit draw={rgb,255: red,255; green,225; blue,193}, line cap=round, line join=round]
\tikzstyle{black arrow}=[tikzit draw=, ->]
\tikzstyle{black hook arrow}=[tikzit draw=,left hook ->]

\usetikzlibrary{decorations.pathreplacing,decorations.markings}
\tikzset{
	on each segment/.style={
		decorate,
		decoration={
			show path construction,
			moveto code={},
			lineto code={
				\path [#1]
				(\tikzinputsegmentfirst) -- (\tikzinputsegmentlast);
			},
			curveto code={
				\path [#1] (\tikzinputsegmentfirst)
				.. controls
				(\tikzinputsegmentsupporta) and (\tikzinputsegmentsupportb)
				..
				(\tikzinputsegmentlast);
			},
			closepath code={
				\path [#1]
				(\tikzinputsegmentfirst) -- (\tikzinputsegmentlast); 
			},
		},
	},
	mid-arrow/.style={postaction={decorate,decoration={
				markings,
				mark=at position .5 with {\arrow[#1]{stealth}}
	}}},
	red semicircle/.style={postaction={decorate,decoration={
				markings,
				mark=at position .65 with {
					\arrow[#1]{Circle[left,fill=red,length=6pt,width=6pt]}
				}
	}}},
}

\vfuzz \hfuzz
\makeatletter
\newcommand{\raisemath}[1]{\mathpalette{\raisem@th{#1}}}
\newcommand{\raisem@th}[3]{\raisebox{#1}{$#2#3$}}
\makeatother

\newcommand{\btimes}{\mathbin{\square}} 
\newcommand*{\longhookleftarrow}{\ensuremath{\leftarrow\joinrel\relbar\joinrel\rhook}}
\newcommand*{\longhookrightarrow}{\ensuremath{\lhook\joinrel\relbar\joinrel\rightarrow}}

\newcommand{\I}{\text{i}}
\newcommand{\B}{\mathcal{B}}
\newcommand{\Bfd}{\mathcal{B}^{\textrm{fd}}}

\newcommand{\Beq}{\B_{\mathrm{eq}}}

\newcommand{\Cc}{\mathcal{C}}
\newcommand{\D}{\mathds{D}}

\newcommand{\Q}{\mathds{Q}}
\newcommand{\R}{\mathds{R}}
\newcommand{\Z}{\mathds{Z}}

\newcommand{\Zc}{\mathcal{Z}}

\def\1{\ifmmode\mathrm{1\!l}\else\mbox{\(\mathrm{1\!l}\)}\fi}
\newcommand{\one}{\mathbbm{1}}
\newcommand{\be}{\begin{equation}}
  \newcommand{\ee}{\end{equation}}
\newcommand{\bes}{\begin{equation*}}
  \newcommand{\ees}{\end{equation*}}

\newcommand{\id}{\operatorname{id}}
\newcommand{\Id}{\operatorname{Id}}

\newcommand{\Hom}{\operatorname{Hom}}
\newcommand{\Aut}{\operatorname{Aut}}
\newcommand{\End}{\operatorname{End}}

\def\LG{\mathcal{LG}}
\def\LGgr{\mathcal{LG}^{\mathrm{gr}}}

\def\LGs{\mathcal{LG}'}

\def\LGs{\mathcal{LG}^{\textrm{{\tiny $\bullet/2$}}}}

\newcommand{\ev}{\operatorname{ev}}

\newcommand{\tev}{\widetilde{\operatorname{ev}}}
\newcommand{\coev}{\operatorname{coev}}
\newcommand{\tcoev}{\widetilde{\operatorname{coev}}}
\def\lra{\longrightarrow}

\def\lmt{\longmapsto}

\DeclareMathOperator{\Jac}{Jac}

\newcommand{\FD}[1]{\textrm{2D}^{#1}}
\newcommand{\Fun}{\textrm{Fun}^{\mathrm{sm}}}
\newcommand{\Nat}{\textrm{Nat}}

\newcommand{\Spin}{\textrm{Spin}}
\newcommand{\GL}{\textrm{GL}}
\newcommand{\SO}{\textrm{SO}}

\newcommand{\Bord}{\textrm{Bord}}
\newcommand{\BBord}{\mathds{B}\mathbbold{o}\mathbb{r}\mathbbold{d}}
\newcommand{\Bordfr}{\Bord_{2,1,0}^{\textrm{fr}}}
\newcommand{\Bordor}{\Bord_{2,1,0}^{\textrm{or}}}
\newcommand{\Bordrspin}{\Bord_{2,1,0}^{r\textrm{-spin}}}

\newcommand{\zz}{\mathcal{Z}}

\newcommand{\Vect}{\operatorname{Vect}}
\newcommand{\Vectk}{\operatorname{Vect}_\Bbbk}

\newcommand{\eps}{\varepsilon}

\newcommand{\sta}{\boxempty}
\newcommand{\fus}{\otimes}

\newcommand{\dual}{\#}
\newcommand{\Ae}{A^{\textrm{e}}}

\newcommand{\dX}{{}^\dagger\hspace{-1.8pt}X}
\newcommand{\Xd}{X^\dagger}

\newcommand{\deqX}{{}^\star\hspace{-1.8pt}X} 
 
\newcommand{\dY}{{}^\dagger\hspace{-0.3pt}Y}

\newcommand\arxiv[2]      {\href{https://arXiv.org/abs/#1}{#2}}
\newcommand\doi[2]        {\href{https://dx.doi.org/#1}{#2}}

\allowdisplaybreaks

\deffootnote[1em]{1em}{1em}{\textsuperscript{\thefootnotemark}}

\theoremstyle{definition} 
\newtheorem{definition}{Definition}
\newtheorem{proposition}[definition]{Proposition}
\newtheorem{theorem}[definition]{Theorem}

\newtheorem{lemma}[definition]{Lemma}
\newtheorem{corollary}[definition]{Corollary}
\newtheorem{remark}[definition]{Remark}

\newtheorem{convention}[definition]{Convention}
\newtheorem{example}[definition]{Example}

\numberwithin{equation}{section}
\numberwithin{definition}{section}
\numberwithin{figure}{section}

\newcommand\void[1]{}

\begin{document}

\title{%
Fully extended $\boldsymbol{r}$-spin TQFTs%
}

\author{%
	Nils Carqueville$^*$ \quad
	L\'or\'ant Szegedy$^\#$ 
	\\[0.5cm]
	\normalsize{\texttt{\href{mailto:nils.carqueville@univie.ac.at}{nils.carqueville@univie.ac.at}}} \\  %
	\normalsize{\texttt{\href{mailto:lorant.szegedy@ist.ac.at}{lorant.szegedy@ist.ac.at}}}
	\\[0.3cm]  %
	\hspace{-1.2cm} {\normalsize\slshape $^*$Universit\"at Wien, Fakult\"at f\"ur Physik, Boltzmanngasse 5, Wien, \"{O}sterreich}\\[-0.1cm]
	\hspace{-1.2cm} {\normalsize\slshape $^\#$Institute for Science and Technology Austria, Am Campus 1, Klosterneuburg, \"{O}sterreich}
}

\date{}
\maketitle

\begin{abstract} 
We prove the $r$-spin cobordism hypothesis in the setting of (weak) 2-categories for every positive integer $r$: The 2-groupoid of 2-dimensional fully extended $r$-spin TQFTs with given target is equivalent to the homotopy fixed points of an induced $\textrm{Spin}_2^r$-action. In particular, such TQFTs are classified by fully dualisable objects together with a trivialisation of the $r$-th power of their Serre automorphisms. For $r=1$ we recover the oriented case (on which our proof builds), while ordinary spin structures correspond to $r=2$. 

To construct examples, we explicitly describe $\textrm{Spin}_2^r$-homotopy fixed points in the equivariant completion of any symmetric monoidal 2-category. We also show that every object in a 2-category of Landau--Ginzburg models gives rise to fully extended spin TQFTs, and that half of these do not factor through the oriented bordism 2-category. 
\end{abstract}

\newpage

\tableofcontents

\newpage

\section{Introduction and summary}

The spin group $\Spin_n$ in dimension~$n$ is by definition the double cover of the group of rotations $\SO_n$ in Euclidean space~$\R^n$. 
A spin structure on an $n$-dimensional oriented manifold is a lift of its tangent bundle along the covering $\Spin_n \lra \SO_n$. 
Such geometric structures and their close cousins in Lorentzian geometry are fundamental in theoretical physics, since e.\,g.\ electrons are classically modelled as sections of spin bundles. 

More generally, for any continuous group homomorphism $\xi\colon G \lra \GL_n$, a tangential structure on an $n$-dimensional manifold~$M$ is a principal $G$-bundle on~$M$ together with a bundle map to the frame bundle of~$M$ that is compatible with~$\xi$ (see Section~\ref{subsec:ClosedrSpinTQFTs} for details). 
The case of spin structures is precisely when~$\xi$ is the covering map $\Spin_n\lra \SO_n$ post-composed with the inclusion $\SO_n\subset\GL_n$; in the case of orientations~$\xi$ is just that inclusion, while in the case of framings~$\xi$ is the inclusion of the trivial group into~$\GL_n$. 

Given the relevance of spin structures in physics, and the motivation to study functorial topological quantum field theories (TQFTs) as a means to gain insight into physics, it is natural to consider spin TQFTs. 
These are (higher) symmetric monoidal functors on (higher) categories of bordisms with prescribed spin structures. 
The case of closed spin TQFTs in dimension $n=2$ was first considered in \cite{MooreSegal, BarrettTavares, NovakRunkel, RunkelSzegedyArf}, and in \cite{MooreSegal, SternSzegedy} they were classified\footnote{In fact \cite[Thm.\,5.2.1]{SternSzegedy} provides a classification of open/closed $r$-spin TQFTs, of which the closed case for $r=2$ discussed here is a special case.} 
in terms of ``closed $\Lambda_2$-Frobenius algebras'' (see Section~\ref{subsec:ClassificationThroughClosedLambdaFrob} for the definition). 
Such algebraic structures formalise the relation between topological Neveu--Schwarz and Ramond sectors, examples of which can be obtained as a $\Z_2$-graded version of the centre construction of \cite{lp0602047}. 
In particular, there is a $(1|1)$-dimensional example in $\Vect_{\mathds{C}}^{\Z_2}$ whose associated TQFT computes the Arf invariant of spin surfaces. 
Not many other explicit examples have been studied in the literature, and all previously known classes of examples are constructed from semisimple algebraic data. 

In the setting of symmetric monoidal $(\infty,n)$-categories, fully extended TQFTs with $G$-structure are widely believed to be classified by homotopy fixed points of a $G$-action (induced from the $G$-action on framed bordisms) on the maximal $\infty$-subgroupoids of fully dualisable objects in the target $(\infty,n)$-categories. 
This is described in significant, yet non-exhaustive, detail in \cite{l0905.0465}. 
To our knowledge, this general version of the cobordism hypothesis, originally put forward in \cite{BDpaper}, is established as a theorem only up to a completion of the extended proof sketch in \cite{l0905.0465}, or up 
to a conjecture on the relation between factorisation homology and adjoints, see \cite[Conj.\,1.2]{AyalaFrancis2017CH}. 

On the other hand, in dimension $n=2$ and in the setting of (weak) 2-categories the cobordism hypothesis for the framed and oriented case was proved explicitly in \cite{Pstragowski} and \cite{HSV, HV, Hesse}, respectively: 
For any symmetric monoidal 2-category~$\B$ the 2-groupoid of fully extended framed TQFTs $\Bordfr \lra \B$ is equivalent to the maximal sub-2-groupoid $(\Bfd)^\times$ of fully dualisable objects in~$\B$, while fully extended oriented TQFTs $\Bordor\lra\B$ are described by $\SO_2$-homotopy fixed points. 
The latter are objects of a 2-groupoid $[(\Bfd)^\times]^{\SO_2}$ and correspond to pairs $(\alpha,\lambda)$, where $\alpha\in\Bfd$ and $\lambda\colon S_\alpha \cong 1_\alpha$ is a trivialisation of the Serre automorphism of~$\alpha$. 
In Sections~\ref{subsubsec:DualeAndSerre} and \ref{sec:framed-CH}--\ref{subsubsec:CHTangential} we recall the notions just mentioned, in particular how the Serre automorphism $S_\alpha\colon \alpha \lra \alpha$, defined in~\eqref{eq:def-serre}, corresponds to one full rotation of frames.

\subsubsection*{\texorpdfstring{$\boldsymbol{r}$}{r}-spin cobordism hypothesis}

In the present paper we classify fully extended spin TQFTs valued in an arbitrary symmetric monoidal 2-category~$\B$ (Section~\ref{sec:FullyExtendedRSpinTQFTs}), and we construct a number of examples (Section~\ref{sec:Examples}). 
More precisely, we consider $r$-spin TQFTs for any positive integer~$r$. 
Recall that while for $n\geqslant 3$, the double cover $\Spin_n\lra \SO_n$ is also the universal cover, this is not true for $n\leqslant 2$. 
Hence there is less reason to single out double covers of $\SO_2$ and instead consider the $r$-fold cover $\Spin_2^r \lra \SO_2$ for all $r\in\Z_{\geqslant 1}$.\footnote{For $n=1$, we have $\SO_1 = \{1\}$, and its $r$-fold cover is the unique map $\{1,\dots,r\}\lra \{1\}$.} 
Note that necessarily $\Spin_2^r \cong \SO_2$ as groups, and that 
by definition $\Spin_2 = \Spin_2^2$ and $\Spin_2^1 = \SO_2$. 

Following \cite{spthesis}, in Section~\ref{subsec:bord-2-cat} we describe a 2-category $\Bordrspin$ of bordisms with $r$-spin structure related to $\xi\colon \Spin^r_2 \lra \SO_2 \longhookrightarrow \GL_2$, and in Section~\ref{subsubsec:CHTangential} we construct a 2-category $\FD{r}(\Bfd)$ whose objects are pairs $(\alpha,\theta)$, where $\alpha\in\Bfd$ and  $\theta\colon S_\alpha^r \cong 1_\alpha$. 
Then we prove (Lemma~\ref{lem:r-spin-fixed-points} and Theorem~\ref{thm:r-spin-CH}): 

\vspace{0.3cm} 
\noindent
\textbf{Theorem ($\boldsymbol{r}$-spin cobordism hypothesis). } 
Let~$\B$ be a symmetric monoidal 2-category and let $r\in\Z_{\geqslant 1}$. 
The 2-groupoid of fully extended $r$-spin TQFTs valued in~$\B$ is equivalent to the homotopy fixed points $[(\Bfd)^\times]^{\Spin_2^r}$. 
This in turn is equivalent to $\FD{r}((\Bfd)^\times)$, and under these equivalences we have 
\be 
\begin{tikzpicture}[
baseline=(current bounding box.base),
descr/.style={fill=white,inner sep=3.5pt},
normal line/.style={->}, 
baseline=-0.57cm
]
\matrix (m) [matrix of math nodes, row sep=0.8em, column sep=2em, text height=1.5ex, text depth=0.1ex] {%
	\Fun\big(\Bordrspin,\B\big)   &    \big[ (\Bfd)^{\times}\big]^{\Spin_2^r}   &    \FD{r}\big((\Bfd)^\times\big) 
	\\
	\zz & & 
	\big( \zz(+), S_{\zz(+)}^r \cong 1_{\zz(+)} \big) \, . 
	\\
};
\path[font=\footnotesize] (m-1-1) edge[->] node[above] {$ \cong $} (m-1-2);
\path[font=\footnotesize] (m-1-2) edge[->] node[above] {$ \cong $} (m-1-3);
\path[font=\footnotesize] (m-2-1) edge[|->] node[above] {} (m-2-3);
\end{tikzpicture}
\ee 
\vspace{0.3cm}  

Put differently, (fully) extended $r$-spin TQFTs are classified by what they assign to the positively framed point $+\in \Bordrspin$ together with a trivialisation of the $r$-th power of the associated Serre automorphism. 
The main ingredients of the proof are a generators-and-relations presentation of $\Bordrspin$, inspired by the work \cite{HV}, and an explicit description of $r$-spin bordisms in terms of holonomies, following \cite{RandalWilliams}.

\subsubsection*{Examples}

The choice of target 2-category is essential for extended TQFTs. 
To broaden the class of known $r$-spin TQFTs, in Section~\ref{sec:Examples} we explicitly describe $\Spin_2^r$-homotopy fixed points in the ``equivariant completion'' $\Beq$ of any given symmetric monoidal 2-category~$\B$. 
As introduced in \cite{cr1210.6363} and reviewed in Section~\ref{subsec:EquivariantCompletion}, objects in~$\Beq$ are pairs $(\alpha,A)$, where $\alpha\in\B$ and $A\in\B(\alpha,\alpha)$ is endowed with the structure of a $\Delta$-separable Frobenius algebra, while 1- and 2-morphisms are bimodules and bimodule maps. 
We show (see Corollary~\ref{cor:SpinTQFTwithBeq}, and~\eqref{eq:NakayamaA} for the definition of the Nakayama automorphism $\gamma_A\colon A \lra A$): 

\vspace{0.3cm} 
\noindent
\textbf{Proposition. } 
Let $(\alpha,A)\in\Beq$ be such that $\alpha\in\Bfd$ as well as $S_\alpha^r\cong 1_\alpha$ and $\gamma_A^r = 1_A$ in~$\B$. 
Then there is an $r$-spin TQFT
\begin{align}
\zz\colon \Bordrspin & \lra \Beq \nonumber
\\ 
+ & \lmt (\alpha,A) \, .  
\end{align} 
\vspace{0.01cm} 

Moreover, in Section~\ref{subsubsec:FrobenisAlgebraInBeq} we explain how to compute the invariants such TQFTs associate to $r$-spin surfaces, by explicitly constructing the closed $\Lambda_r$-Frobenius algebras which classify the underlying non-extended TQFTs. 

An advantage of considering $\Beq$-valued (as opposed to $\B$-valued, for a given~$\B$) TQFTs is as follows. 
As explained in Remark~\ref{rem:PivotalNotGoodForSpin}, $r$-spin TQFTs valued in a pivotal 2-category~$\B$ cannot detect all $r$-spin structures if $r\geqslant 3$. 
However, the equivariant completion~$\Beq$ of a pivotal 2-category~$\B$ is itself not pivotal. 

\medskip 

As a specific example of a target~$\B$, in Section~\ref{subsec:LGmodels} we consider the symmetric monoidal 2-category~$\LG$ of Landau--Ginzburg models, constructed in \cite{cm1208.1481, CMM}. 
(Examples of extended 2-spin TQFTs were first considered in \cite{GunninghamSpinHurwitz}.) 
Objects of~$\LG$ are ``potentials'' $W\in\Bbbk[x_1,\dots,x_n]$ that describe isolated singularities, and Hom categories are homotopy categories of matrix factorisations. 
In \cite{CMM} it was observed that every object in~$\LG$ is fully dualisable, and that precisely those potentials $W(x_1,\dots,x_n)$ that depend on an even number of variables give rise to fully extended oriented TQFTs. 
Moreover, these oriented TQFTs indeed extend the closed TQFTs associated to the (generically non-semisimple) Jacobi algebras $\Jac_W$ to the point. 
In light of the $r$-spin cobordism hypothesis proved in Section~\ref{sec:FullyExtendedRSpinTQFTs}, it is straightforward to extend these results as follows (Theorem~\ref{thm:LGspin}):\footnote{In the case of Landau--Ginzburg models the vector space of automorphisms of the identity 1-morphism on any object is 1-dimensional. Hence the choice of trivialisation of the square of the Serre automorphism is unique up to a non-zero scalar.}

\vspace{0.3cm} 
\noindent
\textbf{Theorem. } 
Every object $W(x_1,\dots,x_n) \in \LG$ gives rise to an extended 2-spin TQFT valued in $\LG$. 
These TQFTs factor through the oriented bordism 2-category iff~$n$ is even.
\vspace{0.3cm} 

Explicitly, the 2-spin TQFT associated to an object $W\in\LG$ with an odd number of variables consists of the even Neveu--Schwarz sector $\Jac_W \in \Vect_\Bbbk \subset \Vect_\Bbbk^{\Z_2}$ and the odd Ramond sector $\Jac_W[1] \in \Vect_\Bbbk^{\Z_2}$, together with the structure maps described in general in Section~\ref{subsubsec:closed-Lr-Fa-B}. 
Moreover, in Example~\ref{exa:LGeqSimplesExample} we illustrate how to apply our results on equivariant completion (Section~\ref{subsec:EquivariantCompletion}) to a variant of $\B=\LG$ and explicitly compute the invariants of $r$-spin tori in the simplest non-trivial (and novel) example.

\subsubsection*{Examples not treated in this paper}

We close this introductory section with a few comments on potential further applications of the $r$-spin cobordism hypothesis. 
Besides the 2-categories $\textrm{Alg}_\Bbbk$ and $\LG$ (as well as their variants with additional $\Z_2$-, $\Z$- or $\Q$-gradings), it is natural to consider the 2-category $\mathcal{V}\textrm{ar}$ of \cite{cw1007.2679} of smooth projective varieties and derived categories (see also Example~\ref{exa:SymmetricMonoidal2Categories}), which appears in the study of B-twisted sigma models. 
The 2-category $\mathcal{V}\textrm{ar}$ has a natural symmetric monoidal structure \cite{BanksOnRozanskyWitten}. 
As explained in \cite{l0905.0465, cw1007.2679}, the Serre automorphism~$S_X$ of $X\in \mathcal{V}\textrm{ar}$ can be identified with the Serre functor of the derived category associated to~$X$. 

In \cite{KuznetsovCalabiYau}, Kuznetsov constructs ``fractional Calabi--Yau categories'' $\mathcal A_X$ as the admissible subcategories of semiorthogonal decompositions of derived categories of certain varieties $X\in \mathcal{V}\textrm{ar}$. 
This means in particular that~$\mathcal A_X$ is a triangulated category with suspension functor~$\Sigma$, such that~$\mathcal A_X$ has a Serre functor~$S$ which satisfies $S^q \cong \Sigma^p$ for some $p,q\in\Z$ with $q\neq 0$. 
It follows that the orbit category $\mathcal A_X/\Z$ has a Serre functor whose $(p-q)$-th power is trivialisable, see e.\,g.\ \cite[Thm.\,5.14]{GrantSerreFunctors}. 

It is tempting to expect that some of the fractional Calabi--Yau categories constructed in \cite[Sect.\,4]{KuznetsovCalabiYau} classify $(p-q)$-spin TQFTs whose target is $\mathcal{V}\textrm{ar}$ up to the $\Z$-action quotiented out in orbit categories. 
This is possible only if one can identify the Serre functor of $\mathcal A_X/\Z$ with the Serre automorphism of some other object in the target 2-category. 
More generally, we could work in the 2-category of smooth and proper triangulated differential graded categories described in \cite[App.\,A]{bfk1105.3177v3}. 
In this setting, both the geometric constructions of \cite{KuznetsovCalabiYau} and the representation theoretic examples of fractional Calabi--Yau categories in \cite[Sect.\,6]{GrantSerreFunctors} may lead to interesting $r$-spin TQFTs.

\subsubsection*{Acknowledgements} 

We are grateful to  
Tobias Dyckerhoff, 
Lukas M\"uller, 
Ingo Runkel, 
and
Christopher Schommer-Pries 
for helpful discussions. 
N.\,C.\ is supported by the DFG Heisenberg Programme. 

\newpage

\section[Non-extended \texorpdfstring{$\boldsymbol{r}$}{r}-spin TQFTs]{Non-extended \texorpdfstring{$\boldsymbol{r}$}{r}-spin TQFTs}
\label{sec:open-closed}

In this section we review the classification of non-extended closed 
$r$-spin and framed TQFTs following \cite{SternSzegedy}, to which we refer for details. 
We recall the relevant categories of 2-dimensional bordisms as well as 
the notion of a closed $\Lambda_r$-Frobenius algebra,
and we state the main classification result: 2-dimensional $r$-spin and framed 
($r=0$) TQFTs are equivalent to closed $\Lambda_r$-Frobenius algebras in the 
target category.

\subsection[Framed and \texorpdfstring{$r$}{r}-spin TQFTs]{Framed and \texorpdfstring{$\boldsymbol{r}$}{r}-spin TQFTs}
\label{subsec:ClosedrSpinTQFTs}

By a \textsl{surface} we mean a 2-dimensional compact smooth manifold.
Let $G$ be a topological group, let  
\begin{equation}
\xi\colon G\lra \GL_2
\label{eq:def-xi-G}
\end{equation}
be a continuous group homomorphism, 
and recall that the frame bundle 
$F\Sigma\lra\Sigma$ of a surface $\Sigma$ is a principal $\GL_2$-bundle.
A \textsl{$G$-structure} (more precisely: a \textsl{tangential structure for $\xi\colon G\lra \GL_2$}) on $\Sigma$ is a
principal $G$-bundle $\pi \colon P \lra \Sigma$ together with a bundle map~$q$ intertwining the group actions via~$\xi$:
\begin{equation}
\begin{tikzpicture}[
baseline=(current bounding box.base),
descr/.style={fill=white,inner sep=3.5pt},
normal line/.style={->}
]
\matrix (m) [matrix of math nodes, row sep=2em, column sep=2em, text height=1.5ex, text depth=0.1ex] {%
	P&&F\Sigma
	\\
	&\Sigma&\\
};
\path[font=\footnotesize] (m-1-1) edge[->] node[below] {$ \pi\phantom{i} $} (m-2-2);
\path[font=\footnotesize] (m-1-1) edge[->] node[above] {$ q $} (m-1-3);
\path[font=\footnotesize] (m-1-3) edge[->] node[above] { } (m-2-2);
\end{tikzpicture}
\label{eq:def-tangential-str}
\end{equation}
A \textsl{map of surfaces with $G$-structure} is
a bundle map which is a local diffeomorphism of the underlying surfaces.
Such a map is called a \textsl{diffeomorphism} 
if its underlying map of surfaces is a diffeomorphism,
and an \textsl{isomorphism of $G$-structures} if the underlying map of surfaces is the identity.

We will consider the following tangential structures:
\begin{itemize}
	\item 
	A \textsl{framing} is a tangential structure for the trivial group: 
	\begin{equation}
	\star\lra \GL_2 \, .
	\label{eq:include-1}
	\end{equation}
	\item 
	An \textsl{orientation} is a tangential structure for the inclusion
	\begin{equation}
	\textrm{SO}_2 \simeq \GL_2^+\longhookrightarrow \GL_2\,, 
	\label{eq:include-GL2p}
	\end{equation}
	where $\GL_2^+$ is the subgroup of elements in $\GL_2$ with positive determinant. 
	\item 
	For $r\in\Z_{\geqslant 0}$, an \textsl{$r$-spin structure} is a tangential structure 
	\begin{equation}
	\begin{tikzpicture}[
	baseline=(current bounding box.base),
	descr/.style={fill=white,inner sep=3.5pt},
	normal line/.style={->}
	]
	\matrix (m) [matrix of math nodes, row sep=2em, column sep=2em, text height=1.5ex, text depth=0.1ex] {%
		\widetilde{\GL_2^+}^r&\GL_2^+&\GL_2 \, , 
		\\
	};
	\path[font=\footnotesize] (m-1-1) edge[->>] node[above] {$ p^r $} (m-1-2);
	\path[font=\footnotesize] (m-1-2) edge[right hook->] node[above] {} (m-1-3);
	\end{tikzpicture}
	\label{eq:def-r-spin-group}
	\end{equation}
	where 
	$p^r \colon \widetilde{\GL_2^+}^r \lra \GL_2^+$ is the $r$-fold covering for $r>0$, 
	while for $r=0$ it is the universal cover.
\end{itemize}

By a trivial $r$-spin structure on a surface $\Sigma$ 
we mean an $r$-spin structure isomorphic to the $r$-spin structure with
trivial bundles $P=\Sigma\times\widetilde{\GL_2^+}^r$, 
$F\Sigma=\Sigma\times \GL_2$ and trivial bundle map 
$q^{(+)}=\id_{\Sigma}\times p^r$ (positive orientation) or
$q^{(-)}=\id_{\Sigma}\times (T\circ p^r)$ (negative orientation),
where $T$ is composition with the matrix 
$(\begin{smallmatrix}+1&0\\0&-1\end{smallmatrix})\in \GL_2$.

\begin{remark}
	\label{rem:rSpinFramingsOrientations}
	A 1-spin structure is the same as an orientation, and a 2-spin structure is usually called a \textsl{spin structure}. 
	Moreover, we can identify framings with 0-spin structures by noting that the fibres of a $0$-spin bundle are contractible,
	see \cite[Prop.\,2.2]{RunkelSzegedyArf}. 
	This is consistent with the fact that
	for any $r\in\Z_{\geqslant 0}$, an $r$-spin structure is a $\Z_r$-bundle 
	over the oriented frame bundle. 
\end{remark}

\medskip 

Let $r\in\Z_{\geqslant 0}$. 
There is a symmetric monoidal category of 2-dimensional $r$-spin bordisms 
$\Bord_{2,1}^{r\text{-spin}}$ as follows. 
An object~$S$ is a 1-dimensional closed manifold~$s$ embedded in a cylinder 
$s\times (-1,1)$, together with an $r$-spin structure on the cylinder.
For an object $S$ we write $S^{(+)}:=s\times [0,1)$ and $S^{(-)}:=s\times (-1,0]$
with the restricted $r$-spin structures.
The morphisms of $\Bord_{2,1}^{r\text{-spin}}$ are diffeomorphism classes of $r$-spin bordisms: 
For $S,S'\in \Bord_{2,1}^{r\text{-spin}}$, an $r$-spin bordism 
$S\lra S'$ is a compact surface~$\Sigma$ 
with $r$-spin structure, together with a boundary parametrisation map
$S^{(+)}\sqcup {S'}^{(-)}\longhookrightarrow \Sigma$,
i.\,e.\ a map of $r$-spin surfaces that identifies the
boundary of~$\Sigma$ with the 1-dimensional embedded manifolds 
$s \times \{0\} \subset S$ and $s'\times \{0\} \subset S'$.
Finally, a diffeomorphism of $r$-spin bordisms is a diffeomorphism of $r$-spin surfaces 
compatible with the boundary parametrisations. 
We usually refer to a morphism in $\Bord_{2,1}^{r\text{-spin}}$ by a bordism that represents it. 

A particular class of $r$-spin bordisms are \textsl{deck transformation bordisms}. 
These are cylinders whose boundary parametrisations are given by
deck transformations of the $r$-spin bundle on the source or target object. 

The composition of morphisms in $\Bord_{2,1}^{r\text{-spin}}$ is given by glueing along boundary parametrisations, hence the unit morphisms are given by cylinders with trivial boundary parametrisations. 
Taking disjoint unions endows $\Bord_{2,1}^{r\text{-spin}}$ with its standard symmetric monoidal structure.
In light of Remark~\ref{rem:rSpinFramingsOrientations} we write 
\be 
\Bord_{2,1}^\mathrm{fr} = \Bord_{2,1}^{0\text{-spin}}
\, , \quad 
\Bord_{2,1}^\mathrm{or} = \Bord_{2,1}^{1\text{-spin}} \, . 
\ee 

\begin{definition}
	Let~${\mathds{C}}c$ be a symmetric monoidal category. 
	A \textsl{(closed) $r$-spin TQFT valued in~${\mathds{C}}c$} is a symmetric monoidal functor
	\begin{equation}
	\Zc\colon\Bord_{2,1}^{r\text{-spin}}\lra\Cc\,.
	\label{eq:next-r-spin-tqft}
	\end{equation}
\end{definition}

The case of 2-dimensional closed spin TQFTs ($r=2$) was first described and classified in \cite{MooreSegal}, 
including concrete examples in terms of Clifford algebras viewed as objects in the category of super vector spaces $\mathcal C = \Vect_{\mathds{C}}^{\Z_2}$. 
Spin TQFTs were further discussed from the perspective of extended TQFTs in \cite{GunninghamSpinHurwitz}, and spin state sum constructions were given in \cite{BarrettTavares, NovakRunkel}. 
TQFTs with $r$-spin structure for arbitrary~$r$ were introduced in \cite{Novak} and further studied in \cite{RunkelSzegedyArf}. 
The classification of general $r$-spin TQFTs appears in \cite{SternSzegedy}, in terms of the algebraic structures we review next.

\subsection[Classification in terms of closed \texorpdfstring{$\Lambda_r$}{Lambda\_r}-Frobenius algebras]{Classification in terms of closed \texorpdfstring{$\boldsymbol{\Lambda_r}$}{Lambda\_r}-Frobenius algebras}
\label{subsec:ClassificationThroughClosedLambdaFrob}

A \textsl{closed $\Lambda_r$-Frobenius algebra}~$C$ 
in a symmetric monoidal category~$\Cc$ consists of a collection of objects $C_a\in\Cc$ for all $a\in\Z_r$ as well as morphisms 
\begin{align}
\mu_{a,b}=
\begin{tikzpicture}[very thick,scale=0.53,color=blue!50!black, baseline=0.65cm]
\draw[-dot-] (2.5,0.75) .. controls +(0,1) and +(0,1) .. (3.5,0.75);
\draw (3,1.5) -- (3,2.25); 
\fill[color=black!80] (3,1.5) circle (0pt) node[below] (0up) {{\tiny $a,b$}};
\end{tikzpicture} 
&\colon C_a\otimes C_b\lra C_{a+b-1}\,,
&\eta_1= 
\begin{tikzpicture}[very thick,scale=0.4,color=blue!50!black, baseline=0]
\draw (-0.5,-0.5) node[Odot] (unit) {}; 
\draw (unit) -- (-0.5,0.5);
\end{tikzpicture} 
&\colon\one\lra C_1\,,\\
\Delta_{a,b}= 
\begin{tikzpicture}[very thick,scale=0.53,color=blue!50!black, rotate=180, baseline=-0.9cm]
\draw[-dot-] (2.5,0.75) .. controls +(0,1) and +(0,1) .. (3.5,0.75);
\draw (3,1.5) -- (3,2.25); 
\fill[color=black!80] (3,1.5) circle (0pt) node[above] (0up) {{\tiny $a,b$}};
\end{tikzpicture} 
&\colon C_{a+b+1}\lra C_{a}\otimes C_b\,,
&\varepsilon_{-1}= 
\begin{tikzpicture}[very thick,scale=0.4,color=blue!50!black, baseline=0cm, rotate=180]
\draw (-0.5,-0.5) node[Odot] (unit) {}; 
\draw (unit) -- (-0.5,0.5);
\end{tikzpicture} 
&\colon C_{-1}\lra \one
\label{eq:cLrFa-mor}
\end{align}
for all $a,b\in\Z_r$.
The \textsl{Nakayama automorphisms} of~$C$ are 
\begin{equation}
N_a:=
\tikzzbox{%
	\begin{tikzpicture}[very thick,scale=0.53,color=blue!50!black, baseline=-0.4cm]
	\draw[-dot-] (0,0) .. controls +(0,1) and +(0,1) .. (1,0);
	\draw (0.5,1.2) node[Odot] (unit) {}; 
	\draw (unit) -- (0.5,0.7);
	\draw (0,0) .. controls +(0,-0.5) and +(0,0.5) .. (-1,-1);
	\draw[-dot-] (0,-1) .. controls +(0,-1) and +(0,-1) .. (1,-1);
	\draw (0.5,-2.2) node[Odot] (unit2) {}; 
	\draw (unit2) -- (0.5,-1.7);
	\draw (0,-1) .. controls +(0,0.5) and +(0,-0.5) .. (-1,0);
	\draw (1,0) -- (1,-1);
	\draw (-1,0) -- (-1,2);
	\draw (-1,-1) -- (-1,-3);
	\fill[color=black!80] (0.55,-1.75) circle (0pt) node[right] (0up) {{\tiny $a,- a$}};
	\fill[color=black!80] (0.55,0.75) circle (0pt) node[right] (0up) {{\tiny $a,- a$}};
	\end{tikzpicture} 
}%
\colon C_a\lra C_a
\label{eq:Nakayama-Ca}
\end{equation}
for all $a\in\Z_r$. 
These data by definition satisfy the following conditions:
\begin{align}
\textrm{(co)associativity:} &&  
\begin{tikzpicture}[very thick,scale=0.53,color=blue!50!black, baseline=0.4cm]
\draw[-dot-] (3,0) .. controls +(0,1) and +(0,1) .. (2,0);
\draw[-dot-] (2.5,0.75) .. controls +(0,1) and +(0,1) .. (3.5,0.75);
\draw (3.5,0.75) -- (3.5,0); 
\draw (3,1.5) -- (3,2.25); 
\fill[color=black!80] (2.5,0.75) circle (0pt) node[below] (0up) {{\tiny $a,b$}};
\fill[color=black!80] (3,1.5) circle (0pt) node[left] (0up) {{\tiny $a\!+\!b\!-\!1,c$}};
\end{tikzpicture} 
=
\begin{tikzpicture}[very thick,scale=0.53,color=blue!50!black, baseline=0.4cm]
\draw[-dot-] (3,0) .. controls +(0,1) and +(0,1) .. (2,0);
\draw[-dot-] (2.5,0.75) .. controls +(0,1) and +(0,1) .. (1.5,0.75);
\draw (1.5,0.75) -- (1.5,0); 
\draw (2,1.5) -- (2,2.25); 
\fill[color=black!80] (2.5,0.75) circle (0pt) node[below] (0up) {{\tiny $b,c$}};
\fill[color=black!80] (2,1.5) circle (0pt) node[right] (0up) {{\tiny $a,b\!+\!c\!-\!1$}};
\end{tikzpicture} 
\, , 
\quad
\begin{tikzpicture}[very thick,scale=0.53,color=blue!50!black, baseline=-0.8cm, rotate=180]
\draw[-dot-] (3,0) .. controls +(0,1) and +(0,1) .. (2,0);
\draw[-dot-] (2.5,0.75) .. controls +(0,1) and +(0,1) .. (1.5,0.75);
\draw (1.5,0.75) -- (1.5,0); 
\draw (2,1.5) -- (2,2.25); 
\fill[color=black!80] (2.5,0.75) circle (0pt) node[above] (0up) {{\tiny $a,b$}};
\fill[color=black!80] (2,1.5) circle (0pt) node[left] (0up) {{\tiny $a\!+\!b\!+\!1,c$}};
\end{tikzpicture} 
=
\begin{tikzpicture}[very thick,scale=0.53,color=blue!50!black, baseline=-0.8cm, rotate=180]
\draw[-dot-] (3,0) .. controls +(0,1) and +(0,1) .. (2,0);
\draw[-dot-] (2.5,0.75) .. controls +(0,1) and +(0,1) .. (3.5,0.75);
\draw (3.5,0.75) -- (3.5,0); 
\draw (3,1.5) -- (3,2.25); 
\fill[color=black!80] (2.5,0.75) circle (0pt) node[above] (0up) {{\tiny $b,c$}};
\fill[color=black!80] (3,1.5) circle (0pt) node[right] (0up) {{\tiny $a,b\!+\!c\!+\!1$}};
\end{tikzpicture} 
\, ,&
\label{eq:clLrFa-co-associativity}
\\
\textrm{(co)unitality:} &&  
\begin{tikzpicture}[very thick,scale=0.4,color=blue!50!black, baseline]
\draw (-0.5,-0.5) node[Odot] (unit) {}; 
\fill (0,0.6) circle (5.0pt) node (meet) {};
\fill[color=black!80] (meet) circle (0pt) node[left] (0up) {{\tiny $1,a$}};
\draw (unit) .. controls +(0,0.5) and +(-0.5,-0.5) .. (0,0.6);
\draw (0,-1.5) -- (0,1.5); 
\end{tikzpicture} 
=
\begin{tikzpicture}[very thick,scale=0.4,color=blue!50!black, baseline]
\draw (0,-1.5) -- (0,1.5); 
\end{tikzpicture} 
=
\begin{tikzpicture}[very thick,scale=0.4,color=blue!50!black, baseline]
\draw (0.5,-0.5) node[Odot] (unit) {}; 
\fill (0,0.6) circle (5.0pt) node (meet) {};
\draw (unit) .. controls +(0,0.5) and +(0.5,-0.5) .. (0,0.6);
\draw (0,-1.5) -- (0,1.5); 
\fill[color=black!80] (meet) circle (0pt) node[right] (0up) {{\tiny $a,1$}};
\end{tikzpicture} 
\, , 
\quad
\quad
\begin{tikzpicture}[very thick,scale=0.4,color=blue!50!black, baseline=0, rotate=180]
\draw (0.5,-0.5) node[Odot] (unit) {}; 
\fill (0,0.6) circle (5.0pt) node (meet) {};
\draw (unit) .. controls +(0,0.5) and +(0.5,-0.5) .. (0,0.6);
\draw (0,-1.5) -- (0,1.5); 
\fill[color=black!80] (meet) circle (0pt) node[left] (0up) {{\tiny $-1,a$}};
\end{tikzpicture} 
=
\begin{tikzpicture}[very thick,scale=0.4,color=blue!50!black, baseline=0, rotate=180]
\draw (0,-1.5) -- (0,1.5); 
\end{tikzpicture} 
=
\begin{tikzpicture}[very thick,scale=0.4,color=blue!50!black, baseline=0cm, rotate=180]
\draw (-0.5,-0.5) node[Odot] (unit) {}; 
\fill (0,0.6) circle (5.0pt) node (meet) {};
\draw (unit) .. controls +(0,0.5) and +(-0.5,-0.5) .. (0,0.6);
\draw (0,-1.5) -- (0,1.5); 
\fill[color=black!80] (meet) circle (0pt) node[right] (0up) {{\tiny $a,-1$}};
\end{tikzpicture} 
\, ,&
\label{eq:clLrFa-co-unitality}
\\
\textrm{Frobenius relation:} &&  
\begin{tikzpicture}[very thick,scale=0.4,color=blue!50!black, baseline=0cm]
\draw[-dot-] (0,0) .. controls +(0,-1) and +(0,-1) .. (-1,0);
\draw[-dot-] (1,0) .. controls +(0,1) and +(0,1) .. (0,0);
\draw (-1,0) -- (-1,1.5); 
\draw (1,0) -- (1,-1.5); 
\draw (0.5,0.8) -- (0.5,1.5); 
\draw (-0.5,-0.8) -- (-0.5,-1.5); 
\fill[color=black!80] (-0.5,-0.8) circle (0pt) node[left] (0up) {{\tiny $c,a\!-\!c\!-\!1$}};
\fill[color=black!80] (0.5,0.8) circle (0pt) node[right] (0up) {{\tiny $d\!-\!b\!+\!1,b$}};
\end{tikzpicture}
=
\begin{tikzpicture}[very thick,scale=0.4,color=blue!50!black, baseline=0cm]
\draw[-dot-] (0,0) .. controls +(0,1) and +(0,1) .. (-1,0);
\draw[-dot-] (1,0) .. controls +(0,-1) and +(0,-1) .. (0,0);
\draw (-1,0) -- (-1,-1.5); 
\draw (1,0) -- (1,1.5); 
\draw (0.5,-0.8) -- (0.5,-1.5); 
\draw (-0.5,0.8) -- (-0.5,1.5); 
\fill[color=black!80] (0.5,-0.8) circle (0pt) node[right] (0up) {{\tiny $b\!-\!d\!-\!1,d$}};
\fill[color=black!80] (-0.5,0.8) circle (0pt) node[left] (0up) {{\tiny $a,c\!-\!a\!+\!1$}};
\end{tikzpicture}
\, ,&
\label{eq:clLrFa-Frobenius}
\\
\textrm{commutativity:} &&  
\tikzzbox{%
	\begin{tikzpicture}[very thick,scale=0.53,color=blue!50!black, baseline=0cm]
	\draw[-dot-] (-0.5,0) .. controls +(0,1) and +(0,1) .. (0.5,0);
	\draw (0,0.7) -- (0,1.5);
	\draw (-0.5,0) -- (-0.5,-0.75);
	\draw (0.5,0) -- (0.5,-0.75);
	\fill (-0.5,-0.25) circle (3.5pt) node[left] (mult1) {{\tiny $N_b^{1-a}$}};
	\fill[color=black!80] (0,0.7) circle (0pt) node[below] (0up) {{\tiny $b,\!a$}};
	\end{tikzpicture} 
}%
= 
\tikzzbox{%
	\begin{tikzpicture}[very thick,scale=0.53,color=blue!50!black, baseline=0cm]
	\draw[-dot-] (-0.5,0) .. controls +(0,1) and +(0,1) .. (0.5,0);
	\draw (0,0.7) -- (0,1.5);
	\draw (-0.5,0) .. controls +(0,-0.25) and +(0,0.25) .. (0.5,-0.75);
	\draw (0.5,0) .. controls +(0,-0.25) and +(0,0.25) .. (-0.5,-0.75);
	%
	\fill[color=black!80] (0,0.7) circle (0pt) node[below] (0up) {{\tiny $a,\!b$}};
	\end{tikzpicture} 
}%
= 
\tikzzbox{%
	\begin{tikzpicture}[very thick,scale=0.53,color=blue!50!black, baseline=0cm]
	\draw[-dot-] (-0.5,0) .. controls +(0,1) and +(0,1) .. (0.5,0);
	\draw (0,0.7) -- (0,1.5);
	\draw (-0.5,0) -- (-0.5,-0.75);
	\draw (0.5,0) -- (0.5,-0.75);
	\fill (0.5,-0.25) circle (3.5pt) node[right] (mult1) {{\tiny $N_a^{b-1}$}};
	\fill[color=black!80] (0,0.7) circle (0pt) node[below] (0up) {{\tiny $b,\!a$}};
	\end{tikzpicture} 
}%
\,,&
\label{eq:clLrFa-commutativity}
\\
\textrm{twist relations:} &&  
N_a^a = 1_{C_a}\,,
\quad
\quad
\tikzzbox{%
	\begin{tikzpicture}[very thick,scale=0.53,color=blue!50!black, baseline=0cm]
	\draw[-dot-] (-0.5,0) .. controls +(0,1) and +(0,1) .. (0.5,0);
	\draw (0,0.7) -- (0,1.5);
	\draw[-dot-] (-0.5,0) .. controls +(0,-1) and +(0,-1) .. (0.5,0);
	\draw (0,-1.2) node[Odot] (unit) {}; 
	\draw (unit) -- (0,-0.8);
	\fill (-0.5,0) circle (3.5pt) node[left] (mult1) {{\tiny $N_a^{b}$}};
	\fill[color=black!80] (0,0.8) circle (0pt) node[left] (0up) {{\tiny $a,\!-a$}};
	\fill[color=black!80] (0,-0.8) circle (0pt) node[left] (0up) {{\tiny $a,\!-a$}};
	\end{tikzpicture} 
}%
=
\tikzzbox{%
	\begin{tikzpicture}[very thick,scale=0.53,color=blue!50!black, baseline=0cm]
	\draw[-dot-] (-0.5,0) .. controls +(0,1) and +(0,1) .. (0.5,0);
	\draw (0,0.7) -- (0,1.5);
	\draw[-dot-] (-0.5,0) .. controls +(0,-1) and +(0,-1) .. (0.5,0);
	\draw (0,-1.2) node[Odot] (unit) {}; 
	\draw (unit) -- (0,-0.8);
	\fill (-0.5,0) circle (3.5pt) node[left] (mult1) {{\tiny $N_{a+b-1}^{b}$}};
	\fill[color=black!80] (0,0.8) circle (0pt) node[left] (0up) {{\tiny $a+b-1,\!-a-b+1$}};
	\fill[color=black!80] (0,-0.8) circle (0pt) node[left] (0up) {{\tiny $a+b-1,\!-a-b+1$}};
	\end{tikzpicture} 
}%
\,,&
\label{eq:clLrFa-twist}\\
\textrm{deck transformation relations:} &&  
N_a^r = 1_{C_a}
\,.&
\label{eq:clLrFa-deck-trf}
\end{align}

A \textsl{map of closed $\Lambda_r$-Frobenius algebras} $\varphi\colon C\lra D$ is 
a collection of morphisms $\varphi_a\colon C_a\lra D_a$ preserving the structure morphisms.
Analogously to the case of ordinary Frobenius algebras, maps of
closed $\Lambda_r$-Frobenius algebras are always isomorphisms.

\begin{example}
	\label{exa:ClosedLambdaRFrobeniusAlgebras}
	\begin{enumerate}[label={(\roman*)}]
		\item 
		One class of closed $\Lambda_r$-Frobenius algebras in a given symmetric monoidal category~$\Cc$ can be constructed from ordinary Frobenius algebras~$A$ in~$\mathcal C$ whose ordinary Nakayama automorphism~$\gamma_A$ satisfies $\gamma_A^r = 1_A$ (see Section~\ref{subsubsec:EquivariantCompletion} and~\eqref{eq:NakayamaA} below for details). 
		Indeed, as explained in \cite{RunkelSzegedyArf} and \cite[Sect.\,4.2]{SternSzegedy}, the construction of commutative Frobenius algebras as the centres of certain types of non-commutative Frobenius algebras in \cite[Sect.\,2.7]{lp0602047} is naturally the special case of $r=1$ of a construction of ``$\Z_r$-graded centre'' for any $r\in\Z_{\geqslant 0}$. 
		\item 
		\label{item:ClosedLambdaRAlgebraBordRspin}
		In the category $\Bord_{2,1}^{r\text{-spin}}$, $r$-spin circles, pair-of-pants, cups, and caps naturally assemble into a closed $\Lambda_r$-Frobenius algebra~$C$. 
		The precise presentation is given in \cite[Sect.\,5.1]{SternSzegedy} in terms of a combinatorial description of $r$-spin structures. 
		In particular, it follows from \cite[Eq.\,(5.2)]{SternSzegedy} that the Nakayama automorphisms of~$C$ are deck transformation bordisms.  
	\end{enumerate}
\end{example}

The closed $\Lambda_r$-Frobenius algebra~$C$ of Example~\ref{exa:ClosedLambdaRFrobeniusAlgebras}\ref{item:ClosedLambdaRAlgebraBordRspin} is not just any example. 
As proven in \cite[Thm.\,5.2.1]{SternSzegedy}, $\Bord_{2,1}^{r\text{-spin}}$ is generated as a symmetric monoidal category by the data of~$C$, subject to relations given by the defining properties \eqref{eq:clLrFa-co-associativity}--\eqref{eq:clLrFa-deck-trf}.
This implies: 

\begin{theorem}[{\cite[Cor.\,5.2.2]{SternSzegedy}}]
	There is an equivalence of symmetric monoidal groupoids between 
	the groupoid of $r$-spin TQFTs with target~$\Cc$ and the
	groupoid of closed $\Lambda_r$-Frobenius algebras in $\Cc$.
	\label{thm:closed-r-spin-TQFT-classification}
\end{theorem}

For this reason we will refer to the objects~$C_a$ of a closed $\Lambda_r$-Frobenius algebra in
any given symmetric monoidal category $\Cc$ (not necessarily equivalent to $\Bord_{2,1}^{r\text{-spin}}$)
as the \textsl{$a$-th circle spaces}. 
The $a$-th circle space in $\Bord_{2,1}^{r\text{-spin}}$ is simply the circle with ``framing number''~$a$, and we denote it $S_a^1$. 
Below in Sections~\ref{subsubsec:closed-Lr-Fa-B}, \ref{subsubsec:FrobenisAlgebraInBeq} and~\ref{subsec:LGmodels} we will use Theorem~\ref{thm:closed-r-spin-TQFT-classification} to construct examples of closed $r$-spin TQFTs beyond those mentioned in Section~\ref{subsec:ClosedrSpinTQFTs}.

\subsection{Computing invariants}
\label{subsec:ComputingInvariants}

The above classification theorem provides a way to compute invariants of 
$r$-spin surfaces from $r$-spin TQFTs in terms of the algebraic data 
of a closed $\Lambda_r$-Frobenius algebra. 
As the number of diffeomorphism classes of $r$-spin structures on a connected oriented 
surface of genus $g\geqslant 2$ is, if non-zero, either one ($r$ odd) or two ($r$ even), 
we are mainly interested in surfaces of genus~0 and~1; see e.\,g.\ 
\cite[Sect.\,3]{SzegedyPhD} and the references therein for a detailed account.

The sphere~$S^2$ admits an $r$-spin structure only if $r\in\{ 1,2 \}$, 
in which case it is unique up to isomorphism, 
hence the torus~$T^2$ is of most interest. 
Any torus with $r$-spin structure can be presented in terms of the 
closed $\Lambda_r$-Frobenius algebra in $\Bord_{2,1}^{r\text{-spin}}$ as  
\begin{equation}
T(a,b):= 
\tikzzbox{%
	\begin{tikzpicture}[very thick,scale=0.53,color=blue!50!black, baseline=-0.1cm]
	\draw[-dot-] (-0.5,0) .. controls +(0,1) and +(0,1) .. (0.5,0);
	\draw (0,1.2) node[Odot] (unit2) {}; 
	\draw (unit2) -- (0,0.7);
	\draw[-dot-] (-0.5,0) .. controls +(0,-1) and +(0,-1) .. (0.5,0);
	\draw (0,-1.2) node[Odot] (unit) {}; 
	\draw (unit) -- (0,-0.8);
	\fill (-0.5,0) circle (3.5pt) node[left] (mult1) {{\tiny $N_{-a}^{1-b}$}};
	\fill[color=black!80] (0,0.8) circle (0pt) node[left] (0up) {{\tiny $-a,\!a$}};
	\fill[color=black!80] (0,-0.8) circle (0pt) node[left] (0up) {{\tiny $-a,\!a$}};
	\end{tikzpicture} 
}%
\,\in \End_{\Bord_{2,1}^{r\text{-spin}}}(\varnothing) 
\label{eq:r-spin-torus-decomp}
\end{equation}
for some $a,b\in\Z_r$. 
Moreover, as shown in \cite[Prop.\,4.1.4]{SternSzegedy}, the $r$-spin torus $T(a,b)$ is diffeomorphic to the $r$-spin torus
$T(\mathrm{gcd}(a,b,r),0)$, and in fact diffeomorphism classes of $r$-spin tori are in bijection with divisors of~$r$. 
Hence we write 
\be 
T(d)=T(d,0) 
\ee 
for the class of $r$-spin tori 
corresponding to the divisor~$d$.

\begin{proposition}[{\cite[Prop.\,4.1.4]{SternSzegedy}}]
	The invariant of the $r$-spin torus $T(d)$ computed by a $\Cc$-valued closed $r$-spin TQFT $\Zc$ classified by a
	closed $\Lambda_r$-Frobenius algebra~$C$ is the quantum dimension of~$C_d$, 
	\begin{equation}
	\Zc\big(T(d)\big)
	=
	\dim(C_d)
	=
	\ev_{C_d} \circ \, b_{C_{-d}, C_d} \circ\coev_{C_d} \, , 
	\label{eq:torus-inv-dim-Cd}
	\end{equation}
	where~$b$ is the braiding of~$\Cc$. 
	\label{prop:torus-inv-dim-Cd}
\end{proposition}

\section[Fully extended \texorpdfstring{$\boldsymbol{r}$}{r}-spin TQFTs]{Fully extended \texorpdfstring{$\boldsymbol{r}$}{r}-spin TQFTs}
\label{sec:FullyExtendedRSpinTQFTs}

In this section we describe fully extended $r$-spin TQFTs and prove the corresponding cobordism hypothesis in the 2-categorical setting. 
In Section~\ref{subsec:DualisabilityInSM2Categories} we recall some aspects of symmetric monoidal 2-categories~$\B$, their Serre automorphisms, and we construct canonical closed $\Lambda_0$-Frobenius algebras. 
The short Section~\ref{subsec:bord-2-cat} describes the 2-category $\Bordrspin$ of $r$-spin bordisms. 
Then in Section~\ref{sec:bord-2-cat} we define the 
	 2-groupoid 
of fully extended $r$-spin TQFTs $\Bordrspin \lra \B$ and explain that it is equivalent to the 
	 2-groupoid 
	of $\Spin_2^r$-homotopy fixed points.

\subsection{Dualisability in symmetric monoidal 2-categories}
\label{subsec:DualisabilityInSM2Categories}

In this section we present our notational conventions for dualisability in symmetric monoidal 2-categories. 
Moreover, we construct a closed $\Lambda_0$-Frobenius algebra (in the sense of Section~\ref{subsec:ClassificationThroughClosedLambdaFrob}) for every fully dualisable object. 

For complete definitions we refer to \cite{benabou, spthesis, Pstragowski} and references therein; 
with an eye towards examples in Section~\ref{sec:Examples}, below we mostly use the same conventions as in \cite[Sect.\,2]{CMM}.

\subsubsection{Conventions for 2-categories}

By a 2-category we mean a (possibly non-strict) bicategory~$\B$ in the sense of \cite[App.\,A.1]{spthesis}. 
For objects $\alpha, \beta \in \B$, we denote the category of 1-morphisms $\alpha \lra \beta$ by $\B(\alpha,\beta)$; for 1-morphisms $X,Y \in \B(\alpha,\beta)$, we write $\Hom_\B(X,Y)$, or simply $\Hom(X,Y)$, for the set of 2-morphisms $X\lra Y$. 
Horizontal and vertical composition are denoted by~$\otimes$ and~$\circ$, respectively: 
\begin{align}
\otimes \colon \B(\beta, \gamma) \times \B(\alpha,\beta) & \lra \B(\alpha,\gamma) \nonumber
\\
(X',X) & \lmt X'\otimes X \, , 
\\ 
\circ\colon \Hom(Y,Z) \times \Hom(X,Y) & \lra \Hom(X,Z) \nonumber 
\\ 
(\psi,\varphi) & \lmt \psi \circ \varphi \, . 
\end{align}
We read string diagrams from right to left and from bottom to top. 
For instance, for 1-morphisms $X\in \B(\alpha,\beta)$, $X'\in\B(\beta,\gamma)$ and $V\in\B(\alpha,\gamma)$, a 2-morphism $\varphi \in \Hom(X'\otimes X,V)$ is represented by 
\be 
\tikzzbox{%
	\begin{tikzpicture}[thick,scale=3.5,color=blue!50!black, baseline=1.5cm
	]
	%
	\coordinate (rho) at (0.5, 1);
	\coordinate (kappa) at (0.4, 1);
	\coordinate (eta) at (0.9, 1);
	%
	\coordinate (Psi) at (0.5, 0.5);
	%
	\coordinate (nu) at (0.25, 0);
	\coordinate (mu) at (0.75, 0);
	\fill [red!30,opacity=0.7] (Psi) -- (rho) -- (1, 1) -- (1, 0) -- (mu);
	\fill [orange!50!white, opacity=0.8] (Psi) -- (mu) -- (nu);
	\fill [orange!30!white, opacity=0.8] (Psi) -- (nu) -- (0, 0) -- (0, 1) -- (rho);
	\draw[color=blue!50!black, ultra thick] (Psi) -- node[pos=0.6, color=blue!50!black, left] {\footnotesize $V$} (rho);  
	\draw[color=blue!50!black, ultra thick] (Psi) -- node[pos=0.6, color=blue!50!black, right] {\footnotesize $X\vphantom{X'}$} (mu);  
	\draw[color=blue!50!black, ultra thick] (Psi) -- node[pos=0.6, color=blue!50!black, left] {\footnotesize $X'$} (nu);  
	\fill[color=black!80] (Psi) circle (0.9pt) node[left] (0up) {\footnotesize $\varphi$};
	%
	\draw[
	color=gray, 
	opacity=0.9, 
	thick, 
	->
	] 
	(1, -0.1) -- node[pos=0.5, color=black, below, sloped] {\footnotesize $\fus\text{-composition}$} (0, -0.1); 
	\draw[
	color=gray, 
	opacity=0.9, 
	thick, 
	->
	] 
	(-0.1, 0) -- node[pos=0.5, color=black, above, sloped] {\footnotesize $\circ\text{-composition}$} (-0.1, 1); 
	\draw[line width=1] (0.8, 0.45) node[line width=0pt] (beta) {{\footnotesize $\alpha$}};
	\draw[line width=1] (0.5, 0.15) node[line width=0pt] (beta) {{\footnotesize $\beta$}};
	\draw[line width=1] (0.2, 0.45) node[line width=0pt] (beta) {{\footnotesize $\gamma$}};
	\end{tikzpicture}
}%
\ee 
though sometimes we will suppress object labels in string diagrams, as e.\,g.\ in~\eqref{eq:LeftAdjointFirstMention} below.

\subsubsection{Adjoints}
\label{subsubsec:Adjoints}

Let~$\B$ be a 2-category. 
A 1-morphism $X\in\B(\alpha,\beta)$ has a \textsl{left adjoint} if there exists a 1-morphism $\dX\in\B(\beta,\alpha)$ together with \textsl{adjunction 2-morphisms} 
\be 
\label{eq:LeftAdjointFirstMention}
\tikzzbox{%
	\begin{tikzpicture}[very thick,scale=1.0,color=blue!50!black, baseline=.4cm]
	\draw[line width=0pt] 
	(3,0) node[line width=0pt] (D) {{\small$X\vphantom{X^\dagger}$}}
	(2,0) node[line width=0pt] (s) {{\small$\dX$}}; 
	\draw[directed] (D) .. controls +(0,1) and +(0,1) .. (s);
	\end{tikzpicture}
}
= \ev_X \colon \dX \otimes X \lra 1_\alpha 
\, , \quad 
\tikzzbox{
	\begin{tikzpicture}[very thick,scale=1.0,color=blue!50!black, baseline=-.4cm,rotate=180]
	\draw[line width=0pt] 
	(3,0) node[line width=0pt] (D) {{\small$X\vphantom{X^\dagger}$}}
	(2,0) node[line width=0pt] (s) {{\small$\dX$}}; 
	\draw[redirected] (D) .. controls +(0,1) and +(0,1) .. (s);
	\end{tikzpicture}
}
= \coev_X \colon 1_\beta \lra X\otimes \dX 
\ee 
such that the Zorro moves 
\be 
\label{eq:ZorroSStringDiagrams}
\tikzzbox{
	\begin{tikzpicture}[very thick,scale=0.85,color=blue!50!black, baseline=0cm]
	\draw[line width=0] 
	(-1,1.25) node[line width=0pt] (A) {}
	(1,-1.25) node[line width=0pt] (A2) {}; 
	\draw[directed] (0,0) .. controls +(0,-1) and +(0,-1) .. (-1,0);
	\draw[directed] (1,0) .. controls +(0,1) and +(0,1) .. (0,0);
	\draw (-1,0) -- (A); 
	\draw (1,0) -- (A2); 
	\end{tikzpicture}
}
=
\tikzzbox{
	\begin{tikzpicture}[very thick,scale=0.85,color=blue!50!black, baseline=0cm]
	\draw[line width=0] 
	(0,1.25) node[line width=0pt] (A) {}
	(0,-1.25) node[line width=0pt] (A2) {}; 
	\draw
	(A2) -- (A); 
	\end{tikzpicture}
}
= 1_X 
\, , \quad
\tikzzbox{
	\begin{tikzpicture}[very thick,scale=0.85,color=blue!50!black, baseline=0cm]
	\draw[line width=0] 
	(1,1.25) node[line width=0pt] (A) {}
	(-1,-1.25) node[line width=0pt] (A2) {}; 
	\draw[directed] (0,0) .. controls +(0,1) and +(0,1) .. (-1,0);
	\draw[directed] (1,0) .. controls +(0,-1) and +(0,-1) .. (0,0);
	\draw (-1,0) -- (A2); 
	\draw (1,0) -- (A); 
	\end{tikzpicture}
}
=
\tikzzbox{
	\begin{tikzpicture}[very thick,scale=0.85,color=blue!50!black, baseline=0cm]
	\draw[line width=0] 
	(0,1.25) node[line width=0pt] (A) {}
	(0,-1.25) node[line width=0pt] (A2) {}; 
	\draw (A) -- (A2); 
	\end{tikzpicture}
}
= 1_{\dX} 
\ee 
are satisfied. 
Similarly, a \textsl{right adjoint} for~$X$ consists of $X^\dagger \in \B(\beta,\alpha)$ with 
\be 
\tikzzbox{
	\begin{tikzpicture}[very thick,scale=1.0,color=blue!50!black, baseline=.4cm]
	\draw[line width=0pt] 
	(3,0) node[line width=0pt] (D) {{\small$\Xd$}}
	(2,0) node[line width=0pt] (s) {{\small$X\vphantom{\dX}$}}; 
	\draw[redirected] (D) .. controls +(0,1) and +(0,1) .. (s);
	\end{tikzpicture}
}
= \tev_X \colon X\otimes \Xd \lra 1_\beta 
\, , \quad 
\tikzzbox{
	\begin{tikzpicture}[very thick,scale=1.0,color=blue!50!black, baseline=-.4cm,rotate=180]
	\draw[line width=0pt] 
	(3,0) node[line width=0pt] (D) {{\small$\Xd$}}
	(2,0) node[line width=0pt] (s) {{\small$X\vphantom{\dX}$}}; 
	\draw[directed] (D) .. controls +(0,1) and +(0,1) .. (s);
	\end{tikzpicture}
}
= \tcoev_X \colon 1_\alpha \lra \Xd \otimes X
\ee 
that satisfy analogous Zorro moves. 

If $X,Y \in \B(\alpha,\beta)$ have left and right adjoints (with chosen adjunction maps), then we write 
\be 
{}^\dagger \varphi := 
\tikzzbox{
	\begin{tikzpicture}[very thick,scale=0.85,color=blue!50!black, baseline=0cm]
	\draw[line width=0] 
	(1,1.25) node[line width=0pt] (A) {{\small $\dX$}}
	(-1,-1.25) node[line width=0pt] (A2) {{\small $\dY$}}; 
	\draw[directed] (0,0) .. controls +(0,1) and +(0,1) .. (-1,0);
	\draw[directed] (1,0) .. controls +(0,-1) and +(0,-1) .. (0,0);
	\fill (0,0) circle (3.0pt) node[right] {{\small $\varphi$}};
	\draw (-1,0) -- (A2); 
	\draw (1,0) -- (A); 
	\end{tikzpicture}
}
, \quad 
\varphi^\dagger := 
\tikzzbox{
	\begin{tikzpicture}[very thick,scale=0.85,color=blue!50!black, baseline=0cm]
	\draw[line width=0] 
	(-1,1.25) node[line width=0pt] (A) {{\small $X^\dagger$}}
	(1,-1.25) node[line width=0pt] (A2) {{\small $Y^\dagger$}}; 
	\fill (0,0) circle (3.0pt) node[left] {{\small $\varphi$}};
	\draw[redirected] (0,0) .. controls +(0,-1) and +(0,-1) .. (-1,0);
	\draw[redirected] (1,0) .. controls +(0,1) and +(0,1) .. (0,0);
	\draw (-1,0) -- (A); 
	\draw (1,0) -- (A2); 
	\end{tikzpicture}
}
\ee 
for the left and right adjoint of $\varphi \in \Hom(X,Y)$, respectively. 
We call~$\B$ \textsl{pivotal} if every 1-morphism~$X$ comes with chosen left and right adjunction data such that $\dX = \Xd$, ${}^\dagger\!\varphi = \varphi^\dagger$ for all 2-morphisms~$\varphi$, and  
\be 
\begin{tikzpicture}[very thick,scale=0.65,color=blue!50!black, baseline=-0.2cm, rotate=180]
\draw[line width=1pt] 
(2,-1.5) node[line width=0pt] (Y) {{\small $Y^\dagger$}}
(3,-1.5) node[line width=0pt] (X) {{\small $X^\dagger$}}
(-1,2) node[line width=0pt] (XY) {{\small $(Y\otimes X)^\dagger$}}; 
\draw[redirected] (1,0) .. controls +(0,1) and +(0,1) .. (2,0);
\draw[redirected] (0,0) .. controls +(0,2) and +(0,2) .. (3,0);
\draw[redirected] (-1,0) .. controls +(0,-1) and +(0,-1) .. (0.5,0);
\draw (2,0) -- (Y);
\draw (3,0) -- (X);
\draw[dotted] (0,0) -- (1,0);
\draw (-1,0) -- (XY);
\end{tikzpicture}
\! =  \!
\begin{tikzpicture}[very thick,scale=0.65,color=blue!50!black, baseline=-0.2cm, rotate=180]
\draw[line width=1pt] 
(-2,-1.5) node[line width=0pt] (X) {{\small $Y^\dagger$}}
(-1,-1.5) node[line width=0pt] (Y) {{\small $X^\dagger$}}
(2,2) node[line width=0pt] (XY) {{\small $(Y\otimes X)^\dagger$}}; 
\draw[redirected] (0,0) .. controls +(0,1) and +(0,1) .. (-1,0);
\draw[redirected] (1,0) .. controls +(0,2) and +(0,2) .. (-2,0);
\draw[redirected] (2,0) .. controls +(0,-1) and +(0,-1) .. (0.5,0);
\draw (-1,0) -- (Y);
\draw (-2,0) -- (X);
\draw[dotted] (0,0) -- (1,0);
\draw (2,0) -- (XY);
\end{tikzpicture}
\ee 
for all composable 1-morphisms $X,Y$.

\subsubsection{Symmetric monoidal structure}

Let~$\B$ be a 2-category. 
A monoidal structure on~$\B$ consists of a 2-functor 
\be 
\btimes \colon \B \times \B \lra \B 
\ee 
called monoidal product, a unit object 
\be 
\one\in\B \, , 
\ee 
a pseudonatural transformation $a \colon \btimes \circ \, (\btimes \times  \textrm{Id}_\B) \lra \btimes \circ (\textrm{Id}_\B \times \btimes)$ called associator, a weak inverse~$a^-$ for~$a$, as well as unitors, 2-unitors and a pentagonator (which we will usually suppress), subject to the coherence axioms in \cite[Sect.\,2.3]{spthesis}. 

Viewing a monoidal 2-category~$\B$ as a 3-category with a single object and using the strictification results of \cite{Gurskibook, Guthmann2019}, we can use the 3-dimensional graphical calculus of \cite{BMS, TrimbleSurfaceDiagrams}. 
For this we extend our diagrammatic conventions by reading 3-dimensional diagrams from front to back. 
For instance, for 1-morphisms $X \in \B(\varepsilon\btimes \delta, \alpha)$, $Y \in \B(\gamma\btimes \zeta, \varepsilon)$, $Z\in \B(\one, \zeta\btimes \delta)$, $X'\in \B(\beta,\alpha)$ and $Y'\in \B(\gamma,\beta)$, the diagram
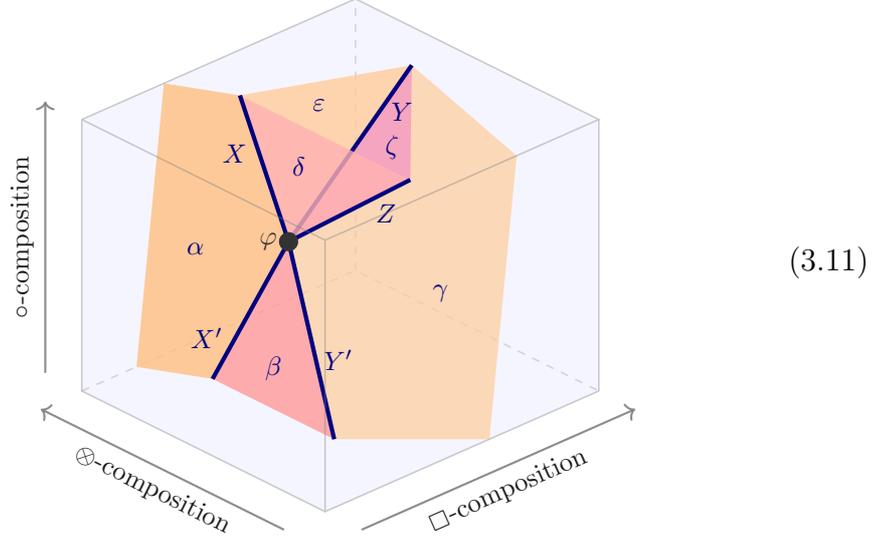
\be 
\tikzzbox{%
	\begin{tikzpicture}[thick,scale=4.0,color=blue!50!black, baseline=0.0cm, 
	style={x={(-0.9cm,-0.4cm)},y={(0.8cm,-0.4cm)},z={(0cm,0.9cm)}}]
	\draw[
	color=gray, 
	opacity=0.3, 
	semithick,
	dashed
	] 
	(1,0,0) -- (0,0,0) -- (0,1,0)
	(0,0,0) -- (0,0,1);
	\fill [blue!40,opacity=0.1] (0,0,0) -- (1,0,0) -- (1,1,0) -- (0,1,0) -- (0,1,1) -- (0,0,1) -- (1,0,1) -- (1,0,0);
	%
	\coordinate (rho) at (0.6, 0.2, 1);
	\coordinate (kappa) at (0.15, 0.4, 1);
	\coordinate (eta) at (0.6, 0.9, 1);
	%
	\coordinate (Psi) at (0.6, 0.4, 0.55);
	%
	\coordinate (nu) at (0.7, 0.2, 0);
	\coordinate (mu) at (0.7, 0.7, 0);
	%
	\fill [orange!40!white, opacity=0.8] (Psi) -- (kappa) -- (rho);
	\draw[line width=1] (0.4, 0.3, 0.92) node[line width=0pt] (beta) {{\footnotesize $\eps$}};
	\fill [orange!40,opacity=0.7] (Psi) -- (kappa) -- (0.3, 1, 1) -- (0.4, 1, 0) -- (mu);
	\draw[line width=1] (0.4, 0.8, 0.45) node[line width=0pt] (beta) {{\footnotesize $\gamma$}};
	\fill [magenta!50!white, opacity=0.8] (Psi) -- (kappa) -- (eta);
	\draw[line width=1] (0.4, 0.6, 0.9) node[line width=0pt] (beta) {{\footnotesize $\zeta$}};
	\draw[color=blue!50!black, ultra thick] (Psi) -- node[pos=0.86, color=blue!50!black, below] {\footnotesize $\;\;Y$} (kappa);  
	\fill [red!30!white, opacity=0.8] (Psi) -- (rho) -- (eta);
	\draw[line width=1] (0.5, 0.33, 0.75) node[line width=0pt] (beta) {{\footnotesize $\delta$}};
	\draw[color=blue!50!black, ultra thick] (Psi) -- node[pos=0.8, color=blue!50!black, below] {\footnotesize $Z$} (eta);  
	\fill [red!40!white, opacity=0.8] (Psi) -- (mu) -- (nu);
	\draw[line width=1] (0.7, 0.45, 0.15) node[line width=0pt] (beta) {{\footnotesize $\beta$}};
	\fill [orange!50!white, opacity=0.8] (Psi) -- (nu) -- (0.8, 0, 0) -- (0.7, 0., 1) -- (rho);
	\draw[line width=1] (0.7, 0.13, 0.45) node[line width=0pt] (beta) {{\footnotesize $\alpha$}};
	\draw[color=blue!50!black, ultra thick] (Psi) -- node[pos=0.6, color=blue!50!black, left] {\footnotesize $X$} (rho);  
	\draw[color=blue!50!black, ultra thick] (Psi) -- node[pos=0.6, color=blue!50!black, right] {\footnotesize $\!Y'$} (mu);  
	\draw[color=blue!50!black, ultra thick] (Psi) -- node[pos=0.7, color=blue!50!black, left] {\footnotesize ${}\;{}X'$} (nu);  
	%
	\fill[color=black!80] (Psi) circle (0.9pt) node[left] (0up) {\footnotesize $\varphi$};
	%
	\draw[
	color=gray, 
	opacity=0.4, 
	semithick
	] 
	(0,1,1) -- (0,1,0) -- (1,1,0) -- (1,1,1) -- (0,1,1) -- (0,0,1) -- (1,0,1) -- (1,0,0) -- (1,1,0)
	(1,0,1) -- (1,1,1);
	%
	\draw[
	color=gray, 
	opacity=0.9, 
	thick, 
	->
	] 
	(1, 1.15, 0) -- node[pos=0.5, color=black, below, sloped] {\footnotesize $\sta\text{-composition}$} (0, 1.15, 0); 
	\draw[
	color=gray, 
	opacity=0.9, 
	thick, 
	->
	] 
	(1.15, 1, 0) -- node[pos=0.5, color=black, below, sloped] {\footnotesize $\fus\text{-composition}$} (1.15, 0, 0); 
	\draw[
	color=gray, 
	opacity=0.9, 
	thick, 
	->
	] 
	(1, -0.15, 0) -- node[pos=0.5, color=black, above, sloped] {\footnotesize $\circ\text{-composition}$} (1, -0.15, 1); 
	\end{tikzpicture}
}
\ee 
represents a 2-morphism $\varphi \in \Hom(X'\otimes Y', X \otimes (Y\btimes 1_\delta) \otimes (1_\gamma \btimes Z))$, compare \cite[Sect.\,3.1.2]{CMS}. 

\medskip 

Let~$\B$ be a monoidal 2-category. 
Writing $\tau \colon \B \times \B \lra \B \times \B$ for the strict 2-functor that acts as $(\zeta,\xi) \lmt (\xi,\zeta)$ on objects, 1- and 2-morphisms, a symmetric braided structure on~$\B$ consists of a pseudonatural transformation 
\be 
b \colon \btimes \lra \btimes \circ \tau 
\ee 
called braiding, a weak inverse~$b^-$ for~$b$, and an invertible modification $\sigma \colon 1_{\btimes} \lra b^- \circ b$, as well as two further invertible modifications between compositions of $a,a^-, b,b^-$, subject to the coherence axioms of \cite[Sect.\,2.3]{spthesis}. 

The braiding~$b$ consists of 1-morphism components
\be 
\tikzzbox{%
	\begin{tikzpicture}[thick,scale=1.0,color=black, baseline=2cm]
	\coordinate (p1) at (0,0);
	\coordinate (p2) at (2,-0.5);
	\coordinate (p3) at (2,0.5);
	\coordinate (p4) at (4,0);
	\coordinate (u1) at (0,3);
	\coordinate (u2) at (2,2.5);
	\coordinate (u3) at (2,3.5);
	\coordinate (u4) at (4,3);
	\coordinate (ld) at (-1,1);
	\coordinate (lu) at (-1,4);
	\coordinate (rd) at (5,1);
	\coordinate (ru) at (5,4);
	\coordinate (rd2) at (5.5,0);
	\coordinate (ru2) at (5.5,3);
	\coordinate (ld2) at (-1.5,0);
	\coordinate (lu2) at (-1.5,3);
	%
	\fill [orange!25!white, opacity=0.8] (p3) -- (ld) -- (lu) -- (u3); 
	\fill [orange!25!white, opacity=0.8] (p3) -- (rd) -- (ru) -- (u3); 
	\draw[thin] (lu) --  (ld); 
	\draw[thin] (ru) --  (rd); 
	%
	\draw[very thick, red!80!black] (p3) -- (ld);
	\draw[very thick, red!80!black] (u3) -- (lu); 
	\draw[very thick, red!80!black] (p3) -- (rd);
	\draw[very thick, red!80!black] (ru) -- (u3); 
	%
	\fill [orange!30!white, opacity=0.8] (p3) -- (ld2) -- (lu2) -- (u3); 
	\fill [orange!30!white, opacity=0.8] (p3) -- (rd2) -- (ru2) -- (u3); 
	%
	\draw[ultra thick] (p3) --  (u3); 
	\fill (2.5,2) circle (0pt) node {{\small $b_{\alpha,\alpha'}$}};
	%
	\draw[very thick, red!80!black] (ru2) -- (u3); 
	\draw[very thick, red!80!black] (rd2) -- (p3); 
	\draw[very thick, red!80!black] (lu2) -- (u3); 
	\draw[very thick, red!80!black] (ld2) -- (p3); 
	\draw[thin] (ru2) --  (rd2); 
	\draw[thin] (lu2) --  (ld2); 
	%
	%
	\fill (-0.5,3.6) circle (0pt) node {{\small $\alpha'$}};
	\fill (4.5,3.6) circle (0pt) node {{\small $\alpha\vphantom{\alpha'}$}};
	\fill (-1,0.5) circle (0pt) node {{\small $\alpha\vphantom{\alpha'}$}};
	\fill (5,0.5) circle (0pt) node {{\small $\alpha'$}};
	%
	\end{tikzpicture}
} 
\;\; \widehat{=} \;\;
b_{\alpha,\alpha'} 
\colon \alpha \btimes \alpha' \lra \alpha' \btimes \alpha 
\ee 
for all $\alpha,\alpha'\in\B$, and of 2-morphism components 
\be 
b_{X,Y}  
\colon 
(Y\btimes X) \otimes b_{\alpha,\beta}\lra b_{\alpha',\beta'} \otimes (X \btimes Y) 
\ee 
for all $X\in \B(\alpha,\alpha')$ and $Y\in\B(\beta,\beta')$. 
Graphically, the 2-morphism components are depicted as 
\be 
\tikzzbox{%
	\begin{tikzpicture}[thick,scale=1.0,color=black, baseline=2cm]
	\coordinate (p3) at (2,0.5);
	\coordinate (u3) at (2,3.5);
	\coordinate (ld) at (-1,1);
	\coordinate (lu) at (-1,4);
	\coordinate (rd) at (5,1);
	\coordinate (ru) at (5,4);
	\coordinate (rd2) at (5.5,0);
	\coordinate (ru2) at (5.5,3);
	\coordinate (ld2) at (-1.5,0);
	\coordinate (lu2) at (-1.5,3);
	%
	\fill [orange!25!white, opacity=0.8] (p3) -- (ld) -- (lu) -- (u3); 
	\fill [orange!25!white, opacity=0.8] (p3) -- (rd) -- (ru) -- (u3); 
	\draw[thin] (lu) --  (ld); 
	\draw[thin] (ru) --  (rd); 
	%
	\coordinate (xlb) at (0,0.85);
	\coordinate (xrb) at (3.5,3.75);
	\draw[ultra thick, color=blue!50!black] (xlb) .. controls +(0,0.25) and +(-0.2,-0.5) .. (2,2); 
	\draw[ultra thick, color=blue!50!black] (xrb) .. controls +(0,-0.75) and +(0.2,0.5) .. (2,2); 
	\fill[color=blue!50!black] (0.4,1.4) circle (0pt) node {{\small $Y$}};
	\fill[color=blue!50!black] (3.7,3.5) circle (0pt) node {{\small $X$}};
	%
	\draw[very thick, red!80!black] (p3) -- (ld);
	\draw[very thick, red!80!black] (u3) -- (lu); 
	\draw[very thick, red!80!black] (p3) -- (rd);
	\draw[very thick, red!80!black] (ru) -- (u3); 
	%
	\fill [orange!30!white, opacity=0.8] (p3) -- (ld2) -- (lu2) -- (u3); 
	\fill [orange!30!white, opacity=0.8] (p3) -- (rd2) -- (ru2) -- (u3); 
	%
	\draw[ultra thick] (p3) --  (u3); 
	\fill (2.5,1) circle (0pt) node {{\small $b_{\alpha,\beta}$}};
	\fill (1.5,3) circle (0pt) node {{\small $b_{\alpha',\beta'}$}};\\
	%
	\coordinate (xlf) at (0.75,0.32);
	\coordinate (xrf) at (4,3.2);
	\draw[ultra thick, color=blue!50!black] (xlf) .. controls +(0,0.5) and +(-0.2,-0.5) .. (2,2); 
	\draw[ultra thick, color=blue!50!black] (xrf) .. controls +(0,-0.75) and +(0.2,0.3) .. (2,2); 
	\fill[color=blue!50!black] (0.55,0.5) circle (0pt) node {{\small $X$}};
	\fill[color=blue!50!black] (4.2,2.9) circle (0pt) node {{\small $Y$}};
	\fill[color=blue!50!black] (2,2) circle (2.5pt) node[left] {{\small $b_{X,Y}$}};
	%
	\draw[very thick, red!80!black] (ru2) -- (u3); 
	\draw[very thick, red!80!black] (rd2) -- (p3); 
	\draw[very thick, red!80!black] (lu2) -- (u3); 
	\draw[very thick, red!80!black] (ld2) -- (p3); 
	\draw[thin] (ru2) --  (rd2); 
	\draw[thin] (lu2) --  (ld2); 
	%
	%
	\fill (-0.5,3.5) circle (0pt) node {{\small $\beta'$}};
	\fill (4.5,3.5) circle (0pt) node {{\small $\alpha\vphantom{\beta'}$}};
	\fill (-1,0.5) circle (0pt) node {{\small $\alpha'\vphantom{\beta}$}};
	\fill (5,0.5) circle (0pt) node {{\small $\beta$}};
	%
	\end{tikzpicture}
}  
\, . 
\ee

\subsubsection{Duality and Serre automorphism}
\label{subsubsec:DualeAndSerre}

Let~$\B$ be a symmetric monoidal 2-category. 
An object $\alpha \in \B$ has a \textsl{dual} if there exists an object $\alpha^\dual$ together with \textsl{adjunction 1-morphisms} 
\begin{align}
\tikzzbox{%
	\begin{tikzpicture}[thick,scale=1.0,color=black, baseline=1.5cm]
	\coordinate (p1) at (0,0);
	\coordinate (p2) at (2,-0.5);
	\coordinate (p3) at (2.5,0.5);
	\coordinate (u1) at (0,3);
	\coordinate (u2) at (2,2.5);
	\coordinate (u3) at (2.5,3.5);
	%
	\fill [orange!20!white, opacity=0.8] 
	(p1) .. controls +(0,0.25) and +(-1,0) ..  (p3)
	-- (p3) --  (u3)
	-- (u3) .. controls +(-1,0) and +(0,0.25) ..  (u1)
	;
	%
	\draw[very thick, red!80!black] (p1) .. controls +(0,0.25) and +(-1,0) ..  (p3); 
	%
	\fill [orange!30!white, opacity=0.8] 
	(p1) .. controls +(0,-0.25) and +(-1,0) ..  (p2)
	-- (p2) --  (u2)
	-- (u2) .. controls +(-1,0) and +(0,-0.25) ..  (u1)
	;
	\draw[thin] (p1) --  (u1); 
	\draw[thin] (p2) --  (u2); 
	\draw[thin] (p3) --  (u3); 
	%
	\draw[very thick, red!80!black] (p1) .. controls +(0,-0.25) and +(-1,0) ..  (p2); 
	\draw[very thick, red!80!black] (u1) .. controls +(0,-0.25) and +(-1,0) ..  (u2); 
	\draw[very thick, red!80!black] (u1) .. controls +(0,0.25) and +(-1,0) ..  (u3); 
	%
	\fill[color=blue!50!black] (1.5,0) circle (0pt) node {{\small $\alpha^\dual$}};
	\fill[color=blue!50!black] (2,3) circle (0pt) node {{\small $\alpha\vphantom{\alpha^\dual}$}};
	%
	\end{tikzpicture}
}%
\;
& 
\;
\widehat{=} \; \tev_\alpha \colon \alpha \btimes \alpha^\dual \lra \one 
\, 
\\ 
\tikzzbox{%
	\begin{tikzpicture}[thick,scale=1.0,color=black, baseline=1.5cm]
	\coordinate (p1) at (0,0);
	\coordinate (p2) at (-2,-0.5);
	\coordinate (p3) at (-2.5,0.5);
	\coordinate (u1) at (0,3);
	\coordinate (u2) at (-2,2.5);
	\coordinate (u3) at (-2.5,3.5);
	%
	\fill [orange!20!white, opacity=0.8] 
	(p1) .. controls +(0,0.25) and +(1,0) ..  (p3)
	-- (p3) --  (u3)
	-- (u3) .. controls +(1,0) and +(0,0.25) ..  (u1)
	;
	%
	\draw[very thick, red!80!black] (p1) .. controls +(0,0.25) and +(1,0) ..  (p3); 
	%
	\fill [orange!30!white, opacity=0.8] 
	(p1) .. controls +(0,-0.25) and +(1,0) ..  (p2)
	-- (p2) --  (u2)
	-- (u2) .. controls +(1,0) and +(0,-0.25) ..  (u1)
	;
	\draw[thin] (p1) --  (u1); 
	\draw[thin] (p2) --  (u2); 
	\draw[thin] (p3) --  (u3); 
	%
	\draw[very thick, red!80!black] (p1) .. controls +(0,-0.25) and +(1,0) ..  (p2); 
	\draw[very thick, red!80!black] (u1) .. controls +(0,-0.25) and +(1,0) ..  (u2); 
	\draw[very thick, red!80!black] (u1) .. controls +(0,0.25) and +(1,0) ..  (u3); 
	%
	\fill[color=blue!50!black] (-1.5,0) circle (0pt) node {{\small $\alpha\vphantom{\alpha^\dual}$}};
	\fill[color=blue!50!black] (-2,3) circle (0pt) node {{\small $\alpha^\dual$}};
	%
	\end{tikzpicture}
} 
\;
& 
\;
\widehat{=} \; \tcoev_\alpha \colon \one \lra \alpha^\dual \btimes \alpha 
\end{align}
and \textsl{cusp 2-isomorphisms} 
\begin{align}
\label{eq:cuspl}
\tikzzbox{%
	\begin{tikzpicture}[thick,scale=0.75,color=black, baseline=1.8cm, xscale=-1]
	\coordinate (p1) at (0,0);
	\coordinate (p2) at (1.5,-0.5);
	\coordinate (p2s) at (4,-0.5);
	\coordinate (p3) at (1.5,0.5);
	\coordinate (p4) at (3,1);
	\coordinate (p5) at (1.5,1.5);
	\coordinate (p6) at (-1,1.5);
	\coordinate (u1) at (0,2.5);
	\coordinate (u2) at (1.5,2);
	\coordinate (u2s) at (4,3.5);
	\coordinate (u3) at (1.5,3);
	\coordinate (u4) at (3,3.5);
	\coordinate (u5) at (1.5,4);
	\coordinate (u6) at (-1,3.5);
	%
	\fill [orange!90!white, opacity=1] 
	(p4) .. controls +(0,0.25) and +(1,0) ..  (p5) 
	-- (p5) -- (0.79,1.5) -- (1.5,2.5) 
	-- (1.5,2.5) .. controls +(0,0) and +(0,0.5) ..  (p4);
	%
	\fill [orange!30!white, opacity=0.8] 
	(1.5,2.5) .. controls +(0,0) and +(0,0.5) ..  (p4)
	-- 
	(p4) .. controls +(0,-0.25) and +(1,0) ..  (p3)
	-- (p3) .. controls +(-1,0) and +(0,0.25) ..  (p1)
	-- 
	(p1) .. controls +(0,0.5) and +(0,0) ..  (1.5,2.5);
	%
	\draw[very thick, red!80!black] (p4) .. controls +(0,0.25) and +(1,0) ..  (p5) -- (p6); 
	\draw[very thick, red!80!black] (p1) .. controls +(0,0.25) and +(-1,0) ..  (p3)
	-- (p3) .. controls +(1,0) and +(0,-0.25) ..  (p4); 
	%
	\draw[thin] (p2s) --  (u2s); 
	\draw[thin] (p6) --  (u6); 
	\fill [orange!20!white, opacity=0.8] 
	(p2s) -- (p2)
	-- (p2) .. controls +(-1,0) and +(0,-0.25) ..  (p1) 
	-- (p1) .. controls +(0,0.5) and +(0,0) ..  (1.5,2.5)
	-- (1.5,2.5) .. controls +(0,0) and +(0,0.5) .. (p4)
	-- (p4) .. controls +(0,0.25) and +(-1,0) ..  (p5) 
	-- (p5) -- (p6) -- (u6) -- (u2s)
	;
	%
	\draw[thin, dotted] (p4) .. controls +(0,0.5) and +(0,0) ..  (1.5,2.5); 
	\draw[thin] (p1) .. controls +(0,0.5) and +(0,0) ..  (1.5,2.5); 
	%
	%
	\draw[very thick, red!80!black] (u6) -- (u2s); 
	\draw[very thick, red!80!black] (p1) .. controls +(0,-0.25) and +(-1,0) ..  (p2) -- (p2s); 
	\draw[very thick, red!80!black] (0.79,1.5) -- (p6); 
	%
	\fill[color=blue!50!black] (3.5,0.2) circle (0pt) node {{\small $\alpha$}};
	\fill (1.5,2.5) circle (2.5pt) node[right] {{\small $c_{\textrm{l}}$}};
	%
	\end{tikzpicture}
}     
& = c_{\textrm{l}} \colon \big(\tev_\alpha \btimes 1_\alpha\big) \otimes \big(1_\alpha \btimes \tcoev_\alpha\big) \lra 1_\alpha \, , 
\\ 
\tikzzbox{%
	\begin{tikzpicture}[thick,scale=0.75,color=black, baseline=1.8cm]
	\coordinate (p1) at (0,0);
	\coordinate (p2) at (1.5,-0.5);
	\coordinate (p2s) at (4,-0.5);
	\coordinate (p3) at (1.5,0.5);
	\coordinate (p4) at (3,1);
	\coordinate (p5) at (1.5,1.5);
	\coordinate (p6) at (-1,1.5);
	\coordinate (u1) at (0,2.5);
	\coordinate (u2) at (1.5,2);
	\coordinate (u2s) at (4,3.5);
	\coordinate (u3) at (1.5,3);
	\coordinate (u4) at (3,3.5);
	\coordinate (u5) at (1.5,4);
	\coordinate (u6) at (-1,3.5);
	%
	\fill [orange!90!white, opacity=1] 
	(p4) .. controls +(0,0.25) and +(1,0) ..  (p5) 
	-- (p5) -- (0.79,1.5) -- (1.5,2.5) 
	-- (1.5,2.5) .. controls +(0,0) and +(0,0.5) ..  (p4);
	%
	\fill [orange!30!white, opacity=0.8] 
	(1.5,2.5) .. controls +(0,0) and +(0,0.5) ..  (p4)
	-- 
	(p4) .. controls +(0,-0.25) and +(1,0) ..  (p3)
	-- (p3) .. controls +(-1,0) and +(0,0.25) ..  (p1)
	-- 
	(p1) .. controls +(0,0.5) and +(0,0) ..  (1.5,2.5);
	%
	\draw[very thick, red!80!black] (p4) .. controls +(0,0.25) and +(1,0) ..  (p5) -- (p6); 
	\draw[very thick, red!80!black] (p1) .. controls +(0,0.25) and +(-1,0) ..  (p3)
	-- (p3) .. controls +(1,0) and +(0,-0.25) ..  (p4); 
	%
	\draw[thin] (p2s) --  (u2s); 
	\draw[thin] (p6) --  (u6); 
	\fill [orange!20!white, opacity=0.8] 
	(p2s) -- (p2)
	-- (p2) .. controls +(-1,0) and +(0,-0.25) ..  (p1) 
	-- (p1) .. controls +(0,0.5) and +(0,0) ..  (1.5,2.5)
	-- (1.5,2.5) .. controls +(0,0) and +(0,0.5) .. (p4)
	-- (p4) .. controls +(0,0.25) and +(-1,0) ..  (p5) 
	-- (p5) -- (p6) -- (u6) -- (u2s)
	;
	%
	\draw[thin, dotted] (p4) .. controls +(0,0.5) and +(0,0) ..  (1.5,2.5); 
	\draw[thin] (p1) .. controls +(0,0.5) and +(0,0) ..  (1.5,2.5); 
	%
	%
	\draw[very thick, red!80!black] (u6) -- (u2s); 
	\draw[very thick, red!80!black] (p1) .. controls +(0,-0.25) and +(-1,0) ..  (p2) -- (p2s); 
	\draw[very thick, red!80!black] (0.79,1.5) -- (p6); 
	%
	\fill[color=blue!50!black] (3.5,0.2) circle (0pt) node {{\small $\alpha^\dual$}};
	\fill (1.5,2.5) circle (2.5pt) node[right] {{\small $c_{\textrm{r}}$}};
	%
	\end{tikzpicture}
}   
& = c_{\textrm{r}} \colon \big( 1_{\alpha^\dual} \btimes \tev_\alpha \big) \otimes \big( \tcoev_\alpha \btimes 1_{\alpha^\dual} \big) \lra 1_{\alpha^\dual} \, . 
\end{align}
More precisely, these data witness~$\alpha^\dual$ as the \textsl{right dual} of~$\alpha$. 
Using the symmetric braiding of~$\B$, the object~$\alpha^\dual$ is also the left dual of~$\alpha$, with adjunction maps 
\be 
\ev_{\alpha} = \tev_\alpha \otimes b_{\alpha^\dual,\alpha} 
\, , \quad 
\coev_\alpha = b_{\alpha^\dual,\alpha} \otimes \tcoev_\alpha \, . 
\ee 

If $\alpha, \beta \in \B$ have duals $\alpha^\dual, \beta^\dual$ with chosen adjunction 1-morphisms $\tev_\alpha$, $\tcoev_\alpha$, $\tev_\beta$, $\tcoev_\beta$, the associated dual $X^\dual$ of a 1-morphism $X\in\B(\alpha,\beta)$ is 
\be 
\tikzzbox{%
	\begin{tikzpicture}[thick,scale=1.0,color=black, baseline=2cm]
	\coordinate (p1) at (0,0);
	\coordinate (p2) at (1.5,-0.5);
	\coordinate (p2s) at (5,-0.5);
	\coordinate (p3) at (1.5,0.5);
	\coordinate (p4) at (3,1);
	\coordinate (p5) at (1.5,1.5);
	\coordinate (p6) at (-2,1.5);
	\coordinate (u1) at (0,3);
	\coordinate (u2) at (1.5,2.5);
	\coordinate (u2s) at (5,2.5);
	\coordinate (u3) at (1.5,3.5);
	\coordinate (u4) at (3,4);
	\coordinate (u5) at (1.5,4.5);
	\coordinate (u6) at (-2,4.5);
	%
	\fill [orange!20!white, opacity=0.8] 
	(p4) .. controls +(0,0.25) and +(1,0) ..  (p5) -- (p6)
	-- (u6) -- (u5)
	-- (u5) .. controls +(1,0) and +(0,0.25) .. (u4);
	%
	\fill [orange!25!white, opacity=0.8] 
	(p1) .. controls +(0,0.25) and +(-1,0) ..  (p3)
	-- (p3) .. controls +(1,0) and +(0,-0.25) ..  (p4)
	-- (p4) --  (u4)
	-- (u4) .. controls +(0,-0.25) and +(1,0) ..  (u3)
	-- (u3) .. controls +(-1,0) and +(0,0.25) ..  (u1)
	;
	%
	\draw[thin] (p1) --  (u1); 
	\draw[thin] (p2s) --  (u2s); 
	\draw[thin] (p4) --  (u4); 
	\draw[thin] (p6) --  (u6); 
	%
	%
	\draw[very thick, red!80!black] (p1) .. controls +(0,-0.25) and +(-1,0) ..  (p2) -- (p2s); 
	\draw[very thick, red!80!black] (p4) .. controls +(0,0.25) and +(1,0) ..  (p5) -- (p6); 
	%
	\draw[ultra thick, blue!50!black] (p3) --  (u3); 
	\fill[color=blue!50!black] (1.75,3) circle (0pt) node {{\small $X$}};
	%
	\draw[very thick, red!80!black] (p1) .. controls +(0,0.25) and +(-1,0) ..  (p3)
	-- (p3) .. controls +(1,0) and +(0,-0.25) ..  (p4); 
	%
	\fill [orange!30!white, opacity=0.8] 
	(p1) .. controls +(0,-0.25) and +(-1,0) ..  (p2)
	-- (p2) -- (p2s)
	-- (p2s) --  (u2s)
	-- (u2s) -- (u2)
	-- (u2) .. controls +(-1,0) and +(0,-0.25) ..  (u1)
	;
	%
	\draw[very thick, red!80!black] (p1) .. controls +(0,-0.25) and +(-1,0) ..  (p2) -- (p2s); 
	\draw[very thick, red!80!black] (u1) .. controls +(0,-0.25) and +(-1,0) ..  (u2) -- (u2s); 
	\draw[very thick, red!80!black] (u1) .. controls +(0,0.25) and +(-1,0) ..  (u3)
	-- (u3) .. controls +(1,0) and +(0,-0.25) ..  (u4); 
	\draw[very thick, red!80!black] (u4) .. controls +(0,0.25) and +(1,0) ..  (u5) -- (u6); 
	%
	\fill[color=blue!50!black] (2.5,3) circle (0pt) node {{\small $\alpha$}};
	\fill[color=blue!50!black] (1,3) circle (0pt) node {{\small $\beta$}};
	\fill[color=blue!50!black] (-1.5,2) circle (0pt) node {{\small $\alpha^\dual$}};
	\fill[color=blue!50!black] (4.5,0) circle (0pt) node {{\small $\beta^\dual$}};
	%
	\end{tikzpicture}
}    
\;\;\widehat{=}\;\; 
X^\dual \in \B(\beta^\dual,\alpha^\dual) \, . 
\ee 
For another 1-morphism $Y\in\B(\alpha,\beta)$ and a 2-morphism $\varphi \in \Hom(X,Y)$, its dual $\varphi^\dual$ is 
\be 
\tikzzbox{%
	\begin{tikzpicture}[thick,scale=1.0,color=black, baseline=2cm]
	\coordinate (p1) at (0,0);
	\coordinate (p2) at (1.5,-0.5);
	\coordinate (p2s) at (5,-0.5);
	\coordinate (p3) at (1.5,0.5);
	\coordinate (p4) at (3,1);
	\coordinate (p5) at (1.5,1.5);
	\coordinate (p6) at (-2,1.5);
	\coordinate (u1) at (0,3);
	\coordinate (u2) at (1.5,2.5);
	\coordinate (u2s) at (5,2.5);
	\coordinate (u3) at (1.5,3.5);
	\coordinate (u4) at (3,4);
	\coordinate (u5) at (1.5,4.5);
	\coordinate (u6) at (-2,4.5);
	%
	\fill [orange!20!white, opacity=0.8] 
	(p4) .. controls +(0,0.25) and +(1,0) ..  (p5) -- (p6)
	-- (u6) -- (u5)
	-- (u5) .. controls +(1,0) and +(0,0.25) .. (u4);
	%
	\fill [orange!25!white, opacity=0.8] 
	(p1) .. controls +(0,0.25) and +(-1,0) ..  (p3)
	-- (p3) .. controls +(1,0) and +(0,-0.25) ..  (p4)
	-- (p4) --  (u4)
	-- (u4) .. controls +(0,-0.25) and +(1,0) ..  (u3)
	-- (u3) .. controls +(-1,0) and +(0,0.25) ..  (u1)
	;
	%
	\draw[thin] (p1) --  (u1); 
	\draw[thin] (p2s) --  (u2s); 
	\draw[thin] (p4) --  (u4); 
	\draw[thin] (p6) --  (u6); 
	%
	%
	\draw[very thick, red!80!black] (p1) .. controls +(0,-0.25) and +(-1,0) ..  (p2) -- (p2s); 
	\draw[very thick, red!80!black] (p4) .. controls +(0,0.25) and +(1,0) ..  (p5) -- (p6); 
	%
	\draw[ultra thick, blue!50!black] (p3) --  (u3); 
	\fill[color=blue!50!black] (1.75,3) circle (0pt) node {{\small $Y$}};
	\fill[color=blue!50!black] (1.2,1) circle (0pt) node {{\small $X$}};
	\fill[color=blue!50!black] (1.5,2) circle (2.5pt) node[right] {{\small $\varphi$}};
	%
	\draw[very thick, red!80!black] (p1) .. controls +(0,0.25) and +(-1,0) ..  (p3)
	-- (p3) .. controls +(1,0) and +(0,-0.25) ..  (p4); 
	%
	\fill [orange!30!white, opacity=0.8] 
	(p1) .. controls +(0,-0.25) and +(-1,0) ..  (p2)
	-- (p2) -- (p2s)
	-- (p2s) --  (u2s)
	-- (u2s) -- (u2)
	-- (u2) .. controls +(-1,0) and +(0,-0.25) ..  (u1)
	;
	%
	\draw[very thick, red!80!black] (p1) .. controls +(0,-0.25) and +(-1,0) ..  (p2) -- (p2s); 
	\draw[very thick, red!80!black] (u1) .. controls +(0,-0.25) and +(-1,0) ..  (u2) -- (u2s); 
	\draw[very thick, red!80!black] (u1) .. controls +(0,0.25) and +(-1,0) ..  (u3)
	-- (u3) .. controls +(1,0) and +(0,-0.25) ..  (u4); 
	\draw[very thick, red!80!black] (u4) .. controls +(0,0.25) and +(1,0) ..  (u5) -- (u6); 
	%
	\fill[color=blue!50!black] (2.5,3) circle (0pt) node {{\small $\alpha$}};
	\fill[color=blue!50!black] (1,3) circle (0pt) node {{\small $\beta$}};
	\fill[color=blue!50!black] (-1.5,2) circle (0pt) node {{\small $\alpha^\dual$}};
	\fill[color=blue!50!black] (4.5,0) circle (0pt) node {{\small $\beta^\dual$}};
	%
	\end{tikzpicture}
}    
=
\varphi^\dual \in \Hom(X^\dual, Y^\dual) \, . 
\ee 

An object~$\alpha$ in a symmetric monoidal 2-category~$\B$ is called \textsl{fully dualisable} if it has a dual~$\alpha^\dual$ such that the 1-morphisms $\tev_\alpha, \tcoev_\alpha$ have both left and right adjoints (as in Section~\ref{subsubsec:Adjoints}). 
The full sub-2-category of fully dualisable objects is denoted~$\Bfd$, and we call~$\B$ fully dualisable if $\B\cong\Bfd$. 

\begin{convention}
	Whether or not a symmetric monoidal 2-category~$\B$ is fully dualisable is a property of~$\B$. 
	If it is fully dualisable, we will assume that we have chosen explicit duality data $(\alpha^\dual, \tev_\alpha, \tcoev_\alpha)$ and adjunction data 
	$({}^\dagger \tev_\alpha, \ev_{\tev_\alpha}, \coev_{\tev_\alpha})$, 
	$(\tev^\dagger_\alpha, \tev_{\tev_\alpha}, \tcoev_{\tev_\alpha})$, 
	$({}^\dagger \tcoev_\alpha, \ev_{\tcoev_\alpha}, \coev_{\tcoev_\alpha})$, 
	$(\tcoev^\dagger_\alpha, \tev_{\tcoev_\alpha}, \tcoev_{\tcoev_\alpha})$ 
	for all $\alpha\in\B$. 
	Put differently, we then view~$\B$ as ``fully dualised''. 
\end{convention}

As shown in \cite{Pstragowski}, the adjunction 1-morphisms $\tev_\alpha, \tcoev_\alpha$ of a fully dualisable object~$\alpha$ do not only have left and right adjoints, but these again have left and right adjoints, and so on infinitely. 
The relations between multiple adjoints are negotiated by the \textsl{Serre automorphism} 
\begin{align}
S_\alpha 
& = \big( 1_\alpha \btimes 
	 \tev_\alpha 
\big) \otimes \big( b_{\alpha,\alpha} \btimes 1_{\alpha^\dual} \big) \otimes \big( 1_\alpha \btimes \tev_\alpha^\dagger \big) 
\nonumber
\\
& \widehat{=} \;\;
\tikzzbox{%
	\begin{tikzpicture}[thick,scale=1.0,color=black, baseline=2cm]
	\coordinate (p1) at (0,0);
	\coordinate (p2) at (2,-0.5);
	\coordinate (p3) at (2,0.5);
	\coordinate (p4) at (4,0);
	\coordinate (u1) at (0,3);
	\coordinate (u2) at (2,2.5);
	\coordinate (u3) at (2,3.5);
	\coordinate (u4) at (4,3);
	\coordinate (ld) at (-2,1);
	\coordinate (lu) at (-2,4);
	\coordinate (rd) at (6,1);
	\coordinate (ru) at (6,4);
	%
	\fill [orange!20!white, opacity=0.8] (p3) -- (ld) -- (lu) -- (u3); 
	\fill [orange!20!white, opacity=0.8] (p3) -- (rd) -- (ru) -- (u3); 
	%
	\fill [orange!25!white, opacity=0.8] 
	(p1) .. controls +(0,0.25) and +(-1,-0.2) ..  (p3)
	-- (p3) .. controls +(1,-0.2) and +(0,0.25) ..  (p4)
	-- (p4) --  (u4)
	-- (u4) .. controls +(0,0.25) and +(1,-0.2) ..  (u3)
	-- (u3) .. controls +(-1,-0.2) and +(0,0.25) ..  (u1)
	;
	%
	\draw[ultra thick] (p3) --  (u3); 
	\fill (2.5,3) circle (0pt) node {{\small $b_{\alpha,\alpha}$}};
	%
	\draw[very thick, red!80!black] (p1) .. controls +(0,0.25) and +(-1,-0.2) ..  (p3); 
	\draw[very thick, red!80!black] (p4) .. controls +(0,0.25) and +(1,-0.2) ..  (p3); 
	\draw[very thick, red!80!black] (p3) -- (ld);
	\draw[very thick, red!80!black] (u3) -- (lu); 
	\draw[very thick, red!80!black] (p3) -- (rd);
	\draw[very thick, red!80!black] (ru) -- (u3); 
	\draw[thin] (lu) --  (ld); 
	\draw[thin] (ru) --  (rd); 
	%
	\fill [orange!30!white, opacity=0.8] 
	(p1) .. controls +(0,-0.25) and +(-1,0) ..  (p2)
	-- (p2) .. controls +(1,0) and +(0,-0.25) ..  (p4)
	-- (p4) --  (u4)
	-- (u4) .. controls +(0,-0.25) and +(1,0) ..  (u2)
	-- (u2) .. controls +(-1,0) and +(0,-0.25) ..  (u1)
	;
	\draw[thin] (p1) --  (u1); 
	\draw[thin] (p4) --  (u4); 
	%
	\draw[very thick, red!80!black] (p1) .. controls +(0,-0.25) and +(-1,0) ..  (p2)
	--(p2) .. controls +(1,0) and +(0,-0.25) .. (p4); 
	\draw[very thick, red!80!black] (u1) .. controls +(0,-0.25) and +(-1,0) ..  (u2) 
	--(u2) .. controls +(1,0) and +(0,-0.25) .. (u4); 
	\draw[very thick, red!80!black] (u4) .. controls +(0,0.25) and +(1,-0.2) ..  (u3); 
	\draw[very thick, red!80!black] (u1) .. controls +(0,0.25) and +(-1,-0.2) .. (u3); 
	%
	\fill (2,0) circle (0pt) node {{\small $\alpha^\dual$}};
	\fill (0.7,3) circle (0pt) node {{\small $\alpha$}};
	\fill (3.3,3) circle (0pt) node {{\small $\alpha$}};
	\fill (5.5,1.5) circle (0pt) node {{\small $\alpha$}};
	\fill (-1.5,1.5) circle (0pt) node {{\small $\alpha$}};
	\fill (-0.3,1.8) circle (0pt) node {{\small $\tev_{\alpha}\vphantom{\tev_{\alpha}^\dagger}$}};
	\fill (4.35,1.8) circle (0pt) node {{\small $\tev_{\alpha}^\dagger$}};
	%
	\end{tikzpicture}
}   
\label{eq:def-serre}
\end{align}
with inverse 
\begin{align}
S^{-1}_\alpha 
& = \big( 1_\alpha \btimes 
	 \tev_\alpha 
\big) \otimes \big( b_{\alpha,\alpha} \btimes 1_{\alpha^\dual} \big) \otimes \big( 1_\alpha \btimes {}^\dagger\tev_\alpha \big) 
\nonumber
\\ 
& \widehat{=} \;\;
\tikzzbox{%
	\begin{tikzpicture}[thick,scale=1.0,color=black, baseline=2cm]
	\coordinate (p1) at (0,0);
	\coordinate (p2) at (2,-0.5);
	\coordinate (p3) at (2,0.5);
	\coordinate (p4) at (4,0);
	\coordinate (u1) at (0,3);
	\coordinate (u2) at (2,2.5);
	\coordinate (u3) at (2,3.5);
	\coordinate (u4) at (4,3);
	\coordinate (ld) at (-2,1);
	\coordinate (lu) at (-2,4);
	\coordinate (rd) at (6,1);
	\coordinate (ru) at (6,4);
	%
	\fill [orange!20!white, opacity=0.8] (p3) -- (ld) -- (lu) -- (u3); 
	\fill [orange!20!white, opacity=0.8] (p3) -- (rd) -- (ru) -- (u3); 
	%
	\fill [orange!25!white, opacity=0.8] 
	(p1) .. controls +(0,0.25) and +(-1,-0.2) ..  (p3)
	-- (p3) .. controls +(1,-0.2) and +(0,0.25) ..  (p4)
	-- (p4) --  (u4)
	-- (u4) .. controls +(0,0.25) and +(1,-0.2) ..  (u3)
	-- (u3) .. controls +(-1,-0.2) and +(0,0.25) ..  (u1)
	;
	%
	\draw[ultra thick] (p3) --  (u3); 
	\fill (2.5,3) circle (0pt) node {{\small $b_{\alpha,\alpha}$}};
	%
	\draw[very thick, red!80!black] (p1) .. controls +(0,0.25) and +(-1,-0.2) ..  (p3); 
	\draw[very thick, red!80!black] (p4) .. controls +(0,0.25) and +(1,-0.2) ..  (p3); 
	\draw[very thick, red!80!black] (p3) -- (ld);
	\draw[very thick, red!80!black] (u3) -- (lu); 
	\draw[very thick, red!80!black] (p3) -- (rd);
	\draw[very thick, red!80!black] (ru) -- (u3); 
	\draw[thin] (lu) --  (ld); 
	\draw[thin] (ru) --  (rd); 
	%
	\fill [orange!30!white, opacity=0.8] 
	(p1) .. controls +(0,-0.25) and +(-1,0) ..  (p2)
	-- (p2) .. controls +(1,0) and +(0,-0.25) ..  (p4)
	-- (p4) --  (u4)
	-- (u4) .. controls +(0,-0.25) and +(1,0) ..  (u2)
	-- (u2) .. controls +(-1,0) and +(0,-0.25) ..  (u1)
	;
	\draw[thin] (p1) --  (u1); 
	\draw[thin] (p4) --  (u4); 
	%
	\draw[very thick, red!80!black] (p1) .. controls +(0,-0.25) and +(-1,0) ..  (p2)
	--(p2) .. controls +(1,0) and +(0,-0.25) .. (p4); 
	\draw[very thick, red!80!black] (u1) .. controls +(0,-0.25) and +(-1,0) ..  (u2) 
	--(u2) .. controls +(1,0) and +(0,-0.25) .. (u4); 
	\draw[very thick, red!80!black] (u4) .. controls +(0,0.25) and +(1,-0.2) ..  (u3); 
	\draw[very thick, red!80!black] (u1) .. controls +(0,0.25) and +(-1,-0.2) .. (u3); 
	%
	\fill (2,0) circle (0pt) node {{\small $\alpha^\dual$}};
	\fill (0.7,3) circle (0pt) node {{\small $\alpha$}};
	\fill (3.3,3) circle (0pt) node {{\small $\alpha$}};
	\fill (5.5,1.5) circle (0pt) node {{\small $\alpha$}};
	\fill (-1.5,1.5) circle (0pt) node {{\small $\alpha$}};
	\fill (-0.3,1.8) circle (0pt) node {{\small $\tev_{\alpha}\vphantom{\tev_{\alpha}^\dagger}$}};
	\fill (4.4,1.8) circle (0pt) node {{\small ${}^\dagger\!\tev_{\alpha}$}};
	%
	\end{tikzpicture}
}  
\, . 
\label{eq:def-serre-inv}
\end{align}
The general result \cite[Thm.\,3.9]{Pstragowski} on multiple adjoints in~$\Bfd$ implies in particular
\be 
\label{eq:evalphadagger} 
\tev_\alpha^\dagger \cong (S_\alpha \btimes 1_{\alpha^\dual} ) \otimes b_{\alpha^\dual,\alpha} \otimes \tcoev_\alpha 
\, , \quad 
{}^\dagger\tev_\alpha \cong (S^{-1}_\alpha \btimes 1_{\alpha^\dual} ) \otimes b_{\alpha^\dual,\alpha} \otimes \tcoev_\alpha \, . 
\ee 

Let $(\Bfd)^\times$ be the maximal sub-2-groupoid of~$\Bfd$. 
Then as shown in \cite[Prop.\,2.8]{HV}, for all $X \in (\Bfd)^\times(\alpha,\beta)$, there are 2-morphisms 
\be 
S_X \colon X \otimes S_\alpha \lra S_\beta \otimes X \, , 
\ee 
which together with the components~$S_\alpha$ assemble into a pseudonatural transformation $\textrm{Id}_{(\Bfd)^\times} \lra \textrm{Id}_{(\Bfd)^\times}$. 
This can be slightly generalised: 

\begin{proposition}
	\label{prop:SonBfd}
	Let~$\B$ be a symmetric monoidal 
		pivotal 
	2-category such that~$\Bfd$ has adjoints for all 1-morphisms. 
	Then the Serre automorphisms~$S_\alpha$ together with the 2-morphisms
		(expressed in terms of the graphical calculus of \cite{BMS} for symmetric monoidal pivotal 2-categories)
	\be 
	\label{eq:SXdiagram}
	S_X = 
	\tikzzbox{%
		\begin{tikzpicture}[thick,scale=1.0,color=black, baseline=2cm]
		\coordinate (p1) at (0,0);
		\coordinate (p2) at (2,-0.5);
		\coordinate (p3) at (2,0.5);
		\coordinate (p4) at (4,0);
		\coordinate (u1) at (0,3);
		\coordinate (u2) at (2,2.5);
		\coordinate (u3) at (2,3.5);
		\coordinate (u4) at (4,3);
		\coordinate (ld) at (-2,1);
		\coordinate (lu) at (-2,4);
		\coordinate (rd) at (6,1);
		\coordinate (ru) at (6,4);
		%
		\fill [orange!20!white, opacity=0.8] (p3) -- (ld) -- (lu) -- (u3); 
		\fill [orange!20!white, opacity=0.8] (p3) -- (rd) -- (ru) -- (u3); 
		%
		\coordinate (x2) at (2,2);
		\coordinate (x5) at (5.5,3.94);
		\draw[ultra thick, blue!50!black] (x2) .. controls +(0,0.75) and +(0,-0.75) .. (x5);
		%
		\fill [orange!25!white, opacity=0.8] 
		(p1) .. controls +(0,0.25) and +(-1,-0.2) ..  (p3)
		-- (p3) .. controls +(1,-0.2) and +(0,0.25) ..  (p4)
		-- (p4) --  (u4)
		-- (u4) .. controls +(0,0.25) and +(1,-0.2) ..  (u3)
		-- (u3) .. controls +(-1,-0.2) and +(0,0.25) ..  (u1)
		;
		%
		\draw[ultra thick] (p3) --  (u3); 
		\fill (1.55,2.2) circle (0pt) node {{\small $b_{X,X}$}};
		\fill[color=blue!50!black] (x2) circle (2.5pt) node {};
		%
		\coordinate (x1) at (-1.5,0.93);
		\coordinate (x3) at (4,2);
		\coordinate (x4) at (0,1);
		\draw[ultra thick, blue!50!black] (x1) .. controls +(0,0.5) and +(-0.2,-0.25) .. (x2); 
		\draw[ultra thick, color=blue!50!black, postaction={decorate}, decoration={markings,mark=at position .59 with {\arrow[draw=blue!50!black]{>}}}] (x2) .. controls +(0.2,0.25) and +(-0.2,0.25) .. (x3); 
		\draw[ultra thick, color=blue!50!black, postaction={decorate}, decoration={markings,mark=at position .45 with {\arrow[draw=blue!50!black]{>}}}] (x4) .. controls +(0.2,-0.65) and +(-0.2,-0.85) .. (x2); 
		\fill[color=blue!50!black] (-1.6,1.2) circle (0pt) node {{\small $X$}};
		\fill[color=blue!50!black] (5.7,3.7) circle (0pt) node {{\small $X$}};
		%
		\draw[very thick, red!80!black] (p1) .. controls +(0,0.25) and +(-1,-0.2) ..  (p3); 
		\draw[very thick, red!80!black] (p4) .. controls +(0,0.25) and +(1,-0.2) ..  (p3); 
		\draw[very thick, red!80!black] (p3) -- (ld);
		\draw[very thick, red!80!black] (u3) -- (lu); 
		\draw[very thick, red!80!black] (p3) -- (rd);
		\draw[very thick, red!80!black] (ru) -- (u3); 
		\draw[thin] (lu) --  (ld); 
		\draw[thin] (ru) --  (rd); 
		%
		\fill [orange!30!white, opacity=0.8] 
		(p1) .. controls +(0,-0.25) and +(-1,0) ..  (p2)
		-- (p2) .. controls +(1,0) and +(0,-0.25) ..  (p4)
		-- (p4) --  (u4)
		-- (u4) .. controls +(0,-0.25) and +(1,0) ..  (u2)
		-- (u2) .. controls +(-1,0) and +(0,-0.25) ..  (u1)
		;
		\draw[thin] (p1) --  (u1); 
		\draw[thin] (p4) --  (u4); 
		%
		\draw[ultra thick, blue!50!black] (x3) .. controls +(0.2,-0.5) and +(0.2,0.5) .. (x4); 
		%
		\draw[very thick, red!80!black] (p1) .. controls +(0,-0.25) and +(-1,0) ..  (p2)
		--(p2) .. controls +(1,0) and +(0,-0.25) .. (p4); 
		\draw[very thick, red!80!black] (u1) .. controls +(0,-0.25) and +(-1,0) ..  (u2) 
		--(u2) .. controls +(1,0) and +(0,-0.25) .. (u4); 
		\draw[very thick, red!80!black] (u4) .. controls +(0,0.25) and +(1,-0.2) ..  (u3); 
		\draw[very thick, red!80!black] (u1) .. controls +(0,0.25) and +(-1,-0.2) .. (u3); 
		%
		\fill (2,0) circle (0pt) node {{\small $\alpha^\dual$}};
		\fill (1.3,2.9) circle (0pt) node {{\small $\beta$}};
		\fill (2.6,2.9) circle (0pt) node {{\small $\beta$}};
		\fill (5.5,1.5) circle (0pt) node {{\small $\alpha$}};
		\fill (-0.5,1.1) circle (0pt) node {{\small $\alpha$}};
		\fill (4.2,3.45) circle (0pt) node {{\small $\beta$}};
		\fill (-1.5,3.5) circle (0pt) node {{\small $\beta$}};
		%
		%
		\end{tikzpicture}
	}  
	\quad \textrm{for all } X \in \Bfd(\alpha,\beta) 
	\ee 
	form a pseudonatural transformation $S\colon \textrm{Id}_{\Bfd} \lra \textrm{Id}_{\Bfd}$. 
\end{proposition}
\begin{proof}
	If $X\in \Bfd(\alpha,\beta)$ has a quasi-inverse~$X^{-1}$, then~$X^{-1}$ is isomorphic to the (chosen) adjoint~$\Xd$, and we have $(X^\dual)^{-1} \cong (X^{-1})^\dual \cong (\Xd)^\dual$. 
	Substituting this into the proof of \cite[Prop.\,2.8]{HV}, we find that specifying~$S_X$ amounts to filling the diagram 
	\be 
	\begin{tikzpicture}[
	baseline=(current bounding box.base), 
	descr/.style={fill=white,inner sep=3.5pt}, 
	normal line/.style={->}
	] 
	\matrix (m) [matrix of math nodes, row sep=6em, column sep=4.2em, text height=1.5ex, text depth=0.1ex] {%
		\alpha \cong \alpha \btimes \one & \alpha \btimes \alpha \btimes \alpha^\dual & \alpha \btimes \alpha \btimes \alpha^\dual & \alpha \btimes \one \cong \alpha 
		\\
		\beta \cong \beta \btimes \one & \beta \btimes \beta \btimes \beta^\dual & \beta \btimes \beta \btimes \beta^\dual & \beta \btimes \one \cong \beta 
		\\
	};
	\path[font=\footnotesize] (m-1-1) edge[->] node[above,sloped] {$ 1_\alpha \btimes \tev_\alpha^\dagger $} (m-1-2);
	\path[font=\footnotesize] (m-1-2) edge[->] node[above,sloped] {$ b_{(\alpha,\alpha)} \btimes 1_{\alpha^\dual}$} (m-1-3);
	\path[font=\footnotesize] (m-1-3) edge[->] node[above,sloped] {$ 1_\alpha \btimes \tev_\alpha $} (m-1-4);
	\path[font=\footnotesize] (m-2-1) edge[->] node[below,sloped] {$ 1_\beta \btimes \tev_\beta^\dagger $} (m-2-2);
	\path[font=\footnotesize] (m-2-2) edge[->] node[below,sloped] {$ b_{(\beta,\beta)} \btimes 1_{\beta^\dual}$} (m-2-3);
	\path[font=\footnotesize] (m-2-3) edge[->] node[below,sloped] {$ 1_\beta \btimes \tev_\beta $} (m-2-4);
	\path[font=\footnotesize, transform canvas={xshift=-8mm}] (m-1-1) edge[->] node[left] {$ X $} (m-2-1);
	\path[font=\footnotesize] (m-1-2) edge[->] node[left] {$ X \btimes X \btimes (X^\dagger)^\dual $} (m-2-2);
	\path[font=\footnotesize] (m-1-3) edge[->] node[right] {$ X \btimes X \btimes (X^\dagger)^\dual $} (m-2-3);
	\path[font=\footnotesize, transform canvas={xshift=+8mm}] (m-1-4) edge[->] node[right] {$ X $} (m-2-4);		
	\end{tikzpicture}
	\ee 
	This is precisely what the expression of~$S_X$ in~\eqref{eq:SXdiagram} does. 
\end{proof}

\begin{example}
	\label{exa:SymmetricMonoidal2Categories}
	We sketch a few fully dualisable symmetric monoidal 2-categories that appear in connection with 2-dimensional TQFT: 
	\begin{enumerate}[label={(\roman*)}]
		\item 
		There is a 2-category $\Bordfr$ of 2-framed points, 1-dimensional bordisms and 2-dimensional bordism classes that we review below in Section~\ref{subsec:bord-2-cat}. 
		The Serre automorphism~$S_+$ of the positively framed point $+ \in \Bordfr$ generates an action of $\pi_1(\textrm{SO}_2) \cong \Z$ and corresponds to a twist of the interval over the point~$+$. 
		\item 
		\label{item:StateSumModels}
		State sum models: 
		For~$\Bbbk$ a field, there is a 2-category $\textrm{Alg}_\Bbbk^{\textrm{fd}}$ of separable $\Bbbk$-algebras, bimodules and bimodule maps \cite{l0905.0465, spthesis}. 
		The Serre automorphism of $A \in \textrm{Alg}_\Bbbk^{\textrm{fd}}$ is the $A$-$A$-bimodule $\Hom_\Bbbk(A,\Bbbk)$. 
		\item 
		\label{item:LGmodels}
		Landau--Ginzburg models: 
		There is a 2-category $\LG$ of isolated singularities, matrix factorisations and their maps up to homotopy \cite{cm1208.1481, CMM}, which we briefly review in Section~\ref{subsec:LGmodels} below. 
		The Serre automorphism of $W \in \LG$ is isomorphic to~$1_W$ up to a shift. 
		\item 
		\label{item:BtwistedSigmaModels}
		B-twisted sigma models: 
		There is a 2-category $\mathcal V\textrm{ar}$ of smooth projective varieties, Fourier--Mukai kernels and Ext groups \cite{cw1007.2679, BanksOnRozanskyWitten}. 
		The Serre automorphism of $X\in\mathcal V\textrm{ar}$ is given by tensoring with the canonical line bundle of~$X$ shifted by $-\dim(X)$. 
		\item 
		Topologically twisted models: 
		There is a 2-category $\textrm{DGSat}_\Bbbk$ of essentially small, smooth, proper and triangulated differential graded $\Bbbk$-categories and their derived categories of bimodules \cite[App.\,A]{bfk1105.3177v3}. 
		The Serre automorphism of $\mathcal C \in \textrm{DGSat}_\Bbbk$ is given in terms of the $\Bbbk$-linear dual composed with the canonical trace functor associated to~$\mathcal C$. 
		The 2-categories of parts~\ref{item:LGmodels} and~\ref{item:BtwistedSigmaModels} are equivalent to sub-2-categories of $\textrm{DGSat}_\Bbbk$. 
	\end{enumerate}
\end{example}

\subsubsection{A Frobenius algebra} 
\label{subsubsec:closed-Lr-Fa-B}

Let~$\B$ be a symmetric monoidal 2-category. 
For a fixed fully dualisable object $\alpha \in \Bfd$ with $\alpha^{\dual\dual}=\alpha$, we now consider the \textsl{$a$-th circle spaces}
\begin{align}
C^\alpha_a  := & \; 
\tev_\alpha \otimes \big( 1_\alpha \btimes S_{\alpha^\dual}^{1-a} \big) \otimes {}^\dagger \tev_\alpha 
\nonumber 
\\ 
\cong & \; 
\tev_\alpha \otimes \big( S_\alpha^{1-a} \btimes 1_{\alpha^\dual} \big) \otimes {}^\dagger \tev_\alpha 
\;\in\; \B(\one,\one) 
\quad \textrm{ for all } a \in \Z 	
\label{eq:Ca} 
\end{align}
where the isomorphism in~\eqref{eq:Ca} is induced by 
\be 
\tikzzbox{%
	\begin{tikzpicture}[thick,scale=1.0,color=black, baseline=1.5cm]
	\coordinate (p1) at (0,0);
	\coordinate (p2) at (2,-0.5);
	\coordinate (p3) at (2.5,0.5);
	\coordinate (u1) at (0,3);
	\coordinate (u2) at (2,2.5);
	\coordinate (u3) at (2.5,3.5);
	%
	\fill [orange!20!white, opacity=0.8] 
	(p1) .. controls +(0,0.25) and +(-1,0) ..  (p3)
	-- (p3) --  (u3)
	-- (u3) .. controls +(-1,0) and +(0,0.25) ..  (u1)
	;
	\draw[ultra thick, blue!50!black] (0.5,0.25) --  (0.5,3.25); 
	\fill[color=blue!50!black] (0.8,3) circle (0pt) node {{\small $S_{\alpha}$}};
	%
	\draw[very thick, red!80!black] (p1) .. controls +(0,0.25) and +(-1,0) ..  (p3); 
	%
	\fill [orange!30!white, opacity=0.8] 
	(p1) .. controls +(0,-0.25) and +(-1,0) ..  (p2)
	-- (p2) --  (u2)
	-- (u2) .. controls +(-1,0) and +(0,-0.25) ..  (u1)
	;
	\draw[thin] (p1) --  (u1); 
	\draw[thin] (p2) --  (u2); 
	\draw[thin] (p3) --  (u3); 
	%
	\draw[very thick, red!80!black] (p1) .. controls +(0,-0.25) and +(-1,0) ..  (p2); 
	\draw[very thick, red!80!black] (u1) .. controls +(0,-0.25) and +(-1,0) ..  (u2); 
	\draw[very thick, red!80!black] (u1) .. controls +(0,0.25) and +(-1,0) ..  (u3); 
	%
	\fill (1.5,0) circle (0pt) node {{\small $\alpha^\dual$}};
	\fill (2,3) circle (0pt) node {{\small $\alpha\vphantom{\alpha^\dual}$}};
	%
	\end{tikzpicture}
}  
\;\; \cong \;\;
\tikzzbox{%
	\begin{tikzpicture}[thick,scale=1.0,color=black, baseline=1.5cm]
	\coordinate (p1) at (0,0);
	\coordinate (p2) at (2,-0.5);
	\coordinate (p3) at (2.5,0.5);
	\coordinate (u1) at (0,3);
	\coordinate (u2) at (2,2.5);
	\coordinate (u3) at (2.5,3.5);
	%
	\fill [orange!20!white, opacity=0.8] 
	(p1) .. controls +(0,0.25) and +(-1,0) ..  (p3)
	-- (p3) --  (u3)
	-- (u3) .. controls +(-1,0) and +(0,0.25) ..  (u1)
	;
	%
	\draw[very thick, red!80!black] (p1) .. controls +(0,0.25) and +(-1,0) ..  (p3); 
	%
	\fill [orange!30!white, opacity=0.8] 
	(p1) .. controls +(0,-0.25) and +(-1,0) ..  (p2)
	-- (p2) --  (u2)
	-- (u2) .. controls +(-1,0) and +(0,-0.25) ..  (u1)
	;
	\draw[ultra thick, blue!50!black] (0.5,-0.3) --  (0.5,2.7); 
	\fill[color=blue!50!black] (0.9,1.3) circle (0pt) node {{\small $S_{\alpha^\dual}$}};
	%
	\draw[thin] (p1) --  (u1); 
	\draw[thin] (p2) --  (u2); 
	\draw[thin] (p3) --  (u3); 
	%
	\draw[very thick, red!80!black] (p1) .. controls +(0,-0.25) and +(-1,0) ..  (p2); 
	\draw[very thick, red!80!black] (u1) .. controls +(0,-0.25) and +(-1,0) ..  (u2); 
	\draw[very thick, red!80!black] (u1) .. controls +(0,0.25) and +(-1,0) ..  (u3); 
	%
	\fill (1.5,0) circle (0pt) node {{\small $\alpha^\dual$}};
	\fill (2,3) circle (0pt) node {{\small $\alpha\vphantom{\alpha^\dual}$}};
	%
	\end{tikzpicture}
}   
\ee 
which in turn is the cusp isomorphism~\eqref{eq:cuspl} combined with $S_\alpha^\dual \cong S_{\alpha^\dual}$. 

If~$\B$ is the 2-category $\textrm{Alg}_\Bbbk^{\textrm{fd}}$ of Example~\ref{exa:SymmetricMonoidal2Categories}\ref{item:StateSumModels}, then for an algebra $A\in \textrm{Alg}_\Bbbk^{\textrm{fd}}$ the zeroth and first circle spaces are the zeroth Hochschild homology and cohomology of~$A$, respectively: $C^A_0 \cong \textrm{HH}_0(A)$ and $C^A_1 \cong \textrm{HH}^0(A)$. 
If~$\B$ is the 2-category $\LG$ of Example~\ref{exa:SymmetricMonoidal2Categories}\ref{item:LGmodels}, then the circle spaces of a given object are (shifts of) the associated Jacobi algebra, as we will explain in Section~\ref{subsec:LGmodels} below. 
In the following we will sometimes treat the isomorphism in~\eqref{eq:Ca} as an identity, and we usually drop the index ``$\alpha$'' in $C^\alpha_a$. 

Next we set 
\begin{align}
\eta_1 & := 
\tikzzbox{%
	\begin{tikzpicture}[thick,scale=1.0,color=black, baseline=-0.9cm]
	\coordinate (p1) at (-2.75,0);
	\coordinate (p2) at (-1,0);
	\fill [orange!20!white, opacity=0.8] 
	(p1) .. controls +(0,-0.5) and +(0,-0.5) ..  (p2)
	-- (p2) .. controls +(0,0.5) and +(0,0.5) ..  (p1)
	;
	\fill [orange!30!white, opacity=0.8] 
	(p1) .. controls +(0,-0.5) and +(0,-0.5) ..  (p2)
	-- (p2) .. controls +(0,-2) and +(0,-2) ..  (p1)
	;
	\draw (p1) .. controls +(0,-2) and +(0,-2) ..  (p2); 
	\draw[very thick, red!80!black] (p1) .. controls +(0,0.5) and +(0,0.5) ..  (p2); 
	\draw[very thick, red!80!black] (p1) .. controls +(0,-0.5) and +(0,-0.5) ..  (p2); 
	%
	\fill ($(p1)+(0.1,0)$) circle (0pt) node[left] {{\small $\tev_\alpha$}};
	\fill ($(p2)+(-0.05,0)$) circle (0pt) node[right] {{\small ${}^\dagger\!\tev_\alpha$}};
	\fill (-1.875,-1.8) circle (0pt) node {{\small $\coev_{\tev_\alpha}$}};
	\end{tikzpicture}
}%
\colon 1_\one \lra C_1 \, , 
\label{eq:Lambda_eta}
\\
\varepsilon_{-1} & := 
\tikzzbox{%
	\begin{tikzpicture}[thick,scale=1.0,color=black, rotate=180, baseline=0.5cm]
	\coordinate (p1) at (-2.75,0);
	\coordinate (p2) at (-1,0);
	%
	\draw[very thick, red!80!black] (p1) .. controls +(0,-0.5) and +(0,-0.5) ..  (p2); 
	%
	\fill [orange!20!white, opacity=0.8] 
	(p1) .. controls +(0,-0.5) and +(0,-0.5) ..  (p2)
	-- (p2) .. controls +(0,0.5) and +(0,0.5) ..  (p1)
	;
	\fill [orange!30!white, opacity=0.8] 
	(p1) .. controls +(0,-0.5) and +(0,-0.5) ..  (p2)
	-- (p2) .. controls +(0,-2) and +(0,-2) ..  (p1)
	;
	\coordinate (q1) at (-1.875,0.38);
	\coordinate (q2) at (-2.67,-0.7);
	\draw (p1) .. controls +(0,-2) and +(0,-2) ..  (p2); 
	\draw[color=blue!50!black] (q1) .. controls +(0,-0.5) and +(0,0.5) ..  (q2); 
	%
	\fill (-2.35,-1.3) circle (0pt) node[right] {{\small $\tev_\alpha^\dagger$}};
	\fill ($(p1)+(0.05,0)$) circle (0pt) node[right] {{\small ${}^\dagger\!\tev_\alpha$}};
	\fill ($(p2)+(-0.1,0)$) circle (0pt) node[left] {{\small $\tev_\alpha$}};
	\fill (-1.875,-1.75) circle (0pt) node {{\small $\tev_{\tev_\alpha}$}};
	\fill[color=blue!50!black] (q1) circle (0pt) node[below] {{\small $S_{\alpha^\dual}^{2}$}};
	\fill[color=blue!50!black] (q2) circle (2pt) node {};
	%
	\draw[very thick, red!80!black] (p1) .. controls +(0,0.5) and +(0,0.5) ..  (p2); 
	\end{tikzpicture}
}%
\colon C_{-1} \lra 1_\one \, , 
\label{eq:Lambda_epsilon}
\\ 
\mu_{a,b} & := 
\tikzzbox{%
	\begin{tikzpicture}[thick,scale=1.0,color=black, baseline=1.9cm]
	\coordinate (p1) at (-2.75,0);
	\coordinate (p2) at (-1,0);
	\coordinate (p3) at (1,0);
	\coordinate (p4) at (2.75,0);
	\coordinate (p5) at (0.875,4);
	\coordinate (p6) at (-0.875,4);
	%
	\draw[very thick, red!80!black] (p1) .. controls +(0,0.5) and +(0,0.5) ..  (p2); 
	\draw[very thick, red!80!black] (p3) .. controls +(0,0.5) and +(0,0.5) ..  (p4); 
	%
	\fill [orange!30!white, opacity=0.8] 
	(p1) .. controls +(0,0.5) and +(0,0.5) ..  (p2)
	-- (p2) .. controls +(0,1.5) and +(0,1.5) ..  (p3)
	-- (p3) .. controls +(0,0.5) and +(0,0.5) ..  (p4)
	-- (p4) .. controls +(0,2.5) and +(0,-2.5) ..  (p5)
	-- (p5) .. controls +(0,-0.5) and +(0,-0.5) ..  (p6)
	-- (p6) .. controls +(0,-2.5) and +(0,2.5) ..  (p1)
	;
	\fill [orange!20!white, opacity=0.8] 
	(p5) .. controls +(0,-0.5) and +(0,-0.5) ..  (p6)
	-- (p6) .. controls +(0,0.5) and +(0,0.5) ..  (p5)
	;
	\fill [orange!20!white, opacity=0.8] 
	(p1) .. controls +(0,-0.5) and +(0,-0.5) ..  (p2)
	-- (p2) .. controls +(0,0.5) and +(0,0.5) ..  (p1)
	;
	\fill [orange!20!white, opacity=0.8] 
	(p3) .. controls +(0,-0.5) and +(0,-0.5) ..  (p4)
	-- (p4) .. controls +(0,0.5) and +(0,0.5) ..  (p3)
	;
	\draw (p2) .. controls +(0,1.5) and +(0,1.5) ..  (p3); 
	\draw (p4) .. controls +(0,2.5) and +(0,-2.5) ..  (p5); 
	\draw (p6) .. controls +(0,-2.5) and +(0,2.5) ..  (p1); 
	\coordinate (q1) at (-1.875,-0.37);
	\coordinate (q2) at (1.875,-0.37);
	\coordinate (q3) at (0,2);
	\coordinate (q4) at (0,3.63);
	\draw[color=blue!50!black] (q1) .. controls +(0,1) and +(-0.3,0) ..  (q3); 
	\draw[color=blue!50!black] (q2) .. controls +(0,1) and +(0.3,0) ..  (q3); 
	\draw[color=blue!50!black] (q4) -- (q3); 
	%
	\fill[color=blue!50!black] (q1) circle (0pt) node[below] {{\small $S_{\alpha^\dual}^{1-a}$}};
	\fill[color=blue!50!black] (q2) circle (0pt) node[below] {{\small $S_{\alpha^\dual}^{1-b}$}};
	\fill[color=blue!50!black] (q4) circle (0pt) node[above] {{\small $S_{\alpha^\dual}^{2-a-b}$}};
	\fill ($(p1)+(0.1,0)$) circle (0pt) node[left] {{\small $\tev_\alpha$}};
	\fill ($(p3)+(0.1,0)$) circle (0pt) node[left] {{\small $\tev_\alpha$}};
	\fill ($(p6)+(0.1,0)$) circle (0pt) node[left] {{\small $\tev_\alpha$}};
	\fill ($(p2)+(-0.05,0)$) circle (0pt) node[right] {{\small ${}^\dagger\!\tev_\alpha$}};
	\fill ($(p4)+(-0.05,0)$) circle (0pt) node[right] {{\small ${}^\dagger\!\tev_\alpha$}};
	\fill ($(p5)+(-0.05,0)$) circle (0pt) node[right] {{\small ${}^\dagger\!\tev_\alpha$}};
	\fill (0,0.75) circle (0pt) node {{\small $\ev_{\tev_\alpha}$}};
	\fill[color=blue!50!black] (q3) circle (2pt) node[above] {};
	%
	\draw[very thick, red!80!black] (p1) .. controls +(0,-0.5) and +(0,-0.5) ..  (p2); 
	\draw[very thick, red!80!black] (p3) .. controls +(0,-0.5) and +(0,-0.5) ..  (p4); 
	\draw[very thick, red!80!black] (p5) .. controls +(0,0.5) and +(0,0.5) ..  (p6); 
	\draw[very thick, red!80!black] (p5) .. controls +(0,-0.5) and +(0,-0.5) ..  (p6); 
	\end{tikzpicture}
}%
\colon C_a \otimes C_b \lra C_{a+b-1} \, , 
\label{eq:Lambda_mu}
\\ 
\Delta_{a,b} & := 
\tikzzbox{%
	\begin{tikzpicture}[thick,scale=1.0,color=black, rotate=180, baseline=-2.5cm]
	\coordinate (p1) at (-2.75,0);
	\coordinate (p2) at (-1,0);
	\coordinate (p3) at (1,0);
	\coordinate (p4) at (2.75,0);
	\coordinate (p5) at (0.875,4);
	\coordinate (p6) at (-0.875,4);
	%
	\draw[very thick, red!80!black] (p5) .. controls +(0,-0.5) and +(0,-0.5) ..  (p6); 
	%
	\fill [orange!30!white, opacity=0.8] 
	(p1) .. controls +(0,0.5) and +(0,0.5) ..  (p2)
	-- (p2) .. controls +(0,1.5) and +(0,1.5) ..  (p3)
	-- (p3) .. controls +(0,0.5) and +(0,0.5) ..  (p4)
	-- (p4) .. controls +(0,2.5) and +(0,-2.5) ..  (p5)
	-- (p5) .. controls +(0,-0.5) and +(0,-0.5) ..  (p6)
	-- (p6) .. controls +(0,-2.5) and +(0,2.5) ..  (p1)
	;
	\fill [orange!20!white, opacity=0.8] 
	(p1) .. controls +(0,-0.5) and +(0,-0.5) ..  (p2)
	-- (p2) .. controls +(0,0.5) and +(0,0.5) ..  (p1)
	;
	\fill [orange!20!white, opacity=0.8] 
	(p3) .. controls +(0,-0.5) and +(0,-0.5) ..  (p4)
	-- (p4) .. controls +(0,0.5) and +(0,0.5) ..  (p3)
	;
	\fill [orange!20!white, opacity=0.8] 
	(p5) .. controls +(0,-0.5) and +(0,-0.5) ..  (p6)
	-- (p6) .. controls +(0,0.5) and +(0,0.5) ..  (p5)
	;
	\draw (p2) .. controls +(0,1.5) and +(0,1.5) ..  (p3); 
	\draw (p4) .. controls +(0,2.5) and +(0,-2.5) ..  (p5); 
	\draw (p6) .. controls +(0,-2.5) and +(0,2.5) ..  (p1); 
	\coordinate (q1) at (-1.875,0.37);
	\coordinate (q2) at (1.875,0.37);
	\coordinate (q3) at (0,2.5);
	\coordinate (q4) at (0,4.37);
	\coordinate (q5) at (1.25,1.6);
	\coordinate (q6) at (0.94,0.43);
	\draw[color=blue!50!black] (q1) .. controls +(0,1) and +(-0.3,0) ..  (q3); 
	\draw[color=blue!50!black] (q2) .. controls +(0,1) and +(0.3,0) ..  (q3); 
	\draw[color=blue!50!black] (q4) -- (q3); 
	\draw[color=blue!50!black] (q5) .. controls +(0,-0.5) and +(0.3,0.4) ..  (q6); 
	%
	\fill[color=blue!50!black] (q1) circle (0pt) node[above] {{\small $S_{\alpha^\dual}^{1-b}$}};
	\fill[color=blue!50!black] (q2) circle (0pt) node[above] {{\small $S_{\alpha^\dual}^{1-a}$}};
	\fill[color=blue!50!black] (q4) circle (0pt) node[below] {{\small $S_{\alpha^\dual}^{-a-b}$}};
	\fill[color=blue!50!black] (0.2,2) circle (0pt) node {{\scriptsize $S_{\alpha^\dual}^{-1-a}$}};
	\fill[color=blue!50!black] (0.9,1.25) circle (0pt) node {{\scriptsize $S_{\alpha^\dual}^{-2}$}};
	\fill ($(p2)+(-0.1,0)$) circle (0pt) node[left] {{\small $\tev_\alpha$}};
	\fill ($(p4)+(-0.1,0)$) circle (0pt) node[left] {{\small $\tev_\alpha$}};
	\fill ($(p5)+(-0.1,0)$) circle (0pt) node[left] {{\small $\tev_\alpha$}};
	\fill ($(p1)+(0.05,0)$) circle (0pt) node[right] {{\small ${}^\dagger\!\tev_\alpha$}};
	\fill ($(p3)+(0.05,0)$) circle (0pt) node[right] {{\small ${}^\dagger\!\tev_\alpha$}};
	\fill ($(p6)+(0.05,0)$) circle (0pt) node[right] {{\small ${}^\dagger\!\tev_\alpha$}};
	\fill (-0.3,1.37) circle (0pt) node {{\small $\tcoev_{\tev_\alpha}$}};
	\fill (0.4,0.7) circle (0pt) node {{\small $\tev_\alpha^\dagger$}};
	\fill[color=blue!50!black] (q3) circle (2pt) node[above] {};
	\fill[color=blue!50!black] (q5) circle (2pt) node[above] {};
	\fill[color=blue!50!black] (q6) circle (2pt) node[above] {};
	%
	%
	\draw[very thick, red!80!black] (p1) .. controls +(0,0.5) and +(0,0.5) ..  (p2); 
	\draw[very thick, red!80!black] (p1) .. controls +(0,-0.5) and +(0,-0.5) ..  (p2); 
	\draw[very thick, red!80!black] (p3) .. controls +(0,0.5) and +(0,0.5) ..  (p4); 
	\draw[very thick, red!80!black] (p3) .. controls +(0,-0.5) and +(0,-0.5) ..  (p4); 
	\draw[very thick, red!80!black] (p5) .. controls +(0,0.5) and +(0,0.5) ..  (p6); 
	\end{tikzpicture}
}%
\colon C_{a+b+1} \lra C_a \otimes C_b 
\label{eq:Lambda_Delta}
\end{align}
for all $a,b \in \Z$, where in the expressions for~$\varepsilon_{-1}$ and $\Delta_{a,b}$ we use the isomorphisms 
\be 
\label{eq:IsosWithSsquared}
\big(S_\alpha^2 \btimes 1_{\alpha^\dual} \big) \otimes {}^\dagger \tev_\alpha \lra \tev_\alpha^\dagger 
\, , \quad 
\big(1_\alpha \btimes S_{\alpha^\dual}^{-2} \big) \otimes \tev_\alpha^\dagger \lra {}^\dagger \tev_\alpha 
\ee 
obtained from~\eqref{eq:evalphadagger}. 
The above data have a familiar structure (recall the definition of closed $\Lambda_0$-Frobenius algebras in Section~\ref{subsec:ClassificationThroughClosedLambdaFrob}):  

\begin{proposition}
	\label{prop:Lambda0FrobAlgInB11}
	The data $\{ C_a \}_{a\in\Z}$ and $\eta_1, \varepsilon_{-1}, \{ \mu_{a,b} , \Delta_{a,b} \}_{a,b\in\Z}$ have the properties of a closed $\Lambda_0$-Frobenius algebra in the symmetric monoidal category $\B(\one, \one)$. 
\end{proposition}
\begin{proof}
	That $\eta_1, \varepsilon_{-1}, \{ \mu_{a,b} , \Delta_{a,b} \}_{a,b\in\Z}$ satisfy the (co)associativity, (co)unitality and Frobenius conditions is straightforward to check using the diagrammatic calculus for monoidal 2-categories. 
	The remaining defining relations~\eqref{eq:clLrFa-commutativity}--\eqref{eq:clLrFa-deck-trf} of a closed $\Lambda_0$-Frobenius algebra are more difficult to verify directly. 
	Instead, we use the framed cobordism hypothesis (Theorem~\ref{thm:framed-CH} below) to argue indirectly: 
	To $\alpha\in\Bfd$ corresponds a symmetric monoidal functor $\zz \colon \Bordfr \lra \B$ with $\zz(+)=\alpha$, such that the data $\eta_1, \varepsilon_{-1}, \{ \mu_{a,b} , \Delta_{a,b} \}_{a,b\in\Z}$ are the images under~$\zz$ of 2-morphisms in $\Bordfr$. 
	The latter in turn are generators of the framed bordism 1-category and satisfy all the relations of a closed $\Lambda_0$-Frobenius algebras, which follows from the special case $r=0$ of \cite[Thm.\,5.2.1]{SternSzegedy}. 
	Hence also their images $\eta_1, \varepsilon_{-1}, \{ \mu_{a,b} , \Delta_{a,b} \}_{a,b\in\Z}$ in~$\B$ satisfy these relations. 
\end{proof}

Together with the $r$-spin cobordism hypothesis proved in Section~\ref{sec:bord-2-cat} below this implies: 

\begin{corollary}
	\label{cor:LambdarFrobAlgInB11}
	For $r\in\Z_{\geqslant 1}$, every isomorphism $S_\alpha^r \cong 1_\alpha$ induces a closed $\Lambda_r$-Frobenius algebra structure on $\{C_a \}_{a\in\{ 0,1,\dots,r-1\}}$. 
\end{corollary}
\begin{proof}
	Combine \cite[Thm.\,5.2.1]{SternSzegedy} for $r\in\Z_{\geqslant 1}$ with Theorem~\ref{thm:r-spin-CH} below. 
	This in particular guarantees the existence of isomorphisms $C_a \cong C_{a+r}$ for all $a\in\Z$. 
\end{proof}

\subsection[The 2-category of \texorpdfstring{$r$}{r}-spin bordisms]{The 2-category of \texorpdfstring{$\boldsymbol{r}$}{r}-spin bordisms}
\label{subsec:bord-2-cat}

\subsubsection{2-categories of bordisms with tangential structure}

Here we briefly recall 2-categories of bordisms with $G$-structure. 
For more background and details we refer to \cite[Sections 3.1--3.3]{spthesis}. 

\medskip 

We begin by fixing conventions for double categories.
A \textsl{double category} $\D$ consists of 
a category of objects $\D_0$,
a category of horizontal morphisms $\D_1$,
unit horizontal morphisms,
a composition functor, 
and natural transformations that implement associativity and unitality of the composition.
The morphisms of $\D_1$ are called 2-morphisms.
The \textsl{horizontal 2-category} of a double category~$\D$ is the 2-category consisting of
the objects of $\D_0$, horizontal 1-morphisms and 2-morphisms between parallel 1-morphisms.

\begin{figure}[tb]
	\centering
	\begin{tikzpicture}
	\begin{pgfonlayer}{nodelayer}
		\node [style=black dot] (0) at (-3, -3) {};
		\node [style=black dot] (1) at (-9, 2) {};
		\node [style=black dot] (3) at (-9, -3) {};
		\node [style=none] (4) at (-6, 2) {$\longhookrightarrow$};
		\node [style=none] (5) at (-9, 4.5) {$T$};
		\node [style=black dot] (7) at (9, -3) {};
		\node [style=black dot] (9) at (9, 2) {};
		\node [style=black dot] (10) at (3, -3) {};
		\node [style=none] (11) at (6, 2) {$\longhookleftarrow$};
		\node [style=none] (12) at (9, 4.5) {$S$};
		\node [style=none, rotate=90] (13) at (-9, -0.5) {$\longhookrightarrow$};
		\node [style=none] (15) at (-10, 2) {};
		\node [style=none] (16) at (-8, 2) {};
		\node [style=none] (17) at (8, 2) {};
		\node [style=none] (18) at (10, 2) {};
		\node [style=none] (21) at (0, 4.5) {$M$};
		\node [style=none, rotate=90] (22) at (9, -0.5) {$\longhookrightarrow$};
		\node [style=none, rotate=90] (23) at (0, -0.5) {$\longhookrightarrow$};
		\node [style=black dot] (2) at (-3, 2) {};
		\node [style=black dot] (8) at (3, 2) {};
		\node [style=none] (19) at (4, 2) {};
		\node [style=none] (20) at (-4, 2) {};
	\end{pgfonlayer}
	\begin{pgfonlayer}{edgelayer}
		\draw [style=orange front] (17.center) to (18.center);
		\draw [style=orange front] (15.center) to (16.center);
		\draw [style=red] (15.center) to (1);
		\draw [style=blue] (18.center) to (9);
		\draw [style=red] (9) to (17.center);
		\draw [style=blue] (1) to (16.center);
		\draw [style=red, in=150, out=75, looseness=0.75] (0) to (10);
		\path [style=orange front]
			(20.center) 
			to [bend right] (2)
			to [in=150, out=60, looseness=0.75] (8)
			to [bend right] (19.center);
		\draw [style=blue, bend right] (20.center) to (2);
		\draw [style=red, in=150, out=60, looseness=0.75] (2) to (8);
		\draw [style=blue, bend left] (19.center) to (8);
	\end{pgfonlayer}
\end{tikzpicture}
	\caption{0- and 1-dimensional manifolds (below) and their (vertical) inclusions into 2-haloes (above).
		The horizontal embeddings give a 
			 1-bordism 
		$M\colon S\lra T$.}
	\label{fig:halos-bordisms}
\end{figure}
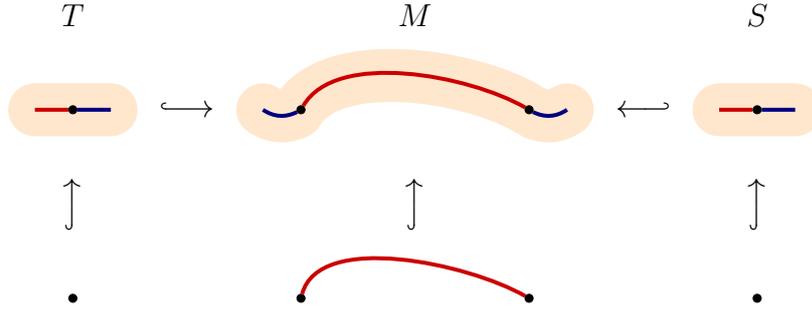

We continue with a sketch of the double category of bordisms with
tangential structure for a chosen group homomorphism $\xi\colon G\lra \GL_2$, which we denote $\BBord^{G}$.
A \textsl{(2-)halo} of a $d$-dimensional manifold ($d\leqslant2$)
is loosely speaking 
a stratified 2-dimensional manifold in which the $d$-manifold is embedded.
We will not need the precise definition, but 
we refer to Figure~\ref{fig:halos-bordisms} for illustrative examples of the cases $d=0$ and $d=1$.
A \textsl{$d$-bordism} between two $(d-1)$-dimensional manifolds with haloes~$S$ and~$T$ 
is a $d$-dimensional compact manifold $M$ together with an embedding $S\sqcup T\longhookrightarrow M$ that identifies the boundary of~$M$ with the $(d-1)$-dimensional manifold underlying the halo $S\sqcup T$. 
We write $M\colon S\lra T$.
A \textsl{diffeomorphism} of such a bordism is a diffeomorphism which
is compatible with the boundary parametrisation maps.

The objects of $(\BBord^{G})_0$ are compact 0-dimensional manifolds with haloes
with $G$-structure, and morphisms are diffeomorphisms of these haloes with $G$-structure.
The objects of $(\BBord^{G})_1$ are 1-bordisms with $G$-structure (recall Section~\ref{subsec:ClosedrSpinTQFTs}),
and its morphisms are diffeomorphism classes of 2-bordisms with
$G$-structure.
The composition functor is given by glueing of bordisms.

The double categories $\BBord^{G}$ are symmetric monoidal via the disjoint union. 
The \textsl{2-category of bordisms with $G$-structure}
$\Bord_{2,1,0}^{G}$ is defined to be the horizontal 2-category of $\BBord^{G}$.
We will use the notation
\begin{equation}
\Bordfr\,,\quad\Bordor\,,\quad\Bordrspin
\label{eq:notation-bord-2-cat}
\end{equation}
for the 2-categories of framed, oriented and $r$-spin bordisms, respectively.
By \cite[Thm.\,1.1]{WesterHansenShulman} these categories inherit 
a symmetric monoidal structure from the respective double categories.

\subsubsection{Functors from group homomorphisms}

Consider the following commutative diagram of homomorphisms of topological groups: 
\begin{equation}
\begin{tikzpicture}[
baseline=(current bounding box.base),
descr/.style={fill=white,inner sep=3.5pt},
normal line/.style={->}
]
\matrix (m) [matrix of math nodes, row sep=2em, column sep=2em, text height=1.5ex, text depth=0.1ex] {%
	G&&G'
	\\
	&\GL_2&\\
};
\path[font=\footnotesize] (m-1-1) edge[->] node[below] {$ \xi $} (m-2-2);
\path[font=\footnotesize] (m-1-1) edge[->] node[above] {$\lambda $} (m-1-3);
\path[font=\footnotesize] (m-1-3) edge[->] node[below] {$\xi'$ } (m-2-2);
\end{tikzpicture}
\label{eq:group-hom-lambda}
\end{equation}
For a $G$-structure $(P,q)$ on a surface $\Sigma$, 
the group homomorphism $\lambda$ induces a $G'$-structure 
on $\Sigma$ via the associated bundle construction: 
\begin{equation}
\begin{tikzpicture}[
baseline=(current bounding box.base),
descr/.style={fill=white,inner sep=3.5pt},
normal line/.style={->}
]
\matrix (m) [matrix of math nodes, row sep=2em, column sep=2em, text height=1.5ex, text depth=0.1ex] {%
	P\times_{\lambda}G'&&F\Sigma
	\\
	&\Sigma&\\
};
\path[font=\footnotesize] (m-1-1) edge[->] node[above] {$q_{\lambda} $} (m-1-3);
\path[font=\footnotesize] (m-1-1) edge[->] node[below] {$ \pi_{\lambda} $} (m-2-2);
\path[font=\footnotesize] (m-1-3) edge[->] node[below] {$ $ } (m-2-2);
\end{tikzpicture}  
\label{eq:induced-tangential-str}
\end{equation}
This construction is compatible with glueing of bordisms with tangential structure,
and with disjoint union. 
Hence it gives rise to symmetric monoidal functors of double categories and of 2-categories:
\begin{equation}
\begin{tikzpicture}[ baseline=(current bounding box.base), 
descr/.style={fill=white,inner sep=3.5pt}, normal line/.style={->} ]
\matrix (m) [matrix of math nodes, row sep=2em, column sep=2em, text height=1.5ex, text depth=0.1ex] {%
	\BBord^{G}&\BBord^{G'} \, , 
	\\
};
\path[font=\footnotesize] (m-1-1) edge[->] node[above] {$\mathbbold{\Lambda}_\lambda$} (m-1-2);
\end{tikzpicture}  
\quad 
\begin{tikzpicture}[ baseline=(current bounding box.base), 
descr/.style={fill=white,inner sep=3.5pt}, normal line/.style={->} ]
\matrix (m) [matrix of math nodes, row sep=2em, column sep=2em, text height=1.5ex, text depth=0.1ex] {%
	\Bord_{2,1,0}^{G}&\Bord_{2,1,0}^{G'} \,. 
	\\
};
\path[font=\footnotesize] (m-1-1) edge[->] node[above] {$\Lambda_\lambda$} (m-1-2);
\end{tikzpicture}  
\label{eq:induced-functor-lambda}
\end{equation}

The group homomorphisms in \eqref{eq:include-1}--\eqref{eq:def-r-spin-group} 
fit into the commutative diagram 
\begin{equation}
\begin{tikzpicture}[ baseline=(current bounding box.base), 
descr/.style={fill=white,inner sep=3.5pt}, normal line/.style={->} ]
\matrix (m) [matrix of math nodes, row sep=4em, column sep=5em, text height=1.5ex, text depth=0.1ex] {%
	\star&\widetilde{\GL_2^+}^r&\GL_2^+ 
	\\
	&\GL_2&\\
};
\path[font=\footnotesize] (m-1-1) edge[right hook->] node[above] {$\widetilde{\lambda}$} (m-1-2);
\path[font=\footnotesize] (m-1-2) edge[->] node[above] {$\lambda:=p^r$} (m-1-3);
\path[font=\footnotesize] (m-1-1) edge[right hook->] node[below] {$ $} (m-2-2);
\path[font=\footnotesize] (m-1-2) edge[->] node[right] {$\iota\circ p^r $ } (m-2-2);
\path[font=\footnotesize] (m-1-3) edge[right hook->] node[below] {$\iota $} (m-2-2);
\end{tikzpicture}  
\label{eq:ex-group-hom-comp}
\end{equation}
and induce symmetric monoidal functors
\begin{equation}
\begin{tikzpicture}[ baseline=(current bounding box.base), 
descr/.style={fill=white,inner sep=3.5pt}, normal line/.style={->} ]
\matrix (m) [matrix of math nodes, row sep=2em, column sep=5em, text height=1.5ex, text depth=0.1ex] {%
	\Bordfr&\Bordrspin&\Bordor
	\, , 
	\\
};
\path[font=\footnotesize] (m-1-1) edge[->] node[above] {$\widetilde{\Lambda} := \Lambda_{\widetilde\lambda}$} (m-1-2);
\path[font=\footnotesize] (m-1-2) edge[->] node[above] {$\Lambda := \Lambda_{p^r}$} (m-1-3);
\end{tikzpicture}  
\label{eq:ex-tlambda-lambda}
\end{equation}
where we use the notation of~\eqref{eq:induced-functor-lambda} and~\eqref{eq:ex-group-hom-comp}. 
The functor $\widetilde{\Lambda}$ assigns to a framed manifold the
manifold with the trivial $r$-spin structure corresponding to the orientation
induced by the framing, and $\Lambda$ assigns 
to a haloed $r$-spin surface the haloed surface with the underlying orientation.

\subsection[Fully extended \texorpdfstring{$r$}{r}-spin TQFTs]{Fully extended \texorpdfstring{$\boldsymbol{r}$}{r}-spin TQFTs}
\label{sec:bord-2-cat}

In this section we consider 2-dimensional extended TQFTs with tangential structure 
and the cobordism hypothesis for $r$-spin structures, $r\in\Z_{\geqslant 0}$.
For this we first recall the framed cobordism hypothesis, homotopy group actions on 2-categories, 
and their homotopy fixed points.
The latter are expected to describe TQFTs with tangential structures,
as is known to be the case for oriented (or equivalently: 1-spin) TQFTs.
After a review of earlier results in the oriented case,
we give a presentation of all $r$-spin bordism 2-categories
in terms of fully dualisable objects,
and prove the $r$-spin cobordism hypothesis 
(for 2-categories, not for $(\infty,2)$-categories).

\begin{definition}
	Let~$\B$ be a symmetric monoidal 2-category. 
	A \textsl{fully extended 2-dimensional TQFT with $G$-structure valued in~$\B$} 
	is a symmetric monoidal functor
	\begin{equation}
	\zz\colon\Bord_{2,1,0}^G\lra\B \, . 
	\label{eq:def:TQFT}
	\end{equation}
	\label{def:TQFT}
\end{definition}
We write 
$ 
\Fun(\Bord_{2,1,0}^G,\B)
$ 
for the {symmetric monoidal 2-groupoid of fully extended TQFTs with $G$-structure and values in~$\B$}.

\subsubsection{The framed cobordism hypothesis}
\label{sec:framed-CH}

Denote with $\FD{0}$ the symmetric monoidal 2-category freely generated by
a single 2-dualisable object $+$, cf.\ \cite{spthesis,Pstragowski}. 
Our slightly ambiguous notation for the generating object in $\FD{0}$ draws justification from the following fact: 

\begin{theorem}[{\cite[Thm.\,7.1]{Pstragowski}}]
	There is an equivalence of symmetric monoidal 2-categories
	\begin{equation}
	\begin{aligned}
	\iota^{0}\colon\FD{0}&\lra\Bordfr\\
	+&\longmapsto +
	\end{aligned}
	\label{eq:iota-0}
	\end{equation}
	sending the object $+\in\FD{0}$ to the positively framed (halo of a) point $+$.
	\label{thm:bord-framed-pres}
\end{theorem}

The framed cobordism hypothesis classifies framed fully extended TQFTs in terms of fully dualisable objects: 

\begin{theorem}[{\cite[Thm.\,8.1]{Pstragowski}}]
	The 2-groupoid of
	framed fully extended TQFTs with target~$\B$ is equivalent to the 
	core of the
	2-category of fully dualisable objects in $\B$ as a symmetric monoidal 2-groupoid:
	\begin{equation}
	\Fun\big(\Bordfr,\B\big)\cong\big(\B^\mathrm{fd}\big)\!{}^{\times}\,.
	\label{eq:framed-CH}
	\end{equation}
	\label{thm:framed-CH}
\end{theorem}

\subsubsection[Homotopy \texorpdfstring{$G$}{G}-actions on 2-categories]{Homotopy \texorpdfstring{$\boldsymbol{G}$}{G}-actions on 2-categories}

In order to state the cobordism hypothesis with 
orientation and more generally with $r$-spin structure, we will need the notion of homotopy action of a group 
on a 2-category, as well as its fixed points.

Let $G$ be a topological group. The \textsl{homotopy action of $G$} 
on a symmetric monoidal 2-category $\B$ is a monoidal functor
\begin{equation}
\rho\colon\Pi_{\leqslant 2}(G)\lra\Aut^{\textrm{sm}}(\B)
\label{eq:htpy-action}
\end{equation}
from the fundamental 2-groupoid of $G$ to the 2-category of
symmetric monoidal autoequivalences of $\B$. 
On $\Pi_{\leqslant 2}(G)$ the monoidal structure comes from the group structure on $G$,
on $\Aut^{\textrm{sm}}$ it is the composition of functors.
Equivalently, a homotopy $G$-action on~$\B$ is a functor 
\begin{equation}
\begin{aligned}
\textrm{B}\Pi_{\leqslant 2}(G)&\lra\text{\{symmetric monoidal 2-categories\}}\\
\star &\longmapsto \B
\end{aligned}
\label{eq:htpy-action-delooping}
\end{equation}
from  the delooping of $\Pi_{\leqslant 2}(G)$ to the 3-category of
symmetric monoidal 2-categories.

For $G=\GL_2^+\simeq \SO_2$ the fundamental 2-groupoid is equivalent to the
2-groupoid $\textrm{B}\underline{\Z}$ with a single object~$\star$ with automorphism group the free abelian group 
on a single generator~$\Z$, and only identity 2-morphisms.
Below we will identify $\Pi_{\leqslant 2}(\GL_2^+)$ with the 2-groupoid $\textrm{B}\underline{\Z}$. 

To define an action $\rho$ of $\GL_2^+$ on a 2-category we only need to specify 
the value of $\rho$ on the generator $-1\in\Z$.
Recall from Proposition~\ref{prop:SonBfd} the Serre automorphism 
$S\colon \textrm{Id}_{\Bfd}\lra \textrm{Id}_{\Bfd}$.
We define the homotopy action of $\GL_2^+$ on fully dualisable objects as follows:
\begin{equation}
\begin{aligned}
\rho\colon \Pi_{\leqslant 2}(\GL_2^+)&\lra\Aut^{\textrm{sm}}(\Bfd)\\
\star&\longmapsto \textrm{Id}_{\Bfd}\\
\Z\ni -1&\longmapsto S\\
1&\longmapsto 1 \, . 
\end{aligned}
\label{eq:rho-SO2}
\end{equation}
Similarly, the homotopy action of the $r$-spin group is defined through the $r$-th power of the 
Serre automorphism:
\begin{equation}
\begin{aligned}
\rho^r\colon \Pi_{\leqslant 2}(\widetilde{\GL_2^+}^r)&\lra\Aut^{\textrm{sm}}(\Bfd)\\
\star&\longmapsto \textrm{Id}_{\Bfd}\\
\Z\ni -1&\longmapsto S^r\\
1&\longmapsto 1 \, . 
\end{aligned}
\label{eq:rho-r-spin}
\end{equation}
Note that this action is the $\GL_2^+$-action \eqref{eq:rho-SO2} composed with 
the functor induced from the covering map $p^r\colon \widetilde{\GL_2^+}^r\lra \GL_2^+$
in \eqref{eq:def-r-spin-group}.

\subsubsection[Presentations of \texorpdfstring{$r$}{r}-spin bordism 2-categories]{Presentations of \texorpdfstring{$\boldsymbol{r}$}{r}-spin bordism 2-categories}
\label{subsubsec:CHTangential}

The 2-category of \textsl{homotopy fixed points} $\B^G$ of a homotopy action~$\rho$ as in~\eqref{eq:htpy-action-delooping} is defined to be the 2-category of natural transformations of functors of 3-categories 
\begin{equation}
\B^G  = \Nat(\Delta_{\star},\rho) \,,
\label{eq:co-invariants-co-limit}
\end{equation}
where the constant functor 
$\Delta_{\star}\colon \textrm{B}\Pi_{\leqslant 2}(G)\lra\text{\{sym.\ mon.\ 2-cat.\}}$
sends the unique object in $\textrm{B}\Pi_{\leqslant 2}(G)$ to the 2-category $\star$ 
with a single object and only identity morphisms, see \cite[Rem.\,3.11--3.14]{HSV}. 
It is expected that $\B^G$ is the 3-limit of the functor \eqref{eq:htpy-action-delooping}, 
but we are not aware of a rigorous development of the theory of 3-limits. 

By the cobordism hypothesis it is expected that 2-dimensional fully extended TQFTs
with $G$-structure and target~$\B$ are classified by homotopy fixed points of a $G$-action on
$(\Bfd)^{\times}$, originating from the $G$-action on $\Bordfr$. 
To our knowledge there is no complete proof for arbitrary~$G$ 
available in the literature, but in the case of orientations this is a known theorem:

\begin{theorem}[{\cite[Cor.\,5.9]{Hesse}}]
	The 2-groupoid of oriented fully extended TQFTs with target~$\B$ 
	is equivalent to the 2-groupoid of homotopy fixed points of the $\SO_2$-action on
	the core of fully dualisable objects in $\B$, 
	\begin{equation}
	\Fun\big(\Bordor,\B\big)\cong \big[\big(\Bfd\big)\!{}^{\times}\big]^{\SO_2}\,.
	\label{eq:oriented-CH}
	\end{equation}
	\label{thm:oriented-CH}
\end{theorem}

The proof in \cite{Hesse} of this uses the presentation of $\Bordor$ 
from \cite{spthesis}, which is not in terms of 2-dualisability data.
We also mention that the equivalence as stated is one of
2-groupoids, but later we will see that this can be extended to
an equivalence of symmetric monoidal 2-groupoids.

\medskip

Let $n\in\Z_{\geqslant 1}$. 
Given a symmetric monoidal 2-category~$\B$, we define a 2-category $\FD{n}(\Bfd)$.
For $n=1$ this reduces to the 2-category in \cite[Thm.\,4.3]{HV} which
is equivalent to the $\SO_2$-homotopy fixed points of $(\Bfd)^\times$. 
Later we will consider the case $n=r$ for $r$-spin TQFTs with $r\geqslant 2$. 
\begin{itemize}
	\item 
	Objects of $\FD{n}(\Bfd)$ are pairs $(\alpha,\theta)$,
	where $\alpha\in\Bfd$, and $\theta\colon S_{\alpha}^n\lra 1_{\alpha}$
	is a 2-isomorphism in $\Bfd$. 
	\item 
	A 1-morphism $(\alpha,\theta) \lra(\alpha',\theta')$ in $\FD{n}(\Bfd)$ 
	is a 1-morphism $X\colon \alpha\lra\alpha'$ in $\Bfd$ such that the following diagram commutes: 
	\begin{equation}
	\begin{tikzpicture}[
	baseline=(current bounding box.base),
	descr/.style={fill=white,inner sep=3.5pt},
	normal line/.style={->}
	]
	\matrix (m) [matrix of math nodes, row sep=1em, column sep=5em, text height=1.5ex, text depth=0.1ex] {%
		X\otimes S^n_{\alpha}&X\otimes 1_{\alpha}&
		\\
		&&X
		\\
		S^n_{\alpha'}\otimes X&1_{\alpha'}\otimes X&
		\\
	};
	\path[font=\footnotesize] (m-1-1) edge[->] node[above,sloped] {$1_{X}\otimes \theta$} (m-1-2);
	\path[font=\footnotesize] (m-3-1) edge[->] node[below,sloped] {$\theta' \otimes 1_{X}$} (m-3-2);
	\path[font=\footnotesize] (m-1-1) edge[->] node[left] {$ S^n_{X} $} (m-3-1);
	\path[font=\footnotesize] (m-1-2) edge[->] node[above,sloped] {$ \cong $} (m-2-3);
	\path[font=\footnotesize] (m-3-2) edge[->] node[below,sloped] {$ \cong $} (m-2-3);
	\end{tikzpicture}
	\label{eq:cond-fixed-pt-morph-B}
	\end{equation}
	\item 
	A 2-morphism $X\lra Y$ in $\FD{n}(\Bfd)$ is a 2-morphism $X\lra Y$ in $\Bfd$. 
	\item 
	Composition and units of $\FD{n}(\Bfd)$ are induced from~$\Bfd$. 
\end{itemize}

To keep the cases $n=1$ and $n\neq 1$ separate, for a given object $(\alpha,\theta) \in \FD{n}(\Bfd)$ we write $\lambda:=\theta$ if $n=1$, and for $n=r\notin \{0,1\}$ we write $\vartheta:=\theta$.

\begin{theorem}[{\cite[Thm.\,4.3]{HV}}]
	There is an equivalence of 2-categories 
	\begin{equation}
	\big[\Bfd\big]^{\SO_2}\cong  \FD{1}(\Bfd)\,.
	\label{eq:SO2-fixed-points}
	\end{equation}
	\label{thm:SO2-fixed-points}
\end{theorem}

\medskip

In the following we will determine a presentation of $\Bordor$ and $\Bordrspin$
in terms of 2-dualisability data. 
The results are collected in Theorems~\ref{thm:iota-1} and~\ref{thm:bord-r-spin-pres}, 
but first we need some preparation.

\begin{lemma}
	Let~$G$ be a topological group, and let $\xi\colon G \lra \GL_2$ be a continuous group homomorphism. 
	\begin{enumerate}[label={(\roman*)}]
		\item Every object in $\Bord_{2,1,0}^G$ is isomorphic
		to a disjoint union of points with trivial $G$-structure.
		\label{lem:bord-fd:part-1}
		\item Every object in $\Bord_{2,1,0}^G$ is fully dualisable.
		\label{lem:bord-fd:part-2}
	\end{enumerate}
	\label{lem:bord-fd}
\end{lemma}
\begin{proof}
	Every connected component~$c$ of the underlying manifold of an object in $\Bord_{2,1,0}^G$ is contractible,
	hence the $G$-structure on~$c$ is trivialisable.
	The mapping cylinder for a trivialisation gives an isomorphism in $\Bord_{2,1,0}^G$.
	This proves part~\ref{lem:bord-fd:part-1}. 
	
	To prove part~\ref{lem:bord-fd:part-2}, consider the commutative diagram of group homomorphisms 
	\begin{equation}
	\begin{tikzpicture}[
	baseline=(current bounding box.base),
	descr/.style={fill=white,inner sep=3.5pt},
	normal line/.style={->}
	]
	\matrix (m) [matrix of math nodes, row sep=2em, column sep=2em, text height=1.5ex, text depth=0.1ex] {%
		\star&&G
		\\
		&\GL_2&
		\\
	};
	\path[font=\footnotesize] (m-1-1) edge[right hook->] node[above,sloped] {incl.} (m-1-3);
	\path[font=\footnotesize] (m-1-1) edge[right hook->] node[above,sloped] {$  $} (m-2-2);
	\path[font=\footnotesize] (m-1-3) edge[->] node[below] {$ \xi $} (m-2-2);
	\end{tikzpicture}
	\label{eq:fr-xi}
	\end{equation}
	from which we get the induced functor $\Lambda_{\text{incl}}$ as in \eqref{eq:induced-functor-lambda}.
	Composing this with the functor in \eqref{eq:iota-0} provides a symmetric monoidal functor
	\begin{equation}
	\begin{tikzpicture}[
	baseline=(current bounding box.base),
	descr/.style={fill=white,inner sep=3.5pt},
	normal line/.style={->}
	]
	\matrix (m) [matrix of math nodes, row sep=1em, column sep=2em, text height=1.5ex, text depth=0.1ex] {%
		\FD{0}&\Bordfr&\Bord_{2,1,0}^G\,.
		\\
	};
	\path[font=\footnotesize] (m-1-1) edge[->] node[above,sloped] {$ \iota^0 $} (m-1-2);
	\path[font=\footnotesize] (m-1-2) edge[->] node[above,sloped] {$ \Lambda_{\text{incl.}} $} (m-1-3);
	\end{tikzpicture}
	\label{eq:iota-lambda-fr-xi}
	\end{equation}
	This composition sends $+\in \FD{0}$ to the haloed point with trivial $G$-structure,
	and symmetric monoidality implies that the image of~$+$ is fully dualisable.
	The claim of part~\ref{lem:bord-fd:part-1} then completes the proof.
\end{proof}

A deck transformation on an $r$-spin surface $(P,q,\Sigma)$ 
is an automorphism of the $r$-spin structure $(P,q)$
which permutes the elements of each fibre of the 
$\Z_r$-bundle $q\colon P\lra F\Sigma$.
We also refer to an $r$-spin bordism as a deck transformation if it is
a mapping cylinder of a deck transformation.
The 1-morphism components $S_{p}$ for $p\in\Bordrspin$ 
of the Serre functor on $\Bordrspin$ are isomorphic to deck transformations (\cite[Rem.\,1.3.1]{dsps1312.7188}):
\begin{equation}
\begin{tikzpicture}
	\begin{pgfonlayer}{nodelayer}
		\node [style=none] (4) at (-9, -2) {$\longhookrightarrow$};
		\node [style=none] (5) at (-11.5, -3.5) {$+$};
		\node [style=none] (50) at (-11, -2) {};
		\node [style=black dot] (51) at (-11.5, -2) {};
		\node [style=none] (52) at (-12, -2) {};
		\node [style=none] (11) at (9, -2) {$\longhookleftarrow$};
		\node [style=none] (21) at (0, -6.5) {$S_{+}$};
		\node [style=none] (23) at (0, 0) {$\cong$};
		\node [style=none] (24) at (11.5, -3.5) {$+$};
		\node [style=none] (27) at (-9, 2) {$\longhookrightarrow$};
		\node [style=none] (28) at (-11.5, 3.5) {$+$};
		\node [style=none] (280) at (-11, 2) {};
		\node [style=black dot] (281) at (-11.5, 2) {};
		\node [style=none] (282) at (-12, 2) {};
		\node [style=none] (31) at (-2.5, 2) {$\longhookleftarrow$};
		\node [style=none] (32) at (-5.75, 3.5) {$1_{+}$};
		\node [style=none] (33) at (0, 3.5) {$+$};
		\node [style=none] (35) at (11.5, 3.5) {$+$};
		\node [style=none] (43) at (1.75, 2) {};
		\node [style=none] (44) at (9.75, 2) {};
		\node [style=none] (45) at (5.75, 2.75) {\scriptsize{deck transformation}};
		\node [style=none] (283) at (11, 2) {};
		\node [style=black dot] (284) at (11.5, 2) {};
		\node [style=none] (285) at (12, 2) {};
		\node [style=none] (295) at (11, -2) {};
		\node [style=black dot] (296) at (11.5, -2) {};
		\node [style=none] (297) at (12, -2) {};
		\node [style=none] (298) at (-0.5, 2) {};
		\node [style=black dot] (299) at (0, 2) {};
		\node [style=none] (300) at (0.5, 2) {};
		\node [style=none] (301) at (-4.5, 2) {};
		\node [style=black dot] (3011) at (-5, 2) {};
		\node [style=none] (302) at (-4.5, 2) {};
		\node [style=none] (303) at (-7, 2) {};
		\node [style=black dot] (3031) at (-6.5, 2) {};
		\node [style=none] (304) at (-7, 2) {};
		\node [style=none] (306) at (7, -2) {};
		\node [style=none] (308) at (-7, -2) {};
		\node [style=none] (37) at (-1.5, -2) {};
		\node [style=none] (38) at (1.5, -2) {};
		\node [style=none] (39) at (-1.5, -3.5) {};
		\node [style=none] (40) at (1.5, -3.5) {};
		\node [style=none] (41) at (-1.5, -5) {};
		\node [style=none] (42) at (1.5, -5) {};
		\node [style=none] (305) at (7, -2) {};
		\node [style=black dot] (3051) at (6.5, -2) {};
		\node [style=none] (307) at (-7, -2) {};
		\node [style=black dot] (3071) at (-6.5, -2) {};
	\end{pgfonlayer}
	\begin{pgfonlayer}{edgelayer}
		\path [style=orange front]
		(305) to (38.center)
		to [in=0, out=-180, looseness=0.75] (39.center)
		to [bend right=90, looseness=1.25]  (41.center)
		to (42.center)
		to [bend right=90, looseness=1.25] (40.center)
		to [in=0, out=180, looseness=0.75] (37.center)
		to [in=360, out=180] (307);
		\draw [style=orange front] (303.center) to (301.center);
		\draw [style=black arrow] (44.center) to (43.center);
		\draw [style=orange front] (50.center) to (52.center);
		\draw [style=blue] (50.center) to (51);
		\draw [style=red] (51) to (52.center);
		\draw [style=orange front] (280.center) to (282.center);
		\draw [style=blue] (280.center) to (281);
		\draw [style=red] (281) to (282.center);
		\draw [style=orange front] (283.center) to (285.center);
		\draw [style=red] (283.center) to (284);
		\draw [style=blue] (284) to (285.center);
		\draw [style=orange front] (295.center) to (297.center);
		\draw [style=red] (295.center) to (296);
		\draw [style=blue] (296) to (297.center);
		\draw [style=orange front] (298.center) to (300.center);
		\draw [style=red] (298.center) to (299);
		\draw [style=blue] (299) to (300.center);
		\draw [style=blue] (3011) to (302.center);
		\draw [style=blue] (3031) to (304.center);
		\draw [style=blue] (3051) to (306.center);
		\draw [style=blue] (3071) to (308.center);
		\draw [style=red mid arrow] (3011) to (3031);
		\draw [style=red mid arrow] (3051) to (38.center);
		\draw [style=red, in=0, out=-180, looseness=0.75] (38.center) to (39.center);
		\draw [style=red mid arrow, bend right=90, looseness=1.25] (39.center) to (41.center);
		\draw [style=red] (41.center) to (42.center);
		\draw [style=red mid arrow, bend right=90, looseness=1.25] (42.center) to (40.center);
		\draw [style=red, in=0, out=180, looseness=0.75] (40.center) to (37.center);
		\draw [style=red mid arrow, in=360, out=180] (37.center) to (3071);
	\end{pgfonlayer}
\end{tikzpicture}
\label{eq:Serre-deck-trf-iso}
\end{equation}
For later use we recall from Example~\ref{exa:ClosedLambdaRFrobeniusAlgebras}\ref{item:ClosedLambdaRAlgebraBordRspin} the relation between deck transformations and Nakayama automorphisms: 
\begin{lemma}
	The Nakayama automorphisms $N_a\colon C_a\lra C_a$ of the closed 
	$\Lambda_r$-Frobenius algebra in $\Bordrspin(\varnothing,\varnothing)$ 
	are deck transformations.
	\label{lem:Nakayama-deck-trf}
\end{lemma}

Another way to express this relation is as follows:
\begin{lemma}
	There are invertible modifications
	\begin{equation}
	\lambda\colon S \stackrel{\cong}{\lra} 1_{\Id_{\Bordor}}
	\quad\text{and}\quad
	\vartheta\colon S^r \stackrel{\cong}{\lra} 1_{\Id_{\Bordrspin}}
	\, . 
	\label{eq:Serre-triv-bord-or-r-spin}
	\end{equation}
	\label{lem:Serre-triv}
\end{lemma}
\begin{proof}
	In $\Bordor$ the 1-morphism components of the Serre automorphism~$S$ are diffeomorphic to the identity, 
	and mapping cylinders of these diffeomorphisms assemble into the modification~$\lambda$.
	In $\Bordrspin$ the $r$-th power of the 1-morphism components of~$S$ 
	are diffeomorphic to the $r$-th power of a deck transformation,
	which in turn is isomorphic to the identity, thus providing~$\vartheta$. 
\end{proof}

This motivates the following definition of a symmetric monoidal 2-category $\FD{n}$ 
via generators and relations, for every $n\in\Z_{\geqslant 1}$. 
The generators of $\FD{n}$ are the objects, 1- and 2-morphisms of $\FD{0}$ (cf.\ Section~\ref{sec:framed-CH}) together with
additional 2-morphisms
\begin{equation}
\theta_{\alpha}\colon S_{\alpha}^n\lra 1_{\alpha}\,, \quad 
\theta_{\alpha}^{-1}\colon 1_{\alpha} \lra S_{\alpha}^n
\label{eq:theta-2D-n}
\end{equation}
for all $\alpha\in\FD{n}$.
The relations of $\FD{n}$ are
\begin{itemize}
	\item 
	the relations of $\FD{0}$, 
	\item $\theta_{\alpha}\circ\theta_{\alpha}^{-1}=1_{1_{\alpha}}$ and 
	$\theta_{\alpha}^{-1}\circ\theta_{\alpha}=1_{S_{\alpha}^n}$,
	\item the commutativity of the diagram 
	\begin{equation}
	\begin{tikzpicture}[
	baseline=(current bounding box.base),
	descr/.style={fill=white,inner sep=3.5pt},
	normal line/.style={->}
	]
	\matrix (m) [matrix of math nodes, row sep=1em, column sep=5em, text height=1.5ex, text depth=0.1ex] {%
		X\otimes S^n_{\alpha}&X\otimes 1_{\alpha}&
		\\
		&&X
		\\
		S^n_{\alpha'}\otimes X&1_{\alpha'}\otimes X&
		\\
	};
	\path[font=\footnotesize] (m-1-1) edge[->] node[above,sloped] {$1_{X}\otimes \theta_{\alpha}$} (m-1-2);
	\path[font=\footnotesize] (m-3-1) edge[->] node[below,sloped] {$\theta_{\alpha'} \otimes 1_{X}$} (m-3-2);
	\path[font=\footnotesize] (m-1-1) edge[->] node[left] {$ S^n_{X} $} (m-3-1);
	\path[font=\footnotesize] (m-1-2) edge[->] node[above,sloped] {$ \cong $} (m-2-3);
	\path[font=\footnotesize] (m-3-2) edge[->] node[below,sloped] {$ \cong $} (m-2-3);
	\end{tikzpicture}
	\label{eq:cond-fixed-pt-morph}
	\end{equation}
	for all $\alpha,\alpha'\in\FD{0}$ and $X\in\FD{0}(\alpha,\alpha')$.
\end{itemize}
We note that the condition in \eqref{eq:cond-fixed-pt-morph} 
expresses the naturality of $S$.
Furthermore the~$\theta_a$ are components of 
an invertible modification $\theta \colon S^n\lra 1_{\Id_{\FD{n}}}$.
For $n=1$ we write $\lambda_{\alpha}:=\theta_{\alpha}$, and 
for $n=r$ we write $\vartheta_{\alpha}:=\theta_{\alpha}$.

Theorems~\ref{thm:oriented-CH}~and~\ref{thm:SO2-fixed-points} together with
the 3-categorical Yoneda lemma \cite[Thm.\,2.12]{BuhnePhD} imply:
\begin{theorem}
	There is an equivalence of symmetric monoidal 2-categories
	\begin{equation}
	\iota^{1}\colon \FD{1} \stackrel{\cong}{\lra} \Bordor\,.
	\label{eq:iota-1}
	\end{equation}
	\label{thm:iota-1}
\end{theorem}
\begin{proof}
	We have a chain of equivalences
	\begin{equation}
	\Fun\big(\Bordor,\B\big) 
	\cong 
	\big[(\Bfd)^{\times}\big]^{\SO_2}
	\cong 
	\FD{1}\big((\Bfd)^{\times}\big)
	\cong 
	\Fun\big(\FD{1},\B\big)
	\, , 
	\label{eq:chain-of-equiv}
	\end{equation}
	all natural in $\B$.
	The first equivalence is from Theorem~\ref{thm:oriented-CH},
	the second is from Theorem~\ref{thm:SO2-fixed-points}. 
	To explain the last equivalence, we will define functors 
	\begin{equation}
	\begin{tikzpicture}[
	baseline=(current bounding box.base),
	descr/.style={fill=white,inner sep=3.5pt},
	normal line/.style={->}
	]
	\matrix (m) [matrix of math nodes, row sep=3em, column sep=3em, text height=1.5ex, text depth=0.1ex] {%
		\Fun\big(\FD{1},\B\big)&\FD{1}\big( (\Bfd)^{\times}\big)
		\\
	};
	\path[font=\footnotesize] (m-1-1) edge[bend left,->] node[above,sloped] {$F$} (m-1-2);
	\path[font=\footnotesize] (m-1-2) edge[bend left,->] node[below,sloped] {$F^{-1}$} (m-1-1);
	\end{tikzpicture}
	\label{eq:Y-Z-diag}
	\end{equation}
	and then show that $F^{-1}$ is in fact a weak inverse of $F$.
	
	The functor $F$ is the evaluation on the generating object $+\in\FD{1}$
	and the corresponding generating 2-morphism $\lambda_+$. In detail:
	\begin{itemize}
		\item for a functor $Y\in\Fun(\FD{1},\B)$ we set 
		$ 
		F(Y):=( Y(+),Y(\lambda_+))
		$, 
		\item for a natural transformation $f\colon Y\lra Y'$ we set 
		\begin{equation}
		F(f):= \Big( f_+\colon \big( Y(+),Y(\lambda_+)\big)\lra\big( Y'(+),Y'(\lambda_+)\big) \Big)\,,
		\label{eq:F-on-1-mor}
		\end{equation}
		\item for a modification $\varphi\colon f\lra f'$ we set
		$
		F(\varphi):=(\varphi_+\colon f_+\lra f'_+ )
		$. 
	\end{itemize}
	We need to show that $F$ indeed lands in $\FD{1}( (\Bfd)^{\times})$,
	i.\,e.\ the diagram~\eqref{eq:cond-fixed-pt-morph-B} commutes for $X=f_+$ and $n=1$.
	The 2-morphism component of $f$ for $U\colon\alpha\lra\alpha'$ is
	$
	f_{U}\colon f_{\alpha'}\otimes Y(U) \lra Y'(U)\otimes f_{\alpha}
	$.
	Substituting $\alpha=\alpha'$ and $U=S_{\alpha}$ we get
	\begin{equation}
	f_{S_{\alpha}}\colon 
	f_{\alpha}\otimes S_{Y(\alpha)}
	\lra 
	S_{Y'(\alpha)}\otimes f_{\alpha}\,,
	\label{eq:f-S-alpha}
	\end{equation}
	where we used monoidality of $Y$ to obtain 
	$Y(S_{\alpha}) \cong S_{Y(\alpha)}$, etc.
	The key observation is that by functoriality and monoidality of $f$
	we have $f_{S_{\alpha}}=S_{f_{\alpha}}$.
	Naturality of $f$ implies that \eqref{eq:cond-fixed-pt-morph-B} then indeed commutes:
	\begin{equation}
	\begin{tikzpicture}[
	baseline=(current bounding box.base),
	descr/.style={fill=white,inner sep=3.5pt},
	normal line/.style={->}
	]
	\matrix (m) [matrix of math nodes, row sep=1em, column sep=5em, text height=1.5ex, text depth=0.1ex] {%
		f_{\alpha}\otimes S_{Y(\alpha)}&f_{\alpha}\otimes 1_{Y(\alpha)}&
		\\
		&&f_{\alpha}
		\\
		S_{Y'(\alpha)}\otimes f_{\alpha}&1_{Y'(\alpha)}\otimes f_{\alpha}&
		\\
	};
	\path[font=\footnotesize] (m-1-1) edge[->] node[above,sloped] {$1_{f_{\alpha}}\otimes Y(\lambda_{\alpha})$} (m-1-2);
	\path[font=\footnotesize] (m-3-1) edge[->] node[below,sloped] {$Y'(\lambda_{\alpha})\otimes 1_{f_{\alpha}}$} (m-3-2);
	\path[font=\footnotesize] (m-1-1) edge[->] node[left] {$ f_{S_{\alpha}}=S_{f_{\alpha}} $} (m-3-1);
	\path[font=\footnotesize] (m-1-2) edge[->] node[above,sloped] {$ \cong $} (m-2-3);
	\path[font=\footnotesize] (m-3-2) edge[->] node[below,sloped] {$ \cong $} (m-2-3);
	\end{tikzpicture}
	\label{eq:f-natural}
	\end{equation}
	
	Now we construct the functor $F^{-1}$.
	Since $\FD{1}$ is defined in terms of generators and relations,
	in order to define a symmetric monoidal functor $Y\colon\FD{1}\lra \B$ it is enough to specify the
	value of $Y$ on the generating objects and morphisms. 
	Then one needs to check that relations in $\FD{1}$ are sent to relations in $\B$.
	The same holds for defining natural transformations and modifications, 
	where we only need to specify their components on generators.
	\begin{itemize}
		\item For $(\alpha,\theta)\in\FD{1}( (\Bfd)^{\times})$ we let
		$F^{-1}(\alpha,\theta)$ be the functor $Y$ with values $Y(+)=\alpha$
		and $Y(\lambda_+)=\theta$.
		This 
		determines the value of $Y$ on all of $\FD{1}$.
		\item For $X\colon(\alpha,\theta)\lra(\alpha',\theta')$ we let
		$F^{-1}(X)$ to be the monoidal natural transformation $f\colon Y\lra Y'$
		with 1-morphism component $f_+ =X$.
		This 
		determines all 1- and 2-morphism components of $f$. 
		\item For $\phi\colon X\lra X'$ we let $F^{-1}(\phi)$ be the modification 
		$\varphi\colon f\lra f'$ with 2-morphism component $\varphi_+=\phi$.
		This determines all 2-morphism components 
		of $\varphi$. 
	\end{itemize}
	We need to check that $F^{-1}$ is well-defined.
	First we note that $Y=F^{-1}(\alpha,\theta)$ is indeed a functor $\FD{1}\lra\B$, 
	because the relations in $\FD{0}$ are satisfied by definition of~$Y$, and 
	so are the relations on~$\lambda$ and its inverse.
	By symmetric monoidality, $Y$ sends $S_+$ in $\FD{1}$ to $S_\alpha$ in 
	$\FD{1}( (\Bfd)^{\times})$, hence naturality is satisfied as well.
	$F^{-1}(X)$ is a natural transformation, since its components satisfy
	\eqref{eq:f-natural}.
	For modifications there are no further conditions to check.
	
	By construction we have $F\circ F^{-1}=\Id$. Moreover, we also have $F^{-1}\circ F\cong \Id$,
	since two functors out of $\FD{1}$ are isomorphic if they agree on generators
	(\cite[Thm.\,2.78]{spthesis}). 
	This shows that~$F$ is an equivalence. 
	
	Finally we observe that the functor $\iota^1$ respects the
	symmetric monoidal structures on $\FD{1}$ and $\Bordor$
	and hence it can canonically be promoted to a symmetric monoidal functor. 
	Thus the claim follows from the 3-categorical Yoneda lemma of \cite[Thm.\,2.12]{BuhnePhD}. 
\end{proof}

By looking at the chain of equivalences in \eqref{eq:chain-of-equiv} we can read off
the value of~$\iota^1$ on generators of the 2-category $\FD{1}$: 
\begin{itemize}
	\item $\iota^1$ sends $+\in\FD{1}$ together with its 2-dualisability data
	to the positively oriented point with its 2-dualisability data
	(cf.\ Lemma~\ref{lem:bord-fd}),
	\item $\iota^1(\lambda_{+})=\lambda_{+}$ 
	from Lemma~\ref{lem:Serre-triv}.
\end{itemize}

\begin{remark}
	Using the symmetric monoidal equivalence $\iota^1$ one also obtains a 
	symmetric monoidal structure on the equivalence in
	Theorem~\ref{thm:oriented-CH}.
	\label{rem:transport-sm-str-to-2D1}
\end{remark}

\medskip

In line with the above presentation of~$\iota^1$, we now define a symmetric monoidal functor 
\begin{equation}
\iota^{r}\colon \FD{r} \lra \Bordrspin
\label{eq:iota-r}
\end{equation}
on generators of $\FD{r}$ for $r\in\Z_{\geqslant 2}$ as follows:
\begin{itemize}
	\item $\iota^r$ sends 
	$+\in\FD{r}$, together with its 2-dualisability data,
	to the point with trivial $r$-spin structure in $\Bordrspin$, 
	with its 2-dualisability data as in Lemma~\ref{lem:bord-fd},
	\item $\iota^r(\vartheta_{+})=\vartheta_{+}$
	from Lemma~\ref{lem:Serre-triv}.
\end{itemize}

\begin{proposition}
	\label{prop:iota-lambda-diag}
	We have a strictly commutative diagram of symmetric monoidal functors 
	\begin{equation}
	\begin{tikzpicture}[
	baseline=(current bounding box.base),
	descr/.style={fill=white,inner sep=3.5pt},
	normal line/.style={->}
	]
	\matrix (m) [matrix of math nodes, row sep=3em, column sep=3em, text height=1.5ex, text depth=0.1ex] {%
		\Bordfr &\Bordrspin &\Bordor
		\\
		\FD{0} &
		\FD{r} &
		\FD{1} 
		\\
	};
	\path[font=\footnotesize] (m-1-1) edge[->] node[above,sloped] {$ \widetilde{\Lambda} $} (m-1-2);
	\path[font=\footnotesize] (m-1-2) edge[->] node[above,sloped] {$ {\Lambda} $} (m-1-3);
	\path[font=\footnotesize] (m-2-1) edge[->] node[below,sloped] {$ \widetilde{K} $} (m-2-2);
	\path[font=\footnotesize] (m-2-2) edge[->] node[below,sloped] {$ {K} \vphantom{\widetilde{K}}$} (m-2-3);
	\path[font=\footnotesize] (m-2-1) edge[->] node[left] {$ \iota^0 $} (m-1-1);
	\path[font=\footnotesize] (m-2-2) edge[->] node[right] {$ \iota^r $} (m-1-2);
	\path[font=\footnotesize] (m-2-3) edge[->] node[right] {$ \iota^1 $} (m-1-3);
	\end{tikzpicture}
	\label{eq:iota-lambda-diag}
	\end{equation}
	where the functors $\widetilde{K}$ and $K$ by definition 
	act as the identity on objects as well as on 1- and  2-morphism generators of $\FD{0}$, 
	while on the other generators we have 
	\begin{equation}
	\begin{tikzpicture}[
	baseline=(current bounding box.base), 
	descr/.style={fill=white,inner sep=3.5pt}, 
	normal line/.style={->}, 
	baseline=-0.1cm
	] 
	\matrix (m) [matrix of math nodes, row sep=3em, column sep=4.0em, text height=1.5ex, text depth=0.1ex] {%
		K\Big( S_{+}^r  & 1_{+} \Big)
		\\
	};
	\path[font=\footnotesize, transform canvas={yshift=+0.8mm}] (m-1-1) edge[->] node[above] { $ \vartheta_{+} $ } (m-1-2);
	\path[font=\footnotesize, transform canvas={yshift=-0.8mm}] (m-1-2) edge[->] node[below] { $ \vartheta_{+}^{-1} $} (m-1-1);
	\end{tikzpicture}
	\!\!=\!\! 
	\begin{tikzpicture}[
	baseline=(current bounding box.base), 
	descr/.style={fill=white,inner sep=3.5pt}, 
	normal line/.style={->}, 
	baseline=-0.1cm
	] 
	\matrix (m) [matrix of math nodes, row sep=3em, column sep=4.0em, text height=1.5ex, text depth=0.1ex] {%
		\Big( S_{+}^r  & 1_+^r & 1_{+} \Big) .
		\\
	};
	\path[font=\footnotesize, transform canvas={yshift=+0.8mm}] (m-1-1) edge[->] node[above] { $ \lambda_{+}^r $ } (m-1-2);
	\path[font=\footnotesize, transform canvas={yshift=-0.8mm}] (m-1-2) edge[->] node[below] { $ \lambda_{+}^{-r} $ } (m-1-1);
	\path[font=\footnotesize, transform canvas={yshift=+0.8mm}] (m-1-2) edge[->] node[above] { $ \cong $ } (m-1-3);
	\path[font=\footnotesize, transform canvas={yshift=-0.8mm}] (m-1-3) edge[->] node[below] { $ \cong $ } (m-1-2);
	\end{tikzpicture}
	\label{eq:K-theta}
	\end{equation}
\end{proposition}

By Theorems~\ref{thm:bord-framed-pres} and~\ref{thm:iota-1}, the functors~$\iota^0$ and~$\iota^1$ are equivalences.
In Section~\ref{subsubsec:ProofrSpinCH} below we will prove:
\begin{theorem}
	The functor in \eqref{eq:iota-r} is an equivalence for all $r\in\Z_{\geqslant 1}$, 
	\begin{equation}
	\iota^{r}\colon \FD{r} \stackrel{\cong}{\lra} \Bordrspin \, . 
	\end{equation}
	\label{thm:bord-r-spin-pres}
\end{theorem}

We also have the analogous statement of Theorem~\ref{thm:SO2-fixed-points}, 
which follows immediately by applying \cite[Thm.\,4.3]{HV} to the natural transformation~$S^r$:
\begin{lemma}
	The homotopy fixed points of the $r$-spin action on~$\Bfd$ are 
	\begin{equation}
	\big[\Bfd\big]^{\Spin_2^r} \cong  \FD{r}(\Bfd) \,.
	\label{eq:r-spin-fixed-points}
	\end{equation}
	\label{lem:r-spin-fixed-points}
\end{lemma}

Theorem~\ref{thm:bord-r-spin-pres} and Lemma~\ref{lem:r-spin-fixed-points}
imply the 2-categorical cobordism hypothesis with $r$-spin structure:

\begin{theorem}
	The 2-groupoid of fully extended $r$-spin TQFTs with target~$\B$ is 
	equivalent to the homotopy fixed points of the $r$-spin action, 
	\begin{equation}
	\Fun\big(\Bordrspin,\B\big)
	\cong
	\big[ (\Bfd)^{\times}\big]^{\Spin_2^r} \, . 
	\label{eq:r-spin-CH}
	\end{equation}
	\label{thm:r-spin-CH}
\end{theorem}
\begin{proof}
	We have a chain of equivalences
	\begin{equation}
	\Fun\big(\Bordrspin,\B\big)
	\cong
	\Fun\big(\FD{r},\B\big)
	\cong
	\FD{r}\big((\Bfd)^{\times}\big)
	\cong
	\big[ (\Bfd)^{\times} \big]^{\Spin_2^r}\,.
	\label{eq:pf-thm-r-spin-CH-2}
	\end{equation}
	The first equivalence is from Theorem~\ref{thm:bord-r-spin-pres},
	the last equivalence is Lemma~\ref{lem:r-spin-fixed-points}, and the 
	proof of the second equivalence is completely analogous to the $n=1$ case in
	the proof of Theorem~\ref{thm:iota-1}.
\end{proof}

\begin{remark}
	The proof of the oriented cobordism hypothesis in \cite{Hesse} 
	(Theorem~\ref{thm:oriented-CH}) 
	uses the presentation of $\Bordor$ of \cite{spthesis}, which is not
	in terms of 2-dualisability data, and a direct computation of the
	$\SO_2$-homotopy fixed points 
	$[(\Bfd)^{\times}]^{\SO_2}$ (Theorem~\ref{thm:SO2-fixed-points}).
	In order to prove the $r$-spin cobordism hypothesis (Theorem~\ref{thm:r-spin-CH})
	we need a presentation of $\Bordrspin$ (Theorem~\ref{thm:bord-r-spin-pres})
	in terms of 2-dualisability data, 
	and a direct computation of $\Spin_2^r$-homotopy fixed points 
	$[ (\Bfd)^{\times} ]^{\Spin_2^r}$ (Lemma~\ref{lem:r-spin-fixed-points}).
	\label{rem:or-vs-r-spin-CH}
\end{remark}

It is straightforward to check the following factorisation of $r$-spin TQFTs:
\begin{proposition}
	Let $\alpha\in\Bfd$ and $k<r$ be such that there are 2-isomorphisms
	$\varphi\colon S_{\alpha}^k\lra 1_{\alpha}$
	and $\psi\colon S_{\alpha}^r\lra 1_{\alpha}$.
	Then
	\begin{enumerate}
		\item there is a 2-isomorphism $\phi\colon S_{\alpha}^g\lra 1_{\alpha}$,
		where $g=\mathrm{gcd}(k,r)$, and
		\item the diagram
		\begin{equation}
		\begin{tikzpicture}[
		baseline=(current bounding box.base),
		descr/.style={fill=white,inner sep=3.5pt},
		normal line/.style={->}
		]
		\matrix (m) [matrix of math nodes, row sep=3em, column sep=3em, text height=1.5ex, text depth=0.1ex] {%
			\Bordrspin
			&
			&
			\Bord_{2,1,0}^{g\text{-spin}}
			\\
			&
			\B 
			&
			\\
		};
		\path[font=\footnotesize] (m-1-1) edge[->] node[above,sloped] {$ {\Lambda_{r,g}} $} (m-1-3);
		\path[font=\footnotesize] (m-1-1) edge[->] node[below] {$ {\zz_{\psi}} $} (m-2-2);
		\path[font=\footnotesize] (m-1-3) edge[->] node[below] {$ \phantom{ii} \zz_{\phi}  $} (m-2-2);
		\end{tikzpicture}
		\label{eq:r-spin-TQFT-factorise}
		\end{equation}
		commutes up to a natural isomorphism, where $\Lambda_{r,g}$ is the functor 
		from \eqref{eq:induced-functor-lambda} for the group homomorphism 
		$\widetilde{\mathrm{GL}_2^+}^r\lra\widetilde{\mathrm{GL}_2^+}^g$, 
		while~$\zz_{\psi}$ and~$\zz_{\phi}$ denote the $r$-spin and $g$-spin TQFTs from 
		Theorem~\ref{thm:r-spin-CH} corresponding to $(\alpha,\psi)\in\FD{r}(\Bfd)$
		and $(\alpha,\phi)\in\FD{g}(\Bfd)$, respectively.
	\end{enumerate}
	\label{prop:r-spin-TQFT-factorise}
\end{proposition}
\begin{remark}
	Let us assume that the adjoints of 1-morphisms in $\B$ satisfy $\Xd=\dX$,
	which is for example the case when~$\B$ is pivotal. 
	Then by the definition of the Serre automorphism \eqref{eq:def-serre} 
	and its inverse \eqref{eq:def-serre-inv}, we have $S=S^{-1}$.
	Hence under this assumption, 
	$r$-spin TQFTs with target $\B$ factorise through oriented TQFTs
	($r$ odd), 
	or through 2-spin TQFTs
	($r$ even).
	\label{rem:pivotal-boring-ext}
\end{remark}

\subsubsection[Proof of the \texorpdfstring{$r$}{r}-spin cobordism hypothesis]{Proof of the \texorpdfstring{$\boldsymbol{r}$}{r}-spin cobordism hypothesis}
\label{subsubsec:ProofrSpinCH}

Here we prove Theorem~\ref{thm:bord-r-spin-pres}.
To do this, we will check the conditions listed in the following Whitehead-type theorem for the functor $\iota^r$ in~\eqref{eq:iota-r}.
\begin{theorem}[{\cite[Thm.\,2.25]{spthesis}}]
	A functor $F\colon \Cc\lra\B$ of 2-categories 
	is an equivalence if and only if it is 
	essentially surjective on objects and 1-morphisms, and
	fully faithful on 2-morphisms.
	\label{thm:equivalence-conditions}
\end{theorem}

\begin{lemma}
	The functor $\iota^r$ is essentially surjective on objects.
	\label{lem:iota-r-ess-surj-obj}
\end{lemma}
\begin{proof}
	This follows from Lemma~\ref{lem:bord-fd}\ref{lem:bord-fd:part-1}. 
\end{proof}

\begin{lemma}
	The functor $\iota^r$ is essentially surjective on 1-morphisms.
	\label{lem:iota-r-ess-surj-1-mor}
\end{lemma}
\begin{proof}
	For every connected component~$c$ of a 1-morphism $P\lra FX$ in $\Bordrspin$, we obtain an element 
	$\delta(c)\in\Z_r$ as follows.
	If the 1-manifold of which~$c$ is a halo is closed, set $\bar{c}:=c$, otherwise let $\bar{c}$ be the 
	$r$-spin surface with closed embedded 1-manifold defined by 
	identifying the two boundary points of the embedded 1-manifold 
	in $c$ via the boundary parametrisation maps.
	This identification is possible by choosing a trivialisation of the 
	$r$-spin structures of the objects parametrising the boundary of $c$.
	Consider a curve $\Gamma\colon S^1 \lra \bar{c}$ parametrising the 1-manifold 
	in $\bar{c}$ and its lift $\widetilde{\Gamma}\colon S^1 \lra F\bar{c}$ to the frame bundle
	defined by picking at every point a tangent vector to $\Gamma$
	and another vector so that the induced orientation agrees with the 
	orientation underlying the $r$-spin structure of $c$.
	This lift is unique up to homotopy.
	There is a unique lift $\widehat{\Gamma}\colon S^1 \lra P|_{\bar{c}}$ of $\widetilde{\Gamma}$ 
	after fixing it at one point, as the fibres are discrete.
	We define $\delta(c)\in \Z_r$ to be the holonomy of $\widehat{\Gamma}$,
	which only depends on $c$; 
	for more details of this construction we refer to \cite{RandalWilliams} or \cite[Sect.\,5.2]{RunkelSzegedyArf}. 
	
	Recall that by \eqref{eq:Serre-deck-trf-iso}, 
	$S_{+}$ is isomorphic to a deck transformation with holonomy $-1\in\Z_r$.
	If $c=\bar{c}$, then 
	\begin{equation}
	c\cong C_{\delta(c)}\in\Bordrspin(\varnothing,\varnothing)
	\label{eq:c-iso-C_a}
	\end{equation}
	from \eqref{eq:Ca}.
	Otherwise $c$ is isomorphic either to the endomorphism $S_{\pm}^{- \delta(c)}$, or to $S_{+}^{-\delta(c)}$ pre- or post-composed with one of the adjunction 1-morphisms of~$+$, for example $c \cong \tev_{+}\circ(S_{+}^{-\delta(c)}\sqcup 1_{-})$. 
\end{proof}

\begin{lemma}
	The functor $\iota^r$ is full on 2-morphisms.
	\label{lem:iota-r-full-2-mor}
\end{lemma}
\begin{proof}
	Let $X,X'\colon \alpha\lra\alpha'$ be parallel 1-morphisms in $\FD{r}$, and let
	$\Sigma\colon \iota^r(X)\lra \iota^r(X')$ be a 2-morphism in $\Bordrspin$.
	Without loss of generality we can assume that~$\Sigma$ is connected of genus~$g$. 
	We write $\Lambda(\Sigma)$ for the oriented surface underlying~$\Sigma$. 
	Hence the $r$-spin structure on~$\Sigma$ is represented by a bundle $P\lra \Lambda(\Sigma)$ and a $\Z_r$-bundle map $q\colon P\lra F\Lambda(\Sigma)$. 
	
	The strategy of our proof is as follows. 
	We describe the $r$-spin structure on~$\Sigma$ 
	up to diffeomorphisms of $r$-spin surfaces with 
	underlying diffeomorphism of surfaces the identity, 
	which we refer to as \textsl{isomorphisms of $r$-spin structures}.
	Then we consider a decomposition of the oriented surface $\Lambda(\Sigma)$
	suitable for our description of $r$-spin structures.
	Finally we lift the oriented generators to $r$-spin generators and
	restore the $r$-spin structure up to isomorphism in the above sense.
	Therefore the $r$-spin surface we build from the generators
	is in particular
	diffeomorphic to~$\Sigma$, thus representing the same 2-morphism in $\Bordrspin$. 
	
	\textsl{Step 1:}
	Following \cite{RandalWilliams}, we describe the $r$-spin structure of~$\Sigma$ in terms of holonomies along curves in the underlying oriented surface $\Lambda(\Sigma)$ in the $\Z_r$-bundle $q\colon P\lra F\Lambda(\Sigma)$. 
	
	\begin{figure}[tb]
		\centering
		\begin{tikzpicture}
	\begin{pgfonlayer}{nodelayer}
		\node [style=none] (64) at (8, 2.5) {};
		\node [style=none] (65) at (6.75, 0.25) {};
		\node [style=none] (66) at (4.75, 0) {};
		\node [style=none] (67) at (8.25, 0) {};
		\node [style=none] (130) at (8.5, 1) {$a_3$};
		\node [style=none] (131) at (6.5, -1.75) {$b_3$};
		\node [style=none] (183) at (5.5, 0.25) {};
		\node [style=none] (184) at (7.5, 0.25) {};
		\node [style=none] (187) at (6.5, 0.25) {};
		\node [style=none] (188) at (-9.5, 0) {};
		\node [style=none] (189) at (-6, 3) {};
		\node [style=none] (190) at (-3.5, 2.5) {};
		\node [style=none] (191) at (3.5, 2.5) {};
		\node [style=none] (192) at (6, 3) {};
		\node [style=none] (193) at (9.5, 0) {};
		\node [style=none] (194) at (6, -3) {};
		\node [style=none] (195) at (3.5, -2.5) {};
		\node [style=none] (196) at (-3.5, -2.5) {};
		\node [style=none] (197) at (-6, -3) {};
		\node [style=none] (210) at (0, -3) {};
		\node [style=none] (211) at (0, 3) {};
		\node [style=none] (213) at (-1.75, 0) {};
		\node [style=none] (214) at (1.75, 0) {};
		\node [style=none] (215) at (1.25, 1.5) {$a_2$};
		\node [style=none] (216) at (0, -1.75) {$b_2$};
		\node [style=none] (217) at (-1, 0.25) {};
		\node [style=none] (218) at (1, 0.25) {};
		\node [style=none] (221) at (0, 0.25) {};
		\node [style=none] (222) at (-8, 2.5) {};
		\node [style=none] (223) at (-6.25, 0) {};
		\node [style=none] (224) at (-7.75, -0.25) {};
		\node [style=none] (225) at (-4.25, -0.25) {};
		\node [style=none] (226) at (-6, 2) {$a_1$};
		\node [style=none] (227) at (-6, -2) {$b_1$};
		\node [style=none] (228) at (-7, 0) {};
		\node [style=none] (229) at (-5, 0) {};
		\node [style=none] (232) at (-6, 0) {};
		\node [style=none] (236) at (6, 0) {};
		\node [style=none] (237) at (7, 0) {};
		\node [style=none] (238) at (-0.5, 0) {};
		\node [style=none] (239) at (0.5, 0) {};
		\node [style=none] (240) at (-6.5, -0.25) {};
		\node [style=none] (241) at (-5.5, -0.25) {};
	\end{pgfonlayer}
	\begin{pgfonlayer}{background}
		\path [fill=orange!20!white, opacity=0.8]
			(188.center) 
			to [in=180, out=90]  (189.center)
			to [in=180, out=0] (190.center)
			to [in=180, out=0] (211.center)
			to [in=180, out=0] (191.center)
			to [in=180, out=0] (192.center)
			to [in=90, out=0] (193.center)
			to [in=0, out=-90] (194.center)
			to [in=0, out=180] (195.center)
			to [in=0, out=180] (210.center)
			to [in=0, out=180] (196.center)
			to [in=0, out=180] (197.center)
			to [in=270, out=180] (188.center);
		\path [fill=white, opacity=1]
			(240.center) to [bend left=45, looseness=1.50] (241.center)
			to [bend left=45] (240.center);
		\path [fill=white, opacity=1]
			(238.center) to [bend left=45, looseness=1.50] (239.center)
			to [bend left=45] (238.center);
		\path [fill=white, opacity=1]
			(236.center) to [bend left=45, looseness=1.50] (237.center)
			to [bend left=45] (236.center);
	\end{pgfonlayer}
	\begin{pgfonlayer}{edgelayer}
		\draw [style=blue mid arrow, bend left=90, looseness=0.75] (65.center) to (64.center);
		\draw [style=blue bg, bend right=90, looseness=0.75] (65.center) to (64.center);
		\draw [style=blue mid arrow, bend right=90] (66.center) to (67.center);
		\draw [style=blue, bend left=90] (66.center) to (67.center);
		\draw [style=black, bend right, looseness=1.50] (183.center) to (184.center);
		\draw [style=blue mid arrow, bend right=90] (213.center) to (214.center);
		\draw [style=blue, bend left=90] (213.center) to (214.center);
		\draw [style=black, bend right, looseness=1.50] (217.center) to (218.center);
		\draw [style=blue mid arrow, bend left=90, looseness=0.75] (223.center) to (222.center);
		\draw [style=blue bg, bend right=90, looseness=0.75] (223.center) to (222.center);
		\draw [style=blue mid arrow, bend right=90] (224.center) to (225.center);
		\draw [style=blue, bend left=90] (224.center) to (225.center);
		\draw [style=black, bend right, looseness=1.50] (228.center) to (229.center);
		\draw [style=black, bend right=45] (236.center) to (237.center);
		\draw [style=black, bend right=45] (238.center) to (239.center);
		\draw [style=black, bend right=45] (240.center) to (241.center);
		\draw [style=black, bend left=45, looseness=1.50] (240.center) to (241.center);
		\draw [style=blue mid arrow, bend left=90, looseness=0.75] (221.center) to (211.center);
		\draw [style=blue bg, bend right=90, looseness=0.75] (221.center) to (211.center);
		\draw [style=black, bend left=45, looseness=1.50] (238.center) to (239.center);
		\draw [style=black, bend left=45, looseness=1.50] (236.center) to (237.center);
		\draw [style=black, in=180, out=90] (188.center) to (189.center);
		\draw [style=black, in=180, out=0] (189.center) to (190.center);
		\draw [style=black, in=180, out=0] (190.center) to (211.center);
		\draw [style=black, in=180, out=0] (211.center) to (191.center);
		\draw [style=black, in=180, out=0] (191.center) to (192.center);
		\draw [style=black, in=90, out=0] (192.center) to (193.center);
		\draw [style=black, in=0, out=-90] (193.center) to (194.center);
		\draw [style=black, in=0, out=180] (194.center) to (195.center);
		\draw [style=black, in=180, out=0] (210.center) to (195.center);
		\draw [style=black, in=180, out=0] (196.center) to (210.center);
		\draw [style=black, in=0, out=180] (196.center) to (197.center);
		\draw [style=black, in=270, out=180] (197.center) to (188.center);
	\end{pgfonlayer}
\end{tikzpicture}
		\caption{Curves on a closed surface.}
		\label{fig:surface-curves}
	\end{figure}
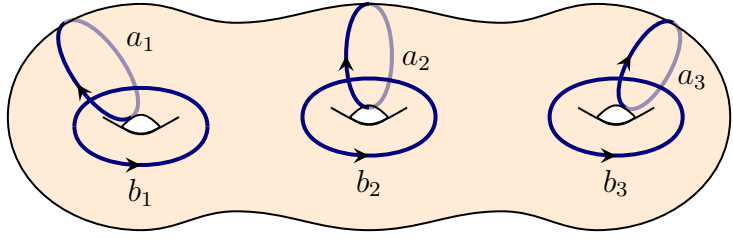
	
	\textsl{Step 1.1:}
	If $\Lambda(\Sigma)$ is a closed surface 
	(i.\,e.\ if $\alpha=\alpha'=\varnothing$ and $X=1_{\varnothing}=X'$), 
	then $r$-spin structures up to isomorphism on $\Lambda(\Sigma)$ are in (non-canonical) bijection with
	$H^1(\Lambda(\Sigma),\Z_r)\cong \Z_r^{2g}$. 
	The latter bijection is given by picking simple closed curves in
	$\Lambda(\Sigma)$ which represent a basis of $H_1(\Lambda(\Sigma))$. 
	For each handle we choose two curves $a_k,b_k$  
	which intersect at precisely one point and 
	which do not intersect with the curves associated to the other handles,
	see Figure~\ref{fig:surface-curves}.
	
	\textsl{Step 1.2:}
	If $\Lambda(\Sigma)$ is not a closed surface we introduce a new surface 
	$\widetilde{\Sigma}$ and additional curves on $\Lambda(\Sigma)$.
	
	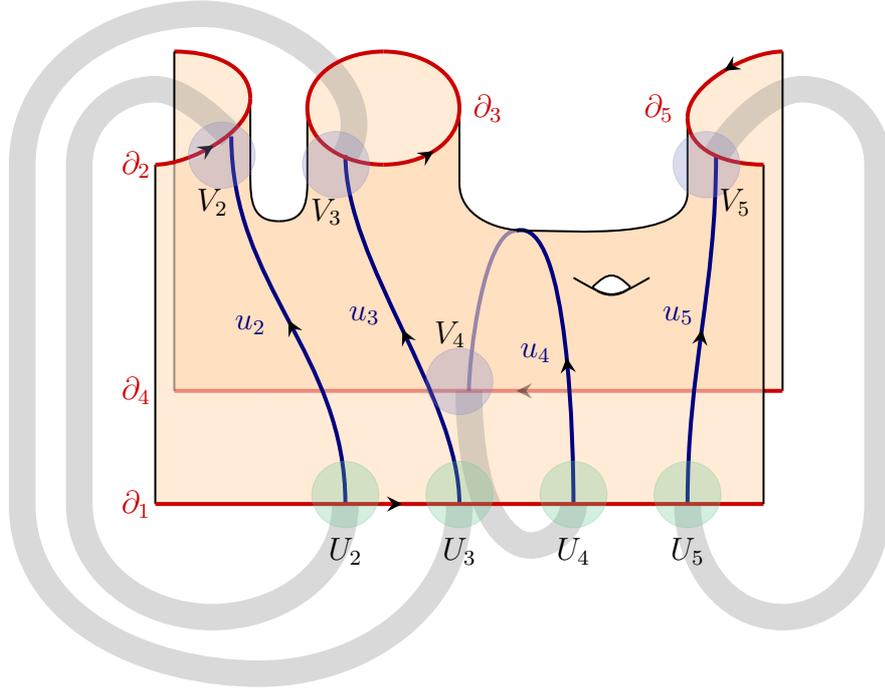
\begin{figure}[tb]
		\centering
		\begin{tikzpicture}
	\begin{pgfonlayer}{nodelayer}
		\node [style=none] (0) at (-6, -2) {};
		\node [style=none] (1) at (10, -2) {};
		\node [style=none] (2) at (-5.5, 1) {};
		\node [style=none] (3) at (10.5, 1) {};
		\node [style=none] (7) at (10.5, 8) {};
		\node [style=none] (8) at (-3.5, 6.75) {};
		\node [style=none] (9) at (8, 6.25) {};
		\node [style=none] (12) at (0, 5) {};
		\node [style=none] (13) at (0, 8) {};
		\node [style=none] (14) at (2, 6.5) {};
		\node [style=none] (15) at (-2, 6.5) {};
		\node [style=none] (17) at (-2.75, 5.5) {};
		\node [style=none] (18) at (3.75, 5.25) {};
		\node [style=none] (19) at (5, 4) {};
		\node [style=none] (20) at (7, 4) {};
		\node [style=none] (23) at (8.75, 7.25) {};
		\node [style=none] (24) at (-1, 7.25) {};
		\node [style=none] (25) at (-4, 7.75) {};
		\node [style=none] (27) at (2.25, 1) {};
		\node [style=none] (29) at (-1, -2) {};
		\node [style=none] (30) at (2, -2) {};
		\node [style=none] (31) at (5, -2) {};
		\node [style=none] (32) at (8, -2) {};
		\node [style=none] (33) at (3.5, 5.25) {};
		\node [style=none] (34) at (10, 7) {};
		\node [style=none] (35) at (10.5, 10) {};
		\node [style=none] (36) at (8, 8.25) {};
		\node [style=none] (38) at (0, 7) {};
		\node [style=none] (39) at (0, 10) {};
		\node [style=none] (40) at (2, 8.5) {};
		\node [style=none] (41) at (-2, 8.5) {};
		\node [style=none] (43) at (-6, 7) {};
		\node [style=none] (44) at (-5.5, 10) {};
		\node [style=none] (45) at (-3.5, 8.75) {};
		\node [style=none, text={rgb,255: red,204; green,0; blue,0}] (51) at (-6.5, -2) {$\partial_1$};
		\node [style=none, text={rgb,255: red,204; green,0; blue,0}] (52) at (-6.5, 1) {$\partial_4$};
		\node [style=none, text={rgb,255: red,204; green,0; blue,0}] (53) at (-6.5, 7) {$\partial_2$};
		\node [style=none, text={rgb,255: red,204; green,0; blue,0}] (54) at (7.25, 8.5) {$\partial_5$};
		\node [style=none, text={rgb,255: red,204; green,0; blue,0}] (55) at (2.75, 8.5) {$\partial_3$};
		\node [style=none, text={rgb,255: red,0; green,0; blue,128}] (56) at (-3.5, 2.75) {$u_2$};
		\node [style=none, text={rgb,255: red,0; green,0; blue,128}] (57) at (-0.5, 3) {$u_3$};
		\node [style=none, text={rgb,255: red,0; green,0; blue,128}] (58) at (4, 2) {$u_4$};
		\node [style=none, text={rgb,255: red,0; green,0; blue,128}] (59) at (7.75, 3) {$u_5$};
		\node [style=blue cloud] (60) at (-1.25, 7) {};
		\node [style=blue cloud] (61) at (8.5, 7) {};
		\node [style=blue cloud] (62) at (-4.25, 7.25) {};
		\node [style=blue cloud] (63) at (2, 1.25) {};
		\node [style=green cloud] (64) at (-1, -1.75) {};
		\node [style=green cloud] (65) at (2, -1.75) {};
		\node [style=green cloud] (66) at (5, -1.75) {};
		\node [style=green cloud] (67) at (8, -1.75) {};
		\node [style=none] (68) at (-1, -3.25) {$U_2$};
		\node [style=none] (69) at (2, -3.25) {$U_3$};
		\node [style=none] (70) at (5, -3.25) {$U_4$};
		\node [style=none] (71) at (8, -3.25) {$U_5$};
		\node [style=none] (72) at (-4.5, 6) {$V_2$};
		\node [style=none] (73) at (-1.5, 5.75) {$V_3$};
		\node [style=none] (74) at (1.75, 2.5) {$V_4$};
		\node [style=none] (75) at (9.25, 6) {$V_5$};
		\node [style=none] (76) at (-4.5, -5) {};
		\node [style=none] (77) at (-8, -2) {};
		\node [style=none] (78) at (-8, 7) {};
		\node [style=none] (79) at (-6, 9) {};
		\node [style=none] (80) at (-3.25, -6.5) {};
		\node [style=none] (81) at (-9.5, -2) {};
		\node [style=none] (82) at (-9.5, 7) {};
		\node [style=none] (83) at (-4.75, 11) {};
		\node [style=none] (84) at (3.5, -3) {};
		\node [style=none] (85) at (10.5, -5) {};
		\node [style=none] (86) at (13, -2) {};
		\node [style=none] (87) at (13, 7) {};
		\node [style=none] (88) at (11, 9) {};
		\node [style=none] (93) at (10, 1) {};
		\node [style=none] (94) at (5.5, 3.75) {};
		\node [style=none] (95) at (6.5, 3.75) {};
		\node [style=none] (96) at (-5.5, 7.01) {};
		\node [style=none] (97) at (7, -2) {};
		\node [style=none] (98) at (-3, 1) {};
	\end{pgfonlayer}
	\begin{pgfonlayer}{background}
		\path [fill=orange!20!white, opacity=0.8]
			(45.center) 
			to [in=0, out=-90, looseness=0.75] (96.center)
			to (44.center)
			to [in=90, out=0] (45.center);
		\path [fill=orange!20!white, opacity=0.8]
			(38.center)
			to [in=-90, out=-180] (41.center)
			to [in=-180, out=90] (39.center)
			to [in=90, out=0] (40.center)
			to [in=0, out=-90] (38.center);
		\path [fill=orange!20!white, opacity=0.8]
			(34.center)
			to [in=-90, out=180] (36.center)
			to [in=180, out=90, looseness=0.75] (35.center)
			to (3.center)
			to (93.center)
			to (34.center);
		\path [fill=orange!30!white, opacity=0.8]
			(45.center) 
			to [in=0, out=-90, looseness=0.75] (96.center)
			to (2.center)
			to (93.center)
			to (34.center)
			to [in=-90, out=180] (36.center)
			to (9.center)
			to [in=0, out=-90, looseness=0.75] (18.center)
			to [in=-90, out=180] (14.center)
			to (40.center)
			to [in=0, out=-90] (38.center)
			to [in=-90, out=-180] (41.center)
			to (15.center)
			to [in=0, out=-90, looseness=1.25] (17.center)
			to [in=-90, out=-180, looseness=1.25] (8.center)
			to (45.center);
		\path [fill=orange!20!white, opacity=0.8]
			(43.center)
			to (0.center)
			to (1.center)
			to (93.center)
			to (2.center)
			to (96.center)
			to (43.center);
		\path [fill=white, opacity=1]
			(94.center) to [bend right=45] (95.center)
			to [bend right=45, looseness=1.50] (94.center);
	\end{pgfonlayer}
	\begin{pgfonlayer}{edgelayer}
		\path [style=grey] 
			(29.center) 
			to [in=0, out=-90] (76.center)
			to [in=270, out=180] (77.center)
			to (78.center)
			to [in=180, out=90] (79.center)
			to [in=135, out=0, looseness=0.75] (25.center);
		\path [style=grey] 
			(24.center) 
			to [in=360, out=60] (83.center)
			to [in=90, out=-180] (82.center)
			to (81.center)
			to [in=180, out=-90] (80.center)
			to [in=270, out=0]  (30.center);
		\path [style=grey] 
			(31.center)
			to [in=-30, out=270] (84.center) 
			to [in=-90, out=150, looseness=0.75] (27.center);
		\path [style=grey] 
			(23.center) 
			to [in=-180, out=60, looseness=0.75] (88.center)
			to [in=90, out=0] (87.center)
			to (86.center)
			to [in=0, out=-90]  (85.center)
			to [in=270, out=180] (32.center);
		\draw [style=black] (7.center) to (3.center);
		\draw [style=black, bend right, looseness=1.50] (19.center) to (20.center);
		\draw [style=black] (35.center) to (7.center);
		\draw [style=red mid arrow, in=90, out=180, looseness=0.75] (35.center) to (36.center);
		\draw [style=red, in=90, out=0] (39.center) to (40.center);
		\draw [style=red, in=-180, out=90] (41.center) to (39.center);
		\draw [style=red, in=90, out=0] (44.center) to (45.center);
		\draw [style=blue bg, in=-180, out=90, looseness=0.50] (27.center) to (33.center);
		\draw [style=red mid arrow bg] (93.center) to (98.center);
		\draw [style=red bg] (98.center) to (2.center);
		\draw [style=red] (3.center) to (93.center);
		\draw [style=black, bend right=45] (94.center) to (95.center);
		\draw [style=black, bend left=45, looseness=1.50] (94.center) to (95.center);
		\draw [style=blue mid arrow, in=-90, out=90, looseness=0.75] (30.center) to (24.center);
		\draw [style=blue mid arrow, in=270, out=90] (32.center) to (23.center);
		\draw [style=blue mid arrow, in=15, out=90, looseness=0.50] (31.center) to (33.center);
		\draw [style=blue mid arrow, in=-90, out=90] (29.center) to (25.center);
		\draw [style=black] (44.center) to (96.center);
		\draw [style=black bg] (96.center) to (2.center);
		\draw [style=black] (43.center) to (0.center);
		\draw [style=red mid arrow] (0.center) to (97.center);
		\draw [style=red] (97.center) to (1.center);
		\draw [style=black] (1.center) to (34.center);
		\draw [style=black] (36.center) to (9.center);
		\draw [style=black, in=0, out=-90, looseness=0.75] (9.center) to (18.center);
		\draw [style=black] (40.center) to (14.center);
		\draw [style=red mid arrow, in=-90, out=0] (38.center) to (40.center);
		\draw [style=black] (41.center) to (15.center);
		\draw [style=black, in=0, out=-90, looseness=1.25] (15.center) to (17.center);
		\draw [style=black, in=-180, out=-90, looseness=1.25] (8.center) to (17.center);
		\draw [style=black] (45.center) to (8.center);
		\draw [style=black, in=180, out=-90] (14.center) to (18.center);
		\draw [style=red, in=-180, out=-90] (41.center) to (38.center);
		\draw [style=red mid arrow, in=-90, out=0, looseness=0.75] (43.center) to (45.center);
		\draw [style=red, in=180, out=-90] (36.center) to (34.center);
	\end{pgfonlayer}
\end{tikzpicture}
		\caption{An example of an oriented surface $\Lambda(\Sigma)$ with non-empty boundary, together with choices of open neighbourhoods near parametrised boundary components ($U_i$ and $V_j$)
			and their identification (from Step 1.2.1) along the boundary indicated by thick grey lines.
			The curves $u_j$ (from Step 1.2.2) are shown in blue.}
		\label{fig:boundary-identifications}
	\end{figure}
	
	\textsl{Step 1.2.1 (definition of $\widetilde{\Sigma}$)}: 
	We define the new oriented surface $\widetilde{\Sigma}$ using the boundary parametrisation maps. 
	Let $\partial_i$, $i\in \{ 1,\dots , |\pi_0(\partial\Sigma)|\}$, denote the parametrised boundary components of $\Lambda(\Sigma)$, 
	which may be circles or intervals. 
	We arbitrarily single out the component $\partial_1$, and we choose a connected subset~$U_j$ 
	of an open neighbourhood of~$\partial_1$ 
	for each remaining boundary component $\partial_j$, $j\in\{2,\dots,|\pi_0(\partial\Sigma)|\}$ 
	so that the~$U_j$ are pairwise disjoint.
	Furthermore we choose a connected subset~$V_j$ of an open neighbourhood of~$\partial_j$ 
	in each remaining boundary component ($j\neq 1$). 
	We illustrate such choices in Figure~\ref{fig:boundary-identifications}.
	
	\begin{figure}[tb]
		\centering
		\resizebox{\textwidth}{!}{
			\input{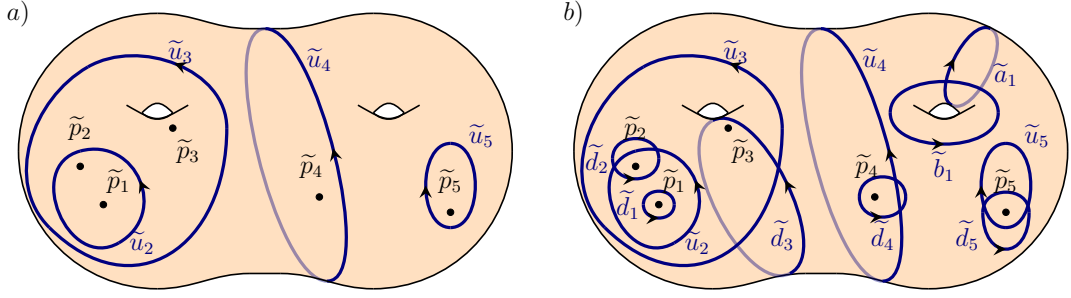}
		}
		\caption{$a)$ The surface $\widetilde{\Sigma}$ obtained from $\Lambda(\Sigma)$ as in Figure~\ref{fig:boundary-identifications} in Step 1.2.1, and the additional curves~$\widetilde{u}_j$ from Step 1.2.2. 
			\newline 
			$b)$ The surface $\widetilde{\Sigma}$ with all the curves 
			$\widetilde{a_k},\widetilde{b_k}, \widetilde{u}_j, \widetilde{d}_i$ from Step 1.2.2.
		}
		\label{fig:curves-on-tilde-Sigma}
	\end{figure}
	Using the boundary parametrisation maps we glue 
	$U_j\cap\partial_1$ to $V_j\cap\partial_j$.
	Finally we retract each remaining boundary component to a single point.
	The surface $\widetilde{\Sigma}$ obtained in this way is 
	a closed surface, whose genus is the sum of~$g$ and the number of closed parametrised boundary components of $\Lambda(\Sigma)$, 
	with a point 
	$\widetilde{p}_i$ 
	removed for each parametrised boundary component 
	$\partial_i$; 
	see Figure~\ref{fig:curves-on-tilde-Sigma}.
	
	\textsl{Step 1.2.2 (additional curves)}: 
	\begin{figure}[tb]
		\centering
		\resizebox{\textwidth}{!}{
			\input{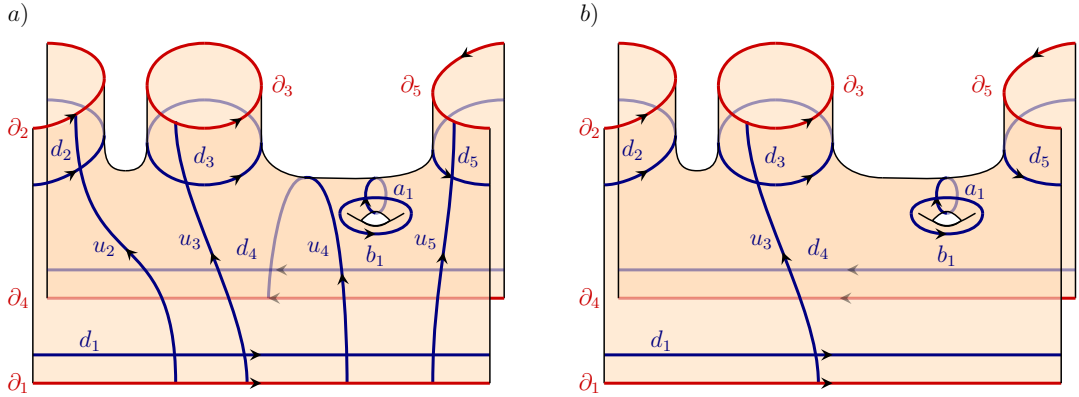}
		}
		\caption{$a)$ All curves on $\Lambda(\Sigma)$ obtained in Step 1.2.2.
			\newline 
			$b)$ Only those curves on $\Lambda(\Sigma)$ whose images in $\widetilde{\Sigma}$ form a minimal generating set of
			$\pi_1(\widetilde{\Sigma})$.}
		\label{fig:curves-on-Sigma}
	\end{figure}
	We extend our collection of curves in $\Lambda(\Sigma)$ 
	by defining curves in $\widetilde{\Sigma}$:
	We pick a simple closed curve~$\widetilde{u}_j$ 
	for each $\partial_j$, $j\in\{2,\dots,|\pi_0(\partial\Sigma)|\}$, encircling $\widetilde{p}_j$,
	and not intersecting with each other, as well as a simple closed curve~$\widetilde{d}_i$
	for each boundary component ``parallel'' along the boundary curve, 
	see Figure~\ref{fig:curves-on-tilde-Sigma}\,$b)$.
	Then we obtain the following curves in $\Lambda(\Sigma)$ corresponding the curves in 
	$\widetilde{\Sigma}$, see Figure~\ref{fig:curves-on-Sigma}\,$a)$:
	\begin{itemize}
		\item two closed curves $a_k,b_k$, $k\in \{1,\dots,g\}$, for each handle, $k \in \{1,\dots,g\}$,
		\item one curve $d_i$ for each boundary component~$\partial_i$, $i\in \{ 1,\dots , |\pi_0(\partial\Sigma)|\}$,
		\item one curve $u_j$ for each boundary component~$\partial_j$, $j\in \{ 2,\dots , |\pi_0(\partial\Sigma)|\}$.
	\end{itemize}
	
	Since near each $U_j$ and $V_j$ the $r$-spin
	surface is trivial, we obtain an $r$-spin structure on $\widetilde{\Sigma}$.
	Also note that the set of $r$-spin structures on $\Lambda(\Sigma)$ with prescribed
	$r$-spin structure near the boundary and the set of $r$-spin structures on 
	$\widetilde{\Sigma}$ with prescribed $r$-spin structure near its punctures
	are in bijection by construction.
	These sets are in bijection with the set
	\begin{equation}
	\left\{\delta\in H^1(\widetilde{\Sigma},\Z_r)\ \Big|\ \delta(d_i)=x_i \text{ for all $i$};\ \delta(u_j)=y_j
	\text{ for all $j$ with $\partial_j\ncong S^1$}\right\},
	\label{eq:H1-tSigma}
	\end{equation}
	where $x_i$ is the holonomy along $d_i$ and $y_j$ is the holonomy along $u_j$, 
	which are fixed by the boundary parametrisation.
	
	In order to describe the $r$-spin structure on $\Sigma$ it is enough to remember 
	the holonomies along a set of curves that generate $\pi_1(\widetilde{\Sigma})$.
	Therefore we reduce the set of curves $a_k,b_k,d_i,u_j$ by discarding the curves 
	$u_j$ with $\partial_j\not\cong S^1$.
	The remaining curves are illustrated in Figure~\ref{fig:curves-on-Sigma}\,$b)$.
	
	\begin{figure}[tb]
		\centering
		\begin{tikzpicture}
	\begin{pgfonlayer}{nodelayer}
		\node [style=none] (35) at (8, 2.25) {};
		\node [style=none] (40) at (8, -2.25) {};

		\node [style=none] (38) at (10, -3.5) {};
		\node [style=none] (39) at (10, -1) {};
		\node [style=none] (33) at (10, 1) {};
		\node [style=none] (34) at (10, 3.5) {};

		\node [style=none] (43) at (11, 1) {};
		\node [style=none] (44) at (11, 3.5) {};
		\node [style=none] (48) at (11, -3.5) {};
		\node [style=none] (49) at (11, -1) {};

		\node [style=none] (46) at (13, 2.25) {};
		\node [style=none] (51) at (13, -2.25) {};

		\node [style=none] (83) at (-6, -3.5) {};
		\node [style=none] (59) at (-6, 0) {};
		\node [style=none] (76) at (-6, 3.5) {};

		\node [style=none] (81) at (-4, -4.75) {};
		\node [style=none] (82) at (-4, -2.25) {};
		\node [style=none] (57) at (-4, -1.25) {};
		\node [style=none] (58) at (-4, 1.25) {};
		\node [style=none] (74) at (-4, 2.25) {};
		\node [style=none] (75) at (-4, 4.75) {};

		\node [style=none] (87) at (-3, -2.25) {};
		\node [style=none] (63) at (-3, -1.25) {};
		\node [style=none] (64) at (-3, 1.25) {};
		\node [style=none] (99) at (-3, 2.25) {};
		\node [style=none] (90) at (-3, 4.75) {};

		\node [style=none] (65) at (-1, 0) {};

		\node [style=none] (69) at (1, 0) {};

		\node [style=none] (91) at (3, -4.75) {};
		\node [style=none] (88) at (3, -2.25) {};
		\node [style=none] (67) at (3, -1.25) {};
		\node [style=none] (68) at (3, 1.25) {};
		\node [style=none] (92) at (3, 2.25) {};

		\node [style=none] (84) at (4, -4.75) {};
		\node [style=none] (85) at (4, -2.25) {};
		\node [style=none] (70) at (4, -1.25) {};
		\node [style=none] (71) at (4, 1.25) {};
		\node [style=none] (77) at (4, 2.25) {};
		\node [style=none] (78) at (4, 4.75) {};

		\node [style=none] (86) at (6, -3.5) {};
		\node [style=none] (72) at (6, 0) {};
		\node [style=none] (79) at (6, 3.5) {};
		\node [style=none] (203) at (-6.5, 5.5) {$a)$};
		\node [style=none] (204) at (7.5, 5.5) {$b)$};`
	\end{pgfonlayer}

	\begin{pgfonlayer}{background}
	\path [fill=orange!20!white, opacity=0.8]
		(34.center) to [in=90, out=180] (35.center)
		to [in=-180, out=-90] (33.center)
		-- (43.center)
		to [in=-90, out=0] (46.center)
		to [in=0, out=90] (44.center)
		-- (34.center);
	\path [fill=orange!30!white, opacity=0.8]
		(35.center)
		to [in=-180, out=-90] (33.center)
		-- (43.center)
		to [in=-90, out=0] (46.center)
		to (51.center)
		to [in=0, out=90] (49.center)
		to (39.center) 
		to [in=90, out=180] (40.center)
		;
	\path [fill=orange!20!white, opacity=0.8]
		(51.center)
		to [in=0, out=-90] (48.center)
		to (38.center) 
		to [in=-90, out=-180] (40.center)
		to [in=180, out=90] (39.center) 
		to (49.center)
		to [in=90, out=0] (51.center)
		;
	\path [fill=orange!20!white, opacity=0.8]
		(76.center) to [in=-180, out=-90] (74.center)
		to (77.center)
		to [in=-90, out=0] (79.center)
		to[in=0, out=90] (78.center)
		to (75.center)
		to [in=90, out=180] (76.center);
	\path [fill=orange!30!white, opacity=0.8]
		(76.center) to [in=-180, out=-90] (74.center)
		to (77.center)
		to [in=-90, out=0] (79.center)
		to (86.center)
		to [in=0, out=90] (85.center)
		to (82.center)
		to [in=90, out=180] (83.center)
		to (76.center);
	\path [fill=orange!20!white, opacity=0.8]
		(86.center)
		to [in=0, out=-90] (84.center)
		to (81.center)
		to [in=-90, out=-180] (83.center)
		to [in=180, out=90] (82.center)
		to (85.center)
		to [in=90, out=0] (86.center);
	\path [fill=white, opacity=1]
		(65.center) to [bend left=90, looseness=1.50] (69.center)
		to [bend left=90, looseness=1.50] (65.center);
	\end{pgfonlayer}
	\begin{pgfonlayer}{edgelayer}
		\draw [style=black bg, in=180, out=90] (40.center) to (39.center);
		\draw [style=black, in=-180, out=-90] (40.center) to (38.center);
		\draw [style=black] (35.center) to (40.center);
		\draw [style=black, in=-90, out=0] (48.center) to (51.center);
		\draw [style=black bg, in=90, out=0] (49.center) to (51.center);
		\draw [style=black] (46.center) to (51.center);
		\draw [style=black bg] (39.center) to (49.center);
		\draw [style=black] (38.center) to (48.center);
		\draw [style=black dashed bg, in=180, out=90] (59.center) to (58.center);
		\draw [style=black dashed, in=-180, out=-90] (59.center) to (57.center);
		\draw [style=black dashed, in=-90, out=0] (63.center) to (65.center);
		\draw [style=black dashed bg, in=90, out=0] (64.center) to (65.center);
		\draw [style=black dashed bg] (58.center) to (64.center);
		\draw [style=black dashed] (57.center) to (63.center);
		\draw [style=black dashed bg, in=180, out=90] (69.center) to (68.center);
		\draw [style=black dashed, in=-180, out=-90] (69.center) to (67.center);
		\draw [style=black dashed, in=-90, out=0] (70.center) to (72.center);
		\draw [style=black dashed bg, in=90, out=0] (71.center) to (72.center);
		\draw [style=black dashed bg] (68.center) to (71.center);
		\draw [style=black dashed] (67.center) to (70.center);
		\draw [style=black, in=180, out=90] (76.center) to (75.center);
		\draw [style=black, in=-180, out=-90] (76.center) to (74.center);
		\draw [style=black, in=-90, out=0] (77.center) to (79.center);
		\draw [style=black, in=90, out=0] (78.center) to (79.center);
		\draw [style=black] (75.center) to (78.center);
		\draw [style=black] (74.center) to (77.center);
		\draw [style=black bg, in=180, out=90] (83.center) to (82.center);
		\draw [style=black, in=-180, out=-90] (83.center) to (81.center);
		\draw [style=black, in=-90, out=0] (84.center) to (86.center);
		\draw [style=black bg, in=90, out=0] (85.center) to (86.center);
		\draw [style=black bg] (82.center) to (85.center);
		\draw [style=black] (81.center) to (84.center);
		\draw [style=black] (76.center) to (59.center);
		\draw [style=black] (59.center) to (83.center);
		\draw [style=black] (72.center) to (86.center);
		\draw [style=black] (79.center) to (72.center);
		\draw [style=black dashed] (74.center) to (57.center);
		\draw [style=black dashed] (57.center) to (81.center);
		\draw [style=black dashed] (90.center) to (99.center);
		\draw [style=black dashed bg] (99.center) to (64.center);
		\draw [style=black dashed bg] (64.center) to (87.center);
		\draw [style=black, bend left=90, looseness=1.50] (65.center) to (69.center);
		\draw [style=black, bend right=90, looseness=1.50] (65.center) to (69.center);
		\draw [style=black dashed] (92.center) to (88.center);
		\draw [style=black dashed] (33.center) to (38.center);
		\draw [style=black dashed] (44.center) to (43.center);
		\draw [style=black dashed bg] (43.center) to (49.center);
		\draw [style=black dashed] (88.center) to (91.center);
		\draw [style=black dashed bg] (77.center) to (71.center);
		\draw [style=black dashed] (78.center) to (77.center);
		\draw [style=black dashed bg] (71.center) to (85.center);
		\draw [style=red mid arrow, in=90, out=180] (34.center) to (35.center);
		\draw [style=red, in=-180, out=-90] (35.center) to (33.center);
		\draw [style=red] (33.center) to (43.center);
		\draw [style=red mid arrow, in=-90, out=0] (43.center) to (46.center);
		\draw [style=red, in=0, out=90] (46.center) to (44.center);
		\draw [style=red] (44.center) to (34.center);
	\end{pgfonlayer}
\end{tikzpicture}
		\caption{Decomposition of a surface into oriented generators.
			\newline
			$a)$ Decomposition of a handle. 
			\newline
			$b)$ Decomposition near a closed parametrised boundary component.
		}
		\label{fig:surf-dec-gen}
	\end{figure}
	
	\textsl{Step 2:}
	We pick a decomposition of $\Lambda(\Sigma)$ into oriented generators, so that
	\begin{itemize}
		\item for each handle we have the decomposition as in 
		Figure~\ref{fig:surf-dec-gen}\,$a)$, 
		\item for each boundary component we have the decomposition as in 
		Figure~\ref{fig:surf-dec-gen}\,$b)$.
	\end{itemize}
	Note that we can require that 1-morphism components of the Serre automorphism only
	appear at the boundary of generating 2-morphisms in $\Bordor$ 
	and not in their interior, as in $\Bordor$ there is a
	trivialisation of the Serre automorphism, cf.\ Lemma~\ref{lem:Serre-triv}.
	We furthermore require that the curves $a_k, b_k, u_j, d_i$ cross the generating 2-morphisms as in Figure~\ref{fig:insert-S}. 
	
	\begin{figure}[tb]
		\centering
		\input{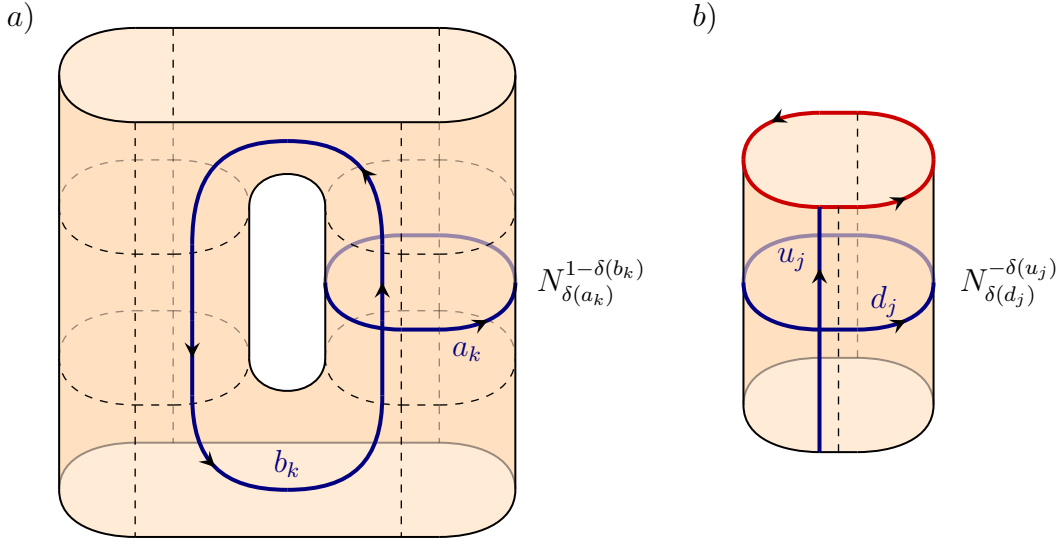}
		\caption{Insertion of the Nakayama automorphism at intersection of curves.}
		\label{fig:insert-S}
	\end{figure}
	
	\textsl{Step 3:}
	Recall from Lemma~\ref{lem:Nakayama-deck-trf} that the Nakayama automorphisms
	$N_a\colon C_a\lra C_a$ are deck transformations.
	We lift the oriented generators to $r$-spin generators by inserting
	(i.\,e.\ by replacing a neighbourhood of $a_k$ and $d_j$ with)
	\begin{itemize}
		\item $N_{\delta(a_k)}^{1-\delta(b_k)}$ at the intersection of $a_k$ and $b_k$,
		see Figure~\ref{fig:insert-S}\,$a)$,
		\item $N_{\delta(d_j)}^{-\delta(u_j)}$ at the intersection of $d_j$ and $u_j$,
		see Figure~\ref{fig:insert-S}\,$b)$,
		\item 
		$\mu_{-\delta(a_k),\delta(a_k)}$ and $\Delta_{-\delta(a_k),\delta(a_k)}$
		from \eqref{eq:Lambda_mu}--\eqref{eq:Lambda_Delta}
		at the saddles crossed by $b_k$,
		\item identity 2-morphisms where no intersections occur.
	\end{itemize}
	
	The $r$-spin structure given by this construction has the same holonomies
	along the above mentioned curves as the $r$-spin structure of $\Sigma$.
	Note that a full circle along a positively oriented simple closed loop, 
	where no insertions of $N_{\delta(a_k)}$ appear, contributes +1 to the holonomy.
	Hence the two $r$-spin structures are isomorphic, and thus the two $r$-spin
	surfaces represent the same 2-morphisms.
\end{proof}

\begin{lemma}
	The functor $\iota^r$ is faithful on 2-morphisms.
	\label{lem:iota-r-faithful-2-mor}
\end{lemma}
\begin{proof}
	Let $\sigma,\sigma'\in\FD{r}(\alpha,\beta)(X,Y)$. 
	Assume that $\iota^r(\sigma)=\iota^r(\sigma')$, and that these are connected bordisms.
	By Proposition~\ref{prop:iota-lambda-diag} we have
	\begin{equation}
	\big( \Lambda\circ\iota^r \big)(\sigma)=
	\big( \iota^1 \circ K \big)(\sigma) \, , 
	\label{eq:L-i-i-L}
	\end{equation}
	and analogously for~$\sigma'$. 
	Hence since~$\iota^1$ is an equivalence (Theorem~\ref{thm:iota-1}), we have
	\begin{equation}
	K(\sigma)=K(\sigma')\,.
	\label{eq:L-sigma}
	\end{equation}
	This means that there is a sequence of relations in $\FD{1}$ relating
	$K(\sigma)$ and $K(\sigma')$.
	
	By \eqref{eq:K-theta} the numbers of $\lambda$ and $\lambda^{-1}$ 
	in $K(\sigma)$ and $K(\sigma')$ are each divisible by $r$.
	Using the coherence theorem for 2-categories, 
	we can bundle together the relations
	in the sequence involving $\lambda$ in tuples of $r$.
	These can be lifted to relations in $\FD{r}$ via \eqref{eq:K-theta}.
	Noting that all other relations in $\FD{1}$ and $\FD{r}$ are the same
	(to wit, those of $\FD{0}$) and that they can hence also be lifted, 
	it follows that $\sigma = \sigma'$. 
\end{proof}

\subsubsection[Computing invariants of \texorpdfstring{$r$}{r}-spin bordisms with closed boundary]{Computing invariants of \texorpdfstring{$\boldsymbol{r}$}{r}-spin bordisms with closed boundary}

With the $r$-spin cobordism hypothesis at hand, we can describe the closed
$\Lambda_r$-Frobenius algebra which classifies the non-extended
$r$-spin TQFT associated to a fully extended $r$-spin TQFT.
In particular we can explicitly describe the values of 
the non-extended $r$-spin TQFT on $r$-spin surfaces with closed boundary.
For convenience we also present the corresponding results for framed TQFTs.

Recall from Sections~\ref{sec:open-closed}~and~\ref{subsec:bord-2-cat} that 
\begin{equation}
\Bord_{2,1}^{\text{fr}} \cong \Bordfr(\varnothing,\varnothing)\quad\text{and}\quad
\Bord_{2,1}^{r\text{-spin}} \cong \Bordrspin(\varnothing,\varnothing)\,.
\label{eq:bord-nonext}
\end{equation}
Consider the
fully extended framed and $r$-spin TQFTs
\begin{equation}
\begin{aligned}
Y\colon\Bordfr&\lra\B &\text{and}&& Z\colon\Bordrspin&\lra\B \\
+&\longmapsto\alpha
&&&
+&\longmapsto\alpha
\\
&&&&\vartheta_+&\longmapsto \big( \vartheta\colon S_{\alpha}^r \stackrel{\cong}{\lra} 1_{\alpha}\big) \, . 
\end{aligned}
\label{eq:Y-Z}
\end{equation}
The corresponding non-extended TQFTs
\begin{equation}
\begin{aligned}
Y|\colon \Bord_{2,1}^{\text{fr}}&\lra \B(\one,\one)&&\text{and}&
Z|\colon \Bord_{2,1}^{r\text{-spin}}&\lra \B(\one,\one)\\
S^1_a&\lmt C^{\alpha}_a&&&
S^1_a&\lmt C^{\alpha}_a
\end{aligned}
\label{eq:Z-restr}
\end{equation}
are classified by the closed $\Lambda_0$- and $\Lambda_r$-Frobenius algebras 
in Proposition~\ref{prop:Lambda0FrobAlgInB11} and Corollary~\ref{cor:LambdarFrobAlgInB11}, respectively.

In particular, the invariants assigned to 
the framed and $r$-spin tori $T(d)$ introduced in Section~\ref{subsec:ComputingInvariants} are 
quantum dimensions of the circle spaces in $\B(\one,\one)$: 
\begin{equation}
T(d)\lmt\dim(C_d^{\alpha})\,.
\label{eq:Z-restr-circle-space-and-tori}
\end{equation}

\begin{remark}
	\label{rem:PivotalNotGoodForSpin}
	\begin{enumerate}[label={(\roman*)}]
		\item 
		Assume that, as in Remark~\ref{rem:pivotal-boring-ext},
		the left and right adjoints of 1-morphisms in the target 2-category~$\B$ agree.
		In this case we effectively have oriented TQFTs ($r$ odd) with all the circle spaces being isomorphic,
		or 2-spin TQFTs ($r$ even) with $C^{\alpha}_a\cong C^{\alpha}_{a+2}$ for every $a\in\Z$.
		Accordingly, the invariants associated to framed and $r$-spin tori may take at most two distinct values if left and right adjoints agree in~$\B$.
		\item 
		\label{item:EulerArf}
		For oriented and 2-spin surfaces there already exist TQFTs which compute 
		complete invariants (the oriented TQFT of \cite{Quinnlectures} with target $\Vectk$ computed from the relative Euler characteristic, and the 2-spin TQFT of \cite{MooreSegal,RunkelSzegedyArf} with target $\Vectk^{\Z_2}$ computing the Arf invariant). 
		For $r>2$, TQFTs with pivotal 2-categories as targets
		cannot distinguish all $r$-spin structures, 
		but other targets may allow for more interesting $r$-spin TQFTs. 
	\end{enumerate}
\end{remark}

\section{Examples}
\label{sec:Examples}

By the main result of the previous section, constructing extended $r$-spin TQFTs amounts to finding fully dualisable objects whose Serre automorphisms are such that their $r$-th power is trivialisable. 
In Section~\ref{subsec:EquivariantCompletion} we increase our chances to find such objects by passing from a given target 2-category to its ``equivariant completion'', where we translate the condition on the Serre automorphism to a condition on the Nakayama automorphism of certain Frobenius algebras (Corollary~\ref{cor:SpinTQFTwithBeq}), and we study the associated circle spaces and hence torus invariants in detail (Section~\ref{subsubsec:FrobenisAlgebraInBeq}). 
Then in Section~\ref{subsec:LGmodels} we show that every object in the 2-category of Landau--Ginzburg models $\LG$ gives rise to an extended 2-spin TQFT, and we illustrate how to do computations in the equivariant completion of $\LG$ and its variants.

\subsection{Equivariant completion} 
\label{subsec:EquivariantCompletion}

In this section, we consider the representation 2-category $\Beq$ of certain Frobenius algebras internal to a given symmetric monoidal 
	pivotal 
2-category~$\B$. 
In particular, we explicitly determine the Serre automorphisms and circle spaces associated to objects in $\Beq$, from which invariants of extended $r$-spin TQFTs with values in $\Beq$ can be computed with the help of Theorem~\ref{thm:equivalence-conditions}. 
We stress that even if the original 2-category~$\B$ is pivotal, its completion $\Beq$ need not be pivotal, which in light of Remark~\ref{rem:PivotalNotGoodForSpin} is a desired feature. 

\medskip 

Throughout this section we fix a symmetric monoidal 
	pivotal 
2-category~$\B$ which satisfies the condition \hyperref[cond:ConditionOnB]{\text{($*$)}} below. 
(The symmetric monoidal structure will not be relevant before Section~\ref{subsubsec:SMstructureOfBeq}.)

\subsubsection{Equivariant completion of a 2-category}
\label{subsubsec:EquivariantCompletion}

A \textsl{$\Delta$-separable Frobenius algebra on an object $\alpha\in\B$} consists of
\be 
A \in \B(\alpha,\alpha) 
\, , \quad 
\mu_A = 
\begin{tikzpicture}[very thick,scale=0.53,color=green!50!black, baseline=0.65cm]
\draw[-dot-] (2.5,0.75) .. controls +(0,1) and +(0,1) .. (3.5,0.75);
\draw (3,1.5) -- (3,2.25); 
\end{tikzpicture} 
\, , \quad 
\eta_A = 
\begin{tikzpicture}[very thick,scale=0.4,color=green!50!black, baseline=0]
\draw (-0.5,-0.5) node[Odot] (unit) {}; 
\draw (unit) -- (-0.5,0.5);
\end{tikzpicture} 
\, , \quad 
\Delta_A = 
\begin{tikzpicture}[very thick,scale=0.53,color=green!50!black, rotate=180, baseline=-0.9cm]
\draw[-dot-] (2.5,0.75) .. controls +(0,1) and +(0,1) .. (3.5,0.75);
\draw (3,1.5) -- (3,2.25); 
\end{tikzpicture} 
\, , \quad 
\varepsilon_A = 
\begin{tikzpicture}[very thick,scale=0.4,color=green!50!black, baseline=0cm, rotate=180]
\draw (-0.5,-0.5) node[Odot] (unit) {}; 
\draw (unit) -- (-0.5,0.5);
\end{tikzpicture} 
\ee 
such that 
\begin{align}
& 
\begin{tikzpicture}[very thick,scale=0.53,color=green!50!black, baseline=0.59cm]
\draw[-dot-] (3,0) .. controls +(0,1) and +(0,1) .. (2,0);
\draw[-dot-] (2.5,0.75) .. controls +(0,1) and +(0,1) .. (3.5,0.75);
\draw (3.5,0.75) -- (3.5,0); 
\draw (3,1.5) -- (3,2.25); 
\end{tikzpicture} 
=
\begin{tikzpicture}[very thick,scale=0.53,color=green!50!black, baseline=0.59cm]
\draw[-dot-] (3,0) .. controls +(0,1) and +(0,1) .. (2,0);
\draw[-dot-] (2.5,0.75) .. controls +(0,1) and +(0,1) .. (1.5,0.75);
\draw (1.5,0.75) -- (1.5,0); 
\draw (2,1.5) -- (2,2.25); 
\end{tikzpicture} 
\, , 
\quad
\begin{tikzpicture}[very thick,scale=0.4,color=green!50!black, baseline]
\draw (-0.5,-0.5) node[Odot] (unit) {}; 
\fill (0,0.6) circle (5.0pt) node (meet) {};
\draw (unit) .. controls +(0,0.5) and +(-0.5,-0.5) .. (0,0.6);
\draw (0,-1.5) -- (0,1.5); 
\end{tikzpicture} 
=
\begin{tikzpicture}[very thick,scale=0.4,color=green!50!black, baseline]
\draw (0,-1.5) -- (0,1.5); 
\end{tikzpicture} 
=
\begin{tikzpicture}[very thick,scale=0.4,color=green!50!black, baseline]
\draw (0.5,-0.5) node[Odot] (unit) {}; 
\fill (0,0.6) circle (5.0pt) node (meet) {};
\draw (unit) .. controls +(0,0.5) and +(0.5,-0.5) .. (0,0.6);
\draw (0,-1.5) -- (0,1.5); 
\end{tikzpicture} 
\, , 
\quad
\begin{tikzpicture}[very thick,scale=0.53,color=green!50!black, baseline=-0.59cm, rotate=180]
\draw[-dot-] (3,0) .. controls +(0,1) and +(0,1) .. (2,0);
\draw[-dot-] (2.5,0.75) .. controls +(0,1) and +(0,1) .. (1.5,0.75);
\draw (1.5,0.75) -- (1.5,0); 
\draw (2,1.5) -- (2,2.25); 
\end{tikzpicture} 
=
\begin{tikzpicture}[very thick,scale=0.53,color=green!50!black, baseline=-0.59cm, rotate=180]
\draw[-dot-] (3,0) .. controls +(0,1) and +(0,1) .. (2,0);
\draw[-dot-] (2.5,0.75) .. controls +(0,1) and +(0,1) .. (3.5,0.75);
\draw (3.5,0.75) -- (3.5,0); 
\draw (3,1.5) -- (3,2.25); 
\end{tikzpicture} 
\, , 
\quad
\begin{tikzpicture}[very thick,scale=0.4,color=green!50!black, baseline=0, rotate=180]
\draw (0.5,-0.5) node[Odot] (unit) {}; 
\fill (0,0.6) circle (5.0pt) node (meet) {};
\draw (unit) .. controls +(0,0.5) and +(0.5,-0.5) .. (0,0.6);
\draw (0,-1.5) -- (0,1.5); 
\end{tikzpicture} 
=
\begin{tikzpicture}[very thick,scale=0.4,color=green!50!black, baseline=0, rotate=180]
\draw (0,-1.5) -- (0,1.5); 
\end{tikzpicture} 
=
\begin{tikzpicture}[very thick,scale=0.4,color=green!50!black, baseline=0cm, rotate=180]
\draw (-0.5,-0.5) node[Odot] (unit) {}; 
\fill (0,0.6) circle (5.0pt) node (meet) {};
\draw (unit) .. controls +(0,0.5) and +(-0.5,-0.5) .. (0,0.6);
\draw (0,-1.5) -- (0,1.5); 
\end{tikzpicture} 
\, , 
\nonumber 
\\ 
& 
\begin{tikzpicture}[very thick,scale=0.4,color=green!50!black, baseline=0cm]
\draw[-dot-] (0,0) .. controls +(0,-1) and +(0,-1) .. (-1,0);
\draw[-dot-] (1,0) .. controls +(0,1) and +(0,1) .. (0,0);
\draw (-1,0) -- (-1,1.5); 
\draw (1,0) -- (1,-1.5); 
\draw (0.5,0.8) -- (0.5,1.5); 
\draw (-0.5,-0.8) -- (-0.5,-1.5); 
\end{tikzpicture}
=
\begin{tikzpicture}[very thick,scale=0.4,color=green!50!black, baseline=0cm]
\draw[-dot-] (0,0) .. controls +(0,1) and +(0,1) .. (-1,0);
\draw[-dot-] (1,0) .. controls +(0,-1) and +(0,-1) .. (0,0);
\draw (-1,0) -- (-1,-1.5); 
\draw (1,0) -- (1,1.5); 
\draw (0.5,-0.8) -- (0.5,-1.5); 
\draw (-0.5,0.8) -- (-0.5,1.5); 
\end{tikzpicture}
\, , 
\quad
\begin{tikzpicture}[very thick,scale=0.4,color=green!50!black, baseline=0cm]
\draw[-dot-] (0,0) .. controls +(0,-1) and +(0,-1) .. (1,0);
\draw[-dot-] (0,0) .. controls +(0,1) and +(0,1) .. (1,0);
\draw (0.5,-0.8) -- (0.5,-1.5); 
\draw (0.5,0.8) -- (0.5,1.5); 
\end{tikzpicture}
\, = \, 
\begin{tikzpicture}[very thick,scale=0.4,color=green!50!black, baseline=0cm]
\draw (0.5,-1.5) -- (0.5,1.5); 
\end{tikzpicture}
\, . 
\end{align}

Recall, e.\,g.\ from \cite[Sect.\,2.2]{cr1210.6363}, the notions of (bi)modules and (bi)module maps over the underlying algebra $(A,\mu_A,\eta_A)$. 
If~$X$ is a right $A$-module and~$Y$ is a left $A$-module, then the \textsl{relative tensor product} $X\otimes_A Y$ is the coequaliser of the canonical maps 
$\!\!\!\!
\begin{tikzpicture}[
baseline=(current bounding box.base), 
descr/.style={fill=white,inner sep=3.5pt}, 
normal line/.style={->}, 
baseline=0cm
] 
\matrix (m) [matrix of math nodes, row sep=0em, column sep=1.6em, text height=0ex, text depth=0ex] {%
	X\otimes A\otimes Y & X\otimes Y
	\\
};
\path[font=\footnotesize, transform canvas={yshift=2.3mm}] (m-1-1) edge[->] node {  } (m-1-2);
\path[font=\footnotesize, transform canvas={yshift=0.7mm}] (m-1-1) edge[->] node {  } (m-1-2);
\end{tikzpicture}
\!\!\!\!$. 
Since~$A$ is a $\Delta$-separable Frobenius algebra, the map 
\be 
\pi_A^{X,Y} = 
\tikzzbox{%
	\begin{tikzpicture}[very thick,scale=0.75,color=blue!50!black, baseline]
	\draw (-1,-1) node[left] (X) {{\small$X$}};
	\draw (1,-1) node[right] (Xu) {{\small$Y$}};
	\draw (-1,-1) -- (-1,1); 
	\draw (1,-1) -- (1,1); 
	\fill[color=green!50!black] (-1,0.6) circle (2.5pt) node (meet) {};
	\fill[color=green!50!black] (1,0.6) circle (2.5pt) node (meet) {};
	\draw[-dot-, color=green!50!black] (0.35,-0.0) .. controls +(0,-0.5) and +(0,-0.5) .. (-0.35,-0.0);
	\draw[color=green!50!black] (0.35,-0.0) .. controls +(0,0.25) and +(-0.25,-0.25) .. (1,0.6);
	\draw[color=green!50!black] (-0.35,-0.0) .. controls +(0,0.25) and +(0.25,-0.25) .. (-1,0.6);
	\draw[color=green!50!black] (0,-0.75) node[Odot] (down) {}; 
	\draw[color=green!50!black] (down) -- (0,-0.35); 
	\end{tikzpicture} 
}
\ee 
is an idempotent. 
If $\pi_A^{X,Y}$ splits, then $X\otimes_A Y$ can be identified with $\textrm{Im}(\pi_A^{X,Y})$, see e.\,g.\ \cite[Lem.\,2.3]{cr1210.6363}. 
Hence we will make the following assumption: 
\begin{itemize}
	\item[($*$)]
	\label{cond:ConditionOnB}
	For all $\Delta$-separable Frobenius algebras~$A$ on all objects of~$\B$, the idempotents $\pi_A^{X,Y}$ split for all modules $X,Y$, and we choose adjunction data for~$A$ such that ${}^\dagger\!A = A^\dagger$ as well as ${}^\dagger\!\mu = \mu^\dagger$ and ${}^\dagger\!\Delta = \Delta^\dagger$. 
\end{itemize}
Thus we have splitting maps 
\be 
\begin{tikzpicture}[
baseline=(current bounding box.base), 
descr/.style={fill=white,inner sep=3.5pt}, 
normal line/.style={->}
] 
\matrix (m) [matrix of math nodes, row sep=3em, column sep=4.0em, text height=1.5ex, text depth=0.1ex] {%
	X\otimes_A Y & X\otimes Y 
	\\
};
\path[font=\footnotesize, transform canvas={yshift=0.8mm}] (m-1-1) edge[->] node[above] { $ \iota_A^{X,Y}$ } (m-1-2);
\path[font=\footnotesize, transform canvas={yshift=-0.8mm}] (m-1-2) edge[->] node[below] { $ \vartheta_A^{X,Y}$ } (m-1-1);
\end{tikzpicture}
\ee 
with $\iota_A^{X,Y} \circ \vartheta_A^{X,Y} = \pi^{X,Y}_A$ and $\vartheta_A^{X,Y}\circ \iota^{X,Y}_A = 1_{X\otimes_A Y}$. 
Note that every Frobenius algebra is self-dual, so there always exist adjunction data such that ${}^\dagger\!A = A^\dagger$. 
The conditions ${}^\dagger\!\mu = \mu^\dagger$ and ${}^\dagger\!\Delta = \Delta^\dagger$ automatically hold if~$\B$ is pivotal, but we do not make this stronger assumption on~$\B$. 

\begin{definition}
	\label{def:Beq}
	The \textsl{equivariant completion} $\Beq$ of~$\B$ is the 2-category whose
	\begin{itemize}
		\item 
		objects are pairs $(\alpha,A)$ with $\alpha\in\B$ and $A\in\B(\alpha,\alpha)$ a $\Delta$-separable Frobenius algebra; 
		\item 
		1-morphisms $(\alpha,A) \lra (\beta,B)$ are 1-morphisms $\alpha\lra\beta$ in~$\B$ together with a $B$-$A$-bimodule structure; 
		\item 
		2-morphisms are bimodule maps in~$\B$; 
		\item 
		horizontal composition is the relative tensor product, and $1_{(\alpha,A)}$ is~$A$ with its canonical $A$-$A$-bimodule structure; 
		\item 
		vertical composition and unit 2-morphisms are induced from~$\B$. 
	\end{itemize}
\end{definition}

Equivariant completion was introduced in \cite{cr1210.6363} in connection with generalised orbifold constructions of oriented TQFTs. 
The attribute ``equivariant'' derives from the fact that an action $\rho \colon \underline{G} \lra \B(\alpha,\alpha)$ of a finite group~$G$ (viewed as a discrete monoidal category~$\underline{G}$) gives rise to a $\Delta$-separable Frobenius algebra structure on $A_G := \bigoplus_{g\in G} \rho(g)$ if $\B(\alpha,\alpha)$ has finite direct sums, and that $G$-equivariantisation can be described in terms of categories of $A_G$-modules. 

The assignment $\B\lmt\Beq$ is a completion in the sense that $(\Beq)_{\textrm{eq}} \cong \Beq$, see \cite[Prop.\,4.2]{cr1210.6363}. 
Equivariant completion is the same as (unital and counital) ``condensation completion'' in dimension~2, as introduced in \cite{GJF} for arbitrary dimension in the context of fully extended framed TQFTs and topological orders.

\subsubsection{Adjoints}
\label{subsubsec:AdjointsInBeq}

Recall that we assume that every $\Delta$-separable Frobenius algebra $A\in\B(\alpha,\alpha)$ comes with chosen adjunction data such that ${}^\dagger\!A = A^\dagger$. 
Hence we can define the \textsl{Nakayama automorphism} and its inverse as follows: 
\be 
\label{eq:NakayamaA}
\gamma_A= 
\tikzzbox{%
	\begin{tikzpicture}[very thick, scale=0.5,color=green!50!black, baseline=-0.35cm]
	\draw (0,0.8) -- (0,2);
	\draw[-dot-] (0,0.8) .. controls +(0,-0.5) and +(0,-0.5) .. (-0.75,0.8);
	\draw (-0.75,0.8) .. controls +(0,0.5) and +(0,0.5) .. (-1.5,0.8);
	\draw[-dot-] (0,-1.8) .. controls +(0,0.5) and +(0,0.5) .. (-0.75,-1.8);
	\draw (-0.75,-1.8) .. controls +(0,-0.5) and +(0,-0.5) .. (-1.5,-1.8);
	\draw (0,-1.8) -- (0,-3);
	\draw[postaction={decorate}, decoration={markings,mark=at position .5 with {\arrow{>}}}] (-1.5,0.8) -- (-1.5,-1.8);
	\draw (-0.375,-0.2) node[Odot] (D) {}; 
	\draw (-0.375,0.4) -- (D);
	\draw (-0.375,-0.8) node[Odot] (E) {}; 
	\draw (-0.375,-1.4) -- (E);
	\end{tikzpicture}
}%
\, , \quad
\gamma_A^{-1} = 
\tikzzbox{%
	\begin{tikzpicture}[very thick, scale=0.5,color=green!50!black, baseline=-0.35cm]
	\draw (0,0.8) -- (0,2);
	\draw[-dot-] (0,0.8) .. controls +(0,-0.5) and +(0,-0.5) .. (0.75,0.8);
	\draw (0.75,0.8) .. controls +(0,0.5) and +(0,0.5) .. (1.5,0.8);
	\draw[-dot-] (0,-1.8) .. controls +(0,0.5) and +(0,0.5) .. (0.75,-1.8);
	\draw (0.75,-1.8) .. controls +(0,-0.5) and +(0,-0.5) .. (1.5,-1.8);
	\draw (0,-1.8) -- (0,-3);
	\draw[postaction={decorate}, decoration={markings,mark=at position .5 with {\arrow{>}}}] (1.5,0.8) -- (1.5,-1.8);
	\draw (0.375,-0.2) node[Odot] (D) {}; 
	\draw (0.375,0.4) -- (D);
	\draw (0.375,-0.8) node[Odot] (E) {}; 
	\draw (0.375,-1.4) -- (E);
	\end{tikzpicture}
}
\, . 
\ee 
Since we also assume ${}^\dagger\!\mu = \mu^\dagger$ and ${}^\dagger\!\Delta = \Delta^\dagger$, the Nakayama automorphism is a map of the underlying algebra and coalgebra structures of~$A$, and~$A$ is a symmetric Frobenius algebra iff $\gamma_A = 1_A$, see e.\,g.\ \cite{FS:FA}. 

\begin{remark}
  Let $\Cc$ be a symmetric monoidal 1-category with 
  left duals, and let $A\in\Cc$ be a Frobenius algebra.
  We can endow $\Cc$ with right duals by setting $\Xd:=\dX$,
  $\tev_X:=\ev_X\circ \, b_{\Xd,X}$ and 
  $\tcoev_X:=b_{X,\Xd}\circ \coev_X$ using the braiding $b$.
  In this case the inverse Nakayama automorphism of $A$ is
  \begin{equation}
    \gamma_A^{-1}=
\tikzzbox{%
	\begin{tikzpicture}[very thick,scale=0.53,color=green!50!black, baseline=-0.4cm]
	\draw[-dot-] (0,0) .. controls +(0,1) and +(0,1) .. (1,0);
	\draw (0.5,1.2) node[Odot] (unit) {}; 
	\draw (unit) -- (0.5,0.7);
	\draw (0,0) .. controls +(0,-0.5) and +(0,0.5) .. (-1,-1);
	\draw[-dot-] (0,-1) .. controls +(0,-1) and +(0,-1) .. (1,-1);
	\draw (0.5,-2.2) node[Odot] (unit2) {}; 
	\draw (unit2) -- (0.5,-1.7);
	\draw (0,-1) .. controls +(0,0.5) and +(0,-0.5) .. (-1,0);
	\draw (1,0) -- (1,-1);
	\draw (-1,0) -- (-1,2);
	\draw (-1,-1) -- (-1,-3);
	\end{tikzpicture} 
}%
\, 
	.
\label{eq:Nakayama-sym-mon-cat}
\end{equation}
Comparing this to the definition~\eqref{eq:Nakayama-Ca} of the Nakayama automorphisms of a closed $\Lambda_0$-Frobenius algebra~$C$, we see that the conventions for these two different Nakayama structures~$\gamma$ and~$N$, for two different algebraic entities~$A$ and~$C$, respectively, are not maximally aligned. 
  \label{rem:Nakayama-sym-mon-cat}
\end{remark}

Given a $B$-$A$-bimodule $X\in\B(\alpha,\beta)$ together with algebra automorphisms $\varphi \colon A\lra A$ and $\psi \colon B\lra B$, the \textsl{$\psi$-$\varphi$-twisted bimodule} ${}_\psi X_\varphi$ is given by 
\be 
\tikzzbox{%
\begin{tikzpicture}[very thick,scale=0.75,color=blue!50!black, baseline]

\draw (0,-1) node[below] (X) {{\small${}_\psi X_\varphi$}};
\draw (0,1) node[above] (Xu) {{\small${}_\psi X_\varphi$}};
\draw[color=green!50!black] (1.0,-1) node[below] (A) {{\small$A$}};
\draw[color=green!50!black] (-0.75,-1) node[below] (B) {{\small$B$}};
\draw (0,-1) -- (0,1); 
\fill[color=green!50!black] (0,-0.3) circle (2.5pt) node (meet) {};
\fill[color=green!50!black] (0,0.3) circle (2.5pt) node (meet) {};
\draw[color=green!50!black] (-0.75,-1) .. controls +(0,0.25) and +(-0.5,-0.25) .. (0,-0.3);
\draw[color=green!50!black] (1.0,-1) .. controls +(0,0.5) and +(0.5,-0.5) .. (0,0.3);
\end{tikzpicture} 
}%
=
\tikzzbox{%
\begin{tikzpicture}[very thick,scale=0.75,color=blue!50!black, baseline]
\draw (0,-1) node[below] (X) {{\small$X$}};
\draw (0,1) node[above] (Xu) {{\small$X$}};
\draw[color=green!50!black] (1.0,-1) node[below] (A) {{\small$A$}};
\draw[color=green!50!black] (-0.75,-1) node[below] (B) {{\small$B$}};
\draw (0,-1) -- (0,1); 
\fill[color=green!50!black] (0,-0.3) circle (2.5pt) node (meet) {};
\fill[color=green!50!black] (0,0.3) circle (2.5pt) node (meet) {};
\fill[color=green!50!black] (-0.5,-0.6) circle (2.5pt) node[left] (meet) {{\small$\psi$}};
\fill[color=green!50!black] (0.7,-0.35) circle (2.5pt) node[right] (meet) {{\small$\varphi$}};
\draw[color=green!50!black] (-0.75,-1) .. controls +(0,0.25) and +(-0.5,-0.25) .. (0,-0.3);\draw[color=green!50!black] (1.0,-1) .. controls +(0,0.5) and +(0.5,-0.5) .. (0,0.3);
\end{tikzpicture} 
}%
\ee 
where the unlabelled vertices on the right-hand side correspond to the original bimodule structure on~$X$. 

If the 2-category~$\B$ has adjoints for 1-morphisms, then, as shown in \cite[Prop.\,4.2]{cr1210.6363}, its equivariant completion~$\Beq$ inherits this property, by twisting with Nakayama automorphisms. 
This implies that even if~$\B$ is pivotal, $\Beq$ typically is not: 
\begin{proposition}
	\label{prop:AdjointsInBeq}
	Let~$\B$ be a 2-category, and let $X\in\Beq((\alpha,A), (\beta,B))$ be such that the underlying 1-morphism $X\colon \alpha\lra \beta$ in~$\B$ has left and right adjoints~$\dX$ and~$\Xd$, respectively. 
	Then~$X$ also has left and right adjoints
	\be 
	\deqX = {}_{\gamma_A^{-1}} (\dX) 
		\, , \quad 
		X^\star = (X^\dagger)_{\gamma_B}
		\label{eq:AdjointsInBeq-1}
	\ee 
	in $\Beq$, witnessed by the adjunction maps 
	\begin{align}
		\ev_X^{\Beq} & = 
		\tikzzbox{%
		\begin{tikzpicture}[very thick,scale=1.0,color=blue!50!black, baseline=.8cm]
		\draw[line width=0pt] 
		(1.75,1.75) node[line width=0pt, color=green!50!black] (A) {{\small$A\vphantom{\dX }$}}
		(1,0) node[line width=0pt] (D) {{\small$X\vphantom{\dX }$}}
		(0,0) node[line width=0pt] (s) {{\small$\dX $}}; 
		\draw[directed] (D) .. controls +(0,1.5) and +(0,1.5) .. (s);
		\draw[color=green!50!black] (1.25,0.55) .. controls +(0.0,0.25) and +(0.25,-0.15) .. (0.86,0.95);
		\draw[-dot-, color=green!50!black] (1.25,0.55) .. controls +(0,-0.5) and +(0,-0.5) .. (1.75,0.55);
		\draw[color=green!50!black] (1.5,-0.1) node[Odot] (unit) {}; 
		\draw[color=green!50!black] (1.5,0.15) -- (unit);
		\fill[color=green!50!black] (0.86,0.95) circle (2pt) node (meet) {};
		\draw[color=green!50!black] (1.75,0.55) -- (A);
		\end{tikzpicture}
		}%
		\circ \iota_B^{\dX,X} \, , 
		&
		\coev_X^{\Beq} & = 
		\vartheta_A^{X,\dX} \circ 
		\tikzzbox{%
		\begin{tikzpicture}[very thick,scale=1.0,color=blue!50!black, baseline=-.8cm,rotate=180]
		\draw[line width=0pt] 
		(3.21,1.85) node[line width=0pt, color=green!50!black] (B) {{\small$B\vphantom{\dX }$}}
		(3,0) node[line width=0pt] (D) {{\small$X\vphantom{\dX }$}}
		(2,0) node[line width=0pt] (s) {{\small$\dX$}}; 
		\draw[redirected] (D) .. controls +(0,1.5) and +(0,1.5) .. (s);
		\fill[color=green!50!black] (2.91,0.85) circle (2pt) node (meet) {};
		\draw[color=green!50!black] (2.91,0.85) .. controls +(0.2,0.25) and +(0,-0.75) .. (B);
		\end{tikzpicture}
		}%
		\, , 
		\nonumber 
		\\ 
		\tev_X^{\Beq} & = 
		\tikzzbox{%
			\begin{tikzpicture}[very thick,scale=1.0,color=blue!50!black, baseline=.8cm, xscale=-1]
			\draw[line width=0pt] 
			(1.75,1.75) node[line width=0pt, color=green!50!black] (A) {{\small$B\vphantom{\dX}$}}
			(1,0) node[line width=0pt] (D) {{\small$X\vphantom{\dX }$}}
			(0,0) node[line width=0pt] (s) {{\small$\Xd $}}; 
			\draw[directed] (D) .. controls +(0,1.5) and +(0,1.5) .. (s);
			\draw[color=green!50!black] (1.25,0.55) .. controls +(0.0,0.25) and +(0.25,-0.15) .. (0.86,0.95);
			\draw[-dot-, color=green!50!black] (1.25,0.55) .. controls +(0,-0.5) and +(0,-0.5) .. (1.75,0.55);
			\draw[color=green!50!black] (1.5,-0.1) node[Odot] (unit) {}; 
			\draw[color=green!50!black] (1.5,0.15) -- (unit);
			\fill[color=green!50!black] (0.86,0.95) circle (2pt) node (meet) {};
			\draw[color=green!50!black] (1.75,0.55) -- (A);
			\end{tikzpicture}
		}%
		\circ \iota_A^{X,\Xd} \, , 
		&
		\tcoev_X^{\Beq} & = 
		\vartheta_B^{\Xd,X} \circ 
		\tikzzbox{%
			\begin{tikzpicture}[very thick,scale=1.0,color=blue!50!black, baseline=-.8cm,rotate=180, xscale=-1]
			\draw[line width=0pt] 
			(3.21,1.85) node[line width=0pt, color=green!50!black] (B) {{\small$A\vphantom{\dX }$}}
			(3,0) node[line width=0pt] (D) {{\small$X\vphantom{\dX }$}}
			(2,0) node[line width=0pt] (s) {{\small$\Xd$}}; 
			\draw[redirected] (D) .. controls +(0,1.5) and +(0,1.5) .. (s);
			\fill[color=green!50!black] (2.91,0.85) circle (2pt) node (meet) {};
			\draw[color=green!50!black] (2.91,0.85) .. controls +(0.2,0.25) and +(0,-0.75) .. (B);
			\end{tikzpicture}
		}%
		\, . 
		\label{eq:AdjointsInBeq}
	\end{align}
\end{proposition}

As a consistency check, we recall that
\be 
\tikzzbox{
	\begin{tikzpicture}[very thick,scale=0.85,color=green!50!black, baseline=0cm]
	\draw[line width=0] 
	(1,1.25) node[line width=0pt] (Ao) {{\small ${}^\dagger\!A$}}
	(-1.5,-1.25) node[line width=0pt] (A1u) {{\small ${}^\dagger\!A$}}
	(-1,-1.25) node[line width=0pt] (A2u) {{\small $A\vphantom{A^\dagger}$}}; 
	\draw[directedgreen] (-0.5,0.7) .. controls +(0,0.75) and +(0,0.75) .. (-1.5,0.7);
	\draw[-dot-] (0,0) .. controls +(0,1) and +(0,1) .. (-1,0);
	\draw[directedgreen] (1,0) .. controls +(0,-1) and +(0,-1) .. (0,0);
	\draw (-1,0) -- (A2u); 
	\draw (1,0) -- (Ao); 
	\draw (-1.5,0.7) -- (A1u); 
	\end{tikzpicture}
}%
, \quad  
\tikzzbox{%
\begin{tikzpicture}[very thick,scale=0.85,color=green!50!black, baseline=0cm]
\draw[line width=0] 
(1,1.25) node[line width=0pt] (A) {{\small ${}^\dagger\!A$}}
(-1,-1.25) node[line width=0pt] (A2) {{\small ${}^\dagger\!A$}}
(-1.5,-1.25) node[line width=0pt, color=green!50!black] (Algebra) {{\small $A\vphantom{A^\dagger}$}}; 
\draw[directedgreen] (0,0) .. controls +(0,1) and +(0,1) .. (-1,0);
\draw[directedgreen] (1,0) .. controls +(0,-1) and +(0,-1) .. (0,0);
\fill[color=green!50!black] (0,0) circle (2.5pt) node[right] {};
\draw[color=green!50!black] (0,0) .. controls +(-0.15,0.15) and +(-0.2,0) .. (0.25,-0.25);
\draw[color=green!50!black] (0.25,-0.25) .. controls +(0.6,0) and +(0.85,0) .. (-0.5,1.25);
\draw[color=green!50!black] (-0.5,1.25) .. controls +(-0.7,0) and +(0,1) .. (-1.5,-0.25);
\draw[<-,color=green!50!black] (0.21,-0.25) -- (0.22,-0.25);
\draw[<-,color=green!50!black] (-0.41,1.25) -- (-0.42,1.25);
\draw[color=green!50!black] (-1.5,-0.25) -- (Algebra); 
\draw (-1,0) -- (A2); 
\draw (1,0) -- (A); 
\end{tikzpicture}
}%
\ee 
are the canonical $A$-actions on~${}^\dagger\!A$, and the $A$-$A$-bimodule structure on~$A^\dagger$ is obtained as the mirror images of the above diagrams. 
These actions agree by assumption on~$A$. 
From this it is straightforward to verify that 
\be 
\label{eq:NakayamaTwistDual}
\tikzzbox{%
	\begin{tikzpicture}[very thick, scale=0.5,color=green!50!black, baseline=-0.9cm]
	\draw[-dot-] (0,-1.8) .. controls +(0,0.5) and +(0,0.5) .. (0.75,-1.8);
	\draw (0.75,-1.8) .. controls +(0,-0.5) and +(0,-0.5) .. (1.5,-1.8);
	\draw (0,-1.8) -- (0,-3);
	\draw[postaction={decorate}, decoration={markings,mark=at position 0.8 with {\arrow{>}}}] (1.5,-0.6) -- (1.5,-1.8);
	\draw (0.375,-0.8) node[Odot] (E) {}; 
	\draw (0.375,-1.4) -- (E);
	\end{tikzpicture}
}%
\colon 
{}_{\gamma_A}A \stackrel{\cong}{\lra} {}^\dagger \! A 
\, , \quad 
\tikzzbox{%
	\begin{tikzpicture}[very thick, scale=0.5,color=green!50!black, baseline=0.2cm]
	\draw (0,0.8) -- (0,2);
	\draw[-dot-] (0,0.8) .. controls +(0,-0.5) and +(0,-0.5) .. (0.75,0.8);
	\draw (0.75,0.8) .. controls +(0,0.5) and +(0,0.5) .. (1.5,0.8);
	\draw[postaction={decorate}, decoration={markings,mark=at position 0.8 with {\arrow{<}}}] (1.5,-0.6) -- (1.5,0.8);
	\draw (0.375,-0.2) node[Odot] (D) {}; 
	\draw (0.375,0.4) -- (D);
	\end{tikzpicture}
}%
\colon 
A^\dagger \stackrel{\cong}{\lra} A_{\gamma_A^{-1}}
\ee 
are bimodule maps. 
Since ${}^\dagger\!A=A^\dagger$ by assumption, it follows that 
\be 
\label{eq:NakayamaLeftRight}
\gamma_A \colon A_{\gamma_A^{-1}} \stackrel{\cong}{\lra} {}_{\gamma_A}A
\ee 
in $\Beq$. 
Moreover, the special case $X=A=B$ in~\eqref{eq:AdjointsInBeq-1} reads 
${}^\star\! A = {}_{\gamma_A^{-1}}({}^\dagger\!A) \cong {}_{\gamma_A^{-1}}({}_{\gamma_A}A) \cong A$ 
and 
$A^\star = (A^\dagger)_{\gamma_A} \cong (A_{\gamma_A^{-1}})_{\gamma_A} \cong A$ 
in $\Beq$, which is consistent with $1_{(\alpha,A)} = A$.

\subsubsection{Symmetric monoidal structure}
\label{subsubsec:SMstructureOfBeq}

As explained in \cite[Cor.\,6.12]{WesterHansenShulman}, the equivariant completion $\Beq$ is the horizontal 2-category of a symmetric monoidal double category $\mathds{B}_{\textrm{eq}}$, which satisfies the conditions under which the symmetric monoidal structure of $\mathds{B}_{\textrm{eq}}$ is passed on to~$\Beq$: 

\begin{proposition}
	$\Beq$ has a symmetric monoidal structure induced from~$\B$. 
\end{proposition}

Here we collect the ingredients of the symmetric monoidal structure on $\Beq$ in graphical presentation, using the conventions of Section~\ref{subsec:DualisabilityInSM2Categories}. 
The monoidal product on objects $(\alpha,A), (\alpha',A') \in \Beq$ is given by 
\be 
(\alpha,A) \btimes^{\Beq} (\alpha',A') 
	= 
	\big( \alpha \btimes^{\B} \alpha' , A \btimes^{\B} A' \big) 
\ee 
where the $\Delta$-separable Frobenius structure on $A \btimes A' \equiv A \btimes^{\B} A'$ is as follows: 
\begin{align}
\mu_{A \btimes A'} & = 
\tikzzbox{%
	\begin{tikzpicture}[very thick,scale=1.0,color=green!50!black, baseline=0.5cm]
	\coordinate (d1) at (-1,0);
	\coordinate (d2) at (+1,0);
	\coordinate (u1) at (-1,1.5);
	\coordinate (u2) at (+1,1.5);
	\coordinate (s) at ($(-0.75,-0.25)$);
	\coordinate (b1) at ($(d1)+(s)$);
	\coordinate (b2) at ($(d2)+(s)$);
	\coordinate (t1) at ($(u1)+(s)$);
	\coordinate (t2) at ($(u2)+(s)$);
	%
	\fill [orange!25!white, opacity=0.8] (u1) -- (u2) -- (d2) -- (d1); 
	%
	\draw[-dot-] ($(d1)+(0.5,0)$) .. controls +(0,1) and +(0,1) .. ($(d1)+(1.5,0)$);
	\draw ($(d1)+(1,0.8)$) -- ($(d1)+(1,1.5)$);
	\fill ($(d1)+(1,0.8)$) circle (0pt) node[below] {{\small $\mu$}};
	%
	\draw[very thick, red!80!black] (u1) -- (u2); 
	\draw[very thick, red!80!black] (d1) -- (d2); 
	\draw[thin, black] (d1) -- (u1); 
	\draw[thin, black] (d2) -- (u2);
	%
	%
	\fill [orange!20!white, opacity=0.7] (b1) -- (b2) -- (t2) -- (t1); 
	%
	\draw[-dot-] ($(d1)+(0.5,0)+(s)$) .. controls +(0,1) and +(0,1) .. ($(d1)+(1.5,0)+(s)$);
	\draw ($(d1)+(1,0.8)+(s)$) -- ($(d1)+(1,1.5)+(s)$);
	\fill ($(d1)+(1,0.8)+(s)$) circle (0pt) node[below] {{\small $\mu'$}};
	%
	\draw[very thick, red!80!black] (t1) -- (t2); 
	\draw[very thick, red!80!black] (b1) -- (b2); 
	\draw[thin, black] (b1) -- (t1); 
	\draw[thin, black] (b2) -- (t2);
	\end{tikzpicture}
}    
\, , & 
\eta_{A \btimes A'} & = 
\tikzzbox{%
	\begin{tikzpicture}[very thick,scale=1.0,color=green!50!black, baseline=0.5cm]
	\coordinate (d1) at (-1,0);
	\coordinate (d2) at (+1,0);
	\coordinate (u1) at (-1,1.5);
	\coordinate (u2) at (+1,1.5);
	\coordinate (s) at ($(-0.75,-0.25)$);
	\coordinate (b1) at ($(d1)+(s)$);
	\coordinate (b2) at ($(d2)+(s)$);
	\coordinate (t1) at ($(u1)+(s)$);
	\coordinate (t2) at ($(u2)+(s)$);
	%
	\fill [orange!25!white, opacity=0.8] (u1) -- (u2) -- (d2) -- (d1); 
	%
	\draw ($(d1)+(1,0.8)$) node[Odot] (unit) {}; 
	\draw (unit) -- ($(d1)+(1,1.5)$);
	\fill ($(d1)+(1,0.8)$) circle (0pt) node[below] {{\small $\eta$}};
	%
	\draw[very thick, red!80!black] (u1) -- (u2); 
	\draw[very thick, red!80!black] (d1) -- (d2); 
	\draw[thin, black] (d1) -- (u1); 
	\draw[thin, black] (d2) -- (u2);
	%
	%
	\fill [orange!20!white, opacity=0.7] (b1) -- (b2) -- (t2) -- (t1); 
	%
	\draw ($(d1)+(1,0.8)+(s)$) node[Odot] (unit) {}; 
	\draw (unit) -- ($(d1)+(1,1.5)+(s)$);
	\fill ($(d1)+(1,0.8)+(s)$) circle (0pt) node[below] {{\small $\eta'$}};
	%
	\draw[very thick, red!80!black] (t1) -- (t2); 
	\draw[very thick, red!80!black] (b1) -- (b2); 
	\draw[thin, black] (b1) -- (t1); 
	\draw[thin, black] (b2) -- (t2);
	\end{tikzpicture}
} 
\, , 
\nonumber
\\
\Delta_{A \btimes A'} & = 
\tikzzbox{%
	\begin{tikzpicture}[very thick,scale=1.0,color=green!50!black, rotate=180, baseline=-0.75cm]
	\coordinate (d1) at (-1,0);
	\coordinate (d2) at (+1,0);
	\coordinate (u1) at (-1,1.5);
	\coordinate (u2) at (+1,1.5);
	\coordinate (s) at ($(-0.75,-0.25)$);
	\coordinate (b1) at ($(d1)+(s)$);
	\coordinate (b2) at ($(d2)+(s)$);
	\coordinate (t1) at ($(u1)+(s)$);
	\coordinate (t2) at ($(u2)+(s)$);
	%
	\fill [orange!25!white, opacity=0.7] (b1) -- (b2) -- (t2) -- (t1); 
	%
	\draw[-dot-] ($(d1)+(0.5,0)+(s)$) .. controls +(0,1) and +(0,1) .. ($(d1)+(1.5,0)+(s)$);
	\draw ($(d1)+(1,0.8)+(s)$) -- ($(d1)+(1,1.5)+(s)$);
	\fill ($(d1)+(1,0.8)+(s)$) circle (0pt) node[above] {{\small $\Delta$}};
	%
	\draw[very thick, red!80!black] (t1) -- (t2); 
	\draw[very thick, red!80!black] (b1) -- (b2); 
	\draw[thin, black] (b1) -- (t1); 
	\draw[thin, black] (b2) -- (t2);
	%
	%
	%
	\fill [orange!20!white, opacity=0.8] (u1) -- (u2) -- (d2) -- (d1); 
	%
	\draw[-dot-] ($(d1)+(0.5,0)$) .. controls +(0,1) and +(0,1) .. ($(d1)+(1.5,0)$);
	\draw ($(d1)+(1,0.8)$) -- ($(d1)+(1,1.5)$);
	\fill ($(d1)+(1,0.8)$) circle (0pt) node[above] {{\small $\Delta'$}};
	%
	\draw[very thick, red!80!black] (u1) -- (u2); 
	\draw[very thick, red!80!black] (d1) -- (d2); 
	\draw[thin, black] (d1) -- (u1); 
	\draw[thin, black] (d2) -- (u2);
	\end{tikzpicture}
}    
\, , & 
\varepsilon_{A \btimes A'} & = 
\tikzzbox{%
	\begin{tikzpicture}[very thick,scale=1.0,color=green!50!black, rotate=180, baseline=-0.75cm]
	\coordinate (d1) at (-1,0);
	\coordinate (d2) at (+1,0);
	\coordinate (u1) at (-1,1.5);
	\coordinate (u2) at (+1,1.5);
	\coordinate (s) at ($(-0.75,-0.25)$);
	\coordinate (b1) at ($(d1)+(s)$);
	\coordinate (b2) at ($(d2)+(s)$);
	\coordinate (t1) at ($(u1)+(s)$);
	\coordinate (t2) at ($(u2)+(s)$);
	%
	\fill [orange!25!white, opacity=0.7] (b1) -- (b2) -- (t2) -- (t1); 
	%
	\draw ($(d1)+(1,0.8)+(s)$) node[Odot] (unit) {}; 
	\draw (unit) -- ($(d1)+(1,1.5)+(s)$);
	\fill ($(d1)+(1,0.8)+(s)$) circle (0pt) node[above] {{\small $\varepsilon$}};
	%
	\draw[very thick, red!80!black] (t1) -- (t2); 
	\draw[very thick, red!80!black] (b1) -- (b2); 
	\draw[thin, black] (b1) -- (t1); 
	\draw[thin, black] (b2) -- (t2);
	%
	%
	%
	\fill [orange!20!white, opacity=0.8] (u1) -- (u2) -- (d2) -- (d1); 
	%
	\draw ($(d1)+(1,0.8)$) node[Odot] (unit) {}; 
	\draw (unit) -- ($(d1)+(1,1.5)$);
	\fill ($(d1)+(1,0.8)$) circle (0pt) node[above] {{\small $\varepsilon'$}};
	%
	\draw[very thick, red!80!black] (u1) -- (u2); 
	\draw[very thick, red!80!black] (d1) -- (d2); 
	\draw[thin, black] (d1) -- (u1); 
	\draw[thin, black] (d2) -- (u2);
	\end{tikzpicture}
}  
\, .
\end{align}
It follows that the Nakayama automorphism of $A\btimes A'$ factorises with respect to~$\btimes$, 
\be 
\gamma_{A\btimes A'} 
= 
\tikzzbox{%
	\begin{tikzpicture}[very thick,scale=0.75,color=green!50!black, baseline=-0.5cm]
	\coordinate (d1) at (-2.5,-3);
	\coordinate (d2) at (+1,-3);
	\coordinate (u1) at (-2.5,2);
	\coordinate (u2) at (+1,2);
	\coordinate (s) at ($(-0.75,-0.5)$);
	\coordinate (b1) at ($(d1)+(s)$);
	\coordinate (b2) at ($(d2)+(s)$);
	\coordinate (t1) at ($(u1)+(s)$);
	\coordinate (t2) at ($(u2)+(s)$);
	%
	\fill [orange!25!white, opacity=0.8] (u1) -- (u2) -- (d2) -- (d1); 
	%
	\draw (0,0.8) -- (0,2);
	\draw[-dot-] (0,0.8) .. controls +(0,-0.5) and +(0,-0.5) .. (-0.75,0.8);
	\draw (-0.75,0.8) .. controls +(0,0.5) and +(0,0.5) .. (-1.5,0.8);
	\draw[-dot-] (0,-1.8) .. controls +(0,0.5) and +(0,0.5) .. (-0.75,-1.8);
	\draw (-0.75,-1.8) .. controls +(0,-0.5) and +(0,-0.5) .. (-1.5,-1.8);
	\draw (0,-1.8) -- (0,-3);
	\draw[postaction={decorate}, decoration={markings,mark=at position .5 with {\arrow{>}}}] (-1.5,0.8) -- (-1.5,-1.8);
	\draw (-0.375,-0.2) node[Odot] (D) {}; 
	\draw (-0.375,0.4) -- (D);
	\draw (-0.375,-0.8) node[Odot] (E) {}; 
	\draw (-0.375,-1.4) -- (E);
	%
	\draw[very thick, red!80!black] (u1) -- (u2); 
	\draw[very thick, red!80!black] (d1) -- (d2); 
	\draw[thin, black] (d1) -- (u1); 
	\draw[thin, black] (d2) -- (u2);
	%
	\fill [orange!20!white, opacity=0.8, shift={($(s)$)}] (t1) -- (t2) -- (b2) -- (b1); 
	%
	\draw[shift={($(s)$)}] (0,0.8) -- (0,2);
	\draw[-dot-,shift={($(s)$)}] (0,0.8) .. controls +(0,-0.5) and +(0,-0.5) .. (-0.75,0.8);
	\draw[shift={($(s)$)}] (-0.75,0.8) .. controls +(0,0.5) and +(0,0.5) .. (-1.5,0.8);
	\draw[-dot-,shift={($(s)$)}] (0,-1.8) .. controls +(0,0.5) and +(0,0.5) .. (-0.75,-1.8);
	\draw[shift={($(s)$)}] (-0.75,-1.8) .. controls +(0,-0.5) and +(0,-0.5) .. (-1.5,-1.8);
	\draw[shift={($(s)$)}] (0,-1.8) -- (0,-3);
	\draw[shift={($(s)$)}, postaction={decorate}, decoration={markings,mark=at position .5 with {\arrow{>}}}] (-1.5,0.8) -- (-1.5,-1.8);
	\draw[shift={($(s)$)}] (-0.375,-0.2) node[Odot] (D) {}; 
	\draw[shift={($(s)$)}] (-0.375,0.4) -- (D);
	\draw[shift={($(s)$)}] (-0.375,-0.8) node[Odot] (E) {}; 
	\draw[shift={($(s)$)}] (-0.375,-1.4) -- (E);
	%
	\draw[very thick, red!80!black] (t1) -- (t2); 
	\draw[very thick, red!80!black] (b1) -- (b2); 
	\draw[thin, black] (b1) -- (t1); 
	\draw[thin, black] (b2) -- (t2);
	%
	\end{tikzpicture}
}
= 
\gamma_A \btimes \gamma_{A'} \, . 
\ee 

On 1-morphisms $X\in\Beq((\alpha,A),(\beta,B))$ and $X'\in\Beq((\alpha',A'),(\beta',B'))$, the monoidal product is $X\btimes^{\B} X'$ with the left $(A\btimes A')$- and right $(B\btimes B')$-action induced from~$\B$. 
On 2-morphisms, we have $\btimes^{\Beq} = \btimes^\B$. 
From now on we will denote the monoidal product of both~$\B$ and $\Beq$ simply as~$\btimes$. 

The braiding in $\Beq$ has 1-morphism components 
\be 
b^{\Beq}_{(\alpha,A), (\alpha',A')} 
	= (A'\btimes A) \otimes b_{\alpha,\alpha'}^\B 
		\cong 
		b^\B_{\alpha,\alpha'} \otimes (A \btimes A') 
\ee 
for $(\alpha, A), (\alpha',A') \in \Beq$, 
with left $(A'\btimes A)$-action given by $\mu_{A'\btimes A}$, and right $(A\btimes A')$-action given by 
\be 
\tikzzbox{%
	\begin{tikzpicture}[thick,scale=1.0,color=black, baseline=2cm]
	\coordinate (p3) at (2,0.5);
	\coordinate (u3) at (2,3.5);
	\coordinate (ld) at (-1,1);
	\coordinate (lu) at (-1,4);
	\coordinate (rd) at (5,1);
	\coordinate (ru) at (5,4);
	\coordinate (rd2) at (5.5,0);
	\coordinate (ru2) at (5.5,3);
	\coordinate (ld2) at (-1.5,0);
	\coordinate (lu2) at (-1.5,3);
	%
	\fill [orange!25!white, opacity=0.8] (p3) -- (ld) -- (lu) -- (u3); 
	\fill [orange!25!white, opacity=0.8] (p3) -- (rd) -- (ru) -- (u3); 
	\draw[thin] (lu) --  (ld); 
	\draw[thin] (ru) --  (rd); 
	%
	\coordinate (xlb) at (4,0.85);
	\coordinate (xrb) at (0.5,3.75);
	\draw[ultra thick, color=green!50!black] (xlb) .. controls +(0,0.25) and +(0.2,-0.5) .. (2,2); 
	\draw[ultra thick, color=green!50!black] (xrb) -- (0.5,0.75); 
	\draw[ultra thick, color=green!50!black] (0.5,3) .. controls +(0,-0.25) and +(-0.2,0.5) .. (2,2); 
	\fill[color=green!50!black] (0.5,3) circle (2.5pt) node {};
	\fill[color=green!50!black] (4.1,1.1) circle (0pt) node {{\small $A$}};
	\fill[color=green!50!black] (0.25,3.5) circle (0pt) node {{\small $A'$}};
	%
	\draw[very thick, red!80!black] (p3) -- (ld);
	\draw[very thick, red!80!black] (u3) -- (lu); 
	\draw[very thick, red!80!black] (p3) -- (rd);
	\draw[very thick, red!80!black] (ru) -- (u3); 
	%
	\fill [orange!30!white, opacity=0.8] (p3) -- (ld2) -- (lu2) -- (u3); 
	\fill [orange!30!white, opacity=0.8] (p3) -- (rd2) -- (ru2) -- (u3); 
	%
	\draw[ultra thick] (p3) --  (u3); 
	%
	\coordinate (xlf) at (3.25,0.32);
	\coordinate (xrf) at (0,3.2);
	\draw[ultra thick, color=green!50!black] (xlf) .. controls +(0,0.5) and +(0.2,-0.5) .. (2,2); 
	\draw[ultra thick, color=green!50!black] (0,2.7) .. controls +(0,-0.5) and +(-0.2,0.3) .. (2,2); 
	\fill[color=green!50!black] (0,2.7) circle (2.5pt) node {};
	\draw[ultra thick, color=green!50!black] (xrf) -- (0,0.2); 
	\fill[color=green!50!black] (3.5,0.5) circle (0pt) node {{\small $A'$}};
	\fill[color=green!50!black] (-0.2,0.5) circle (0pt) node {{\small $A$}};
	\fill[color=blue!50!black] (2,2) circle (2.5pt) node[right] {{\small $b_{A,A'}$}};
	\draw[thin] (ru2) --  (rd2); 
	\draw[thin] (lu2) --  (ld2); 
	%
	\draw[very thick, red!80!black] (ru2) -- (u3); 
	\draw[very thick, red!80!black] (rd2) -- (p3); 
	\draw[very thick, red!80!black] (lu2) -- (u3); 
	\draw[very thick, red!80!black] (ld2) -- (p3); 
	%
	%
	\fill (-0.5,3.6) circle (0pt) node {{\small $\alpha'$}};
	\fill (4.5,3.6) circle (0pt) node {{\small $\alpha\vphantom{\alpha'}$}};
	\fill (-1,0.5) circle (0pt) node {{\small $\alpha\vphantom{\alpha'}$}};
	\fill (5,0.5) circle (0pt) node {{\small $\alpha'$}};
	%
	\end{tikzpicture}
}   
\, .
\ee 
The 2-morphism components of $b^{\Beq}$ are 
\be 
b^{\Beq}_{X,X'} = 
\tikzzbox{%
	\begin{tikzpicture}[thick,scale=1.0,color=black, baseline=1.9cm]
	\coordinate (p3) at (2,0.5);
	\coordinate (u3) at (2,3.5);
	\coordinate (ld) at (-1,1);
	\coordinate (lu) at (-1,4);
	\coordinate (rd) at (5,1);
	\coordinate (ru) at (5,4);
	\coordinate (rd2) at (5.5,0);
	\coordinate (ru2) at (5.5,3);
	\coordinate (ld2) at (-1.5,0);
	\coordinate (lu2) at (-1.5,3);
	%
	\fill [orange!25!white, opacity=0.8] (p3) -- (ld) -- (lu) -- (u3); 
	\fill [orange!25!white, opacity=0.8] (p3) -- (rd) -- (ru) -- (u3); 
	\draw[thin] (lu) --  (ld); 
	\draw[thin] (ru) --  (rd); 
	%
	\coordinate (xlb) at (0,0.85);
	\coordinate (xrb) at (3.5,3.75);
	\draw[ultra thick, color=blue!50!black] (xlb) .. controls +(0,0.25) and +(-0.2,-0.5) .. (2,2); 
	\draw[ultra thick, color=blue!50!black] (xrb) .. controls +(0,-0.75) and +(0.2,0.5) .. (2,2); 
	\fill[color=blue!50!black] (-0.1,1.1) circle (0pt) node {{\small $X'$}};
	\fill[color=blue!50!black] (3.7,3.5) circle (0pt) node {{\small $X$}};
	%
	\coordinate (A'lb) at (0.5,0.75);
	\fill[color=green!50!black] (0.75,0.95) circle (0pt) node {{\small $A'$}};
	\draw[ultra thick, color=green!50!black] (A'lb) .. controls +(0,0.1) and +(0.1,-0.1) .. (0.45,1.14); 
	\fill[color=green!50!black] (0.45,1.14) circle (2.5pt) node[left] {};
	%
	\coordinate (B'ub) at (-0.15,3.85);
	\fill[color=green!50!black] (0.15,3.55) circle (0pt) node {{\small $B'$}};
	\draw[ultra thick, color=green!50!black] (B'ub) .. controls +(0,-0.5) and +(-0.5,0.5) .. (0.8,1.28); 
	\fill[color=green!50!black] (0.8,1.28) circle (2.5pt) node[left] {};
	%
	\draw[very thick, red!80!black] (p3) -- (ld);
	\draw[very thick, red!80!black] (u3) -- (lu); 
	\draw[very thick, red!80!black] (p3) -- (rd);
	\draw[very thick, red!80!black] (ru) -- (u3); 
	%
	\fill [orange!30!white, opacity=0.8] (p3) -- (ld2) -- (lu2) -- (u3); 
	\fill [orange!30!white, opacity=0.8] (p3) -- (rd2) -- (ru2) -- (u3); 
	%
	\draw[ultra thick] (p3) --  (u3); 
	\fill (2.5,1) circle (0pt) node {{\small $b_{\alpha,\alpha'}$}};
	\fill (1.6,3) circle (0pt) node {{\small $b_{\beta,\beta'}$}};\\
	%
	\coordinate (xlf) at (0.75,0.32);
	\coordinate (xrf) at (4,3.2);
	\draw[ultra thick, color=blue!50!black] (xlf) .. controls +(0,0.5) and +(-0.2,-0.5) .. (2,2); 
	\draw[ultra thick, color=blue!50!black] (xrf) .. controls +(0,-0.75) and +(0.2,0.3) .. (2,2); 
	\fill[color=blue!50!black] (0.55,0.5) circle (0pt) node {{\small $X$}};
	\fill[color=blue!50!black] (4.2,2.9) circle (0pt) node {{\small $X'$}};
	\fill[color=blue!50!black] (2,2) circle (2.5pt) node[left] {{\small $b_{X,X'}$}};
	%
	\coordinate (Alf) at (1.3,0.4);
	\draw[ultra thick, color=green!50!black] (Alf) .. controls +(0,0.25) and +(0.1,-0.1) .. (1.12,0.97); 
	\fill[color=green!50!black] (1.45,0.8) circle (0pt) node {{\small $A$}};
	\fill[color=green!50!black] (1.12,0.97) circle (2.5pt) node[left] {};
	%
	\coordinate (B) at (0.3,3.25);
	\draw[ultra thick, color=green!50!black] (B) .. controls +(0,-0.75) and +(-0.5,0.1) .. (1.38,1.25); 
	\fill[color=green!50!black] (0.55,3) circle (0pt) node {{\small $B$}};
	\fill[color=green!50!black] (1.38,1.25) circle (2.5pt) node[left] {};
	%
	\draw[very thick, red!80!black] (ru2) -- (u3); 
	\draw[very thick, red!80!black] (rd2) -- (p3); 
	\draw[very thick, red!80!black] (lu2) -- (u3); 
	\draw[very thick, red!80!black] (ld2) -- (p3); 
	\draw[thin] (ru2) --  (rd2); 
	\draw[thin] (lu2) --  (ld2); 
	%
	%
	\fill (-0.5,3.5) circle (0pt) node {{\small $\beta'$}};
	\fill (4.5,3.5) circle (0pt) node {{\small $\alpha\vphantom{\beta'}$}};
	\fill (-1,0.5) circle (0pt) node {{\small $\beta\vphantom{\alpha'\beta}$}};
	\fill (5,0.5) circle (0pt) node {{\small $\alpha'\vphantom{\alpha'\beta}$}};
	%
	\end{tikzpicture}
}  
\ee 
for all $X\in\Beq((\alpha,A),(\beta,B))$ and $X'\in\Beq((\alpha',A'),(\beta',B'))$. 

\begin{remark}
	 It follows that $\B \subset \Beq$ is an embedding of symmetric monoidal 2-categories. 
	 In particular, the braiding components of 1-morphisms in $\Beq((\one, 1_\one), (\one, 1_\one))$ are those in $\B(\one,\one)$. 
\end{remark}

\subsubsection{Duality and Serre automorphism}

If $\alpha \in \B$ is dualisable with duality data $(\alpha^\dual, \tev_\alpha, \tcoev_\alpha)$, then every object $(\alpha,A) \in \Beq$ is dualisable with 
\begin{align}
(\alpha,A)^\dual = \big( \alpha^\dual, A^\dual \big) 
	\, , \quad 
	\tev_{(\alpha,A)} 
		& = \tev_\alpha \otimes (A \btimes 1_{\alpha^\dual}) 
			\cong \tev_\alpha \otimes (1_\alpha \btimes A^\dual )\, , 
	\nonumber 
	\\ 
	\tcoev_{(\alpha,A)} 
	& = (1_{\alpha^\dual} \btimes A) \otimes \tcoev_\alpha 
	\cong ( A^\dual \btimes 1_\alpha) \otimes \tcoev_\alpha \, , 
\end{align}
where we used the isomorphism 
\be 
\tikzzbox{%
	\begin{tikzpicture}[thick,scale=1.0,color=black, baseline=1.5cm]
	\coordinate (p1) at (0,0);
	\coordinate (p2) at (2,-0.5);
	\coordinate (p3) at (2.5,0.5);
	\coordinate (u1) at (0,3);
	\coordinate (u2) at (2,2.5);
	\coordinate (u3) at (2.5,3.5);
	%
	\fill [orange!20!white, opacity=0.8] 
	(p1) .. controls +(0,0.25) and +(-1,0) ..  (p3)
	-- (p3) --  (u3)
	-- (u3) .. controls +(-1,0) and +(0,0.25) ..  (u1)
	;
	\draw[ultra thick, green!50!black] (0.5,0.25) .. controls +(0,1) and +(0,-0.25) .. (0,1.5); 
	\fill[color=green!50!black] (0.75,0.6) circle (0pt) node {{\small $A$}};
	%
	\draw[very thick, red!80!black] (p1) .. controls +(0,0.25) and +(-1,0) ..  (p3); 
	%
	\fill [orange!30!white, opacity=0.8] 
	(p1) .. controls +(0,-0.25) and +(-1,0) ..  (p2)
	-- (p2) --  (u2)
	-- (u2) .. controls +(-1,0) and +(0,-0.25) ..  (u1)
	;
	\draw[ultra thick, green!50!black] (0.75,2.62) .. controls +(0,-1) and +(0,0.25) .. (0,1.5); 
	\fill[color=green!50!black] (1.05,2.3) circle (0pt) node {{\small $A^\dual$}};
	%
	\draw[thin] (p1) --  (u1); 
	\draw[thin] (p2) --  (u2); 
	\draw[thin] (p3) --  (u3); 
	%
	\draw[very thick, red!80!black] (p1) .. controls +(0,-0.25) and +(-1,0) ..  (p2); 
	\draw[very thick, red!80!black] (u1) .. controls +(0,-0.25) and +(-1,0) ..  (u2); 
	\draw[very thick, red!80!black] (u1) .. controls +(0,0.25) and +(-1,0) ..  (u3); 
	%
	\fill (1.5,0) circle (0pt) node {{\small $\alpha^\dual$}};
	\fill (2,3) circle (0pt) node {{\small $\alpha\vphantom{\alpha^\dual}$}};
	%
	\end{tikzpicture}
}  
\,  
\colon \tev_\alpha \otimes (A \btimes 1_{\alpha^\dual}) \stackrel{\cong}{\lra} \tev_\alpha \otimes (1_\alpha \btimes A^\dual )
\ee 
induced by the inverse cusp isomorphism~$c_{\textrm{l}}^{-1}$ in~\eqref{eq:cuspl}, and similarly with $\tcoev_{(\alpha,A)}$ and~$c_{\textrm{r}}^{-1}$. 
These isomorphisms together with~$\mu_A$ also give the above adjunction morphisms their bimodule structures. 

The Frobenius algebra structure $(A^\dual, \mu_{A^\dual}, \eta_{A^\dual}, \Delta_{A^\dual}, \varepsilon_{A^\dual})$ on $A^\dual$ is by definition $(A^\dual, \mu_A^\dual, \eta_A^\dual, \Delta_A^\dual, \varepsilon_A^\dual)$ up to cusp isomorphisms as needed. 
We illustrate this with the multiplication 
\be 
\label{eq:MuAdual}
\mu_{A^\dual} = 
\tikzzbox{%
	\begin{tikzpicture}[thick,scale=0.75,color=black, baseline=-0.5cm]
	\coordinate (p1) at (0,0);
	\coordinate (p2) at (1.5,-0.5);
	\coordinate (p2s) at (4,-0.5);
	\coordinate (p3) at (1.5,0.5);
	\coordinate (p4) at (3,1);
	\coordinate (p5) at (1.5,1.5);
	\coordinate (p6) at (-1,1.5);
	\coordinate (u1) at (0,2.5);
	\coordinate (u2) at (1.5,2);
	\coordinate (u2s) at (4,3.5);
	\coordinate (u3) at (1.5,3);
	\coordinate (u4) at (3,3.5);
	\coordinate (u5) at (1.5,4);
	\coordinate (u6) at (-1,3.5);
	%
	\draw[thin] ($(p2s)+(0,-4)$) --  (u2s); 
	\draw[thin] ($(p6)+(0,-6)$) --  (u6); 
	\fill [orange!20!white, opacity=0.8] 
	($(p6)+(0,-6)$) -- ($(p2s)+(0,-4)$) -- (u2s) -- (u6);
	%
	\draw[very thick, color=green!50!black] ($(p6)+(1.5,-6)$) .. controls +(0,1) and +(0,-1.5) .. ($(p6)+(3.5,-3)$);
	\draw[-dot-, very thick, color=green!50!black] ($(p6)+(3.5,-3)$) .. controls +(0,1.5) and +(0,1.5) .. ($(p6)+(1.5,-3)$);
	\draw[very thick, color=green!50!black] ($(p6)+(2.5,-1.9)$) .. controls +(0,1) and +(-0.1,-0.25) .. ($(p6)+(1.5,-0.5)$);
	\draw[very thick, color=green!50!black] ($(p6)+(1.5,-3)$) .. controls +(0,-0.5) and +(-0.25,0.25) .. ($(p6)+(1.2,-3.9)$);
	\fill ($(p6)+(2.5,-1.9)$) circle (0pt) node[below] {{\small $\mu_A$}};
	\fill[color=green!50!black] ($(p6)+(1.15,-5.65)$) circle (0pt) node {{\small $A^\dual$}};
	\fill[color=green!50!black] ($(p6)+(3.95,-5.65)$) circle (0pt) node {{\small $A^\dual$}};
	\fill[color=green!50!black] ($(p6)+(3.9,1.65)$) circle (0pt) node {{\small $A^\dual$}};
	%
	%
	\fill [orange!30!white, opacity=0.8] 
	(p1) .. controls +(0,0.5) and +(0,0) ..  (1.5,2.5)
	-- (1.5,2.5) .. controls +(0,0) and +(0,0.5) .. (p4)
	-- (p4) -- ($(p4)+(0,-2)$) 
	-- ($(p4)+(0,-2)$) .. controls +(0,-0.5) and +(0,0) .. (1.5,-3.5)
	-- (1.5,-3.5) .. controls +(0,0) and +(0,-0.5) .. ($(p1)+(0,-2)$)
	;
	%
	\draw[thin, dotted] (p4) .. controls +(0,0.5) and +(0,0) ..  (1.5,2.5); 
	\draw[thin] (p1) .. controls +(0,0.5) and +(0,0) ..  (1.5,2.5); 
	%
	\draw[thin] (p1) -- ($(p1)+(0,-2)$);
	\draw[thin, dotted] (p4) -- ($(p4)+(0,-2)$);
	%
	\draw[thin, dotted] ($(p4)+(0,-2)$) .. controls +(0,-0.5) and +(0,0) ..  ($(1.5,-1.5)+(0,-2)$); 
	\draw[thin] ($(p1)+(0,-2)$) .. controls +(0,-0.5) and +(0,0) ..  ($(1.5,-1.5)+(0,-2)$); 
	%
	\draw[color=green!50!black] ($(p6)+(1.2,-3.9)$) .. controls +(0.25,-0.25) and +(0,1.5) .. ($(p6)+(3.5,-6)$);
	\draw[very thick, color=green!50!black] ($(p6)+(1.5,-0.5)$) .. controls +(0.1,0.25) and +(0,-1.5) .. ($(u6)+(3.5,0)$);
	%
	\draw[very thick, red!80!black] (u6) -- (u2s); 
	\draw[very thick, red!80!black] ($(p6)+(0,-6)$) -- ($(p2s)+(0,-4)$); 
	%
	\fill[color=blue!50!black] ($(u6)+(0.5,-0.5)$) circle (0pt) node {{\small $\alpha^\dual$}};
	\fill (1.5,2.5) circle (2.5pt) node[above] {{\small $c_{\textrm{r}}$}};
	\fill ($(1.5,-1.5)+(0,-2)$) circle (2.5pt) node[below] {{\small $c_{\textrm{r}}^{-1}$}};
	%
	\end{tikzpicture}
}   
\equiv 
\tikzzbox{%
	\begin{tikzpicture}[very thick,scale=0.75,color=green!50!black, baseline]
	\draw (-0.85,0.2) -- (-0.85,1); 
	\draw (-1.2,-0.3) .. controls +(0,0.25) and +(0,0) .. (-0.85,0.2);
	\draw (-0.5,-0.3) .. controls +(0,0.25) and +(0,0) .. (-0.85,0.2);
	\draw (-0.5,-0.3) .. controls +(0,-0.25) and +(0,0.25) .. (-1.2,-1);
	\draw[color=white, line width=4pt] (-1.2,-0.3) .. controls +(0,-0.25) and +(0,0.25) .. (-0.5,-1);
	\draw (-1.2,-0.3) .. controls +(0,-0.25) and +(0,0.25) .. (-0.5,-1);
	\fill (-0.85,0.2) circle (2.5pt) node (mult1) {};
	\end{tikzpicture} 
}%
\ee 
where the last expression is shorthand for the defining Gray diagram in the middle. 
Note that there is a canonical isomorphism $A\cong A^{\dual\dual}$ in~$\B$ (see \cite[Fig.\,32]{BMS}), which we leave implicit. 
Similarly, we denote the other structure maps of $A^\dual$ as 
\be 
\label{eq:RestAdual}
\eta_{A^\dual} = 
\tikzzbox{%
\begin{tikzpicture}[very thick,scale=0.4,color=green!50!black, baseline=0]
\draw (-0.5,-0.5) node[Odot] (unit) {}; 
\draw (unit) -- (-0.5,0.7);
\end{tikzpicture} 
}%
\, , \quad 
\Delta_{A^\dual} = 
\tikzzbox{%
\begin{tikzpicture}[very thick,scale=0.75,color=green!50!black, rotate=180, baseline]
\draw (-0.85,0.2) -- (-0.85,1); 
\fill (-0.85,0.2) circle (2.5pt) node (mult1) {};
\draw (-1.2,-0.3) .. controls +(0,0.25) and +(0,0) .. (-0.85,0.2);
\draw (-1.2,-0.3) .. controls +(0,-0.25) and +(0,0.25) .. (-0.5,-1);
\draw (-0.5,-0.3) .. controls +(0,0.25) and +(0,0) .. (-0.85,0.2);
\draw[color=white, line width=4pt] (-0.5,-0.3) .. controls +(0,-0.25) and +(0,0.25) .. (-1.2,-1);
\draw (-0.5,-0.3) .. controls +(0,-0.25) and +(0,0.25) .. (-1.2,-1);
\end{tikzpicture} 
}%
\, , \quad 
\varepsilon_{A^\dual} = 
\tikzzbox{%
\begin{tikzpicture}[very thick,scale=0.4,color=green!50!black, rotate=180, baseline=0]
\draw (-0.5,-0.5) node[Odot] (unit) {}; 
\draw (unit) -- (-0.5,0.7);
\end{tikzpicture} 
}%
\, . 
\ee 

The \textsl{enveloping algebra of~$A$} is 
\be 
\Ae = A \btimes A^\dual \, . 
\ee 
The cusp isomorphisms together with~$\mu_A$ give a canonical right $\Ae$-module structure on $\tev_{(\alpha,A)}$ and a left $(\Ae)^\dual$-module structure on $\tcoev_{(\alpha,A)}$. 

\begin{lemma}
	If $\alpha \in \B$ is dualisable, then for $(\alpha,A)\in\Beq$ we have that 
	\be 
	\gamma_{A^\dual} = (\gamma_A^\dual)^{-1}
	\ee 
	up to cusp isomorphisms in~$\B$. 
\end{lemma}
\begin{proof}
	Up to cusp isomorphisms, $\tev_{A^\dual} \colon A^\dual \otimes (A^\dual)^\dagger \lra 1_{\alpha^\dual}$ agrees with $\ev_A^\dual$, 
	\be 
	\ev_A^\dual 
		= 
		\tikzzbox{
			\begin{tikzpicture}[very thick,scale=1.0,color=green!50!black, baseline=-.3cm]
			\draw[line width=0pt] 
			(0,-1) node[line width=0pt] (D) {{\small$A^\dual\vphantom{{}^\dagger\!A^\dual}$}}
			(1,-1) node[line width=0pt] (s) {{\small${}^\dagger\!A^\dual$}}; 
			\draw[redirectedgreen] (0,0) .. controls +(0,0.5) and +(0,0.5) .. (1,0);
			\draw (1,0) .. controls +(0,-0.5) and +(0,0.75) .. (D);
			\draw[color=white, line width=4pt] (0,0) .. controls +(0,-0.5) and +(0,0.75) .. (s);
			\draw (0,0) .. controls +(0,-0.5) and +(0,0.75) .. (s);
			\end{tikzpicture}
		} 
		\equiv 
		\tikzzbox{
			\begin{tikzpicture}[very thick,scale=1.0,color=green!50!black, baseline=.4cm]
			\draw[line width=0pt] 
			(3,0) node[line width=0pt] (D) {{\small$(A^\dual)^\dagger$}}
			(2,0) node[line width=0pt] (s) {{\small$A^\dual\vphantom{(A^\dual)^\dagger}$}}; 
			\draw[redirectedgreen] (D) .. controls +(0,1) and +(0,1) .. (s);
			\end{tikzpicture}
		} 
		= 
		\tev_{A^\dual} \, , 
	\ee 
	and similarly for the coevaluations: 
	\be 
	\tcoev_A^\dual 
	= 
	\tikzzbox{
		\begin{tikzpicture}[very thick,scale=1.0,color=green!50!black, rotate=180, baseline=.3cm]
		\draw[line width=0pt] 
		(0,-1) node[line width=0pt] (D) {{\small$(A^\dagger)^\dual$}}
		(1,-1) node[line width=0pt] (s) {{\small$A^\dual\vphantom{{}^\dagger\!(A^\dual)}$}}; 
		\draw[redirectedgreen] (0,0) .. controls +(0,0.5) and +(0,0.5) .. (1,0);
		\draw (0,0) .. controls +(0,-0.5) and +(0,0.75) .. (s);
		\draw[color=white, line width=4pt] (1,0) .. controls +(0,-0.5) and +(0,0.75) .. (D);
		\draw (1,0) .. controls +(0,-0.5) and +(0,0.75) .. (D);
		\end{tikzpicture}
	} 
	\equiv 
	\tikzzbox{
		\begin{tikzpicture}[very thick,scale=1.0,color=green!50!black, baseline=-.4cm,rotate=180]
		\draw[line width=0pt] 
		(3,0) node[line width=0pt] (D) {{\small$A^\dual\vphantom{{}^\dagger\!A^\dual}$}}
		(2,0) node[line width=0pt] (s) {{\small${}^\dagger\!A^\dual$}}; 
		\draw[redirectedgreen] (D) .. controls +(0,1) and +(0,1) .. (s);
		\end{tikzpicture}
	}  
	= 
	\coev_{A^\dual}	 . 
	\ee 
	Hence together with~\eqref{eq:MuAdual} and~\eqref{eq:RestAdual}, we find 
	\be 
	\gamma_{A^\dual}^{-1} 
	= 
	\tikzzbox{%
		\begin{tikzpicture}[very thick, scale=0.5,color=green!50!black, baseline=-0.1cm]
		\draw (-0.375,0.2) node[Odot] (D) {}; 
		\draw (-0.375,0.8) -- (D);
		\draw[-dot-] (0,1.2) .. controls +(0,-0.5) and +(0,-0.5) .. (-0.75,1.2);
		\draw (0,1.2) .. controls +(0,0.5) and +(0,-2) .. (-0.75,4);
		\draw[color=white, line width=4pt] (-0.75,1.2) .. controls +(0,0.5) and +(0,-0.5) .. (1,2.75);
		\draw (-0.75,1.2) .. controls +(0,0.5) and +(0,-0.5) .. (1,2.75);
		\draw (1,2.75) .. controls +(0,0.75) and +(0,0.75) .. (0.25,2.75);
		\draw[color=white, line width=4pt] (0.25,2.75) .. controls +(0,-1) and +(0,1) .. (1,0);
		\draw[postaction={decorate}, decoration={markings,mark=at position 1 with {\arrow{>}}}] (0.25,2.75) .. controls +(0,-1) and +(0,1) .. (1,0);
		\draw[color=white, line width=4pt] (0.25,-2.75) .. controls +(0,1) and +(0,-1) .. (1,0);
		\draw (0.25,-2.75) .. controls +(0,1) and +(0,-1) .. (1,0);
		\draw (1,-2.75) .. controls +(0,-0.75) and +(0,-0.75) .. (0.25,-2.75);
		\draw (0,-1.2) .. controls +(0,-0.5) and +(0,2) .. (-0.75,-4);
		\draw[color=white, line width=4pt] (-0.75,-1.2) .. controls +(0,-0.5) and +(0,0.5) .. (1,-2.75);
		\draw (-0.75,-1.2) .. controls +(0,-0.5) and +(0,0.5) .. (1,-2.75);
		\draw[color=white, line width=4pt] (0.25,-2.75) .. controls +(0,1) and +(0,-1) .. (1,0);
		\draw (0.25,-2.75) .. controls +(0,1) and +(0,-1) .. (1,0);
		\draw (-0.375,-0.2) node[Odot] (E) {}; 
		\draw (-0.375,-0.8) -- (E);
		\draw[-dot-] (0,-1.2) .. controls +(0,0.5) and +(0,0.5) .. (-0.75,-1.2);
		\end{tikzpicture}
	}%
	= 
	\tikzzbox{%
		\begin{tikzpicture}[very thick, scale=0.5,color=green!50!black, baseline=-0.35cm]
		\draw (0,0.8) -- (0,3.5);
		\draw[-dot-] (0,0.8) .. controls +(0,-0.5) and +(0,-0.5) .. (-0.75,0.8);
		\draw (-0.75,0.8) .. controls +(0,0.5) and +(0,0.5) .. (-1.5,0.8);
		\draw[-dot-] (0,-1.8) .. controls +(0,0.5) and +(0,0.5) .. (-0.75,-1.8);
		\draw (-0.75,-1.8) .. controls +(0,-0.5) and +(0,-0.5) .. (-1.5,-1.8);
		\draw (0,-1.8) -- (0,-4.5);
		\draw[postaction={decorate}, decoration={markings,mark=at position .5 with {\arrow{>}}}] (-1.5,0.8) -- (-1.5,-1.8);
		\draw (-0.375,-0.2) node[Odot] (D) {}; 
		\draw (-0.375,0.4) -- (D);
		\draw (-0.375,-0.8) node[Odot] (E) {}; 
		\draw (-0.375,-1.4) -- (E);
		\end{tikzpicture}
	}%
	=
	\gamma_A^\dual \, . 
	\ee 
\end{proof}

\medskip 

We now turn to full dualisability, which is another property that is compatible with equivariant completion: 

\begin{proposition}
	Let $\alpha\in\B$ be fully dualisable. 
	Then every $(\alpha,A)\in\Beq$ is fully dualisable. 
\end{proposition}
\begin{proof}
	If $\tev_\alpha, \tcoev_\alpha$ witness~$\alpha^\dual$ as a dual of~$\alpha$ in~$\B$, then $\tev_{(\alpha,A)} = \tev_\alpha \otimes (A \btimes 1_{\alpha^\dual})$, $\tcoev_{(\alpha,A)} = (1_{\alpha^\dual} \btimes A) \otimes \tcoev_\alpha$ witness $(\alpha^\dual, A^\dual)$ as a dual of $(\alpha, A)$ in $\Beq$. 
	Moreover, by Proposition~\ref{prop:AdjointsInBeq}, these adjunction 1-morphisms have adjoints themselves since $A\in\B(\alpha,\alpha)$ has adjoints thanks to its Frobenius algebra structure, while $\tev_\alpha, \tcoev_\alpha$ have adjoints by assumption. 
\end{proof}

In particular, according to~\eqref{eq:AdjointsInBeq-1} we have 
\be 
\tev_{(\alpha,A)}^\star 
	= 
	\big( \tev_{(\alpha,A)}^\dagger \big)_{\gamma_{1_\one}}
	= 
	\big[ \tev_\alpha \otimes (A \btimes 1_{\alpha^\dual}) \big]^\dagger 
	\cong 
	\big( A^\dagger \btimes 1_{\alpha^\dual} \big) \otimes \tev_\alpha^\dagger \, . 
\ee 
Hence the Serre automorphism of $(\alpha,A)\in\Beq$ is 
\begin{align}
S_{(\alpha,A)} 
	 & = 
	 \big( A \btimes \tev_{(\alpha,A)} \big) \otimes \big( b_{(\alpha,A), (\alpha,A)} \btimes A^\dual \big) \otimes \big( A \btimes \tev_{(\alpha,A)}^\star \big) 
	 \nonumber
	 \\
	 &\widehat{=} \;
	 \tikzzbox{%
	 	\begin{tikzpicture}[thick,scale=1.0,color=black, baseline=2cm]
	 	\coordinate (p1) at (0,0);
	 	\coordinate (p2) at (2,-0.5);
	 	\coordinate (p3) at (2,0.5);
	 	\coordinate (p4) at (4,0);
	 	\coordinate (u1) at (0,3);
	 	\coordinate (u2) at (2,2.5);
	 	\coordinate (u3) at (2,3.5);
	 	\coordinate (u4) at (4,3);
	 	\coordinate (ld) at (-2,1);
	 	\coordinate (lu) at (-2,4);
	 	\coordinate (rd) at (6,1);
	 	\coordinate (ru) at (6,4);
	 	%
	 	\fill [orange!20!white, opacity=0.8] (p3) -- (ld) -- (lu) -- (u3); 
	 	\fill [orange!20!white, opacity=0.8] (p3) -- (rd) -- (ru) -- (u3); 
	 	%
	 	\fill [orange!25!white, opacity=0.8] 
	 	(p1) .. controls +(0,0.25) and +(-1,-0.2) ..  (p3)
	 	-- (p3) .. controls +(1,-0.2) and +(0,0.25) ..  (p4)
	 	-- (p4) --  (u4)
	 	-- (u4) .. controls +(0,0.25) and +(1,-0.2) ..  (u3)
	 	-- (u3) .. controls +(-1,-0.2) and +(0,0.25) ..  (u1)
	 	;
	 	%
	 	\draw[ultra thick] (p3) --  (u3); 
	 	\fill (2.5,3) circle (0pt) node {{\small $b_{\alpha,\alpha}$}};
	 	%
	 	\draw[ultra thick, green!50!black] (-1.2,0.88) --  (-1.2,3.88);
	 	\draw[ultra thick, green!50!black] (5.2,0.88) --  (5.2,3.88);
	 	\fill[color=green!50!black] (0.5,3) circle (0pt) node {{\small $A\vphantom{A^\dagger}$}};
	 	\fill[green!50!black] (3.55,3) circle (0pt) node {{\small $A^\dagger$}};
	 	\fill[green!50!black] (5.4,1.5) circle (0pt) node {{\small $A$}};
	 	\fill[green!50!black] (-1.4,1.5) circle (0pt) node {{\small $A$}};
	 	\draw[very thick, red!80!black] (p3) -- (ld); 										
	 	\draw[very thick, red!80!black] (p3) -- (rd);										
	 	\draw[ultra thick, green!50!black] (0.7,0.27) --  (0.7,3.27);
	 	\draw[ultra thick, green!50!black] (3.3,0.27) --  (3.3,3.27);
	 	%
	 	\draw[very thick, red!80!black] (p1) .. controls +(0,0.25) and +(-1,-0.2) ..  (p3); 
	 	\draw[very thick, red!80!black] (p4) .. controls +(0,0.25) and +(1,-0.2) ..  (p3);
	 	\draw[very thick, red!80!black] (u3) -- (lu); 
	 	\draw[very thick, red!80!black] (ru) -- (u3); 
	 	\draw[thin] (lu) --  (ld); 
	 	\draw[thin] (ru) --  (rd); 
	 	%
	 	\fill [orange!30!white, opacity=0.8] 
	 	(p1) .. controls +(0,-0.25) and +(-1,0) ..  (p2)
	 	-- (p2) .. controls +(1,0) and +(0,-0.25) ..  (p4)
	 	-- (p4) --  (u4)
	 	-- (u4) .. controls +(0,-0.25) and +(1,0) ..  (u2)
	 	-- (u2) .. controls +(-1,0) and +(0,-0.25) ..  (u1)
	 	;
	 	\draw[thin] (p1) --  (u1); 
	 	\draw[thin] (p4) --  (u4); 
	 	%
	 	\draw[very thick, red!80!black] (p1) .. controls +(0,-0.25) and +(-1,0) ..  (p2)
	 	--(p2) .. controls +(1,0) and +(0,-0.25) .. (p4); 
	 	\draw[very thick, red!80!black] (u1) .. controls +(0,-0.25) and +(-1,0) ..  (u2) 
	 	--(u2) .. controls +(1,0) and +(0,-0.25) .. (u4); 
	 	\draw[very thick, red!80!black] (u4) .. controls +(0,0.25) and +(1,-0.2) ..  (u3); 
	 	\draw[very thick, red!80!black] (u1) .. controls +(0,0.25) and +(-1,-0.2) .. (u3); 
	 	%
	 	\fill (-0.3,1.8) circle (0pt) node {{\small $\tev_{\alpha}\vphantom{\tev_{\alpha}^\dagger}$}};
	 	\fill (4.35,1.8) circle (0pt) node {{\small $\tev_{\alpha}^\dagger$}};
	 	%
	 	\end{tikzpicture}
	 } 
	 \, . 
\end{align}
Applying the 2-isomorphisms $b_{A,A}$ and $b_{A^\dagger, 1_\alpha}^{-1}$, the inner $A$-line can be moved to the right and the $A^\dagger$-line can be moved to the left, respectively. 
Alternatively, the $A^\dagger$-line can be moved to the left of the $b_{\alpha,\alpha}$-line with the help of two cusp isomorphisms, and then to the right by $b_{A^\dagger,1_\alpha}$. 
Thus we have shown: 

\begin{proposition}
	\label{prop:SA}
	Let $\alpha\in\Bfd$. 
	Then $(\alpha,A)\in\Beq^{\textrm{fd}}$, and 
	\be 
	\label{eq:SalphaAinBeq}
	S_{(\alpha,A)} 
		\cong 
		A^\dagger \otimes S_\alpha \otimes A 
		\cong 
		A \otimes S_\alpha \otimes A^\dagger \, . 
	\ee 
\end{proposition}

Below we will frequently not display $A=1_{(\alpha,A)}$ and simply write 
$
S_{(\alpha,A)} 
\cong 
A^\dagger \otimes S_\alpha 
\cong 
S_\alpha \otimes A^\dagger
$. 

\begin{corollary}
	\label{cor:SpinTQFTwithBeq}
	Let $r\in\Z_{\geqslant 1}$, $\alpha\in\Bfd$, and $(\alpha,A)\in\Beq$ such that $S_\alpha^r \cong 1_\alpha$ and $\gamma_A^r = 1_A$ in~$\B$. 
	Then there is an $r$-spin TQFT 
	\begin{align}
	\zz\colon \Bordrspin & \lra \Beq \nonumber
		\\ 
		+ & \lmt (\alpha,A) \, .  
	\end{align} 
\end{corollary}
\begin{proof}
	Combining~\eqref{eq:SalphaAinBeq} with $A^\dagger \cong A_{\gamma_A^{-1}}$, 
	Lemma~\ref{lem:r-spin-fixed-points} and Theorem~\ref{thm:r-spin-CH}, we see that any choice of isomorphism $S_\alpha^r \cong 1_\alpha$ determines a $\Spin_2^r$-homotopy fixed point in $(\Bfd)^{\times}$. 
\end{proof}

\subsubsection{A Frobenius algebra}
\label{subsubsec:FrobenisAlgebraInBeq}

For $(\alpha,A)\in\Beq^{\textrm{fd}}$ and $a\in\Z$, the $a$-th circle space (recall~\eqref{eq:Ca}) is 
\be 
\label{eq:CalphaA}
C^{(\alpha,A)}_a 
	\cong 
	\tev_{(\alpha,A)} \otimes_{\Ae} \big( S_{(\alpha,A)}^{1-a} \btimes 1_{(\alpha,A)^\dual} \big) \otimes_{\Ae} {}^\star \tev_{(\alpha,A)} \, . 
\ee 
By Proposition~\ref{prop:AdjointsInBeq} we have 
\be 
{}^\star \tev_{(\alpha,A)}
	\cong 
	{}_{\gamma_{\Ae}^{-1}} \big( {}^\dagger \! A \btimes 1_{\alpha^\dual} \big) \otimes {}^\dagger \tev_\alpha 
	\cong 
	\big( 1_\alpha \btimes {}_{\gamma_A^{-1}} ({}^\dagger\! A^\dual)_{\gamma_A} \big) \otimes \tev_\alpha \, ,
\ee 
and by Proposition~\ref{prop:SA} together with~\eqref{eq:NakayamaTwistDual} and ${}^\dagger\!A = A^\dagger$ we have 
\be 
\label{eq:SApower} 
S_{(\alpha,A)}^{1-a}
	\equiv 
	S_{(\alpha,A)}^{\otimes_A(1-a)} 
	\cong 
	{}_{\gamma_A^{1-a}} A \otimes S_\alpha^{1-a} \otimes A \, . 
\ee 

Our next goal is to explicitly describe the closed $\Lambda_0$-Frobenius structure on the circle spaces $C^{(\alpha,A)}_a$ in $\Beq$. 
This means that we will determine the (co)multiplication and (co)unit of \eqref{eq:Lambda_eta}--\eqref{eq:Lambda_Delta} of $C^{(\alpha,A)}_a$ directly in terms of data in~$\B$. 
In doing so we will frequently use the fact that the 2-morphism
\be 
\label{eq:SASASAlong}
\tikzzbox{
	\begin{tikzpicture}[very thick,scale=1.0,color=green!50!black, baseline=1cm]
	\coordinate (p0) at (0,0);
	\coordinate (p1) at (1,0);
	\coordinate (p1up) at (1,0.8);
	\coordinate (p2) at (2,0);
	\coordinate (p3) at (3,0);
	\coordinate (p35) at (3.5,0.8);
	\coordinate (p4) at (4,0);
	\coordinate (p5) at (5,0);
	\coordinate (p6) at (6,0);
	\coordinate (q0) at (0,2.5);
	\coordinate (q1) at (2.25,2.5);
	\coordinate (q2) at (3.5,2.5);
	\coordinate (q3) at (6,2.5);
	\draw[blue!50!black,-dot-, thick] (p2) .. controls +(0,2) and +(0,2) .. (p5);
	\draw (p1) -- (p1up);
	\draw[-dot-] (p1up) .. controls +(0,1.5) and +(0,1.5) .. (p35);
	\draw[-dot-] (p3) .. controls +(0,1) and +(0,1) .. (p4);
	\draw (2.25,1.9) -- (q1);
	\draw[blue!50!black, thick] (3.5,1.5) -- (q2);
	\draw (p6) -- (q3);
	\draw (p1) node[below] (bla) {${}_{\gamma_A^x}A$}; 
	\draw[blue!50!black] (p2) node[below] (bla2) {$S_\alpha^x$}; 
	\draw (p3) node[below] (bla3) {$A$}; 
	\draw (p4) node[below] (bla4) {${}_{\gamma_A^y}A$}; 
	\draw[blue!50!black] (p5) node[below] (bla2) {$S_\alpha^y$}; 
	\draw (p6) node[below] (bla2) {$A$}; 
	\draw (q1) node[above] (bla) {${}_{\gamma_A^{x+y}}A$}; 
	\draw[blue!50!black] (q2) node[above] (bla2) {$S_\alpha^{x+y}\vphantom{{}_{\gamma_A^{x+y}}A}$}; 
	\draw (q3) node[above] (bla2) {$A\vphantom{{}_{\gamma_A^{x+y}}A}$}; 
	\draw (q0) node[above left] (bla2) {$S_{(\alpha,A)}^{x+y} \cong \vphantom{{}_{\gamma_A^{x+y}}A}$};
	\draw (p0) node[below left] (bla2) {$S_{(\alpha,A)}^{x} \otimes S_{(\alpha,A)}^{y} \cong \vphantom{{}_{\gamma_A^{x+y}}A}$};
	\fill[color=green!50!black] (1,0.8) circle (2pt) node[left] (meet) {{\small$\gamma_A^y$}};
	\end{tikzpicture}
} 
\, , \quad x,y\in\Z \, , 
\ee 
in~$\B$ induces the identity $S_{(\alpha,A)}^x \otimes_A S_{(\alpha,A)}^y \lra S_{(\alpha,A)}^{x+y}$, up to the isomorphism in~\eqref{eq:SApower}. 
From now on we will no longer display all $A$-lines in diagrams that represent 2-morphisms in $\Beq$, since $A = 1_{(\alpha,A)}$. 
Accordingly, we abbreviate~\eqref{eq:SASASAlong} and its inverse as 
\be 
\label{eq:SASASA}
\tikzzbox{%
	\begin{tikzpicture}[very thick, scale=1,color=green!50!black, baseline=0cm]
	\coordinate (q1) at (-0.5,-1.5);
	\coordinate (q2) at (0.5,-1.5);
	\coordinate (q3) at (0,1.5);
	\draw (q1) .. controls +(0,0.5) and +(0,-0.5) .. (-0.1,0);
	\draw (q2) .. controls +(0,0.5) and +(0,-0.5) .. (+0.1,0);
	\draw (0,0) -- (q3);
	\draw (q3) node[above] (bla2) {$S_{(\alpha,A)}^{x+y}$};
	\draw (q1) node[below] (bla2) {$S_{(\alpha,A)}^{x}$};
	\draw (0.7,-1.5) node[below] (bla2) {$S_{(\alpha,A)}^{y}$};
	\draw[line width=0.5pt] (0,0) node[inner sep=0pt,draw,fill=green!50!black] (f) {{\scriptsize \phantom{opp}}};
	\end{tikzpicture}
}%
=
\tikzzbox{%
	\begin{tikzpicture}[very thick, scale=1,color=green!50!black, baseline=0cm]
	\coordinate (q1) at (-1,-1);
	\coordinate (q2) at (+1,-1);
	\coordinate (q3) at (0,1);
	\coordinate (q3u) at (0,1.5);
	\coordinate (q1d) at (-1,-1.5);
	\coordinate (q2d) at (+1,-1.5);
	\draw[blue!50!black, thick,-dot-] (-0.9,-1) .. controls +(0,1) and +(0,1) .. (1.1,-1);
	\draw[-dot-] (-1.1,-1) .. controls +(0,2) and +(0,2) .. (0.9,-1);
	\draw (q1d) node[below] (bla2) {$S_{(\alpha,A)}^{x}$};
	\draw (q1d)-- (q1);
	\draw (q2d) node[below] (bla2) {$S_{(\alpha,A)}^{y}$};
	\draw (q2d)-- (q2);
	\draw (q3u) node[above] (bla2) {$S_{(\alpha,A)}^{x+y}$};
	\draw (q3u)-- (q3);
	\draw (-0.1,1) -- (-0.1,0.5);
	\draw[blue!50!black, thick] (0.1,1) -- (0.1,-0.3);
	\draw[line width=0.5pt] (q1) node[inner sep=0pt,draw=green!50!gray,fill=green!50!gray] (f) {{\scriptsize \phantom{opp}}};
	\draw[line width=0.5pt] (q2) node[inner sep=0pt,draw=green!50!gray,fill=green!50!gray] (f) {{\scriptsize \phantom{opp}}};
	\draw[line width=0.5pt] (q3) node[inner sep=0pt,draw=green!50!gray,fill=green!50!gray] (f) {{\scriptsize \phantom{opp}}};
	\fill[color=green!50!black] (-1,-0.3) circle (2pt) node[left] (meet) {{\small$\gamma_A^y$}};
	\end{tikzpicture}
}%
, \quad 
\tikzzbox{%
	\begin{tikzpicture}[very thick, scale=1,color=green!50!black, baseline=0cm, yscale=-1]
	\coordinate (q1) at (-0.5,-1.5);
	\coordinate (q2) at (0.5,-1.5);
	\coordinate (q3) at (0,1.5);
	\draw (q1) .. controls +(0,0.5) and +(0,-0.5) .. (-0.1,0);
	\draw (q2) .. controls +(0,0.5) and +(0,-0.5) .. (+0.1,0);
	\draw (0,0) -- (q3);
	\draw (q3) node[below] (bla2) {$S_{(\alpha,A)}^{x+y}$};
	\draw (q1) node[above] (bla2) {$S_{(\alpha,A)}^{x}$};
	\draw (0.7,-1.5) node[above] (bla2) {$S_{(\alpha,A)}^{y}$};
	\draw[line width=0.5pt] (0,0) node[inner sep=0pt,draw,fill=green!50!black] (f) {{\scriptsize \phantom{opp}}};
	\end{tikzpicture}
}%
=
\tikzzbox{%
	\begin{tikzpicture}[very thick, scale=1,color=green!50!black, baseline=0cm, yscale=-1]
	\coordinate (q1) at (-1,-1);
	\coordinate (q2) at (+1,-1);
	\coordinate (q3) at (0,1);
	\coordinate (q3u) at (0,1.5);
	\coordinate (q1d) at (-1,-1.5);
	\coordinate (q2d) at (+1,-1.5);
	\draw[blue!50!black, thick,-dot-] (-0.9,-1) .. controls +(0,1) and +(0,1) .. (1.1,-1);
	\draw[-dot-] (-1.1,-1) .. controls +(0,2) and +(0,2) .. (0.9,-1);
	\draw (q1d) node[above] (bla2) {$S_{(\alpha,A)}^{x}$};
	\draw (q1d)-- (q1);
	\draw (q2d) node[above] (bla2) {$S_{(\alpha,A)}^{y}$};
	\draw (q2d)-- (q2);
	\draw (q3u) node[below] (bla2) {$S_{(\alpha,A)}^{x+y}$};
	\draw (q3u)-- (q3);
	\draw (-0.1,1) -- (-0.1,0.5);
	\draw[blue!50!black, thick] (0.1,1) -- (0.1,-0.3);
	\draw[line width=0.5pt] (q1) node[inner sep=0pt,draw=green!50!gray,fill=green!50!gray] (f) {{\scriptsize \phantom{opp}}};
	\draw[line width=0.5pt] (q2) node[inner sep=0pt,draw=green!50!gray,fill=green!50!gray] (f) {{\scriptsize \phantom{opp}}};
	\draw[line width=0.5pt] (q3) node[inner sep=0pt,draw=green!50!gray,fill=green!50!gray] (f) {{\scriptsize \phantom{opp}}};
	\fill[color=green!50!black] (-1,-0.3) circle (2pt) node[left] (meet) {{\small$\gamma_A^{-y}$}};
	\end{tikzpicture}
}%
. 
\ee 

With the above preparations, we can present the isomorphism~\eqref{eq:IsosWithSsquared} in the case of the equivariant completion: 

\begin{lemma}
	Let $\alpha\in\Bfd$. 
	Then for any $(\alpha,A)\in\Beq$, there are mutually inverse isomorphisms
	\be 
	\begin{tikzpicture}[
	baseline=(current bounding box.base), 
	descr/.style={fill=white,inner sep=3.5pt}, 
	normal line/.style={->}
	] 
	\matrix (m) [matrix of math nodes, row sep=3em, column sep=4.0em, text height=1.5ex, text depth=0.1ex] {%
		\big( 1_{(\alpha,A)} \btimes S^2_{(\alpha,A)^\dual} \big) \otimes_{\Ae} {}^\star \tev_{(\alpha,A)}  & \tev_{(\alpha,A)}^\star 
		\\
	};
	\path[font=\footnotesize, transform canvas={yshift=0.8mm}] (m-1-1) edge[->] node[above] { $f_A$ } (m-1-2);
	\path[font=\footnotesize, transform canvas={yshift=-0.8mm}] (m-1-2) edge[->] node[below] { $f'_A$ } (m-1-1);
	\end{tikzpicture}
	\ee 
	given by: 
	\begin{align}
	\label{eq:fAfAprime}
	f_A & = 
	\tikzzbox{%
		\begin{tikzpicture}[thick,scale=1.5,color=black, baseline=1.5cm]
		\coordinate (p1) at (0,0);
		\coordinate (p2) at (-3.5,-0.5);
		\coordinate (p3) at (-2.5,0.5);
		\coordinate (u1) at (0,3);
		\coordinate (u2) at (-3.5,2.5);
		\coordinate (u3) at (-2.5,3.5);
		%
		\fill [orange!20!white, opacity=0.8] 
		(p1) .. controls +(0,0.25) and +(1,0) ..  (p3)
		-- (p3) --  (u3)
		-- (u3) .. controls +(1,0) and +(0,0.25) ..  (u1)
		;
		%
		\draw[very thick, red!80!black] (p1) .. controls +(0,0.25) and +(1,0) ..  (p3); 
		\draw[thin] (p3) --  (u3); 
		%
		\fill [orange!30!white, opacity=0.8] 
		(p1) .. controls +(0,-0.5) and +(1,0) ..  (p2)
		-- (p2) --  (u2)
		-- (u2) .. controls +(1,0) and +(0,-0.5) ..  (u1)
		;
		\draw[thin] (p1) --  (u1); 
		\draw[thin] (p2) --  (u2); 
		%
		\draw[very thick, red!80!black] (u1) .. controls +(0,0.25) and +(1,0) ..  (u3); 
		%
		\fill[color=blue!50!black] (-2.2,3.2) circle (0pt) node {{\small $\alpha\vphantom{\alpha^\dual}$}};
		\fill[color=blue!50!black] (-3.2,2.2) circle (0pt) node {{\small $\alpha^\dual$}};
		\fill[color=green!50!black] (-2.25,-0.75) circle (0pt) node {{\small $S_{(\alpha,A)^\dual}^2$}};
		\fill[color=green!50!black] (-0.5,-0.6) circle (0pt) node {{\small ${}_{\gamma_A^{-1}}({}^\dagger \! A^\dual)_{\gamma_A}$}};
		\fill[color=green!50!black] (-1.5,2.75) circle (0pt) node {{\small $(A^\dual)^\dagger$}};
		\fill[color=blue!50!black] (-0.25,1.5) circle (0pt) node {{\scriptsize $S_{\alpha^\dual}^2$}};
		%
		\draw[blue!50!black] (-2.4,0) .. controls +(0,0.5) and +(0,-2.5) .. (0,2);
		\fill[color=blue!50!black] (0,2) circle (1.5pt) node[above] {};
		\fill (0.25,2.5) circle (0pt) node {{\small $\tev_{\alpha}^\dagger$}};
		\fill (0.27,0.5) circle (0pt) node {{\small ${}^\dagger\!\tev_{\alpha}$}};
		\draw[ultra thick, green!50!black] (-2.5,-0.5) -- (-2.5,0);
		\draw[ultra thick, green!50!black] (-2.6,0) .. controls +(0,0.5) and +(0,-0.25) .. (-3.25,0.5);
		\draw[ultra thick, green!50!black] (-3.25,0.5) -- (-3.25,0.95);
		\draw[ultra thick, green!50!black] (-3.25,0.95) .. controls +(0,0.25) and +(0,0.25) .. (-2.9,0.95);
		\draw[green!50!black] (-2.725,0.6) node[Odot] (unit1) {}; 
		\draw[ultra thick, green!50!black] (unit1) -- (-2.725,0.8);
		\draw[-dot-,ultra thick, green!50!black] (-2.9,0.95) .. controls +(0,-0.25) and +(0,-0.25) .. (-2.55,0.95);
		\draw[ultra thick, green!50!black] (-2.55,1.5) -- (-2.55,0.95);
		\draw[ultra thick, green!50!black] (-2.5,0) .. controls +(0,0.5) and +(0,-0.25) .. (-1.75,0.5);
		\draw[ultra thick, green!50!black] (-1.75,0.5) -- (-1.75,0.95);
		\draw[ultra thick, green!50!black] (-1.75,0.95) .. controls +(0,0.25) and +(0,0.25) .. (-2.1,0.95);
		\draw[green!50!black] (-2.275,0.6) node[Odot] (unit2) {}; 
		\draw[ultra thick, green!50!black] (unit2) -- (-2.275,0.8);
		\draw[-dot-,ultra thick, green!50!black] (-2.1,0.95) .. controls +(0,-0.25) and +(0,-0.25) .. (-2.45,0.95);
		\draw[-dot-,ultra thick, green!50!black] (-2.45,0.95) .. controls +(0,0.75) and +(0,0.75) .. (-1.2,0.95); 
		\draw[-dot-,ultra thick, green!50!black] (-1.2,0.95) .. controls +(0,-0.25) and +(0,-0.25) .. (-0.85,0.95);
		\draw[green!50!black] (-1.025,0.6) node[Odot] (unit3) {}; 
		\draw[ultra thick, green!50!black] (-1.025,0.8) -- (unit3);
		\draw[ultra thick, green!50!black] (-0.85,0.95) .. controls +(0,0.25) and +(0,0.25) .. (-0.5,0.95);
		\draw[ultra thick, green!50!black] (-0.5,-0.32) -- (-0.5,0.95);
		\draw[-dot-,ultra thick, green!50!black] (-2.55,1.5) .. controls +(0,0.4) and +(0,0.4) .. (-1.82,1.5); 
		\draw[ultra thick, green!50!black] (-2.185,1.8) -- (-2.185,2);
		\draw[-dot-, ultra thick, green!50!black] (-2.185,2) .. controls +(0,0.25) and +(0,0.25) .. (-1.835,2);
		\draw[green!50!black] (-2.01,2.35) node[Odot] (unit4) {}; 
		\draw[ultra thick, green!50!black] (-2.01,2.15) -- (unit4);
		\draw[ultra thick, green!50!black] (-1.485,2) .. controls +(0,-0.25) and +(0,-0.25) .. (-1.835,2);
		\draw[ultra thick, green!50!black] (-1.485,2) -- (-1.485,2.56);
		\draw[line width=0.5pt] (-2.5,0) node[inner sep=0pt,draw=green!50!gray,fill=green!50!gray] (f) {{\scriptsize \phantom{opp}}};
		%
		%
		\draw[very thick, red!80!black] (p1) .. controls +(0,-0.5) and +(1,0) ..  (p2); 
		\draw[very thick, red!80!black] (u1) .. controls +(0,-0.5) and +(1,0) ..  (u2); 
		\end{tikzpicture}
	}
	\! ,  
	& 
	f'_A & = 
	\tikzzbox{%
		\begin{tikzpicture}[thick,scale=1.5,color=black, baseline=1.5cm]
		\coordinate (p1) at (0,0);
		\coordinate (p2) at (-3.5,-0.5);
		\coordinate (p3) at (-2.5,0.5);
		\coordinate (u1) at (0,3);
		\coordinate (u2) at (-3.5,2.5);
		\coordinate (u3) at (-2.5,3.5);
		%
		\fill [orange!20!white, opacity=0.8] 
		(p1) .. controls +(0,0.25) and +(1,0) ..  (p3)
		-- (p3) --  (u3)
		-- (u3) .. controls +(1,0) and +(0,0.25) ..  (u1)
		;
		%
		\draw[very thick, red!80!black] (p1) .. controls +(0,0.25) and +(1,0) ..  (p3); 
		\draw[thin] (p3) --  (u3); 
		%
		\fill [orange!30!white, opacity=0.8] 
		(p1) .. controls +(0,-0.5) and +(1,0) ..  (p2)
		-- (p2) --  (u2)
		-- (u2) .. controls +(1,0) and +(0,-0.5) ..  (u1)
		;
		\draw[thin] (p1) --  (u1); 
		\draw[thin] (p2) --  (u2); 
		%
		\draw[very thick, red!80!black] (u1) .. controls +(0,0.25) and +(1,0) ..  (u3); 
		%
		\fill[color=blue!50!black] (-2.2,3.2) circle (0pt) node {{\small $\alpha\vphantom{\alpha^\dual}$}};
		\fill[color=blue!50!black] (-3.2,2.2) circle (0pt) node {{\small $\alpha^\dual$}};
		\fill[color=green!50!black] (-2.75,2.75) circle (0pt) node {{\small $S_{(\alpha,A)^\dual}^2$}};
		\fill[color=green!50!black] (-2.85,-0.7) circle (0pt) node {{\small ${}^\dagger\!A^\dual$}};
		\fill[color=green!50!black] (-1.15,2.8) circle (0pt) node {{\small $(A^\dagger)^\dual$}};
		%
		\draw[blue!50!black] (-2.4,2) .. controls +(0,-0.9) and +(0,3.3) .. (0,1);
		\fill[color=blue!50!black] (0,1) circle (1.5pt) node[above] {};
		\fill (0.25,0.5) circle (0pt) node {{\small $\tev_{\alpha}^\dagger$}};
		\fill (0.27,2.5) circle (0pt) node {{\small ${}^\dagger\!\tev_{\alpha}$}};
		\fill[color=blue!50!black] (-0.25,1.55) circle (0pt) node {{\scriptsize $S_{\alpha^\dual}^{-2}$}};
		\draw[ultra thick, green!50!black] (-2.5,2.5) -- (-2.5,2);
		\draw[ultra thick, green!50!black] (-2.6,2) -- (-2.6,1.5);
		\draw[ultra thick, green!50!black] (-2.6,1.5) .. controls +(0,-0.25) and +(0,-0.25) .. (-2.95,1.5);
		\draw[green!50!black] (-3.125,1.85) node[Odot] (unit1) {}; 
		\draw[ultra thick, green!50!black] (unit1) -- (-3.125,1.7);
		\draw[-dot-,ultra thick, green!50!black] (-2.95,1.5) .. controls +(0,0.25) and +(0,0.25) .. (-3.3,1.5);
		\draw[ultra thick, green!50!black] (-3.3,1.5) -- (-3.3,1);
		\draw[ultra thick, green!50!black] (-2.5,2.5) -- (-2.5,1.5);
		\draw[ultra thick, green!50!black] (-2.5,1.5) .. controls +(0,-0.25) and +(0,-0.25) .. (-2.15,1.5);
		\draw[green!50!black] (-1.975,1.85) node[Odot] (unit2) {}; 
		\draw[ultra thick, green!50!black] (unit2) -- (-1.975,1.7);
		\draw[-dot-,ultra thick, green!50!black] (-2.15,1.5) .. controls +(0,0.25) and +(0,0.25) .. (-1.8,1.5);
		\draw[ultra thick, green!50!black] (-1.8,1.5) -- (-1.8,1);
		\draw[-dot-, ultra thick, green!50!black] (-1.8,1) .. controls +(0,-0.4) and +(0,-0.4) .. (-0.5,1);
		\draw[ultra thick, green!50!black] (-0.5,1) -- (-0.5,2);
		\draw[-dot-, ultra thick, green!50!black] (-0.5,2) .. controls +(0,0.25) and +(0,0.25) .. (-0.85,2);
		\draw[ultra thick, green!50!black] (-0.85,2) .. controls +(0,-0.25) and +(0,-0.25) .. (-1.2,2);
		\draw[green!50!black] (-0.675,2.35) node[Odot] (unit3) {}; 
		\draw[ultra thick, green!50!black] (-0.675,2.2) -- (unit3);
		\draw[ultra thick, green!50!black] (-1.2,2) -- (-1.2,2.58);
		\draw[ultra thick, green!50!black] (-3.3,1.5) -- (-3.3,0.7);
		\draw[-dot-, ultra thick, green!50!black] (-3.3,0.7) .. controls +(0,-0.5) and +(0,-0.5) .. (-1.15,0.7);
		\draw[ultra thick, green!50!black] (-2.225,0.3) -- (-2.225,0);
		\draw[-dot-, ultra thick, green!50!black] (-2.225,0) .. controls +(0,-0.25) and +(0,-0.25) .. (-2.575,0);
		\draw[green!50!black] (-2.4,-0.35) node[Odot] (unit4) {}; 
		\draw[ultra thick, green!50!black] (-2.4,-0.2) -- (unit4);
		\draw[ultra thick, green!50!black] (-2.575,0) .. controls +(0,0.25) and +(0,0.25) .. (-2.9025,0);
		\draw[ultra thick, green!50!black] (-2.9025,0) -- (-2.9025,-0.5);
		\draw[line width=0.5pt] (-2.5,2) node[inner sep=0pt,draw=green!50!gray,fill=green!50!gray] (f) {{\scriptsize \phantom{opp}}};
		%
		%
		\draw[very thick, red!80!black] (p1) .. controls +(0,-0.5) and +(1,0) ..  (p2); 
		\draw[very thick, red!80!black] (u1) .. controls +(0,-0.5) and +(1,0) ..  (u2); 
		\end{tikzpicture}
	}
	\end{align}
\end{lemma}
\begin{proof}
	Repeated use of~\eqref{eq:NakayamaTwistDual} together with standard manipulations of string diagrams for $\Delta$-separable Frobenius algebras 
		 shows 
	that $f_A, f'_A$ are indeed bimodule maps, and that $f_A \circ f'_A = 1_{{}^\dagger\!A^\dual}$. 
	Since their source and target are isomorphic, it follows that~$f_A$ and~$f'_A$ indeed represent mutually inverse 2-morphisms in $\Beq$. 
\end{proof}

We have now expressed all the ingredients of the closed $\Lambda_0$-Frobenius structure on $\{ C^{(\alpha,A)}_a \}$ in $\Beq$ directly in terms of data in~$\B$. 
With this the Nakayama automorphisms $N^{(\alpha,A)}_a$ can be computed: 

\begin{proposition}
	\label{prop:NakayamaPowerBeq}
	Let $\alpha\in\Bfd$. 
	Then for any $(\alpha,A)\in\Beq$, we have 
	\be 
	N^{(\alpha,A)}_a 
	= 
	\tikzzbox{%
		\begin{tikzpicture}[very thick, scale=1,color=green!50!black, baseline=2.5cm]
		\coordinate (d2) at (2,0);
		\fill [orange!20!white, opacity=0.8] 
		($(d2)+(-3.3,1)$) -- ($(d2)+(-3.3,4)$) -- ($(d2)+(+3.3,4)$) -- ($(d2)+(+3.3,1)$)
		;
		%
		\draw[thick, color=black] ($(d2)+(-3.3,1)$) -- ($(d2)+(-3.3,4)$); 
		\fill[color=black] ($(d2)+(-3.65,1.6)$) circle (0pt) node {{\small $\tev_\alpha\vphantom{{}_{\gamma^{2-a}_{A^\dual}}A^\dual}$}};
		%
		\draw[thick, color=black] ($(d2)+(+3.3,1)$) -- ($(d2)+(+3.3,4)$); 
		\fill[color=black] ($(d2)+(3.7,1.6)$) circle (0pt) node {{\small $\tev_\alpha^\dagger\vphantom{{}_{\gamma^{2-a}_{A^\dual}}A^\dual}$}};
		%
		\draw ($(d2)+(-2,1)$) -- ($(d2)+(-2,4)$);
		\fill ($(d2)+(-2.3,1.6)$) circle (0pt) node {{\small $A^\dual\vphantom{{}_{\gamma^{2-a}_{A^\dual}}A^\dual}$}};
		%
		\draw[ultra thick] ($(d2)+(0,1)$) -- ($(d2)+(0,4)$);
		\fill ($(d2)+(0.7,1.6)$) circle (0pt) node {{\small ${}_{\gamma^{2-a}_{A^\dual}}A^\dual$}}; 
		%
		\draw[color=blue!50!black] ($(d2)+(+2,1)$) -- ($(d2)+(+2,4)$);
		\fill[color=blue!50!black] ($(d2)+(2.5,1.6)$) circle (0pt) node {{\small $S_\alpha^{1-a}\vphantom{{}_{\gamma^{2-a}_{A^\dual}}A^\dual}$}}; 
		\draw[ultra thick, color=blue!50!black] (d2) -- ($(d2)+(0,1)$); 
		\draw[ultra thick, color=blue!50!black] ($(d2)+(0,5)$) -- ($(d2)+(0,4)$); 
		%
		\fill ($(d2)+(0,2.5)$) circle (3pt) node[right] {{\small $\gamma_{A^\dual}$}};
		\fill[color=blue!50!black] (d2) circle (0pt) node[below] {{\small $C_a^{(\alpha,A)}$}};
		%
		\draw[line width=0.5pt] ($(d2)+(0,4)$) node[inner sep=1pt,draw=green!50!gray,fill=green!50!gray] (f) {{\scriptsize \phantom{opppppppppppppppppppppppppppppppppppppppppp}}};
		\draw[line width=0.5pt] ($(d2)+(0,1)$) node[inner sep=1pt,draw=green!50!gray,fill=green!50!gray] (f) {{\scriptsize \phantom{opppppppppppppppppppppppppppppppppppppppppp}}};
		\coordinate (q) at ($(d2)+(0,5)$);
		\coordinate (d) at ($(q)-(-0.25,-0.75)$);
		\draw[ultra thick, color=blue!50!black, postaction={decorate}, decoration={markings,mark=at position .59 with {\arrow[draw=blue!50!black]{>}}}] ($(d)+(-0.25,-0.25)$) .. controls +(0,0.75) and +(0,0.5) .. ($(d)+(1,0)$);
		\draw[ultra thick, color=blue!50!black, postaction={decorate}, decoration={markings,mark=at position .59 with {\arrow[draw=blue!50!black]{<}}}] ($(d)+(-0.25,0.25)$) .. controls +(0,-0.75) and +(0,-0.5) .. ($(d)+(1,0)$);
		\draw[ultra thick, color=blue!50!black] ($(d)+(-0.25,-0.25)$) -- ($(d)+(-0.25,-1.25)$);
		\draw[ultra thick, color=blue!50!black] ($(d)+(-0.25,+0.25)$) -- ($(d)+(-0.25,+1.25)$);
		%
		\draw[color=blue!50!black] ($(d)+(-0.25,+1.25)$) node[above] (bla2) {$C_a^{(\alpha,A)}$};
		\end{tikzpicture}
	}%
	.
	\ee 
\end{proposition}
\begin{proof}
	Our task is to compute $N^{(\alpha,A)}_a$ as defined in~\eqref{eq:Nakayama-Ca} in~$\Beq$. 
	We will first compute $\varepsilon_{-1} \circ \mu_{a,-a}$ and $\Delta_{a,-a} \circ \eta_1$ in~$\Beq$, starting with $\varepsilon_{-1} \circ \mu_{a,-a}$. 
	Using the notation introduced in~\eqref{eq:SASASA}, we have 
	\be 
	\mu_{a,-a} = 
	\tikzzbox{%
		\begin{tikzpicture}[thick,scale=1.0,color=black, baseline=1.9cm]
		\coordinate (p1) at (-2.75,0);
		\coordinate (p2) at (-1,0);
		\coordinate (p3) at (1,0);
		\coordinate (p4) at (2.75,0);
		\coordinate (p5) at (0.875,4);
		\coordinate (p6) at (-0.875,4);
		%
		\draw[very thick, red!80!black] (p1) .. controls +(0,0.5) and +(0,0.5) ..  (p2); 
		\draw[very thick, red!80!black] (p3) .. controls +(0,0.5) and +(0,0.5) ..  (p4); 
		%
		\fill [orange!30!white, opacity=0.8] 
		(p1) .. controls +(0,0.5) and +(0,0.5) ..  (p2)
		-- (p2) .. controls +(0,1.5) and +(0,1.5) ..  (p3)
		-- (p3) .. controls +(0,0.5) and +(0,0.5) ..  (p4)
		-- (p4) .. controls +(0,2.5) and +(0,-2.5) ..  (p5)
		-- (p5) .. controls +(0,-0.5) and +(0,-0.5) ..  (p6)
		-- (p6) .. controls +(0,-2.5) and +(0,2.5) ..  (p1)
		;
		\fill [orange!20!white, opacity=0.8] 
		(p5) .. controls +(0,-0.5) and +(0,-0.5) ..  (p6)
		-- (p6) .. controls +(0,0.5) and +(0,0.5) ..  (p5)
		;
		\fill [orange!20!white, opacity=0.8] 
		(p1) .. controls +(0,-0.5) and +(0,-0.5) ..  (p2)
		-- (p2) .. controls +(0,0.5) and +(0,0.5) ..  (p1)
		;
		\fill [orange!20!white, opacity=0.8] 
		(p3) .. controls +(0,-0.5) and +(0,-0.5) ..  (p4)
		-- (p4) .. controls +(0,0.5) and +(0,0.5) ..  (p3)
		;
		\draw (p2) .. controls +(0,1.5) and +(0,1.5) ..  (p3); 
		\draw (p4) .. controls +(0,2.5) and +(0,-2.5) ..  (p5); 
		\draw (p6) .. controls +(0,-2.5) and +(0,2.5) ..  (p1); 
		\coordinate (q1) at (-1.875,-0.37);
		\coordinate (q2) at (1.875,-0.37);
		\coordinate (q3) at (0,2);
		\coordinate (q4) at (0,3.63);
		\coordinate (q5) at (-1.4,-0.3);
		\coordinate (q6) at (+1.4,-0.3);
		\coordinate (q7) at (-2.6,-0.24);
		\coordinate (q8) at (+2.6,-0.24);
		\coordinate (q9) at (-0.7,3.75);
		\coordinate (q10) at (0.7,3.75);
		\draw[color=green!50!black, ultra thick] (q1) .. controls +(0,1.5) and +(0,-0.5) ..  (-0.1,2); 
		\draw[color=green!50!black, ultra thick] (q2) .. controls +(0,1.5) and +(0,-0.5) ..  (+0.1,2); 
		\draw[color=green!50!black, ultra thick] (q4) -- (q3); 
		\draw[directedgreen, color=green!50!black] (q6) .. controls +(0,2.2) and +(0,2.2) ..  (q5); 
		\draw[color=green!50!black] (q7) .. controls +(0,2.2) and +(0,-2.2) ..  (q9); 
		\draw[color=green!50!black] (q8) .. controls +(0,2.2) and +(0,-2.2) ..  (q10); 
		%
		\fill[color=green!50!black] (q1) circle (0pt) node[below] {{\small $S_{(\alpha,A)^\dual}^{-a+1}$}};
		\fill[color=green!50!black] (q2) circle (0pt) node[below] {{\small $S_{(\alpha,A)^\dual}^{a+1}$}};
		\fill[color=green!50!black] (0,3.6) circle (0pt) node[above] {{\small $S_{(\alpha,A)^\dual}^{2}$}};
		\fill[color=green!50!black] (-1.05,-0.45) circle (0pt) node {{\scriptsize ${}^\dagger\! A^\dual$}};
		\fill[color=green!50!black] (+1.1,-0.45) circle (0pt) node {{\scriptsize $A^\dual$}};
		\fill[color=green!50!black] (-0.8,2.45) circle (0pt) node {{\scriptsize $A^\dual\vphantom{{}^\dagger\! A^\dual}$}};
		\fill[color=green!50!black] (+0.7,2.45) circle (0pt) node {{\scriptsize ${}^\dagger\! A^\dual$}};
		\fill ($(p1)+(0.1,0)$) circle (0pt) node[left] {{\small $\tev_\alpha$}};
		\fill ($(p3)+(0.1,0)$) circle (0pt) node[left] {{\small $\tev_\alpha$}};
		\fill ($(p6)+(0.1,0)$) circle (0pt) node[left] {{\small $\tev_\alpha$}};
		\fill ($(p2)+(-0.05,0)$) circle (0pt) node[right] {{\small ${}^\dagger\!\tev_\alpha$}};
		\fill ($(p4)+(-0.05,0)$) circle (0pt) node[right] {{\small ${}^\dagger\!\tev_\alpha$}};
		\fill ($(p5)+(-0.05,0)$) circle (0pt) node[right] {{\small ${}^\dagger\!\tev_\alpha$}};
		\fill (0,0.75) circle (0pt) node {{\small $\ev_{\tev_\alpha}$}};
		\draw[line width=0.5pt] (q3) node[inner sep=0pt,draw,fill=green!50!black] (f) {{\scriptsize \phantom{opp}}};
		%
		\draw[very thick, red!80!black] (p1) .. controls +(0,-0.5) and +(0,-0.5) ..  (p2); 
		\draw[very thick, red!80!black] (p3) .. controls +(0,-0.5) and +(0,-0.5) ..  (p4); 
		\draw[very thick, red!80!black] (p5) .. controls +(0,0.5) and +(0,0.5) ..  (p6); 
		\draw[very thick, red!80!black] (p5) .. controls +(0,-0.5) and +(0,-0.5) ..  (p6); 
		\end{tikzpicture}
	}%
	\, , 
	\ee 
	and inserting the expression for~$f_A$ in~\eqref{eq:fAfAprime} into~\eqref{eq:Lambda_epsilon}, we have 
	\be 
	\varepsilon_{-1}  = 
	\tikzzbox{%
		\begin{tikzpicture}[very thick,scale=1.5,color=black, baseline=2.5cm]
		\coordinate (p1) at (0,0);
		\coordinate (p2) at (-2.5,-0.5);
		\coordinate (p3) at (-2.5,0.5);
		\coordinate (u1) at (0,3);
		\coordinate (u2) at (-2.5,2.5);
		\coordinate (u3) at (-2.5,3.5);
		\coordinate (u4) at (-5,3);
		\coordinate (p4) at (-5,0);
		\coordinate (u) at (-2.5,4.4);
		%
		\fill [orange!20!white, opacity=0.8] 
		(p1) .. controls +(0,0.5) and +(1,0) ..  (p3)
		-- (p3) .. controls +(-1,0) and +(0,0.5) ..  (p4)
		-- (p4) -- (u4) 
		-- (u4) .. controls +(0,1) and +(-1,0) ..  (u)
		-- (u)  .. controls +(1,0) and +(0,1) ..  (u1)
		;
		%
		\draw[very thick, red!80!black] (p1) .. controls +(0,0.5) and +(1,0) ..  (p3)
		(p3) .. controls +(-1,0) and +(0,0.5) ..  (p4); 
		%
		\fill [orange!30!white, opacity=0.8] 
		(p1) .. controls +(0,-0.5) and +(1,0) ..  (p2)
		-- (p2) .. controls +(-1,0) and +(0,-0.5) ..  (p4)
		-- (p4) -- (u4)
		-- (u4) .. controls +(0,1) and +(-1,0) ..  (u)
		-- (u)  .. controls +(1,0) and +(0,1) ..  (u1)
		;
		\draw[thin] (p1) --  (u1); 
		\draw[thin] (p4) --  (u4) -- (u4) .. controls +(0,1) and +(-1,0) ..  (u)
		-- (u)  .. controls +(1,0) and +(0,1) ..  (u1);
		%
		\fill[color=green!50!black] (-2.25,-0.75) circle (0pt) node {{\small $S_{(\alpha,A)^\dual}^2$}};
		\fill[color=green!50!black] (-0.5,-0.62) circle (0pt) node {{\small ${}_{\gamma_A^{-1}}({}^\dagger \! A^\dual)_{\gamma_A}$}};
		\fill[color=green!50!black] (-4.5,-0.55) circle (0pt) node {{\small $A^\dual$}};
		\fill[color=blue!50!black] (-0.25,1.5) circle (0pt) node {{\scriptsize $S_{\alpha^\dual}^2$}};
		%
		\draw[blue!50!black] (-2.4,0) .. controls +(0,0.5) and +(0,-2.5) .. (0,2);
		\fill[color=blue!50!black] (0,2) circle (1.5pt) node[above] {};
		\fill (0.25,2.5) circle (0pt) node {{\small $\tev_{\alpha}^\dagger$}};
		\fill (0.27,0.5) circle (0pt) node {{\small ${}^\dagger\!\tev_{\alpha}$}};
		\draw[ultra thick, green!50!black] (-2.5,-0.5) -- (-2.5,0);
		\draw[green!50!black] (-2.6,0) .. controls +(0,0.5) and +(0,-0.25) .. (-3.25,0.5);
		\draw[green!50!black] (-3.25,0.5) -- (-3.25,0.95);
		\draw[green!50!black] (-3.25,0.95) .. controls +(0,0.25) and +(0,0.25) .. (-2.9,0.95);
		\draw[green!50!black] (-2.725,0.6) node[Odot] (unit1) {}; 
		\draw[green!50!black] (unit1) -- (-2.725,0.8);
		\draw[-dot-,green!50!black] (-2.9,0.95) .. controls +(0,-0.25) and +(0,-0.25) .. (-2.55,0.95);
		\draw[green!50!black] (-2.55,1.5) -- (-2.55,0.95);
		\draw[green!50!black] (-2.5,0) .. controls +(0,0.5) and +(0,-0.25) .. (-1.75,0.5);
		\draw[green!50!black] (-1.75,0.5) -- (-1.75,0.95);
		\draw[green!50!black] (-1.75,0.95) .. controls +(0,0.25) and +(0,0.25) .. (-2.1,0.95);
		\draw[green!50!black] (-2.275,0.6) node[Odot] (unit2) {}; 
		\draw[green!50!black] (unit2) -- (-2.275,0.8);
		\draw[-dot-,green!50!black] (-2.1,0.95) .. controls +(0,-0.25) and +(0,-0.25) .. (-2.45,0.95);
		\draw[-dot-,green!50!black] (-2.45,0.95) .. controls +(0,0.75) and +(0,0.75) .. (-1.2,0.95); 
		\draw[-dot-,green!50!black] (-1.2,0.95) .. controls +(0,-0.25) and +(0,-0.25) .. (-0.85,0.95);
		\draw[green!50!black] (-1.025,0.6) node[Odot] (unit3) {}; 
		\draw[green!50!black] (-1.025,0.8) -- (unit3);
		\draw[green!50!black] (-0.85,0.95) .. controls +(0,0.25) and +(0,0.25) .. (-0.5,0.95);
		\draw[green!50!black] (-0.5,-0.37) -- (-0.5,0.95);
		\draw[-dot-,green!50!black] (-2.55,1.5) .. controls +(0,0.4) and +(0,0.4) .. (-1.82,1.5); 
		\draw[green!50!black] (-2.185,1.8) -- (-2.185,2);
		\draw[-dot-, green!50!black] (-2.185,2) .. controls +(0,0.25) and +(0,0.25) .. (-1.835,2);
		\draw[green!50!black] (-2.01,2.35) node[Odot] (unit4) {}; 
		\draw[green!50!black] (-2.01,2.15) -- (unit4);
		\draw[green!50!black] (-1.485,2) .. controls +(0,-0.25) and +(0,-0.25) .. (-1.835,2);
		\draw[green!50!black] (-1.485,2) .. controls +(0,2) and +(0,2) .. (-4.5,2);
		\draw[green!50!black] (-4.5,2) -- (-4.5,-0.37);
		\draw[line width=0.5pt] (-2.5,0) node[inner sep=0pt,draw=green!50!gray,fill=green!50!gray] (f) {{\scriptsize \phantom{opp}}};
		%
		%
		\draw[very thick, red!80!black] (p1) .. controls +(0,-0.5) and +(1,0) ..  (p2)
		(p2) .. controls +(-1,0) and +(0,-0.5) ..  (p4); 
		\end{tikzpicture}
	}
	\, . 
	\ee 
	In the composition $\varepsilon_{-1} \circ \mu_{a,-a}$, we first use 
	\be 
	\tikzzbox{%
		\begin{tikzpicture}[very thick, scale=1,color=green!50!black, baseline=0.4cm]
		\coordinate (q1) at (-0.5,-1.5);
		\coordinate (q2) at (0.5,-1.5);
		\coordinate (q3) at (0,1.5);
		\coordinate (p1) at (-0.5,2.5);
		\coordinate (p2) at (0.5,2.5);
		\draw[ultra thick] (q1) .. controls +(0,0.5) and +(0,-0.5) .. (-0.1,0);
		\draw[ultra thick] (q2) .. controls +(0,0.5) and +(0,-0.5) .. (+0.1,0);
		\draw[ultra thick] (0,0) -- (0,1);
		\draw (p1) .. controls +(0,-0.5) and +(0,0.5) .. (-0.1,1);
		\draw (p2) .. controls +(0,-0.5) and +(0,0.5) .. (+0.1,1);
		\draw (0.7,0.45) node (bla2) {{\small $S_{(\alpha,A)^\dual}^{2}$}};
		\draw ($(q1)+(-0.1,0)$) node[below] (bla2) {$S_{(\alpha,A)^\dual}^{-a+1}$};
		\draw ($(q2)+(+0.3,0)$) node[below] (bla2) {$S_{(\alpha,A)^\dual}^{+a+1}$};
		\draw (p1) node[above] (bla2) {${}^\dagger\!A^\dual$};
		\draw (p2) node[above] (bla2) {${}^\dagger\!A^\dual$};
		\draw[line width=0.5pt] (0,0) node[inner sep=0pt,draw,fill=green!50!black] (f) {{\scriptsize \phantom{opp}}};
		\draw[line width=0.5pt] (0,1) node[inner sep=0pt,draw,fill=green!50!gray] (f) {{\scriptsize \phantom{opp}}};
		\end{tikzpicture}
	}%
	=
	\tikzzbox{%
		\begin{tikzpicture}[very thick, scale=1,color=green!50!black, baseline=0.4cm]
		\coordinate (q1) at (-0.5,-1.5);
		\coordinate (q2) at (0.5,-1.5);
		\coordinate (q3) at (0,1.5);
		\coordinate (p1) at (0.5,2.5);
		\coordinate (p2) at (2.5,2.5);
		\draw[ultra thick] (q1) -- ($(q1)+(0,0.5)$);
		\draw[ultra thick] (q2) -- ($(q2)+(0,0.5)$);
		\draw[-dot-] ($(q1)+(0,1)$) .. controls +(0,0.75) and +(0,0.75) .. ($(q2)+(0,1)$);
		\draw(0,0) -- (0,0.5);
		\draw ($(q1)+(0,1)$) -- ($(q1)+(0,0.5)$);
		\draw ($(q2)+(0,1)$) -- ($(q2)+(0,0.5)$);
		\fill ($(q1)+(0,1)$) circle (2pt) node[left] {{\small $\gamma_{A^\dual}^{a+1}$}};
		\draw ($(q1)+(0,1)$); 
		%
		\draw[-dot-] (-1,1.5) .. controls +(0,-1.25) and +(0,-1.25) .. (1,1.5);
		\draw (-0.625,2.2) node[Odot] (unit1) {}; 
		\draw (unit1) -- (-0.625,1.9);
		\fill (-1,1.5) circle (2pt) node[left] {{\small $\gamma_{A^\dual}^{-1}$}};
		\draw[-dot-] (-1,1.5) .. controls +(0,0.5) and +(0,0.5) .. (-0.25,1.5);
		\draw[redirectedgreen] (-0.25,1.5) .. controls +(0,-0.5) and +(0,-0.5) .. (0.5,1.5);
		\draw (0.5,1.5) -- (p1);
		%
		\draw (1.375,2.2) node[Odot] (unit2) {}; 
		\draw (unit2) -- (1.375,1.9);
		\draw[-dot-] (1,1.5) .. controls +(0,0.5) and +(0,0.5) .. (1.75,1.5);
		\draw[redirectedgreen] (1.75,1.5) .. controls +(0,-0.5) and +(0,-0.5) .. (2.5,1.5);
		\draw (2.5,1.5) -- (p2);
		\draw ($(q1)+(-0.1,0)$) node[below] (bla2) {$S_{(\alpha,A)^\dual}^{-a+1}$};
		\draw ($(q2)+(+0.3,0)$) node[below] (bla2) {$S_{(\alpha,A)^\dual}^{+a+1}$};
		\draw (p1) node[above] (bla2) {${}^\dagger\!A^\dual$};
		\draw (p2) node[above] (bla2) {${}^\dagger\!A^\dual$};
		\draw[line width=0.5pt] ($(q1)+(0,0.5)$) node[inner sep=0pt,draw,fill=green!50!gray] (f) {{\scriptsize \phantom{opp}}};
		\draw[line width=0.5pt] ($(q2)+(0,0.5)$) node[inner sep=0pt,draw,fill=green!50!gray] (f) {{\scriptsize \phantom{opp}}};
		\end{tikzpicture}
	}%
	\ee 
	where here and below we suppress $S_\alpha$-strands. 
	This expression cancels with another subdiagram of $\varepsilon_{-1} \circ \mu_{a,-a}$, leaving 
	\begin{align}
	\varepsilon_{-1} \circ \mu_{a,-a} 
	& \equiv  
	\tikzzbox{%
		\begin{tikzpicture}[very thick, scale=1,color=green!50!black, baseline=2.4cm]
		\coordinate (q1) at (1,0);
		\coordinate (q2) at (2,0);
		\coordinate (q3) at (3,0);
		\coordinate (q4) at (4,0);
		\coordinate (q5) at (5,0);
		\coordinate (q6) at (7,0);
		%
		\draw[ultra thick] (q2) -- ($(q2)+(0,0.5)$);
		\draw (q3) -- ($(q3)+(0,0.5)$);
		\draw (q4) -- ($(q4)+(0,0.5)$);
		\draw[ultra thick] (q5) -- ($(q5)+(0,0.5)$);
		\draw (q6) -- ($(q6)+(0,0.5)$);
		\draw[redirectedgreen] ($(q3)+(0,0.5)$) .. controls +(0,1) and +(0,1) .. ($(q4)+(0,0.5)$);
		%
		\draw ($(q2)+(0,1)$) -- ($(q2)+(0,0.5)$);
		\draw ($(q5)+(0,1)$) -- ($(q5)+(0,0.5)$); 
		\draw[-dot-] ($(q2)+(0,1)$) .. controls +(0,2) and +(0,2) .. ($(q5)+(0,1)$);
		\fill ($(q2)+(0,1)$) circle (2pt) node[right] {{\small $\gamma_{A^\dual}^{a+1}$}};
		%
		\draw (3.5,2.5) -- (3.5,3);
		\fill (3.5,3) circle (2pt) node[left] {{\small $\gamma_{A^\dual}^{-1}$}};
		\draw[-dot-] (3.5,3) .. controls +(0,1) and +(0,1) .. (5.5,3);
		\draw (5.5,3) -- (5.5,2); 
		\draw[-dot-] (5.5,2) .. controls +(0,-0.75) and +(0,-0.75) .. (6.25,2);
		\draw (5.875,1.1) node[Odot] (unit1) {}; 
		\draw (unit1) -- (5.875,1.4);
		\draw[directedgreen] (6.25,2) .. controls +(0,0.75) and +(0,0.75) .. (7,2);
		\draw (7,2) -- (q6); 
		%
		\draw (4.5,3.75) -- (4.5,4.2);
		\draw[-dot-] (4.5,4.2) .. controls +(0,0.75) and +(0,0.75) .. (5.25,4.2);
		\draw (4.875,5.1) node[Odot] (unit2) {}; 
		\draw (unit2) -- (4.874,4.8);
		\draw[redirectedgreen] (5.25,4.2) .. controls +(0,-0.75) and +(0,-0.75) .. (6,4.2);
		\draw (6,4.2) -- (6,4.5);
		\draw[directedgreen] (1,4.5) .. controls +(0,1.5) and +(0,1.5) .. (6,4.5);
		\draw (1,4.5) -- (q1);
		%
		\draw (q1) node[below] (bla2) {$A^\dual\vphantom{S_{(\alpha,A)^\dual}^{-a+1}}$};
		\draw (q2) node[below] (bla2) {$S_{(\alpha,A)^\dual}^{-a+1}$};
		\draw (q3) node[below] (bla2) {${}^\dagger\!A^\dual\vphantom{S_{(\alpha,A)^\dual}^{-a+1}}$};
		\draw (q4) node[below] (bla2) {$A^\dual\vphantom{S_{(\alpha,A)^\dual}^{-a+1}}$};
		\draw (q5) node[below] (bla2) {$S_{(\alpha,A)^\dual}^{+a+1}$};
		\draw (q6) node[below] (bla2) {$(A^\dual)^\dagger\vphantom{S_{(\alpha,A)^\dual}^{-a+1}}$};
		\draw[line width=0.5pt] ($(q2)+(0,0.5)$) node[inner sep=0pt,draw,fill=green!50!gray] (f) {{\scriptsize \phantom{opp}}};
		\draw[line width=0.5pt] ($(q5)+(0,0.5)$) node[inner sep=0pt,draw,fill=green!50!gray] (f) {{\scriptsize \phantom{opp}}};
		\end{tikzpicture}
	}%
	\nonumber 
	\\
	& =
	\tikzzbox{%
		\begin{tikzpicture}[very thick, scale=1,color=green!50!black, baseline=2.4cm]
		\coordinate (q1) at (1,0);
		\coordinate (q2) at (2,0);
		\coordinate (q3) at (3,0);
		\coordinate (q4) at (4,0);
		\coordinate (q5) at (5,0);
		\coordinate (q6) at (7,0);
		%
		\draw[ultra thick] (q2) -- ($(q2)+(0,0.5)$);
		\draw (q3) -- ($(q3)+(0,0.5)$);
		\draw (q4) -- ($(q4)+(0,0.5)$);
		\draw[ultra thick] (q5) -- ($(q5)+(0,0.5)$);
		\draw (q6) -- ($(q6)+(0,0.5)$);
		\draw[redirectedgreen] ($(q3)+(0,0.5)$) .. controls +(0,1) and +(0,1) .. ($(q4)+(0,0.5)$);
		%
		\draw ($(q2)+(0,1)$) -- ($(q2)+(0,0.5)$);
		\draw ($(q5)+(0,1)$) -- ($(q5)+(0,0.5)$); 
		\draw[-dot-] ($(q2)+(0,1)$) .. controls +(0,2) and +(0,2) .. ($(q5)+(0,1)$);
		\fill ($(q2)+(0,1)$) circle (2pt) node[right] {{\small $\gamma_{A^\dual}^{a+1}$}};
		%
		\draw (3.5,2.5) -- (3.5,3);
		\fill (3.5,3) circle (2pt) node[left] {{\small $\gamma_{A^\dual}^{-1}$}};
		\draw[-dot-] (3.5,3) .. controls +(0,1) and +(0,1) .. (5.5,3);
		\draw (5.5,3) -- (5.5,2); 
		\draw[-dot-] (5.5,2) .. controls +(0,-0.75) and +(0,-0.75) .. (6.25,2);
		\draw (5.875,1.1) node[Odot] (unit1) {}; 
		\draw (unit1) -- (5.875,1.4);
		\draw[directedgreen] (6.25,2) .. controls +(0,0.75) and +(0,0.75) .. (7,2);
		\draw (7,2) -- (q6); 
		%
		\draw (4.5,3.75) -- (4.5,4);
		\fill (4.5,4) circle (2pt) node[right] {{\small $\gamma_{A^\dual}^{-1}$}};
		\draw[-dot-] (1,4) .. controls +(0,1) and +(0,1) .. (4.5,4);
		\draw (1,4) -- (q1);
		\draw (2.75,5.1) node[Odot] (unit2) {}; 
		\draw (unit2) -- (2.75,4.8);
		%
		\draw (q1) node[below] (bla2) {$A^\dual\vphantom{S_{(\alpha,A)^\dual}^{-a+1}}$};
		\draw (q2) node[below] (bla2) {$S_{(\alpha,A)^\dual}^{-a+1}$};
		\draw (q3) node[below] (bla2) {${}^\dagger\!A^\dual\vphantom{S_{(\alpha,A)^\dual}^{-a+1}}$};
		\draw (q4) node[below] (bla2) {$A^\dual\vphantom{S_{(\alpha,A)^\dual}^{-a+1}}$};
		\draw (q5) node[below] (bla2) {$S_{(\alpha,A)^\dual}^{+a+1}$};
		\draw (q6) node[below] (bla2) {$(A^\dual)^\dagger\vphantom{S_{(\alpha,A)^\dual}^{-a+1}}$};
		\draw[line width=0.5pt] ($(q2)+(0,0.5)$) node[inner sep=0pt,draw,fill=green!50!gray] (f) {{\scriptsize \phantom{opp}}};
		\draw[line width=0.5pt] ($(q5)+(0,0.5)$) node[inner sep=0pt,draw,fill=green!50!gray] (f) {{\scriptsize \phantom{opp}}};
		\end{tikzpicture}
	}%
	\label{eq:CAa_eps_mu}
	\, . 
	\end{align}
	The second step uses the properties of $\Delta$-separable Frobenius algebras and~\eqref{eq:NakayamaA}; moreover, here and below we suppress $\tev_\alpha$ and its adjoints (as they are only spectators in our string diagram manipulations). 
	Analogously, we arrive at 
	\be 
	\label{eq:CAa_Delta_eta}
	\Delta_{a,-a} \circ \eta_1 
	\equiv
	\tikzzbox{%
		\begin{tikzpicture}[very thick, scale=1,color=green!50!black, baseline=-2.4cm]
		\coordinate (q1) at (1,0);
		\coordinate (q2) at (2,0);
		\coordinate (q3) at (3,0);
		\coordinate (q4) at (5,0);
		\coordinate (q5) at (6,0);
		\coordinate (q6) at (7,0);
		%
		\draw (q1) -- ($(q1)+(0,-0.5)$);
		\draw[ultra thick] (q2) -- ($(q2)+(0,-0.5)$);
		\draw (q3) -- ($(q3)+(0,-1)$);
		\draw (q4) -- ($(q4)+(0,-0.5)$);
		\draw[ultra thick] (q5) -- ($(q5)+(0,-0.5)$);
		\draw (q6) -- ($(q6)+(0,-0.5)$);
		%
		\draw[directedgreen] ($(q3)+(0,-1)$) .. controls +(0,-0.75) and +(0,-0.75) .. ($(q3)+(0,-1)+(0.75,0)$);
		\draw[-dot-] ($(q3)+(0,-1)+(0.75,0)$) .. controls +(0,0.75) and +(0,0.75) .. ($(q3)+(0,-1)+(1.5,0)$);
		\draw ($(q3)+(0,-1)+(1.125,0.9)$) node[Odot] (unit1) {}; 
		\draw (unit1) -- ($(q3)+(0,-1)+(1.125,0.6)$);
		\draw ($(q2)+(0,-0.5)$) -- ($(q2)+(0,-1.5)$);
		\draw[-dot-] ($(q2)+(0,-2.2)+(-0.5,0)$) .. controls +(0,1) and +(0,1) .. ($(q2)+(0,-2.2)+(+0.5,0)$);
		\fill ($(q2)+(0,-2.2)+(+0.5,0)$) circle (2pt) node[right] {{\small $\gamma_{A^\dual}$}};
		\fill ($(q2)+(0,-2.2)+(-0.5,0)$) circle (2pt) node[right] {{\small $\gamma_{A^\dual}^{1-a}$}};
		\draw[-dot-] ($(q2)+(0,-2.2)+(+0.5,0)$) .. controls +(0,-1) and +(0,-1) .. ($(q3)+(0,-2.2)+(1.5,0)$);
		\draw ($(q3)+(0,-2.2)+(1.5,0)$) -- ($(q3)+(0,-1)+(1.5,0)$);
		\draw (3.5,-3) -- (3.5,-3.5);
		\draw[-dot-] (3.5,-3.5) .. controls +(0,-1) and +(0,-1) .. (5,-3.5);
		\fill (3.5,-3.5) circle (2pt) node[right] {{\small $\gamma_{A^\dual}$}};
		\draw ($(q4)+(0,-0.5)$) -- (5,-3.5); 
		\draw (4.25,-4.6) node[Odot] (unit2) {}; 
		\draw (unit2) -- (4.25,-4.3);
		\draw[-dot-] ($(q2)+(0,-2.2)+(-0.5,0)$) .. controls +(0,-4) and +(0,-4) .. ($(q5)+(0,-2.2)$);
		\draw ($(q5)+(0,-2.2)$) -- ($(q5)+(0,-0.5)$);
		\draw (3.75,-5.55) node[Odot] (unit3) {}; 
		\draw (unit3) -- (3.75,-5.2);
		\draw[redirectedgreen] ($(q1)+(0,-0.5)$) .. controls +(0,-7.5) and +(0,-7.5) .. ($(q6)+(0,-0.5)$);
		%
		\draw (q1) node[above] (bla2) {$A^\dual\vphantom{S_{(\alpha,A)^\dual}^{-a+1}}$};
		\draw (q2) node[above] (bla2) {$S_{(\alpha,A)^\dual}^{-a+1}$};
		\draw (q3) node[above] (bla2) {$(A^\dual)^\dagger\vphantom{S_{(\alpha,A)^\dual}^{-a+1}}$};
		\draw (q4) node[above] (bla2) {$A^\dual\vphantom{S_{(\alpha,A)^\dual}^{-a+1}}$};
		\draw (q5) node[above] (bla2) {$S_{(\alpha,A)^\dual}^{+a+1}$};
		\draw (q6) node[above] (bla2) {${}^\dagger\!A^\dual\vphantom{S_{(\alpha,A)^\dual}^{-a+1}}$};
		\draw[line width=0.5pt] ($(q2)+(0,-0.5)$) node[inner sep=0pt,draw,fill=green!50!gray] (f) {{\scriptsize \phantom{opp}}};
		\draw[line width=0.5pt] ($(q5)+(0,-0.5)$) node[inner sep=0pt,draw,fill=green!50!gray] (f) {{\scriptsize \phantom{opp}}};
		\end{tikzpicture}
	}%
	\, . 
	\ee 
	
	\vspace{-1.6cm}
	
	Combining~\eqref{eq:CAa_eps_mu} with~\eqref{eq:CAa_Delta_eta} into~\eqref{eq:Nakayama-Ca}, we employ the relation 
	\be 
	\tikzzbox{%
		\begin{tikzpicture}[very thick, scale=1,color=blue!50!black, baseline=0cm]
		\draw[postaction={decorate}, decoration={markings,mark=at position .59 with {\arrow[draw=blue!50!black]{>}}}] (-0.25,-0.25) .. controls +(0,0.75) and +(0,0.5) .. (1,0);
		\draw[postaction={decorate}, decoration={markings,mark=at position .59 with {\arrow[draw=blue!50!black]{<}}}] (-0.25,0.25) .. controls +(0,-0.75) and +(0,-0.5) .. (1,0);
		\draw (-0.25,-0.25) -- (-0.25,-1.25);
		\draw (-0.25,+0.25) -- (-0.25,+1.25);
		\fill (1,0) circle (2.5pt) node[right] {{\small $\varphi$}};
		\end{tikzpicture}
	}%
	= 
	\tikzzbox{%
		\begin{tikzpicture}[very thick, scale=1,color=blue!50!black, baseline=0cm]
		\draw[postaction={decorate}, decoration={markings,mark=at position .59 with {\arrow[draw=blue!50!black]{>}}}] (-0.25,-0.25) .. controls +(0,0.75) and +(0,0.5) .. (1,0);
		\draw[postaction={decorate}, decoration={markings,mark=at position .59 with {\arrow[draw=blue!50!black]{<}}}] (-0.25,0.25) .. controls +(0,-0.75) and +(0,-0.5) .. (1,0);
		\draw (-0.25,-0.25) -- (-0.25,-1.25);
		\draw (-0.25,+0.25) -- (-0.25,+1.25);
		\fill (-0.25,-0.75) circle (2.5pt) node[left] {{\small $\varphi^\dagger$}};
		\end{tikzpicture}
	}%
	\ee 
	to see that $N^{(\alpha,A)}_a$ is given by pre-composing the twist 
	\be 
	\label{eq:CaATwist}
	\tikzzbox{%
		\begin{tikzpicture}[very thick, scale=1,color=blue!50!black, baseline=0cm]
		\draw[postaction={decorate}, decoration={markings,mark=at position .59 with {\arrow[draw=blue!50!black]{>}}}] (-0.25,-0.25) .. controls +(0,0.75) and +(0,0.5) .. (1,0);
		\draw[postaction={decorate}, decoration={markings,mark=at position .59 with {\arrow[draw=blue!50!black]{<}}}] (-0.25,0.25) .. controls +(0,-0.75) and +(0,-0.5) .. (1,0);
		\draw (-0.25,-0.25) -- (-0.25,-1.25);
		\draw (-0.25,+0.25) -- (-0.25,+1.25);
		\draw (-0.25,-1.25) node[below] (bla2) {$C_a^{(\alpha,A)}$};
		\draw (-0.25,+1.25) node[above] (bla2) {$C_a^{(\alpha,A)}$};
		\end{tikzpicture}
	}%
	\ee 
	with (where we continue to suppress $\tev_\alpha$)
	\be 
	\tikzzbox
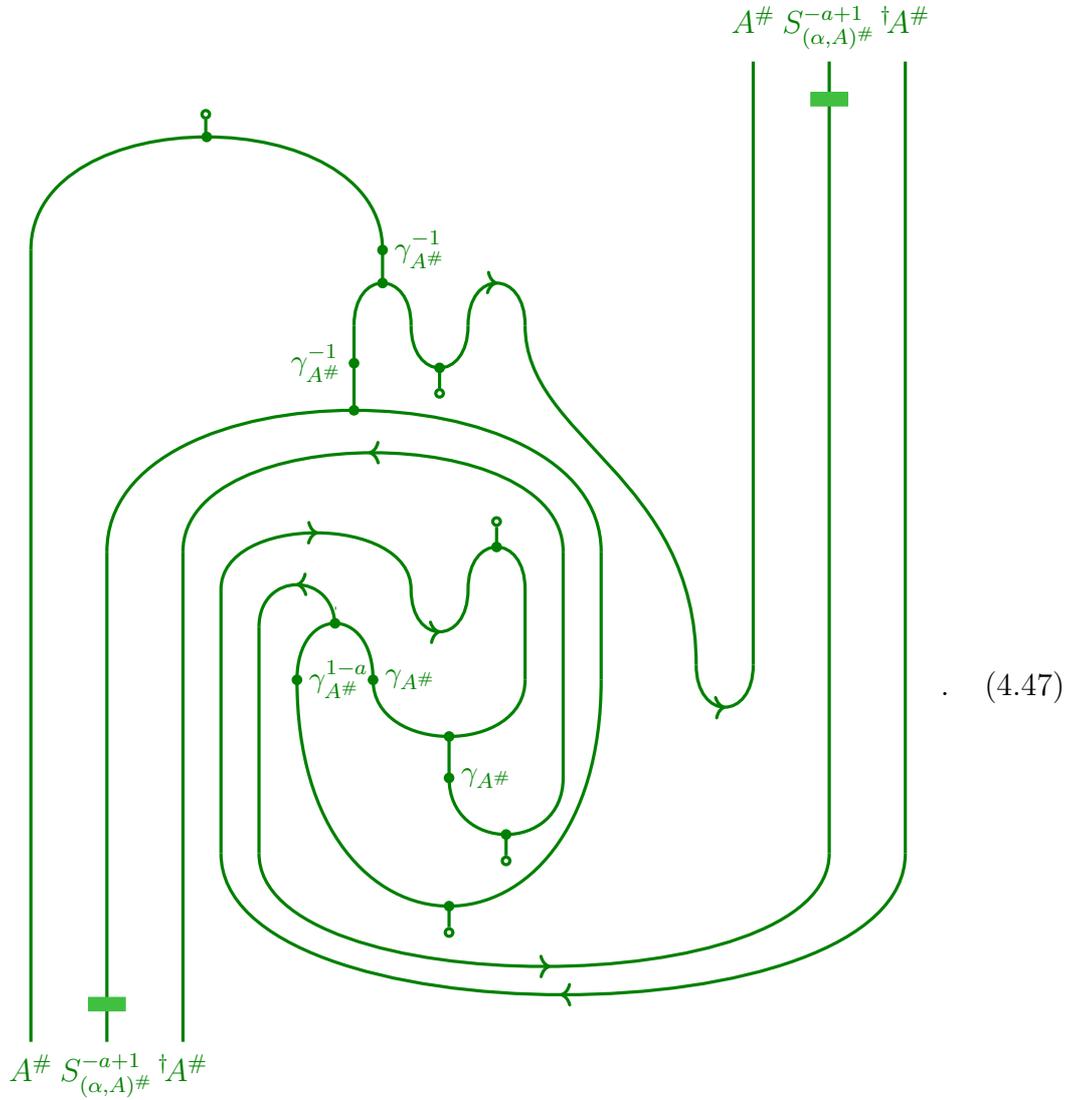
{%
		\begin{tikzpicture}[very thick, scale=1,color=green!50!black, baseline=-2.4cm]
		\coordinate (q1) at (1,0);
		\coordinate (q2) at (2,0);
		\coordinate (q3) at (3,0);
		\coordinate (q4) at (5,0);
		\coordinate (q5) at (6,0);
		\coordinate (q6) at (7,0);
		%
		\draw[directedgreen] ($(q3)+(0,-1)$) .. controls +(0,-0.75) and +(0,-0.75) .. ($(q3)+(0,-1)+(0.75,0)$);
		\draw[-dot-] ($(q3)+(0,-1)+(0.75,0)$) .. controls +(0,0.75) and +(0,0.75) .. ($(q3)+(0,-1)+(1.5,0)$);
		\draw ($(q3)+(0,-1)+(1.125,0.9)$) node[Odot] (unit1) {}; 
		\draw (unit1) -- ($(q3)+(0,-1)+(1.125,0.6)$);
		\draw ($(q2)+(0,-1.25)$) -- ($(q2)+(0,-1.25)$);
		\draw[-dot-] ($(q2)+(0,-2.2)+(-0.5,0)$) .. controls +(0,1) and +(0,1) .. ($(q2)+(0,-2.2)+(+0.5,0)$);
		\fill ($(q2)+(0,-2.2)+(+0.5,0)$) circle (2pt) node[right] {{\small $\gamma_{A^\dual}$}};
		\fill ($(q2)+(0,-2.2)+(-0.5,0)$) circle (2pt) node[right] {{\small $\gamma_{A^\dual}^{1-a}$}};
		\draw[-dot-] ($(q2)+(0,-2.2)+(+0.5,0)$) .. controls +(0,-1) and +(0,-1) .. ($(q3)+(0,-2.2)+(1.5,0)$);
		\draw ($(q3)+(0,-2.2)+(1.5,0)$) -- ($(q3)+(0,-1)+(1.5,0)$);
		\draw (3.5,-3) -- (3.5,-3.5);
		\draw[-dot-] (3.5,-3.5) .. controls +(0,-1) and +(0,-1) .. (5,-3.5);
		\fill (3.5,-3.5) circle (2pt) node[right] {{\small $\gamma_{A^\dual}$}};
		\draw ($(q4)+(0,-0.5)$) -- (5,-3.5); 
		\draw (4.25,-4.6) node[Odot] (unit2) {}; 
		\draw (unit2) -- (4.25,-4.3);
		\draw[-dot-] ($(q2)+(0,-2.2)+(-0.5,0)$) .. controls +(0,-4) and +(0,-4) .. ($(q5)+(-0.5,-2.2)$);
		\draw ($(q5)+(-0.5,-2.2)$) -- ($(q5)+(-0.5,-0.5)$);
		\draw (3.5,-5.55) node[Odot] (unit3) {}; 
		\draw (unit3) -- (3.5,-5.2);
		\draw[directedgreen] ($(q1)+(-0.5,-1)$) .. controls +(0,1) and +(0,1) .. ($(q3)+(0,-1)$);
		\draw[directedgreen] ($(q2)+(0,-1.5)$) .. controls +(0,0.75) and +(0,0.75) .. ($(q2)+(-1,-1.5)$);
		\draw ($(q2)+(-1,-4.5)$) -- ($(q2)+(-1,-1.5)$);
		\draw[directedgreen] ($(q4)+(0,-0.5)$) .. controls +(0,1.75) and +(0,1.75) .. ($(q1)+(-1,-0.5)$);
		\draw[-dot-] ($(q5)+(-0.5,-0.5)$) .. controls +(0,2.5) and +(0,2.5) .. ($(q1)+(-2,-0.5)$);
		\draw ($(q2)+(0.25,1.4)$) -- ($(q2)+(0.25,2.5)$);
		\fill ($(q2)+(0.25,2)$) circle (2pt) node[left] {{\small $\gamma_{A^\dual}^{-1}$}};
		\draw[-dot-] ($(q2)+(0.25,2.5)$) .. controls +(0,0.75) and +(0,0.75) .. ($(q2)+(1,2.5)$);
		\draw ($(q2)+(0.625,3.5)$) -- ($(q2)+(0.625,3)$);
		\fill ($(q2)+(0.625,3.5)$) circle (2pt) node[right] {{\small $\gamma_{A^\dual}^{-1}$}};
		\draw[-dot-] ($(q2)+(0.625,3.5)$) .. controls +(0,2) and +(0,2) .. ($(q1)+(-3,3.5)$);
		\draw ($(q2)+(-1.7,5.3)$) node[Odot] (unit5) {}; 
		\draw (unit5) -- ($(q2)+(-1.7,5)$); 
		\draw ($(q1)+(-3,3.5)$) -- ($(q1)+(-3,-7)$);
		\draw[-dot-] ($(q2)+(1,2.5)$) .. controls +(0,-0.75) and +(0,-0.75) .. ($(q2)+(1.75,2.5)$);
		\draw ($(q2)+(1.375,1.6)$) node[Odot] (unit4) {}; 
		\draw (unit4) -- ($(q2)+(1.375,2)$) ;
		\draw[directedgreen] ($(q2)+(1.75,2.5)$) .. controls +(0,0.75) and +(0,0.75) .. ($(q2)+(2.5,2.5)$);
		\draw ($(q2)+(2.5,2.5)$) .. controls +(0,-1.5) and +(0,2.5) .. ($(q2)+(4.75,-2)$);
		\draw[directedgreen] ($(q2)+(4.75,-2)$) .. controls +(0,-0.75) and +(0,-0.75) .. ($(q2)+(5.5,-2)$);
		\draw ($(q2)+(5.5,-2)$) -- ($(q2)+(5.5,6)$);
		\draw ($(q1)+(-2,-0.5)$) -- ($(q1)+(-2,-7)$);
		\draw ($(q1)+(-1,-0.5)$) -- ($(q1)+(-1,-7)$);
		\draw ($(q1)+(-0.5,-1)$) -- ($(q1)+(-0.5,-4.5)$);
		\draw[redirectedgreen] ($(q1)+(-0.5,-4.5)$) .. controls +(0,-2.5) and +(0,-2.5) .. ($(q6)+(2.5,-4.5)$);
		\draw[directedgreen] ($(q2)+(-1,-4.5)$) .. controls +(0,-2) and +(0,-2) .. ($(q6)+(1.5,-4.5)$);
		\draw ($(q6)+(1.5,-4.5)$) -- ($(q6)+(1.5,6)$);
		\draw ($(q6)+(2.5,-4.5)$) -- ($(q6)+(2.5,6)$);
		%
		\draw ($(q6)+(0.5,6)$) node[above] (bla2) {$A^\dual\vphantom{S_{(\alpha,A)^\dual}^{-a+1}}$};
		\draw ($(q6)+(1.5,6)$) node[above] (bla2) {$S_{(\alpha,A)^\dual}^{-a+1}$};
		\draw ($(q6)+(2.5,6)$) node[above] (bla2) {${}^\dagger\! A^\dual\vphantom{S_{(\alpha,A)^\dual}^{-a+1}}$};
		\draw ($(q1)+(-3,-7)$) node[below] (bla2) {$A^\dual\vphantom{S_{(\alpha,A)^\dual}^{-a+1}}$};
		\draw ($(q1)+(-2,-7)$) node[below] (bla2) {$S_{(\alpha,A)^\dual}^{-a+1}$};
		\draw ($(q1)+(-1,-7)$) node[below] (bla2) {${}^\dagger\! A^\dual\vphantom{S_{(\alpha,A)^\dual}^{-a+1}}$};
		\draw[line width=0.5pt] ($(q1)+(-2,-6.5)$) node[inner sep=0pt,draw=green!50!gray,fill=green!50!gray] (f) {{\scriptsize \phantom{opp}}};
		\draw[line width=0.5pt] ($(q6)+(1.5,5.5)$) node[inner sep=0pt,draw=green!50!gray,fill=green!50!gray] (f) {{\scriptsize \phantom{opp}}};
		\end{tikzpicture}
	}%
	.
	\ee 
	Repeatedly using the defining properties of $\Delta$-separable Frobenius algebras as well the properties of the Nakayama automorphism~$\gamma_{A^\dual}$ collected in Section~\ref{subsubsec:AdjointsInBeq}, a straightforward but lengthy computation shows that the above string diagram is equal to 
	\be 
	\label{eq:LastDiagramInProof}
	\tikzzbox{%
		\begin{tikzpicture}[very thick, scale=1,color=green!50!black, baseline=3cm]
		\coordinate (d1) at (0,0);
		\coordinate (d2) at (2,0);
		\coordinate (d3) at (4,0);
		\coordinate (u1) at (0,6);
		\coordinate (u2) at (2,6);
		\coordinate (u3) at (7,6);%
		\draw (d1) -- (u1);
		\draw (d2) -- (u2);
		\draw[ultra thick] (d2) -- ($(d2)+(0,+0.5)$);
		\draw[ultra thick] (u2) -- ($(u2)+(0,-0.5)$);
		%
		\fill ($(d1)+(0,2)$) circle (2pt) node[right] {{\small $\gamma_{A^\dual}$}};
		\fill ($(d2)+(0,2)$) circle (2pt) node[right] {{\small $\gamma_{A^\dual}$}};
		\fill ($(d3)+(1.5,2)$) circle (2pt) node[right] {{\small $\gamma_{A^\dual}$}};
		\fill ($(d1)+(0.625,4)+(0.75,0)$) circle (2pt) node[left] {};
		\fill ($(d1)+(0.75,4)+(0.75,0)$) circle (0pt) node[left] {{\scriptsize $\gamma_{A^\dual}^{1-a}$}};
		\fill ($(d2)+(1.375,3.5)+(0.75,0)$) circle (2pt) node[left] {};
		\fill ($(d2)+(1.375,3.5)+(0.85,0)$) circle (0pt) node[left] {{\scriptsize $\gamma_{A^\dual}$}};
		\draw (d3) -- ($(d3)+(0,+1)$);
		\draw[redirectedgreen] ($(d3)+(0,+1)$) .. controls +(0,0.75) and +(0,0.75) .. ($(d3)+(0,+1)+(0.75,0)$);
		\draw[-dot-] ($(d3)+(0,+1)+(0.75,0)$) .. controls +(0,-0.75) and +(0,-0.75) .. ($(d3)+(0,+1)+(1.5,0)$);
		\draw ($(d3)+(0,+1)+(1.125,-0.9)$) node[Odot] (unit1) {}; 
		\draw (unit1) -- ($(d3)+(0,+1)+(1.125,-0.5)$);
		\draw (u3) -- ($(u3)+(0,-1)$);
		\draw[directedgreen] ($(u3)+(0,-1)$) .. controls +(0,-0.75) and +(0,-0.75) .. ($(u3)+(0,-1)+(-0.75,0)$);
		\draw[-dot-] ($(u3)+(0,-1)+(-0.75,0)$) .. controls +(0,0.75) and +(0,0.75) .. ($(u3)+(0,-1)+(-1.5,0)$);
		\draw ($(u3)+(0,-1)+(-1.125,0.9)$) node[Odot] (unit2) {}; 
		\draw (unit2) -- ($(u3)+(0,-1)+(-1.125,0.5)$);
		\draw ($(u3)+(0,-1)+(-1.5,0)$) -- ($(d3)+(0,+1)+(1.5,0)$);
		%
		\draw[-dot-] ($(d1)+(0.625,4)$) .. controls +(0,-0.75) and +(0,-0.75) .. ($(d1)+(0.625,4)+(0.75,0)$);
		\draw ($(d1)+(1,3.1)$) node[Odot] (unit3) {}; 
		\draw (unit3) -- ($(d1)+(1,3.4)$);
		\draw ($(d1)+(0.625,4)$) .. controls +(0,0.5) and +(0.2,-0.3) .. ($(d1)+(0,4.8)$);
		\draw ($(d1)+(0.625,4)+(0.75,0)$) .. controls +(0,0.5) and +(-0.2,-0.3) .. ($(d2)+(0,4.8)$);
		\fill ($(d2)+(0,4.8)$) circle (2pt) node[right] {};
		\fill ($(d1)+(0,4.8)$) circle (2pt) node[right] {};
		%
		\draw[-dot-] ($(d2)+(1.375,3.5)$) .. controls +(0,-0.75) and +(0,-0.75) .. ($(d2)+(1.375,3.5)+(0.75,0)$);
		\draw ($(d2)+(1.375,3.5)+(0.375,-0.9)$) node[Odot] (unit4) {}; 
		\draw (unit4) -- ($(d2)+(1.375,3.5)+(0.375,-0.6)$);
		\draw ($(d2)+(1.375,3.5)$) .. controls +(0,0.5) and +(0.2,-0.3) .. ($(d2)+(0,4.3)$);
		\draw ($(d2)+(1.375,3.5)+(0.75,0)$) .. controls +(0,0.5) and +(-0.2,-0.3) .. ($(d2)+(3.5,4.3)$);
		\fill ($(d2)+(0,4.3)$) circle (2pt) node[right] {};
		\fill ($(d2)+(3.5,4.3)$) circle (2pt) node[right] {};
		%
		\draw (d1) node[below] (bla2) {$A^\dual\vphantom{S_{(\alpha,A)^\dual}^{-a+1}}$};
		\draw (d2) node[below] (bla2) {$S_{(\alpha,A)^\dual}^{-a+1}$};
		\draw (d3) node[below] (bla2) {${}^\dagger\! A^\dual\vphantom{S_{(\alpha,A)^\dual}^{-a+1}}$};
		\draw (u1) node[above] (bla2) {$A^\dual\vphantom{S_{(\alpha,A)^\dual}^{-a+1}}$};
		\draw (u2) node[above] (bla2) {$S_{(\alpha,A)^\dual}^{-a+1}$};
		\draw (u3) node[above] (bla2) {${}^\dagger\! A^\dual\vphantom{S_{(\alpha,A)^\dual}^{-a+1}}$};
		\draw[line width=0.5pt] ($(d2)+(0,+0.5)$) node[inner sep=0pt,draw=green!50!gray,fill=green!50!gray] (f) {{\scriptsize \phantom{opp}}};
		\draw[line width=0.5pt] ($(u2)+(0,-0.5)$) node[inner sep=0pt,draw=green!50!gray,fill=green!50!gray] (f) {{\scriptsize \phantom{opp}}};
		\end{tikzpicture}
	}%
	.
	\ee 
	Putting $\tev_{\alpha}, \tev_{\alpha}^\dagger, S_{\alpha}^{1-a}$ back in, a final application of the isomorphisms~\eqref{eq:NakayamaLeftRight} and~\eqref{eq:SApower} allows us to identify~\eqref{eq:LastDiagramInProof} with
	\be 
	\tikzzbox{%
		\begin{tikzpicture}[very thick, scale=1,color=green!50!black, baseline=2.5cm]
		\coordinate (d2) at (2,0);
		\fill [orange!20!white, opacity=0.8] 
		($(d2)+(-3.3,1)$) -- ($(d2)+(-3.3,4)$) -- ($(d2)+(+3.3,4)$) -- ($(d2)+(+3.3,1)$)
		;
		%
		\draw[thick, color=black] ($(d2)+(-3.3,1)$) -- ($(d2)+(-3.3,4)$); 
		\fill[color=black] ($(d2)+(-3.65,1.6)$) circle (0pt) node {{\small $\tev_\alpha\vphantom{{}_{\gamma^{2-a}_{A^\dual}}A^\dual}$}};
		%
		\draw[thick, color=black] ($(d2)+(+3.3,1)$) -- ($(d2)+(+3.3,4)$); 
		\fill[color=black] ($(d2)+(3.7,1.6)$) circle (0pt) node {{\small $\tev_\alpha^\dagger\vphantom{{}_{\gamma^{2-a}_{A^\dual}}A^\dual}$}};
		%
		\draw ($(d2)+(-2,1)$) -- ($(d2)+(-2,4)$);
		\fill ($(d2)+(-2.3,1.6)$) circle (0pt) node {{\small $A^\dual\vphantom{{}_{\gamma^{2-a}_{A^\dual}}A^\dual}$}};
		%
		\draw[ultra thick] ($(d2)+(0,1)$) -- ($(d2)+(0,4)$);
		\fill ($(d2)+(0.7,1.6)$) circle (0pt) node {{\small ${}_{\gamma^{2-a}_{A^\dual}}A^\dual$}}; 
		%
		\draw[color=blue!50!black] ($(d2)+(+2,1)$) -- ($(d2)+(+2,4)$);
		\fill[color=blue!50!black] ($(d2)+(2.5,1.6)$) circle (0pt) node {{\small $S_\alpha^{1-a}\vphantom{{}_{\gamma^{2-a}_{A^\dual}}A^\dual}$}}; 
		\draw[ultra thick, color=blue!50!black] (d2) -- ($(d2)+(0,1)$); 
		\draw[ultra thick, color=blue!50!black] ($(d2)+(0,5)$) -- ($(d2)+(0,4)$); 
		%
		\fill ($(d2)+(0,2.5)$) circle (3pt) node[right] {{\small $\gamma_{A^\dual}$}};
		\fill[color=blue!50!black] (d2) circle (0pt) node[below] {{\small $C_a^{(\alpha,A)}$}};
		\fill[color=blue!50!black] ($(d2)+(0,5)$) circle (0pt) node[above] {{\small $C_a^{(\alpha,A)}$}};
		\draw[line width=0.5pt] ($(d2)+(0,4)$) node[inner sep=1pt,draw=green!50!gray,fill=green!50!gray] (f) {{\scriptsize \phantom{opppppppppppppppppppppppppppppppppppppppppp}}};
		\draw[line width=0.5pt] ($(d2)+(0,1)$) node[inner sep=1pt,draw=green!50!gray,fill=green!50!gray] (f) {{\scriptsize \phantom{opppppppppppppppppppppppppppppppppppppppppp}}};
		\end{tikzpicture}
	}%
	.
	\ee
	Post-composing with~\eqref{eq:CaATwist} thus completes the proof. 
\end{proof}

Combining Proposition~\ref{prop:NakayamaPowerBeq} with the isomorphisms~\eqref{eq:SalphaAinBeq} and~\eqref{eq:SApower}, we obtain closed $\Lambda_r$-Frobenius algebras from $\Delta$-separable Frobenius algebras~$A$ on fully dualisable objects if the $r$-th power of the Nakayama automorphism of~$A$ is the identity: 

\begin{corollary}
	If for $r\in\Z_{\geqslant 1}$ there is an isomorphism $S_\alpha^r \cong 1_\alpha$ in $\Bfd$, and for $(\alpha,A)\in\Beq$ we have $\gamma_A^r = 1_A$, then there is an induced closed $\Lambda_r$-Frobenius algebra structure on $\{ C^{(\alpha,A)}_a \}_{a\in\{ 0,1,\dots,r-1\}}$. 
\end{corollary}

If~$\B$ has internal Homs, as is the case in the examples related to TQFTs of state sum, sigma model and Landau--Ginzburg type, then the computation of both $C^{(\alpha,A)}_a$ and $N^{(\alpha,A)}_a$ can be simplified. 
For ease of presentation, we further assume that $S_\alpha \cong 1_\alpha$; in Section~\ref{subsec:LGmodels} below we will see how this restriction can be lifted in practice. 

\begin{lemma}
	\label{lem:CaalphaAFromProjector}
	Let $\alpha\in\Bfd$ with $S_\alpha \cong 1_\alpha$, and assume that for $(\alpha,A)\in\Beq$ we have $C^{(\alpha,A)}_a \cong \B( 1_\one, C^{(\alpha,A)}_a)$. 
	Then 
	\be 
	\label{eq:CaTwisted}
	C^{(\alpha,A)}_a 
		\cong 
		\Bigg\{ 
		 \varphi \in \B(1_\alpha,A) \;\Bigg|\;  
		 	\tikzzbox{%
		 	\begin{tikzpicture}[very thick,scale=0.65,color=green!50!black, baseline=0cm]
		 	\draw (0,0) -- (0,1.3);
		 	\fill (0,0) circle (2.5pt) node[below] {{\small $\varphi$}};
		 	\fill (-0.67,0) circle (2.5pt) node[left] {{\small $\gamma_A^{1-a}$}};
		 	\draw (0,0.8) .. controls +(-0.9,-0.3) and +(-0.9,0) .. (0,-0.8);
		 	\draw (0,-0.8) .. controls +(0.9,0) and +(0.7,-0.1) .. (0,0.4);
		 	\fill (0,-0.8) circle (2.5pt) node {};
		 	\fill (0,0.4) circle (2.5pt) node {};
		 	\fill (0,0.8) circle (2.5pt) node {};
		 	\draw (0,-1.2) node[Odot] (unit) {};
		 	\draw (0,-0.8) -- (unit);
		 	\end{tikzpicture}
		 	}%
		 	= 
		 	\tikzzbox{%
		 	\begin{tikzpicture}[very thick,scale=0.65,color=green!50!black, baseline]
		 	\fill (0,-0.3) circle (2.5pt) node[right] (D) {{\small $\varphi$}};
		 	\draw (0,-0.3) -- (0,0.8); 
		 	\end{tikzpicture} 
		 	}%
		\Bigg\}
	\ee 
	and $N^{(\alpha,A)}_a$ corresponds to post-composition with~$\gamma_A$. 
\end{lemma}

Note that the above result further elucidates the relation between the two different notions of ``Nakayama morphism''~$N$ and~$\gamma$. 

\begin{proof}[Proof of Lemma~\ref{lem:CaalphaAFromProjector}]
	We have 
	\begin{align}
	C^{(\alpha,A)}_a 
		& \cong \Beq\Big( 1_{(\one,1_\one)}, C^{(\alpha,A)}_a \Big) \nonumber
		\\
		& \cong \Beq\Big( \tev_{(\alpha,A)}, \tev_{(\alpha,A)} \otimes_{\Ae} \big( S_{(\alpha,A)}^{1-a} \btimes 1_{\alpha^\dual} \big)  \Big) \nonumber
		\\
		& \cong \Beq \big( A, S_{(\alpha,A)}^{1-a} \big)\nonumber
		\\
		& \cong \Beq \big( A, {}_{\gamma_A^{1-a}}A \big) \nonumber
		\\ 
		& \cong 
		\Bigg\{ 
		\varphi \in \B(1_\alpha,A) \;\Bigg|\;  
		\tikzzbox{%
		\begin{tikzpicture}[very thick,scale=0.65,color=green!50!black, baseline=0cm]
		\draw (0,0) -- (0,1.3);
		\fill (0,0) circle (2.5pt) node[below] {{\small $\varphi$}};
		\fill (-0.67,0) circle (2.5pt) node[left] {{\small $\gamma_A^{1-a}$}};
		\draw (0,0.8) .. controls +(-0.9,-0.3) and +(-0.9,0) .. (0,-0.8);
		\draw (0,-0.8) .. controls +(0.9,0) and +(0.7,-0.1) .. (0,0.4);
		\fill (0,-0.8) circle (2.5pt) node {};
		\fill (0,0.4) circle (2.5pt) node {};
		\fill (0,0.8) circle (2.5pt) node {};
		\draw (0,-1.2) node[Odot] (unit) {};
		\draw (0,-0.8) -- (unit);
		\end{tikzpicture}
		}%
		= 
		\tikzzbox{%
		\begin{tikzpicture}[very thick,scale=0.65,color=green!50!black, baseline]
		\fill (0,-0.3) circle (2.5pt) node[right] (D) {{\small $\varphi$}};
		\draw (0,-0.3) -- (0,0.8); 
		\end{tikzpicture} 
		}%
		\Bigg\} . 
	\end{align}
	In the second step we used~\eqref{eq:CalphaA} and adjunction for $\tev_{(\alpha,A)}$; 
	the third step is the isomorphism 
	\be 
	\tikzzbox{%
		\begin{tikzpicture}[thick,scale=1.0,color=black, baseline=1.75cm]
		\coordinate (p1) at (0,0);
		\coordinate (p2) at (2,-0.5);
		\coordinate (p3) at (2.5,0.5);
		\coordinate (u1) at (0,3);
		\coordinate (u2) at (2,2.5);
		\coordinate (u3) at (2.5,3.5);
		%
		\fill [orange!20!white, opacity=0.8] 
		(p1) .. controls +(0,0.25) and +(-1,0) ..  (p3)
		-- (p3) --  (u3)
		-- (u3) .. controls +(-1,0) and +(0,0.25) ..  (u1)
		;
		%
		\draw[very thick, red!80!black] (p1) .. controls +(0,0.25) and +(-1,0) ..  (p3); 
		%
		\fill [orange!30!white, opacity=0.8] 
		(p1) .. controls +(0,-0.25) and +(-1,0) ..  (p2)
		-- (p2) --  (u2)
		-- (u2) .. controls +(-1,0) and +(0,-0.25) ..  (u1)
		;
		%
		\draw[ultra thick, green!50!black] (2.25,3.5) -- (2.25,2); 
		\fill[color=green!50!black] (2.25,2) circle (2pt) node[below] {{\small $\varphi$}};
		%
		\draw[thin] (p1) --  (u1); 
		\draw[thin] (p2) --  (u2); 
		\draw[thin] (p3) --  (u3); 
		%
		\draw[very thick, red!80!black] (p1) .. controls +(0,-0.25) and +(-1,0) ..  (p2); 
		\draw[very thick, red!80!black] (u1) .. controls +(0,-0.25) and +(-1,0) ..  (u2); 
		\draw[very thick, red!80!black] (u1) .. controls +(0,0.25) and +(-1,0) ..  (u3); 
		%
		\fill[color=blue!50!black] (1.5,0) circle (0pt) node {{\small $\alpha^\dual$}};
		\fill[color=blue!50!black] (1.1,3) circle (0pt) node {{\small $\alpha$}};
		%
		\end{tikzpicture}
	}%
	\;
	\lmt 
	\;
	\tikzzbox{%
		\begin{tikzpicture}[thick,scale=1.0,color=black, xscale=-1, baseline=2.75cm]
		\coordinate (p1) at (0,0);
		\coordinate (p2) at (1.5,-0.5);
		\coordinate (p2s) at (4,-0.5);
		\coordinate (p3) at (1.5,0.5);
		\coordinate (p4) at (3,1);
		\coordinate (p5) at (1.5,1.5);
		\coordinate (p6) at (-1,1.5);
		\coordinate (u1) at (0,3);
		\coordinate (u2) at (1.5,2.5);
		\coordinate (u2s) at (4,2.5);
		\coordinate (u3) at (1.5,3.5);
		\coordinate (u4) at (3,4);
		\coordinate (u5) at (1.5,4.5);
		\coordinate (u6) at (-1,4.5);
		%
		\fill [orange!20!white, opacity=0.8] 
		(p4) .. controls +(0,0.25) and +(1,0) ..  (p5) -- (p6)
		-- (u6) -- (u5)
		-- (u5) .. controls +(1,0) and +(0,0.25) .. (u4);
		%
		\fill [orange!25!white, opacity=0.8] 
		(p1) .. controls +(0,0.25) and +(-1,0) ..  (p3)
		-- (p3) .. controls +(1,0) and +(0,-0.25) ..  (p4)
		-- (p4) --  (u4)
		-- (u4) .. controls +(0,-0.25) and +(1,0) ..  (u3)
		-- (u3) .. controls +(-1,0) and +(0,0.25) ..  (u1)
		;
		%
		\draw[thin] (p1) --  (u1); 
		\draw[thin] (p2s) --  (u2s); 
		\draw[thin] (p4) --  (u4); 
		\draw[thin] (p6) --  (u6); 
		%
		\draw[ultra thick, green!50!black] (-0.5,4.5) -- (-0.5,3); 
		\fill[color=green!50!black] (-0.5,3) circle (2pt) node[right] {{\small $\varphi$}};
		%
		\draw[very thick, red!80!black] (p1) .. controls +(0,-0.25) and +(-1,0) ..  (p2) -- (p2s); 
		\draw[very thick, red!80!black] (p4) .. controls +(0,0.25) and +(1,0) ..  (p5) -- (p6); 
		%
		\draw[very thick, red!80!black] (p1) .. controls +(0,0.25) and +(-1,0) ..  (p3)
		-- (p3) .. controls +(1,0) and +(0,-0.25) ..  (p4); 
		%
		\fill [orange!30!white, opacity=0.8] 
		(p1) .. controls +(0,-0.25) and +(-1,0) ..  (p2)
		-- (p2) -- (p2s)
		-- (p2s) --  (u2s)
		-- (u2s) -- (u2)
		-- (u2) .. controls +(-1,0) and +(0,-0.25) ..  (u1)
		;
		%
		\draw[very thick, red!80!black] (p1) .. controls +(0,-0.25) and +(-1,0) ..  (p2) -- (p2s); 
		\draw[very thick, red!80!black] (u1) .. controls +(0,-0.25) and +(-1,0) ..  (u2) -- (u2s); 
		\draw[very thick, red!80!black] (u1) .. controls +(0,0.25) and +(-1,0) ..  (u3)
		-- (u3) .. controls +(1,0) and +(0,-0.25) ..  (u4); 
		\draw[very thick, red!80!black] (u4) .. controls +(0,0.25) and +(1,0) ..  (u5) -- (u6); 
		%
		\fill[color=blue!50!black] (1.75,3) circle (0pt) node {{\small $\alpha^\dual$}};
		\fill[color=blue!50!black] (-0.5,2) circle (0pt) node {{\small $\alpha$}};
		\fill[color=blue!50!black] (3.5,0) circle (0pt) node {{\small $\alpha$}};
		%
		\end{tikzpicture}
	}    
	\;
	\stackrel{\cong}{\lmt} 
	\;
	\tikzzbox{%
		\begin{tikzpicture}[thick,scale=1.0,color=black, baseline=2.75cm]
		\coordinate (p5) at (1.5,1.5);
		\coordinate (p6) at (-1,1.5);
		\coordinate (u5) at (1.5,4.5);
		\coordinate (u6) at (-1,4.5);
		%
		\fill [orange!20!white, opacity=0.8] 
		(p5) -- (p6)
		-- (u6) -- (u5);
		%
		\draw[thin] (p5) --  (u5); 
		\draw[thin] (p6) --  (u6); 
		%
		\draw[ultra thick, green!50!black] (0.25,4.5) -- (0.25,3); 
		\fill[color=green!50!black] (0.25,3) circle (2pt) node[right] {{\small $\varphi$}};
		%
		\draw[very thick, red!80!black] (p5) -- (p6); 
		\draw[very thick, red!80!black] (u5) -- (u6); 
		%
		\fill[color=blue!50!black] (-0.5,2) circle (0pt) node {{\small $\alpha$}};
		%
		\end{tikzpicture}
	} 
	\, ; 
	\ee 
	the fourth step is~\eqref{eq:SApower} together with the assumption $S_\alpha \cong 1_\alpha$; 
	the fifth step is a standard computation with $\Delta$-separable Frobenius algebras along the lines of \cite[Sect.\,3.2]{BCP}. 
\end{proof}

\subsection{Landau--Ginzburg models}
\label{subsec:LGmodels}

In this section we briefly review the 2-category of Landau--Ginzburg models $\LG$ and note that every object in $\LG$ gives rise to an extended 2-spin TQFT. 
Then we apply the results of Section~\ref{subsec:EquivariantCompletion} to a closely related 2-category $\LGs$ and consider the simplest non-trivial example. 

\medskip 

Recall from \cite{cm1208.1481} that for every fixed field~$\Bbbk$, there is a 2-category $\LG$ whose objects are pairs $(\Bbbk[x_1,\dots,x_n], W)$, where $n\in \Z_{\geqslant 0}$ and $W=0\in \Bbbk$ if $n=0$, while for $n>0$, $W\in\Bbbk[x] \equiv \Bbbk[x_1,\dots,x_n]$ is such that the \textsl{Jacobi algebra} 
\be 
\Jac_W = \Bbbk[x_1,\dots,x_n]\big/\big( \partial_{x_1}W, \dots, \partial_{x_n}W \big) 
\ee 
is finite-dimensional over~$\Bbbk$. 
We refer to such polynomials~$W$ as \textsl{potentials}. 
The Hom categories of $\LG$ are idempotent completions of homotopy categories of finite-rank matrix factorisations. 
Hence up to technicalities with idempotents (which will not be relevant to our discussions below), a 1-morphism $(\Bbbk[x],W) \lra (\Bbbk[z],V)$ is a free $\Z_2$-graded $\Bbbk[x,z]$-module $X = X^0\oplus X^1$ together with an odd $\Bbbk[x,z]$-linear endomorphism $d_X\colon X \lra X$ such that $d_X^2 = (V-W) \cdot 1_X$. 
The Hom sets of 2-morphisms $(X,d_X) \lra (X',d_{X'})$ consist of the even cohomology classes of the differential defined on $\Z_2$-homogeneous maps as 
\begin{align}
\delta_{X,X'} \colon \Hom_{\Bbbk[x,z]}(X,X') & \lra  \Hom_{\Bbbk[x,z]}(X,X') \nonumber 
	\\ 
	\zeta & \lmt d_{X'} \circ \zeta - (-1)^{|\zeta|} \zeta \circ d_X \, , 
	\label{eq:MFdiff}
\end{align}
and extended linearly to all of $\Hom_{\Bbbk[x,z]}(X,X')$. 

Given $(X,d_X) \colon (\Bbbk[x], W_1) \lra (\Bbbk[y], W_2)$ and $(Y,d_Y) \colon (\Bbbk[y], W_2) \lra (\Bbbk[z], W_3)$, their horizontal composition is $(Y\otimes_{\Bbbk[y]} X, d_Y\otimes 1_X + 1\otimes d_X)$, and the unit 1-morphism of $(\Bbbk[x_1,\dots,x_n], W)$ is $1_W = (I_W,d_{I_W})$ with  
\be 
\label{eq:IW}
I_W = \bigwedge \Big( \bigoplus_{i=1}^n \Bbbk[x,x'] \cdot \theta_i \Big) 
	, \quad  
	d_{I_W} = \sum_{i=1}^n \Big( \partial_{[i]}^{x',x} W \cdot \theta_i + \big( x'_i - x_i \big) \cdot \theta_i^* \Big) , 
\ee 
where $\{\theta_i\}$ is a chosen $\Bbbk[x',x]$-basis of $\Bbbk[x,x']^{\oplus n}$, and 
\be 
\partial_{[i]}^{x',x} W 
= 
\frac{W(x_1,\dots,x_{i-1}, x'_i, \dots x'_n) - W(x_1,\dots,x_i, x'_{i+1}, \dots x'_n)}{x'_i-x_i} \, . 
\ee 
A straightforward computation shows that $\End(1_W) \cong \Jac_W$ in $\LG$. 

Every 1-morphism $X \equiv (X,d_X) \in \LG( (\Bbbk[x_1,\dots,x_n], W), (\Bbbk[z_1,\dots,z_m], V) )$ has a left adjoint~$\dX$ and a right adjoint~$\Xd$, and $\dX \cong \Xd$ iff $m=n \, \textrm{mod} \, 2$. 
The associated adjunction 2-morphisms are explicitly known, see \cite[Thm.\,6.11]{cm1208.1481}, or \cite{cm1303.1389} for a concise review. 

The 2-category $\LG$ has a natural monoidal structure, with the monoidal product on objects given by $(\Bbbk[x],W) \btimes (\Bbbk[z],V) = (\Bbbk[x,z],W+V)$, while on 1- and 2-morphisms it is basically~$\otimes_\Bbbk$, see \cite[Sect.\,2.2]{CMM}. 
Hence $\one := (\Bbbk,0)$ is the unit object. 
Every $(\Bbbk[x],W)\in\LG$ has a (left and right) dual $(\Bbbk[x],W)^\dual = (\Bbbk[x],-W)$ whose associated adjunction 1-morphisms have~$1_W$ as their underlying matrix factorisation, see \cite[Prop.\,2.6]{CMM}. 
Thus, as every 1-morphism in $\LG$ has an adjoint, every object in $\LG$ is fully dualisable. 

The monoidal 2-category $\LG$ has a symmetric braiding, whose 1-morphism components $b_{V,W}$ are given by $1_{V+W}$ (up to a reordering of variables), while the 2-morphism components are compositions of canonical module isomorphisms and structure maps of the underlying 2-category $\LG$. 
For details we refer to \cite[Sect.\,2.3]{CMM}. 
In summary, we have: 

\begin{theorem}[\cite{cm1208.1481, CMM}]
	\label{thm:LGsymmon}
	For every field~$\Bbbk$, the 2-category of Landau--Ginzburg models $\LG$ has a symmetric monoidal structure such that $\LG=\LG^{\textrm{fd}}$. 
\end{theorem}

\begin{remark}
	A variant of $\LG$ is the symmetric monoidal 2-category $\LGs$, which is defined analogously to $\LG$, but 
		 2-morphisms 
	are given by both even and odd cohomology of the differentials~$\delta_{X,X'}$ in~\eqref{eq:MFdiff}, but with classes $-\zeta$ and $+\zeta$ identified. 
	This ad hoc $\Z_2$-quotient allows to stay within the realm of 2-categories (as opposed to super 2-categories) while allowing odd 2-morphisms, compare \cite{kr0401268} and \cite[Rem.\,3.11(ii)]{CMM}. 
	
	Theorem~\ref{thm:LGsymmon} also holds for $\LGs$, i.\,e.\ $\LGs = (\LGs)^{\textrm{fd}}$. 
\end{remark}

The Serre automorphism of $W\equiv (\Bbbk[x_1,\dots,x_n], W)$ was computed in \cite[Lem.\,3.8]{CMM} to be 
\be 
\label{eq:SWLGIW}
S_W \cong 1_W[n] \, , 
\ee 
where $[n]$ denotes the $n$-fold application of the shift functor $[1]$, which sends a matrix factorisation $(X^0\oplus X^1,d_X)$ to $(X^1\oplus X^0, -d_X)$. 
It follows that $[2]$ is the identity functor, and one finds that $\Hom(1_W,1_W[n]) \cong \delta_{n,0\,\textrm{mod}\,2} \cdot \Jac_W$ in $\LG$, as $\LG$ has only even cohomology classes as 2-morphisms, while $\Hom(1_W,1_W[n]) \cong \Jac_W[n]$ is purely odd in $\LGs$ if~$n$ is odd. 
As a consequence, $(\Bbbk[x_1,\dots,x_n],W)$ determines an extended oriented TQFT with values in $\LG$ iff~$n$ is even, and it determines an extended oriented TQFT with values in $\LGs$ for every value of~$n$. 

\begin{remark}
	As shown in \cite[Sect.\,3]{CMM}, fully extended oriented TQFTs with values in $\LG$ are indeed extensions of closed Landau--Ginzburg models to the point: 
	Applying the cobordism hypothesis of \cite{spthesis, Hesse} to the duality data of an object $(\Bbbk[x],W)$ in $\LG$ or in $\LGs$, one recovers the (non-semisimple) commutative Frobenius algebra $\Jac_W$ from the circle, the pair-of-pants, and the disc. 
\end{remark}

As an immediate consequence of Theorem~\ref{thm:r-spin-CH}, Lemma~\ref{lem:r-spin-fixed-points}, \eqref{eq:SWLGIW} and the isomorphism $\textrm{Aut}(1_W) \cong \Bbbk$, we find that every potential depending on an odd number of variables gives rise to a proper extended spin TQFT: 

\begin{theorem}
	\label{thm:LGspin}
	Every object $W\equiv (\Bbbk[x_1,\dots,x_n], W) \in \LG$ gives rise to a unique-up-to-isomorphism extended 2-spin TQFT valued in $\LG$. 
	These TQFTs factor through the oriented bordism 2-category iff~$n$ is even. 
\end{theorem}

\medskip 

It is straightforward to compute that 
\be 
C_a^{W} \equiv 
C_a^{(\Bbbk[x_1,\dots,x_n],W)} \cong \Jac_W [n\cdot (1-a)] 
	\quad \textrm{in } \LG(\one,\one) \cong \textrm{vect}^{\Z_2}_\Bbbk
\ee 
for $a\in\{0,1\}$. 
Hence these circle spaces are the zeroth Hochschild homology and cohomology, respectively, of the differential graded category of matrix factorisations, first computed in \cite{d0904.4713}. 
Moreover, for the Nakayama automorphisms we have $N^W_a = (-1)^{n\cdot (1-a)} \cdot 1_{C_a^{W}}$. 
  
\medskip

We now turn to the equivariant completion of $\LGs$ to look for $r$-spin TQFTs that can detect more $r$-spin structures than oriented TQFTs. 
One type of example+ of $\Delta$-separable Frobenius algebras with trivialisable $r$-th power of its Serre automorphism is the algebra~$A_G$ mentioned in Section~\ref{subsubsec:EquivariantCompletion}, in the case $G = \Z_r$. 

Recall from \cite[Sect.\,7.1]{cr1210.6363} that if a finite group~$G$ acts on $\Bbbk[x_1,\dots,x_n]$ such that a given $W \in \Bbbk[x_1,\dots,x_n]$ is invariant, this induces a $\Delta$-separable Frobenius structure on $A_G := \bigoplus_{g\in G} {}_g(1_W)$, where the $g$-twisted matrix factorisation ${}_g(1_W)$ is obtained from~\eqref{eq:IW} by replacing $x'_i \lmt g^{-1}(x'_i)$. 
Its Nakayama automorphism is 
\be 
\label{eq:NakayamaAG}
\gamma_{A_G} = \sum_{g\in G} \textrm{det}(g)^{-1} \cdot 1_{{}_g(1_W)} \, , 
\ee 
where $\textrm{det}(g)$ is the determinant of the $g$-action on $x_1,\dots,x_n$, cf.\ \cite[Sect.\,3.1]{BCP}. 

\begin{example}
	\label{exa:LGeqSimplesExample}
	For $r\in\Z_{\geqslant 3}$, we consider $(\Bbbk[x],x^r) \in \LGs$ with the $\Z_r$-action $\Z_r \lra \textrm{Aut}_\Bbbk(\Bbbk[x])$, $1\lmt (x\lmt \xi \cdot x)$, where 
	\be 
	\xi := \textrm{e}^{2\pi\I/r} \, . 
	\ee 
	Hence $W := x^r$ is invariant, and we have ${}_g(1_W) = (\Bbbk[x] \cdot 1) \oplus (\Bbbk[x] \cdot \theta)$ with 
	\be
	d_{{}_g(1_W)} 
		= 
		\frac{x'^r-x^r}{\xi^{-g} x' - x} \cdot \theta + \big( \xi^{-g}x' - x \big) \cdot \theta^*
	\ee 
	for $g\in\{ 0,1,\dots,r-1\}$. 
	Setting $A_{\Z_r} = \bigoplus_{g\in G} {}_g(1_W)$, we have an object 
	\begin{equation}
	\big((\Bbbk[x],x^r), A_{\Z_r}\big) \in (\LGs)_{\textrm{eq}}\,.
	\end{equation}
	
	For $g=0$ we have $\Hom(1_W,1_W) \cong \Bbbk \cdot \{ 1,x,\dots, x^{r-2} \}$ as a purely even $\Z_2$-graded vector space, while for $g\neq 0$ one finds that \cite[App.\,2]{BCP} 
	\be 
	\Hom\big( 1_W, {}_g(1_W) \big) 
		\cong
		\Bbbk \cdot \Big( \frac{x'^r-x^r}{(x'-x)(\xi^{-g} x' - x)} \cdot \theta + \theta^* \Big)[1]
	\ee 
	is a purely odd, 1-dimensional $\Z_2$-graded vector space. 
	Moreover, by \eqref{eq:SWLGIW} we have $S_W\cong 1_W$ in $\LGs$, and according to~\eqref{eq:NakayamaAG}, the Nakayama automorphism of $A_{\Z_r}$ is 
	\be 
	\label{eq:gammaAZr}
	\gamma_{A_{\Z_r}} 
		= 
		\sum_{g=0}^{r-1} \xi^{-g} \cdot 1_{{}_g(1_W)} \, . 
	\ee 
	
	We will use Lemma~\ref{lem:CaalphaAFromProjector} to compute the circle spaces $C_a \equiv C_a^{((\Bbbk[x],x^r),A_{\Z_r})}$. 
	Hence we have to identify the image of the projector 
	\begin{align}
	\Hom(1_W, A_{\Z_r}) & \lra \Hom(1_W, A_{\Z_r}) \nonumber 
	\\
	\tikzzbox{%
		\begin{tikzpicture}[very thick,scale=0.65,color=green!50!black, baseline]
		\fill (0,-0.3) circle (2.5pt) node[right] (D) {{\small $\varphi$}};
		\draw (0,-0.3) -- (0,0.8); 
		\end{tikzpicture} 
	}%
	& \lmt 
	\tikzzbox{%
		\begin{tikzpicture}[very thick,scale=0.65,color=green!50!black, baseline=0cm]
		\draw (0,0) -- (0,1.3);
		\fill (0,0) circle (2.5pt) node[below] {{\small $\varphi$}};
		\fill (-0.67,0) circle (2.5pt) node[left] {{\small $\gamma_{A_{\Z_r}}^{1-a}$}};
		\draw (0,0.8) .. controls +(-0.9,-0.3) and +(-0.9,0) .. (0,-0.8);
		\draw (0,-0.8) .. controls +(0.9,0) and +(0.7,-0.1) .. (0,0.4);
		\fill (0,-0.8) circle (2.5pt) node {};
		\fill (0,0.4) circle (2.5pt) node {};
		\fill (0,0.8) circle (2.5pt) node {};
		\draw (0,-1.2) node[Odot] (unit) {};
		\draw (0,-0.8) -- (unit);
		\end{tikzpicture}
	}%
	\, . 	
	\end{align}
	Expanding $\varphi\colon 1_W \lra A_{\Z_r}$ as $\sum_{g=0}^{r-1} \varphi_g$ with $\varphi_g \colon 1_W \lra {}_g(1_W)$ and using~\eqref{eq:gammaAZr}, a variant of~\eqref{eq:CaTwisted} reads
	\be 
	\label{eq:CaLG}
	C_a 
		\cong 
		\bigoplus_{g=0}^{r-1} \Bigg\{ 
		\; 
		\tikzzbox{%
			\begin{tikzpicture}[very thick,scale=0.65,color=green!50!black, baseline]
			\fill (0,-0.3) circle (2.5pt) node[right] (D) {{\small $\varphi_g$}};
			\draw (0,-0.3) -- (0,0.8); 
			\end{tikzpicture} 
		}%
		\; \Big|\; 
		\frac{1}{r}\sum_{h=0}^{r-1} 
		\tikzzbox{%
			\begin{tikzpicture}[very thick,scale=0.65,color=green!50!black, baseline=0cm]
			\draw (0,0) -- (0,1.3);
			\fill (0,0) circle (2.5pt) node[below] {{\small $\varphi_g$}};
			\fill (-0.67,0) circle (2.5pt) node[left] {{\small $\xi^{-h(1-a)}$}};
			\draw (0,0.8) .. controls +(-0.9,-0.3) and +(-0.9,0) .. (0,-0.8);
			\draw (0,-0.8) .. controls +(0.9,0) and +(0.7,-0.1) .. (0,0.4);
			\fill (0,-0.8) circle (2.5pt) node {};
			\fill (0,0.4) circle (2.5pt) node {};
			\fill (0,0.8) circle (2.5pt) node {};
			\draw (0,-1.2) node[Odot] (unit) {};
			\draw (0,-0.8) -- (unit);
			\end{tikzpicture}
		}%
		=
		\tikzzbox{%
			\begin{tikzpicture}[very thick,scale=0.65,color=green!50!black, baseline]
			\fill (0,-0.3) circle (2.5pt) node[right] (D) {{\small $\varphi_g$}};
			\draw (0,-0.3) -- (0,0.8); 
			\end{tikzpicture} 
		}%
		\Bigg \}[1-a]
		\, , 
	\ee 
	where we used~\eqref{eq:SWLGIW} and the isomorphism $\Hom^i(X,Y[1]) \cong \Hom^{i+1}(X,Y)$ in $\LGs$. 
	A direct computation along the lines of \cite[App.\,2]{BCP} then reveals that 
	\be 
	\tikzzbox{%
		\begin{tikzpicture}[very thick,scale=0.65,color=green!50!black, baseline=0cm]
		\draw (0,0) -- (0,1.3);
		\fill (0,0) circle (2.5pt) node[below] {{\small $\varphi_g$}};
		\fill (-0.67,0) circle (2.5pt) node[left] {{\small $\xi^{-h(1-a)}$}};
		\draw (0,0.8) .. controls +(-0.9,-0.3) and +(-0.9,0) .. (0,-0.8);
		\draw (0,-0.8) .. controls +(0.9,0) and +(0.7,-0.1) .. (0,0.4);
		\fill (0,-0.8) circle (2.5pt) node {};
		\fill (0,0.4) circle (2.5pt) node {};
		\fill (0,0.8) circle (2.5pt) node {};
		\draw (0,-1.2) node[Odot] (unit) {};
		\draw (0,-0.8) -- (unit);
		\end{tikzpicture}
	}%
	\cong
	\begin{tikzpicture}[
	baseline=(current bounding box.base),
	descr/.style={fill=white,inner sep=3.5pt},
	normal line/.style={->}
	]
	\matrix (m) [matrix of math nodes, row sep=0.5em, column sep=1em, text height=1.5ex, text depth=0.1ex, cells={anchor=west}] {%
		\xi^{h(a-1)-hj} \cdot \varphi_g & \textrm{if $g=0$ and $\varphi_g=x^j$}
		\\
		{} & {}
		\\
		\xi^{ha} \cdot \varphi_g & \textrm{if $g\neq0$.}
		\\
	};
	\draw [decorate,decoration={brace,amplitude=10pt},xshift=-140pt,yshift=0pt]
	(1.4,-1.2) -- (1.4,1.2) node [midway,xshift=-0.6cm]
	{};
	\end{tikzpicture}
	\ee 
	Hence the summand for $g=0$ in~\eqref{eq:CaLG} is~0 for $a=0$, and equal to the 1-dimensional $\Z_2$-graded vector space $\Bbbk\cdot x^{a-1}[1-a]$ otherwise, while the summands for $g\neq 0$ contribute only if $a=0$, namely a term $\Hom(1_W,{}_g(1_W))[1] \cong \Bbbk[2] = \Bbbk$: 
	\be 
	\label{eq:CaInLG}
	C_a 
		\cong
		\delta_{a\geqslant 1} \cdot \Bbbk[1-a] \oplus \delta_{a,0} \bigoplus_{g=1}^{r-1} \Bbbk \, . 
	\ee 
	It follows that the quantum dimension of~$C_a$ (as an object in $\LG(\one,\one) \cong \textrm{vect}^{\Z_2}$, i.\,e.\ as a super vector space) is $r-1$ for $a=0$, $+1$ for $a\in\Z\setminus\{0\}$ even, and $-1$ for~$a$ odd. 
	However, in $\LGs$ the 2-morphisms $(+1)\cdot 1_{C_a}$ and $(-1)\cdot 1_{C_a}$ are identical. 
	Recalling Proposition~\ref{prop:torus-inv-dim-Cd}, this means that the $(\LGs)_{\textrm{eq}}$-valued TQFT associated to $A_{\Z_r}$ can only distinguish two $r$-spin structures on the torus (for $r\neq 2$), while there are as many non-diffeomorphic $r$-spin structures on~$T^2$ as there are divisors of~$r$. 

	We emphasise that this example works for arbitrary $r\geqslant 3$,
	whereas the example computing the Arf invariant mentioned in Remark~\ref{rem:PivotalNotGoodForSpin}\ref{item:EulerArf} is defined only for~$r$ even.
\end{example}

The computational techniques used in Example~\ref{exa:LGeqSimplesExample} can analogously be applied to more involved examples. 
For instance, there are $\Z_r$-actions on $\Bbbk[x_1,x_2]$ which leave $W=x_1^r + x_2^{2r}$ invariant, and the associated $(\LGs)_{\textrm{eq}}$-valued TQFTs may detect more than two $r$-spin structures on the torus. 
We leave such computations as well as the application of the theory developed in Section~\ref{subsec:EquivariantCompletion} to the 2-category of $\Q$-\textsl{graded} Landau--Ginzburg models $\LGgr$ (see \cite{cm1208.1481} or \cite[Sect.\,2.5]{CMM}) to future work.

\end{document}